\def\real{{\tt I\kern-.2em{R}}}  
\def\nat{{\tt I\kern-.2em{N}}}     
\def\eps{\epsilon}

%This appears to be the correct archive version as of file date. (The date of the v9,tex file that is.) This contains all of the corrections. This is the only correct file as of file data. The EC files have not been corrected. The arxiv files are correct as of the date of the stored arxiv version. This stored version has been broken down into these two files. There is a major change as of 11/22/07 in interpretation of light-clock alterations. This file was downloaded to arxiv on 6/7/09. I have removed, 7/12/09, the nonreversible hyposthesis as unnecessary and have redone the derivation for (13) and (14). Major additional changes made on 8/13/09. Downloaded on 8/22/09. On 8/25/09, and 9/25/09, 9/30/09 a few additions and important corrections and changes have been made. They are in website book. The first sentence of [12] is incorrect and I have changed it. But, I will not be downloaded to arxiv until 2010 or 2011. (This was done). It has been correctly stated on my web site files.  On 17 SEP 2012, I change slightly my statement relative to Hilbert and Einstein. 

%I am dowloading on to arxiv.org on 2/14/2014 all of the improvements I have made since the lst version of 10/26/2010. This includes note 

\def\realp#1{{\tt I\kern-.2em{R}}^#1}
\def\natp#1{{\tt I\kern-.2em{N}}^#1}
\def\hyper#1{\,^*\kern-.2em{#1}}
\def\monad#1{\mu (#1)}

\def\St#1{{\tt st}#1}
\def\st#1{{\tt st}(#1)}
\def\hyperreal{{^*{\real}}}
\def\hyperrealp#1{{\tt ^*{I\kern-.2em{R}}}^#1} 
\def\hypernat{{^*{\nat }}}
\def\hypernatp#1{{{^*{{\tt I\kern-.2em{N}}}}}^#1} 
\def\eskip{\hskip.25em\relax}

\def\Hyper#1{\hyper {\eskip #1}}
\def\leaderfill{\leaders\hbox to 1em{\hss.\hss}\hfill}
\def\srealp#1{{\rm I\kern-.2em{R}}^#1}

\def\power#1{{{\cal P}(#1)}}

\def\pars{\par\smallskip}
\def\parm{\par\medskip}
\def\r#1{{\rm #1}}
\def\b#1{{\bf #1}}
\def\ref#1{$^{#1}$}

\def\m@th{\mathsurround=0pt}
\def\rightarrowfill{$\m@th \mathord- \mkern-6mu \cleaders\hbox{$\mkern-2mu 
\mathord- \mkern-2mu$}\hfil \mkern-6mu \mathord\rightarrow$}
\def\leftarrowfill{$\mathord\leftarrow
\mkern -6mu \m@th \mathord- \mkern-6mu \cleaders\hbox{$\mkern-2mu 
\mathord- \mkern-2mu$}\hfil $}
\def\noarrowfill{$\m@th \mathord- \mkern-6mu \cleaders\hbox{$\mkern-2mu 
\mathord- \mkern-2mu$}\hfil$}
\def\orgate{$\bigcirc \kern-.80em \lor$}
\def\andgate{$\bigcirc \kern-.80em \land$}
\def\inverter{$\bigcirc \kern-.80em \neg$}

 \magnification=1000
\tolerance 10000
\hoffset=0.25in
\hsize 6.00 true in
\vsize 8.75 true in
\font\eightsl=cmsl8
\pageno=1
\headline={\ifnum\pageno= 1\hfil {\quad}
\hfil\else\ifodd\pageno\rightheadline \else\leftheadline\fi\fi} 
\def\rightheadline{\eightsl \hfil Nonstandard Analysis Applied to Relativity \hfil} 
\def\leftheadline{{} \eightsl \hfil Nonstandard Analysis Applied to Relativity  \hfil} 
\voffset=1\baselineskip

\def\id{\par\hangindent2\parindent\textindent}
\def\textindent#1{\indent\llap{#1}}
\baselineskip 14pt
\centerline{}
\vskip 1.00in
 %===SIMPLE TABLE OF ONE BOX WITH MULTIPLE LINES 
%IT IS IN DISPLAYED FORM WHEN PRINTED
$$\vbox{\offinterlineskip
\hrule
\halign{&\vrule#&
\strut\quad\hfil#\hfil\quad\cr
height2pt&\omit&\cr
&\quad&\cr
&{\bf Nonstandard Analysis Applied to}&\cr
&{\bf Special and General Relativity -}&\cr
&{\bf The Theory of Infinitesimal Light-Clocks}&\cr
&\quad&\cr
height2pt&\omit&\cr}
\hrule}$$
\vskip 1.00in
\centerline{\bf Robert A. Herrmann}
\vskip 4.50in
\centerline{SPECIAL ARXIV.ORG EDITION}
\vfil
\eject
{\quad}

\vskip 2.75in
\centerline{\copyright 1995 Robert A Herrmann under the title Einstein Corrected.  This version \copyright 2014 Robert A. Herrmann}\par\bigskip
{\leftskip 1.00in \rightskip 1.00in \noindent Any portion of this monograph may be reproduced, without change and giving proper authorship, by any method without seeking permission and without the payment of any fees. This is the 10 point font version.\par}\par 
\vfil\eject
\centerline{}
\centerline{\bf CONTENTS}
\bigskip
\line{\indent \ \ Preface\leaderfill 5}
\bigskip
\centerline{Article 1}\pars
\centerline{\bf Foundations and Corrections to Einstein's}
\centerline{\bf Special and General Theories of Relativity.}\pars
\line{\indent 1. \ \ Introduction\leaderfill 8}
\line{\indent 2. \ \ Some Ether History\leaderfill 8}
\line{\indent 3. \ \ The Calculus\leaderfill 9}
\line{\indent 4. \ \ Relativity and Logical error\leaderfill 11}
\line{\indent 5. \ \ A Privileged Observer\leaderfill 12}
\line{\indent 6. \ \ The Fock Criticism and Other\hfil} 
\line{\indent \ \ \ \ \ \ \ Significant Matters\leaderfill 14}
\line{\indent 7. \ \ Why Different Derivations?\leaderfill 18}
\line{\indent 8. \ \ A Corrected Derivation\leaderfill 19}
\line{\indent \ \ \ \ \ \ References for part 1\leaderfill 21}
\bigskip
\centerline{ Article 2}\pars
\centerline{\bf A Corrected Derivation for the Special Theory of Relativity}
\line{\indent 1. \ \ The Fundamental Postulates\leaderfill 23}
\line{\indent 2. \ \ Pre-derivation Comments\leaderfill 24}
\line{\indent 3. \ \ The Derivation\leaderfill 25}
\line{\indent 4. \ \ The Time Continuum\leaderfill 33}
\line{\indent 5. \ \ Standard Light-clocks and c\leaderfill 34}
\line{\indent 6. \ \ Infinitesimal Light-clock Analysis\leaderfill 35}
\line{\indent 7. \ \ An Interpretation\leaderfill 39}
\line{\indent 8. \ \ A Speculation and Ambiguous Interpretations\leaderfill 41}
\line{\indent 9. \ \ Reciprocal Relations\leaderfill 44}
\line{\indent \ \ \ \ \ \ NOTES\leaderfill 46}
\line{\indent \ \ \ \ \ \ References\leaderfill 49}
\line{\indent \ \ \ \ \ \ Appendix-A\leaderfill 50}
\bigskip
                                                     
\centerline{Article 3}\pars
\centerline{\bf Foundations and Corrections to Einstein's}
\centerline{\bf Special and General Theories of Relativity.}\pars
\line{\indent 1. \ \ Some Special Theory Effects\leaderfill 55}
\line{\indent 2. \ \ General Effects\leaderfill 58}
\line{\indent 3. \ \ Relativistic Alterations\leaderfill 60}
\line{\indent 4. \ \ Gravitational Alterations\leaderfill 65}
\line{\indent 5. \ \ Nonstandard Particle Medium Analysis\leaderfill 67}
\line{\indent 6. \ \ Minimizing Singularities\leaderfill 69}
\line{\indent 7. \ \ Applications\leaderfill 72}\vfil\eject
\line{\indent 8. \ \ Prior to Expansion$,$ Expansion and \hfil}
\line{\indent \ \ \ \ \ \ Pseudo-White Hole Effects\leaderfill 74}
\line{\indent \ \ \ \ \ \ NOTES\leaderfill 77}
\line{\indent \ \ \ \ \ \ References for part 2\leaderfill 90}
\line{\indent \ \ \ \ \ \ Appendix-B\leaderfill 93}
\bigskip    
\line{\indent{\bf Index}\leaderfill 103}
\vfil\eject
\centerline{\bf Preface}\par\smallskip
It is actually dangerous for me to present the material that appears within 
this book due to the usual misunderstandings. Any scientist who claims that there are fundamental errors within the foundational 
methods used to obtain Einstein's General and Special Theories of 
relativity may be greatly ridiculed by his colleagues who do not read carefully. The reason for this 
has nothing to do with science but has everything to do with scientific 
careers$,$ research grants and the like. Thousands upon thousands of individuals 
have built their entire professional careers upon these two theories and their 
ramifications. The theoretical science produced is claimed to be  ``rational'' 
since it follows the patterns of a mathematical structure. As a mathematician 
who produces such structures$,$ it is particular abhorrent to the scientific 
community if I make such a claim.  Mathematicians seem to have an unsettling effect upon some members of the physical science community$,$ especially when a mathematician delves into a natural science. After all$,$ it was the 
mathematician Hilbert who$,$ by application of the calculus of variations, derived the so-called Einstein gravitational field equations and was actually the first to present the tensor expression in a 
public form. On the other hand, Einstein was very proficient in applying the field expressions to physical situations in order to predict verifiable physical implications. \par

Now please read the following very carefully. The results presented here and in my published papers on this subject are not intended to denigrate those scientists who have$,$ in the past$,$ contributed to these Einsteinian theories or who continue to do so. The corrections I have made are in the foundations for these theories. The corrections are totally related to how the results are interpreted physically. These corrections do not contradict the results obtained when the Einsteinian approach and the language used are considered as {\bf models} for behavior. These corrections are based upon newly discovered rules for rigorous infinitesimal modeling. These results may also be significant to those that hold to the belief that many events within the natural world are produced classically by a zero-point radiation field.\par

Many 
unqualified individuals continue to present their own alternatives to these 
Einstein theories$,$ some claiming that the results are but an exercise in high-school algebra. Certain scientific groups tend to categorize any and all 
criticisms of the Einstein theories as coming from the unqualified and lump 
such criticisms into the same unworthy category. However$,$ many highly 
qualified scientists of the past such as Builder$,$ Fock$,$ Ives and Dingle have 
made such claims relative to the foundations of these two theories.  
For Ives$,$ the fundamental approach was to assume that such a thing as length 
contraction$,$ and not time dilation$,$ is a real natural effect and it is this 
that leads
to 
the Einstein conclusions. In order to eliminate these 
criticisms$,$ Lawden states the ``modern'' interpretation that length 
contraction has no physical meaning$,$ and only ``time dilation'' is of 
significance. This modern assumption is certainly rather ad hoc in character. 
Further$,$ many theory paradoxes still appear within the literature and 
are simply ignored by the scientific community. There is$,$ however$,$ a reason 
for this.\par
 
The actual approach used can now be shown explicitly to contain 
logical error. It was not possible to show this until many years after the 
theory was fully developed. Further$,$ the original approach utilizes a 
``geometric language$,$'' a language that has been criticized by many including 
John Wheeler as the incorrect approach to analyze the fundamental behavior of 
universe in which we dwell. Although Einstein used an explicit operational 
approach in his Special Theory$,$ he was unable to use a mathematical approach 
that encapsules his operational definition since the actual mathematics was 
not discovered until 1961. He used what was available to him at the time. In 
this book$,$ all such errors and paradoxes are removed by use of the modern 
corrected theory of infinitesimals and infinite numbers as discovered by 
Robinson. Moreover$,$ the recently discovered correct rules for infinitesimal 
modeling are used$,$ and this eliminates the need for tensor analysis and 
Riemannian geometry.\par 

The logical errors 
occur when rigorous mathematical procedures are applied but the abstract
mathematical theory uses modeling procedures (i.e. interpretations) that 
specifically contradict mathematical modeling requirements. These errors are 
detailed within this book. Unfortunately$,$ the same confused 
approach often pervades most ``physical'' interpretations for mathematical 
structures.\par

As mentioned$,$ 
in the articles contained in this book$,$ these errors are eliminated by 
application of the corrected theory of infinitesimals as discovered by 
Abraham Robinson. However$,$ in so doing$,$ significant differences in 
``philosophy'' appear necessary. These differences are amply discussed 
in part 1 of the first article. Note that each article begins with an 
extensive abstract the contents of which I will not reproduce in this 
preface.\par

The results in this book eliminate all of the known controversies 
associated with these two theories. 
Indirectly$,$ the results show the logical existence of a 
privileged frame of reference. From very elementary assumptions gleaned from 
laboratory observation$,$ it is shown that there is neither absolute time 
dilation nor length contraction$,$ but there is an alteration in one and only 
one mode of time and length measurement due to relative velocity$,$ (i.e. velocity) 
potential velocities$,$ textual expansion and the like. These alterations have 
a ``physical'' cause. The apparent alteration of 
mass$,$ the gravitational redshift$,$ the transverse Doppler effect and all other 
verified consequences of these two theories are predicted and shown to have  
``physical'' causes -- causes associated with a nonstandard particle medium. \par 

The order in which 
these articles appear in this book is somewhat reversed from the order in 
which they 
were written. Article 2$,$ {\it A Corrected Derivation for the Special Theory of 
Relativity}$,$ was written first and presented first before the 
Mathematical Association of America. Article 1 and Article 3$,$ {\it Foundations and Corrections 
to Einstein's Special and General Theories of Relativity$,$ Article 2} and {\it Article 3}$,$ 
were written from November 1992 -- Sept. 1993. However$,$ Articles  1 and 3
contain$,$ almost exclusively$,$ classical mathematics$,$ (with minor exceptions)
and are more easily 
comprehended by scientists versed in this subject. Article 2 requires the 
additional concepts associated with elementary nonstandard analysis. Article 
3 can be comprehended relative to the results presented. 
The necessary 
formal infinitesimal theory presented in Article 3 should not detract from this 
comprehension. The material in article 3 was partially funded by a research grant from the U. S. Naval Academy Research Council. [Added May 2001. Listed below and elsewhere are four published  
journal articles and a few archived references relative to the methods presented within this book. These references update some of the reference information listed at the end of each article.]\par

It is 
hoped that the conclusions developed throughout these articles that are 
ultimately dependent upon the concepts of nonstandard analysis will motivate 
the scientific community to become more conversant with proper infinitesimal 
modeling techniques. \par 

\vskip 18pt
\leftline{Robert A. Herrmann Ph.D.} 
\leftline{Annapolis$,$ MD}
\leftline{August 8$,$ 1995}\bigskip
\centerline{\bf References }\par\medskip
 Herrmann, R. A. 1992.  A corrected derivation for the Special Theory of 
relativity. Presented before the Mathematical Association of America$,$ Nov. 14$,$ 1992 
at Coppin State College$,$ Baltimore$,$ MD.\par\smallskip
 Herrmann, R. A. 1993. The Theory of Ultralogics. \hfil\break $<$http://www.arXiv.org/abs/math.GM/9903081$>$ and 9903082.\par\smallskip
Herrmann, R. A. 1994. The Special Theory and a nonstandard substratum. {\it Speculations in Science and Technology} 17(1):2-10. http://www.arXiv.org/abs/physics/0005031\par\smallskip
Herrmann, R. A. 1995. Operator equations$,$ separation of variables and relativistic alterations. {\it Intern. J. Math. \& Math. 
Sci.} 18(1):59-62. http://www.arxiv.org/abs/math-ph/0312005
\par\smallskip
Herrmann, R. A. 1996.  An operator equation and relativistic 
alterations in the time for radioactive decay. {\it Intern. J. Math. \& Math. 
Sci.} 19(2):397-402. \hfil\break http://www.arxiv.org/abs/math-ph/0312006
\par\smallskip
Herrmann, R. A. 1997. A hypercontinuous$,$ hypersmooth Schwarzchild line-element transformation. {\it Intern. J. Math. \& Math. 
Sci.} 20(1):201-204. http://www.arxiv.org/abs/math-ph/0312007\par\smallskip
Mehra, J, 1974. Einstein, Hilbert and The Theory of Gravitation, D. Reidel Co., Dondrecht, Holland

\vfil\eject

\centerline{\bf Special and General Theories of Relativity, Article 1} 
\bigskip
\centerline{Robert A. Herrmann}
\bigskip
 Abstract:  In Article 1 of this paper$,$ newly identified logical errors in the 
derivations that yield Einstein's 
Special and General Theories of Relativity are discussed. These errors are 
much more significant than those identified by Fock.
The basic 
philosophy of science used as a foundation for these theories is identified. 
The philosophy of the 
privileged observer is detailed. Artricle 1 concludes with a brief discussion 
of a new derivation for the Lorentz transformation that eliminates all the 
logical errors associated with previous derivations as well as eliminating the 
controversial concepts of ``length contraction'' and ``time dilation.''
\bigskip
\leftline{\bf 1. Introduction}\par\bigskip
Nobel prize winner Max Planck (1932$,$ p. 2) wrote:\pars
{\leftskip=0.5in \rightskip=0.5in \noindent Nature does not allow 
herself to be exhaustively expressed in human thought.\par}\parm
\noindent The D-world mathematical model as discussed by this author (Herrmann 
1990$,$ p. 12) is used to analyze the linguistic methods used by modern 
science and develops the following rational possibility. 
\parm
{\leftskip=0.5in \rightskip=0.5in \noindent  
Human beings do not have the ability 
to comprehend and will not eventually describe in 
human languages all of the true laws or natural events that govern the cosmos. 
This includes the 
laws or natural events that govern the development of individual natural 
systems. \par}\parm
\noindent Hence$,$ the philosophy of science  as espoused by {Planck}$,$ and denoted 
by (A)$,$  
can be 
logically argued for by using mathematical reasoning.\par
On the other hand$,$ the philosophy of {scientism}$,$ which is denoted by (S)$,$ 
assumes the negative 
of Plank's statement. One
 abiding rule of scientism is that nature will allow 
herself to be exhaustively 
expressed in human thought. Nowhere in modern science do these two 
contradictory philosophies of science clash more violently than in 
the {ether or medium concept associated with electromagnetic propagation}.  
 Our concern in
the first sections
 of this paper$,$ will be the \ae{ther} concept$,$ a concept 
that is partially discussed in numerous  journal papers. (See for example Benton$,$ 
1988).
\parm
\leftline{\bf 2. Some \AE{ther} History}\parm
{Newton} attempted to imagine properties of a {medium} for his universal 
theory of gravity for in his letters to Boyle he apparently stated that: 
\pars
{\leftskip=0.5in \rightskip=0.5in \noindent . . . he found he was unable$,$ from 
experiment and observation$,$ to give a satisfactory account of this medium and 
the manner of its operation in producing the chief phenomena of 
nature (Maxwell 1965b$,$ p. 487).\par}\pars  
{Maxwell} (1965b$,$ p. 764) tells us that the only medium that survived as a 
structure and that seems to uphold the propagation of light was that proposed 
by {Huygens}.
It was {Thomson} (1854) who did most of the calculations as to the mechanical 
properties that such a {``luminiferous  \ae{ther}''} should possess. 
Then$,$ in 1864$,$ 
Maxwell (1965a$,$ pp. 526 -- 597)  outlined the general properties he would 
imagine to hold for an electromagnetic field and proposed that the 
luminiferous \ae{ther} and the electromagnetic medium are the same.\par
The mechanical properties of such an electromagnetic \ae{ther} needed to be 
expressed 
in a scientific language and required the basic methods of mathematical 
deduction. Apparently$,$ most of  the believers in such an \ae{ther} were 
following 
philosophy (S) with the additional requirement of {{\it absolute realism}}. Absolute 
realism$,$ in 
this case$,$ means that an individual corresponds specific physical terms 
to a list of selected abstract 
mathematical 
terms 
 and for a ``physical theory'' to be accepted those 
physical entities being named by the physical terms 
are assumed to exist in reality. 
Unfortunately$,$ under the usual 
correspondence$,$ 
various {ether calculations} produce physical behavior that seemingly can 
not exist 
within our universe. One such difficulty was relative to \ae{ther} stresses. \par
{Lorentz} (1952$,$ p. 31) proposed  that the concept of realism be altered. \pars   
{\leftskip=0.5in \rightskip=0.5in \noindent I should add that$,$ while denying 
the real existence of \ae{ther} stresses$,$ we can still avail ourselves of all of 
the mathematical transformations by which the application of the formula (43) 
may be made easier. We need not refrain from reducing the force to a 
surface-integral$,$ and for convenience's sake we may continue to apply to the 
quantities occurring in this integral the name stresses. Only$,$ we must be 
aware that they are only imaginary ones$,$ nothing else than auxiliary 
mathematical quantities.\par}\pars
\noindent However$,$ altering realism in this manner would slightly remove  
human deductive logic$,$ as it is encapsulated within  mathematical reasoning$,$  
as a bases for (S). The burning question would be why the {human mind 
needs these ``imaginary'' entities} to properly describe the behavior of 
the \ae{ther}? Under the (S) philosophy$,$ either the \ae{ther} did not correspond to 
reality or there would need to be  new concepts developed and corresponding 
physical terms defined. But the situation is more complex than this since realism 
always depends upon a theory of correspondence. \parm
\leftline{\bf 3. The Calculus}\parm
Newton's concept of mathematical modeling was firmly rooted in the natural 
world. \pars
{\leftskip 0.5in \rightskip 0.5in \noindent  
Geometry does not teach us to draw lines$,$ but requires them to be drawn$,$ 
for it requires that the learner should first be taught to describe these 
accurately before he enters geometry$,$ then it shows how by these operations 
problems may be solved. To describe right lines and circles are problems$,$ but 
not geometrical problems. The solution of these problems is required from 
mechanics$,$ . . . . therefore geometry is founded in mechanical practice$,$ and is 
nothing but that part of universal mechanics which accurately proposes and 
demonstrates the art of measure (Newton$,$ 1934 p. xvii).\par}\pars
 Newton's claim is that 
our observations  and intuitive comprehension of mechanics  
comes first and then 
these concepts are abstracted to include the vague notion 
that objects have certain ``capacities or potentials to do things.'' 
We are told that it is {\it after} experimentation$,$ observation and 
reflection that the mathematical structure is evoked and these ``easy'' 
capacity concepts are modeled.\par
Newton used the language of {{\it infinitesimals}} within all of his applied 
mathematics. To him$,$ these infinitesimal quantities existed in objective 
reality$,$ they referred to measures of actual natural world behavior.
For example$,$ in 1686$,$ Newton  explains to his critics 
what he claims is the easily 
comprehended notion of the {{\it ultimate velocity}}$,$ or what we now term as
 {{\it instantaneous velocity}}$,$ for an actual real material object.\par
{\leftskip 0.5in \rightskip 0.5in \noindent But by the same argument it may be alleged 
that a body arriving at a certain place$,$ and there stopping$,$ has no ultimate 
velocity; because the velocity$,$ before the body comes to the place$,$ is not its 
ultimate velocity; when it has arrived$,$ there is none. But the answer is easy; 
for by the ultimate velocity is meant that with which the body is moved$,$ 
neither before it arrives at its last place and the motion ceases$,$ nor after$,$ 
but at the very instant it arrives; that is$,$ the velocity with which the body 
arrives at its last place$,$ and with which the motion ceases (Newton$,$ 1934$,$ 
pp. 38--39).\par}\pars  
The abstract notion of instantaneous velocity may have been ``easy'' for Newton 
to grasp$,$ but it was incomprehensible to {Berkeley} and many others who 
believed that such abstractions could not be applied to {actual real material 
objects}. The paramount philosophy of science for Berkeley was a science of the 
material and {directly observed universe}. Any arguments that relied upon 
such abstractions would need to be rejected.\par
Newton's approach created a {schism} in the philosophy of mathematical modeling. 
One group of scientists believed that there exists actual real world entities  
that can be characterized in terms of infinitesimal measures of 
time$,$ mass$,$ 
volume$,$ charge$,$ and the like. Another 
group assumed that such terms are auxiliary in character and do not 
correspond to objective reality. \par
In the early 19th century$,$ {Cauchy} (1821)  using the language of infinitesimals established a 
result that {Able} (1826) showed by counter-example could not be logically 
correct. However$,$ no direct error can be found in the Cauchy's logical 
argument. Hence the intuitive assumptions underlying the behavior of 
Newton's infinitesimals must be logically contradictory. 
Unfortunately$,$ the vast majority of the scientific community  still use  
Newton's infinitesimal concepts. How many know that such a use can lead to
logical contradictions? \par
A few years ago$,$ {Robinson} (1966)  
removed the contradictions from the theory of 
infinitesimals. The new corrected mathematical structure 
allows for rigorously defined modeling processes
(Herrmann$,$ 1991a). This mathematical structure has a great deal to say about 
the mathematical {schism} mentioned above. For example$,$ instantaneous velocity$,$ 
acceleration and even Newton's second law of motion can be derived from  
fundamental properties of the infinitesimal world rather than simply being defined
(Herrmann$,$ 1991a$,$ part 2$,$ pp. 4 --15). Moreover$,$ infinitesimals do not behave
in exactly the same manner as do directly observed real objects. \par
  The schism mentioned above has been ignored by most modern day scientists.
Such scientists simply use the old contradictory infinitesimal language 
to derive 
expressions that are claimed to 
model real world physical behavior without discussing the realism question for 
portions of their derivations. 
The mathematical 
model called {{\it the nonstandard physical world model}} (i.e. NSP-world) 
uses the corrected 
theory of the infinitesimally small and infinitely large$,$ with other 
techniques$,$  
along with 
a new physical language theory of correspondence. Many of these new entities 
need not exist within the natural world since their behavior differs greatly
from any known natural world entity. \par
 The above {schism in mathematical 
modeling} has now become more evident for$,$ from the viewpoint of {{\it indirect 
evidence}}$,$ these new entities might actually exist in an NSP-world in which the 
natural world is specially embedded. Indeed$,$ one could say 
that the NSP-world is omnipresent with respect to the natural world and 
upholds$,$ sustains and guides natural system behavior.  
Notwithstanding this possibility$,$ science cannot eliminate the NSP-world from its 
derivations if mathematical models are used. The NSP-world is always 
lurking in the background. Of recent interest is the possible existence of a 
natural world space medium (Barnes$,$ 1986) that might be operationally obtained 
as 
the standard part  
of an NSP-world nonstandard photon-particle medium (NSPPM), where 
photon particle properties are used. Photon behavior uses portion of the subparticle field. This NSPPM is closely related to the basic Lorentz 
transformation. The NSPPM is not to be taken as a medium that follows the nonstandard extension of the Maxwell field equations.\parm
\vfil\eject
\leftline{\bf 4. Relativity and Logical Error} \parm
Using the (S) philosophy$,$ {Einstein} (1979) wrote that at the age of sixteen$,$ using a {{\it conceptual  
observer}}$,$   
he had concocted a mind 
experiment relative to the known properties of light propagation that seemed 
to imply a paradox. He claimed to have eliminated this paradox some ten years 
later (Einstein$,$ 1905).  When this author was 
twelve$,$ he read with interest Einstein's major book on this subject 
(Einstein$,$ 1945)  and felt that there might be some logical 
error in the basic derivations. However$,$ a resolution of this error had to 
wait until the mathematics itself was correct by Robinson. \par
How did Einstein arrive at his derivation? \pars
{\leftskip 0.5in \rightskip 0.5in \noindent Then I myself wanted to verify the flow of the 
ether with respect to the Earth$,$ in other words$,$ the motion of the Earth. When 
I first thought about this problem$,$ I did not doubt the existence of the \ae{ther} 
or the motion of the Earth through it (Einstein$,$ 1982$,$ p. 46).\par}\pars
Einstein did not state in his 1922 lecture that the \ae{ther} 
did not exist. He said$,$ ``Since then I have come to believe that the motion of 
the Earth cannot be detected by any optical experiment though the Earth is 
revolving about the Sun'' (1982$,$ p. 46). However$,$ 
later he 
and Infeld specifically argue against the \ae{ther} in their 
popularizing book on scientific intuition 
(Einstein and Infeld$,$ 1938$,$ pp. 157--186). 
The argument is based entirely upon 
the (S) philosophy and erroneous hidden assumptions. The most glaring  
assumption is that if they could not describe an \ae{ther} that satisfies the 
experimental conditions$,$ then it does not exist. \par
The {bases for Einstein's original 1905} paper is that  ``Time cannot be absolutely 
defined$,$ and there is an inseparable relation between time and signal 
velocity (Einstein$,$ 1982$,$ p. 46).'' Although Einstein states that an {absolute 
time} is not definable$,$ it has been shown that this claim is false (Herrmann$,$ 
1992). This immediately implies that the Einstein derivation for the Lorentz 
transformation is inconsistent with respect to a basic  premise. 
Einstein claims to be interested in an {operational definition} and 
first uses the term {``clock''}   as meaning  {\it any} measure of time within 
the natural world  without further defining  the apparatus. This does not 
immediately contradict the concept of absolute time as not being definable.                   
But then he restricts the characterization of such clocks by adding 
light propagation terminology relative to their synchronization. Certain 
distances  
are also defined
in terms of these restricted clocks and a property of light propagation. The 
 {predicate}  that is interpreted as ``any time measure within the natural world'' 
has 
been restricted to natural world clocks that are synchronizable 
by light propagation techniques.  
There always remains the possibility that not all identified natural world  
clocks are thus synchronizable. It can be argued that some biological clocks
 fall 
into this category. \par
Einstein next 
extends the domain of  
these ``times'' to include the local absolute {Newtonian time continuum} and  
applies the infinitesimal calculus to these``times.'' The  
infinitesimals represent a modeling {\it concept} used in the calculus to approximate a continuum of real numbers and there is no logical error if this technique is used consistently. The more closely the behavior of the measuring devices is modeled by real number properties$,$ then the better this approximating mathematical device will predict behavior.    
As mentioned, Einstein introduced the ``operational'' time notion by requiring his ``times'' to be restricted by synchronization techniques and for ``proper'' time the ``Einstein measure'' technique. The basic logical error occurs later when these restrictions are dropped and the results are extended to ``time'' as a general concept. [For example, see Einstein 1907.] Such an error occurs when one   
substitutes nonequivalent predicates and is called the  
 {{\it model theoretic error of generalization}}. A 
statement that holds for a specific domain (time restricted by the  
language of light propagation) 
cannot be extended to a
 domain that refers to time as measured by {\it any} device.         
For example$,$ the statement that the usual ordering 
of the integers is a well-ordering  does 
not hold 
when the domain is extended to the rational numbers.\par 
Through use of a standard partial derivative technique$,$ Einstein derives the {Lorentz 
transformation}. On the left-hand side of the equations is the proper 
time which is a measure of time defined in a slightly different manner using  
his synchronized  natural  world time concept. On the  
right-hand side of these equations$,$ time is expressed in the original 
synchronized mode. \par 
In the section of this 1905 paper entitled ``Physical meaning of the equations 
obtained with respect to moving rigid bodies and moving clocks,'' the model theoretic error of generalization occurs. Einstein removes from his clock 
concepts the additional language of light propagation and 
 generalizes the time concept
to any measure of time whether correlated to  
light propagation or 
not. Throughout the remainder of this paper$,$ the time measures, for any clock,  
utilizes the absolute Newtonian 
time continuum with infinitesimals so that the Calculus may be applied. 
Thus, in this paper,     
nonequivalent predicates are assumed equivalent.      
\par
Modern derivations of the Lorentz transformation
 (Bergmann$,$ 1976; 
Lawden$,$ 1982)   
do not apply Einstein's first infinitesimal approach but rather 
start with two coordinate statements $x^2 + y^2 + z^2 = c^2t^2$ and $\overline{x}^2
+ \overline{y}^2 + \overline{z}^2 =c^2 \overline{t}^2.$ The time measures are 
once again restricted by light propagation language and the requirement of 
synchronization. Further$,$ for this previous approach, 
these expressions are 
obtained by application of {spherical wavefront (light) concepts}$,$ 
and the assumed constancy of 
the measure of the velocity of such propagation.  
Notice that these expressions 
require  measurement of the 
propagation velocity to be made in accordance with devices that include 
the restricting light propagation language. 
After the derivation of the Lorentz 
transformation$,$ these modern treatments once again employ the 
model theoretic error of generalization. 
First$,$   
claims are made  that the results 
obtained hold 
not just for the light propagation associated measures of time and distance 
but for all  natural world physical processes (Lawden$,$ 1982$,$ p. 13).  Then 
the domain of application for the  Lorentz 
transformation is further extended to a time continuum with infinitesimals 
(Lawden$,$ 1982$,$ p. 32).  Not only do we 
have the same logical errors$,$ but these modern treatments reject$,$ without 
further thought$,$ the \ae{ther} concept as being a physical entity but use the 
absolute time concept which Einstein claimed also has no physical meaning. 
\par
Relative to the {General Theory}$,$  the same logical error occurs. The square of 
the Minkowski 
space-time interval $\tau$ (Lawden$,$ 1982$,$ p. 14$,$ eqt. (7.4)) 
for  restricted time 
measures is assumed to hold for infinitesimals when $\tau$  is 
 expressed in the famous  differential form 
(Bergman$,$ 1976$,$ p. 44; Lawden$,$ 1982$,$ p. 132). Unless the unfounded extension 
of the restricted time measures 
to a time continuum with infinitesimals is used$,$ the justification for this 
differential form in terms of ``infinitesimal observers'' who can 
synchronize their infinitesimal clocks (Lawden$,$ 1982$,$ p. 132) violates a basic  
requirement of infinitesimal modeling. This tenet states that to pass a time 
related 
property to the classical infinitesimal world it must hold with respect to a 
special approximation process 
for a 
Newtonian time 
interval that is  modeled by a {{\it continuum}} {of real numbers} $\{x\mid a<x <b\}.$ 
The use of any of the highly predictive forms of the classical infinitesimal 
calculus requires this assumption. \parm
\leftline{\bf 5. A Privileged Observer}\parm
Einstein's stated hypotheses$,$ which contradicted all previous modeling assumptions$,$ 
are described as: \par      
{\leftskip 0.5in \rightskip 0.5in \noindent I. The laws of nature are equally 
valid for all inertial frames of reference. II. The velocity of light is 
invariant for all inertial systems$,$ being independent of the velocity of its 
source; more exactly$,$ the measure of this velocity (of light) is a constant.
(Prokhovnik$,$ 1967$,$ pp. 6--7).\par}\pars
These hypotheses$,$ as well as the approach used for the General Theory$,$ 
have significant philosophic implications.  For these hypotheses reject\par
{\leftskip 0.5in \rightskip 0.5in \noindent . . . the Newtonian concept of a 
privileged observer$,$ at rest in absolute space$,$ . . . (Lawden$,$ 1982$,$ p. 
127).\par}\pars
It is interesting how the appliers of the methods of mathematical modeling 
misunderstand even the most basic procedures.  A privileged observer need not 
mean an actual entity within our universe. Whether or not one believes that a 
privileged observer exists in reality is a philosophic question. If a 
privileged observer  
is 
not a real entity$,$ then it may be a member of the class  of {{\it conceptual} 
observers} 
or an NSP-world observer. 
If$,$ as Lawden claims$,$ the privileged observer is to be rejected$,$ then this 
cannot mean the rejection of the privileged conceptual observer or NSP-world 
observer  since
the entire intuitive foundations for the 
Special and General Theories of 
Relativity are based upon Einstein's mind   
experiment and concepts associated with a privileged conceptual observer or 
NSP-world observer.\par 
 The relation of these hypotheses to a Maxwellian type natural \ae{ther} 
is 
that the \ae{ther} is considered a place where  privileged 
observation could take place. But  actual  ``observation'' using 
electromagnetic 
procedures by real human entities within 
our universe seems as if these observations can never be considered as taking
place 
within such an \ae{ther} 
for it appears  
that such entities cannot  determine the basic property of their velocity 
through  it. Hence$,$ how would we know whether or not we are at rest 
with respect to such an \ae{ther}?   
 Further$,$ 
unsatisfactory attempts were made to use perceived natural law to detail 
the behavior of such an \ae{ther}  assuming it is part of the natural world.  
Both of these apparent human inabilities contradict philosophy (S). Thus$,$ to 
remove this contradiction$,$ you simply {postulate away the existence 
of real Maxwellian \ae{ther}} as a foundation. By this process$,$   
\ae{ther} observers are removed if such an \ae{ther} is within the universe. 
Since the theories use infinitesimal 
mathematics$,$ you should not remove$,$ in an  ad hoc manner$,$ NSP-world observation  
for the NSP-world  contains a rigorously defined  
 infinitesimal nonstandard (NS) substratum.   
Without such an NS-substratum$,$ you incorporate into your 
science terms that are claimed to have no content$,$ indeed$,$ no meaning. Hence$,$ 
the view is taken within this new research that privileged observation occurs within the pure NSP-world
and that the fundamental frame of reference is an infinitesimal Cartesian frame 
with an Euclidean styled metric.  
\par
In practice$,$ the idea of absolute length and time is used. Then$,$ it is claimed$,$ that 
such concepts have 
no meaning.  We are told that {\it any} type of measurement
of distance or time$,$ no matter what kind$,$ depends upon  
``relative velocity'' 
without defining such a velocity. We are forbidden from searching for a  
``cause'' for such behavior$,$ a  ``cause'' that may not be fully 
comprehensible. Special and General Theory behavior must simply be accepted 
without further thought.
\par
For nonuniform velocity$,$ the major mathematical structure that was available 
to Einstein and that 
upholds$,$ almost completely$,$ the philosophy of no privileged observer 
within the universe was the pure abstract 
mathematical structure known as {{\it classical Riemannian geometry}} (or the 
absolute or generalized calculus). 
It turns out that Einstein was not a good mathematician but he lived  in 
a region of Germany which was a ``hot bed'' for studies of this structure.
This mathematical structure was the only one available that appeared to 
match his intuition. 
He received considerable help that led to his guess as to a proper expression for a 
new law of gravity. If such a new approach was accepted$,$ it would certainly 
enhance the importance of this mathematician's logical game.\par 
Classical Riemannian geometry is defined entirely by the infinitesimal calculus in terms of a 
required transformation scheme between special types of coordinate systems. 
At least 
locally$,$ each pair of coordinates must satisfy very special  
transformation rules. This means 
that the coordinate systems involved are not all of the possible ones. 
The use of classical tensor analysis actually contradicts the basic philosophy
that  ``Physical space is$,$ then$,$ nothing more than the aggregate of all possible 
coordinate  frames (Lawden$,$ 1982$,$ 127).'' What can only be said is that the 
General Theory of Relativity$,$ in its classical form$,$ applies only to coordinate 
frames that are modeled by this structure$,$ a model that$,$ at that time$,$ 
was thought to correspond correctly to infinitesimal intuition.  
\par
Conceptually$,$ 
the geometry of surfaces requires {``observation'' from a higher dimension}. 
But Einstein's General Theory uses well-known surface concepts where the 
surface has four dimensions. The statement that there can be no privileged 
observer within the universe 
does not preclude  a corresponding type of conceptual observation 
from a required 
``outside'' higher dimension. Although 
higher dimensions can be mathematically introduced  into the theory$,$  
conceptual observation is usually based upon human experience.  
Notwithstanding the balloon analogy$,$ it appears very 
difficult to give a truly meaningful nonmathematical description for the 
appearance of an assumed 4-dimensional ``surface'' within a 5-dimensional 
space.\par   
After the development of our modern approaches to mathematical structures$,$ it 
has become apparent that there is more than one mathematical structure that 
will lead to the same physical consequences. These other structures often 
model a different philosophy of science. For example$,$ a {fractal curve} is 
supposed to be a highly nonsmooth entity.  Yet$,$ it has been shown (Herrmann$,$ 
1989) 
that such fractal effects can be produced by a highly smooth process within the 
NSP-world$,$ an ultrasmooth microeffect. Thus two different 
philosophies  for the physical theory of fractal curves can be
utilized.\par    
We now know that 
the concept of the infinitesimal as viewed  physically prior to 1966 
did not correspond to 
a mathematical structure. Further$,$ in direct opposition to the 
modeling concepts of Newton$,$ the mathematical structure often comes first and 
 {nature 
is required to behave as it dictates}. Once again$,$ we have philosophy (S). \par
The General Theory of Relativity associates pure geometric terms such as 
``intrinsic curvature'' and the associated ``geodesic curvature'' for the 
``force'' concept in Newton's theory and replaces action-at-a-distance with  
propagated gravitational effects. 
But General Relativity 
is a continuum theory. Hence$,$ with a rejection of a space medium for  
radiation or gravitational propagation effects$,$ 
there would be regions surrounding positions within 
our 
universe that over an actual interval of time are totally empty of known or 
imagined entities. 
Yet$,$ the term 
intrinsic curvature would apply to these 4-dimensional regions. Further$,$ 
the rules of the ``game'' again state that 
scientists should not be
allowed to
investigate a more basic ``cause'' for gravitational effects.\par
If other well established mathematical models predict the exact same 
consequences as 
the General Relativity$,$ then substituting the 
nonintuitive concept of the intrinsic curvature of space-time for the 
experiential concept of force would 
be unnecessary.   
Is the concept of action-at-a-distance any less comprehensible to the human mind 
than an intrinsic 
curvature of  ``empty'' 4-dimensional space-time? 
 \parm
\leftline{\bf 6. The Fock Criticism and Significant Other Matters}\parm  
In the 1930s$,$ a major technical criticism of some of the Einsteinian arguments 
was brought forth by the Russian cosmologist {Fock} (1959). His criticism is 
related to the {{\it Equivalence Principle}}  as used by Einstein. Consider$,$ first$,$ the 
global or overall space-time physical law of the equivalence of inertial mass
(i.e. the ``m'' in the Newtonian expression $\vec F = m\vec a$) and 
gravitational mass (i.e. in Newtonian gravity theory$,$ for example$,$ 
the mass stated in this law). Then consider the so-called Equivalence Principle 
states that an acceleration field (i.e. the $\vec a$) and a gravitational 
field cannot be distinguished one from another. There have been many arguments 
presented that such an Equivalence Principle relative to the field concept is false from the global viewpoint. \parm
{\leftskip 0.5in \rightskip 0.5in \noindent  As was mentioned$,$ Einstein 
considered that from the point of view of the Principle of Equivalence it 
is impossible to speak of absolute acceleration just as it is impossible to 
speak of absolute velocity. We consider this conclusion of Einstein's to be 
erroneous. . . . This conclusion is based upon the notion that fields of 
acceleration are indistinguishable from fields of gravitation. But$,$ although
the effects of acceleration and of gravitation may be indistinguishable 
``in the small''$,$ i.e. locally$,$ they are undoubtedly distinguishable ``in the 
large''$,$ i. e. when the  boundary conditions to be imposed on gravitational 
field are taken into account (Fock$,$ 1959$,$ p. 208)\par}\parm
\noindent Fock gives an example of this for a rotating system and many others 
have given examples for accelerating noninfinitesimal structures. Then Fock 
writes:\parm
{\leftskip 0.5in \rightskip 0.5in \noindent In the first place there is here
an incorrect initial assumption. Einstein speaks of arbitrary gravitational 
fields extending as far as one pleases and not limited by boundary conditions. 
Such fields {\it cannot exist.} Boundary conditions or similar conditions which 
characterize space as a whole are absolutely essential and thus the notion 
of ``acceleration relative to space'' retains its significance in some form 
or another.
. . . The essence of the error committed is in the initial assumption 
consists in {\it forgetting that the nature of the equivalence of fields of 
acceleration and of gravitation is strictly local} (1959$,$ p. 369).\par}\parm
The fact that the {effects must be local} is why the infinitesimal calculus$,$ in 
generalized form$,$ is used as a means to model mathematically a general law of 
gravity. Further$,$ Fock criticizes the concept of nonuniqueness relative to 
the general notion that any frame of reference will suffice. As I have pointed out 
the term {\it any} is not correct when the generalized calculus is used since 
the frame of reference$,$ as modeled by a coordinate system$,$ must have certain 
differential properties and be very special locally. Fock claims that the term 
 {\it any} must be restricted greatly. This restriction is somewhat technical 
in character and refers to a claimed wave-like quality that certain solutions to 
the Einstein gravitational field equations seem to require. 
Fock claims that 
the correct coordinate system in which to discuss solutions is an 
 {``harmonic''} system (Fock$,$ 1959$,$ p. 346-352). Further$,$ it is claimed that 
solutions must have an additional boundary condition that at ``infinity'' they 
become the infinitesimal ``chronotopic'' line-element. [Note: this will be 
derived in article 3.] This is also called the {Minkowski-type line-element} 
(i.e. metric)$,$ but Fock  calls it the {Galilean line-element} (metric). 
Sometimes it is termed the Euclidean requirement. Whatever terminology is 
employed$,$ the 
Fock idea is that there are preferred coordinate systems in which measurements 
make sense. \parm
{\leftskip 0.5in \rightskip 0.5in \noindent In the question of an isolated 
system of bodies the question of a coordinate system is answered in the same 
way as in the absence of a gravitational field: there exists a preferred 
system of coordinates (Galilean or harmonic) but it is also possible to use 
any other coordinate system. However$,$ the geometric significance of the latter
can only be established by comparing it with the preferred system (1959$,$ p. 
376).\par}\parm\smallskip
I agree with Fock that application of the Equivalence Principle is only local. 
The necessity for harmonic coordinates has not been established on 
mathematical grounds except that it leads to unique and testable conclusions 
 for  our   universe.   But  Fock  did not   have   our
 present day knowledge for the 
correct rules for infinitesimal modeling and$,$ hence$,$ missed the point entirely. 
The Einstein gravitational equations as they are expressed in terms of the 
general or absolute calculus can be analogue modeled in many ways. {\it One} 
language that can be used is that of Riemannian geometry. But$,$ {Riemannian 
geometry is just that$,$ 
an analogue model}; something that simply represents behavior of something 
else but is not itself reality. This may lead to a certain mental 
visualization that aid in producing conclusions$,$ but it also appears that 
visualizing the Einstein equations in terms of a generalization of the four 
dimensional wave equation also aids in such comprehension. This is precisely
the reason that harmonic coordinates are introduced. Both of these realization 
techniques still remain analogue models for a reality that was not expressible 
in a correct language until now. \par 
Recently$,$ I have discovered that {Fokker} (1955) did guess at the correct 
concepts. Fokker suggested that the name {{\it chronogeometry}} would be more 
appropriate. As it will be established$,$ Fokker was correct with this suggestion.
When one realizes that mathematical coordinate systems are but an abstract 
entity without relation to reality$,$ then one is lead to a theory of 
measurement that would correspond to a particular coordinate system. As will 
be shown in the articles that follow this one$,$ the concept of the privileged 
observer comes prior to the selection of a coordinate system. After such an 
observer is defined$,$ then a definable measuring process may be correlated to 
a compatible coordinate system. The 
privileged observer will essentially be observing infinitesimal light-clocks (defined in article 2) and 
nothing more.\par
For 
comparison purposes$,$ the location$,$ so to speak$,$ of the fundamental observer can be considered as  
fixed in the NSPPM.  The NS-substratum
coordinate system is a four dimensional 
Cartesian (Euclidean) system with the infinitesimal light-clocks being 
oriented along coordinate lines in order to measure the dynamic properties 
of motion. When certain motion occurs 
within the natural world$,$ the infinitesimal light-clocks undergo alterations.
 This 
alteration leads to a physical$,$ not geometric$,$ statement that can  be
considered 
an invariant under certain physical changes that are cause by motion of the 
infinitesimal light-clocks. This statement states that due to certain 
properties of electromagnetic radiation within the natural world$,$ there 
will be a specific relation between coordinate (infinitesimal) light-clocks
that is due to motion with respect to the NSPPM.\par

 Prior to 
any physical alteration in the infinitesimal light-clock counting mechanism$,$ 
an {acceptable coordinate change} (i.e. continuously differentiable with 
nonvanishing Jacobian) has the basic purpose of simply changing the 
``orientation$,$'' 
so to speak$,$ of the infinitesimal light-clocks in order to give a different 
measure of the physical dynamics. On the other hand$,$  
using the General Relativity assumptions associated with 
Riemannian geometry and {\it once a solution is obtained for a particular physical 
scenario}$,$ the {coordinate change 
is interpreted as an 
acceptable alteration in the gravitational field}. It is then claimed that such 
an alteration in the gravitational field will affect the abstract notion of ``time.'' Such an interpretation is in error logically. Physical measures 
altered by a gravitational field are modeled by alterations in 
infinitesimal light-clock measures - a model that yields alterations in physical behavior.\par

The special physical ``invariant'' $dS^2$ is used. The {invariant statement}
says that IF 
one imposes a coordinate change as a means of measuring different dynamic 
properties$,$ then the expression $dS^2$ (i.e. the infinitesimal 
``chronotopic interval'') remains 
fixed. \par
Fock and all previous researchers start their investigation with this 
so-called invariant 
statement$,$ but$,$ of course$,$ never relate the statement to the actual entities 
that are 
being altered by motion$,$ the iinfinitesimal (inf.) light-clocks.  The chronotopic interval can be 
generalized and represents infinitesimal light-clock behavior without the presence of a 
gravitational field of one sort or another as the agent for the motion in 
question.  The major aspect of this 
interval is 
``uniform'' velocity and, when infinitesimalized, nonuniform velocity.
When generalized$,$ the physical invariant (i.e. the fundamental metric of 
space-time$,$ which I term a  ``line-element'' since tensors are not used) looks like $dS^2=g_{11}dx_1dx_1 
+ g_{12} dx_1dx_2 + \cdots + g_{44}dx_4dx_4$ and for such things as 
gravitations fields one can consider various types of ``potential'' 
velocities. However$,$ the various $dx_i$ cannot be uncontrolled infinitesimals 
in that they require an interpretation. This expression comes from classical 
Riemannian geometry and in that discipline each $dx_i$ takes on different 
interpretations. In one case$,$ each $x_i$ is a function of another parameter 
(usually  ``time'') and the expression is used to measure of arc-length along 
a ``curve.'' In the general case$,$ each $dx_i$ is supposed to be interpreted as 
the geometric concept of ``distance between infinitely near points.'' Prior to 
1961 and Robinson's discovery$,$ these interpretations did not follow a 
mathematically rigorous theory. Prior to the late 1980s$,$ the actual method of 
obtaining an interpretation$,$ the ``infinitesimalizing method'' also had not 
been discovered. The discovery of these correct methods shows that the actual 
infinitesimal chronotopic interval should not be interpreted in geometric 
terms$,$ 
but 
it should be interpreted in terms of measures that retain an electromagnetic 
propagation language$,$ and motion or potential motion. Of course$,$ such 
differential models are, usually, intended to give approximations for macroscopic or 
large scale behavior. A correct method to incorporate both of these needed 
requirements will appear later.\par 

For the basic infinitesimal 
chronotopic 
interval$,$ the coefficients are as follows: the variable $dx_i,\ i = 1,2,3,4$ 
are related to types of infinitesimal light-clock measurements.  
The coefficients are $g_{11}=g_{22} = g_{33} = -1,\ g_{44} = 1/c^2,$ where $c$ 
is the measured velocity of light {\it in vacuo.} All other coefficients are 
zero. Fock states$,$ somewhat 
incorrectly$,$ that to  ``understand'' the geometry of space-time one needs to 
compare such changes in space-time geometry with this basic interval. 
In the presence of a gravitational field$,$ the coefficients $g_{ij}$ 
are related to the field's potential. This potential is also related to a
force effect produced by the field$,$ as previously mentioned$,$ and$,$ due to prior 
use of these concepts$,$ leads to a geometric analogue model for what is$,$ in 
reality$,$ gravitational alterations in the behavior of electromagnetic 
radiation. If only the geometric method is used$,$ then for a particular 
physical scenario$,$ the terms of this line-element are appropriately altered
so that the $dx_idx_j$ represent specific geometric coordinate concepts and 
the  corresponding functions $g_{ij}$ satisfy the Einstein gravitational 
equation and also compensate for the $dx_idx_j$ in such a manner that the 
line-element is invariant. In this manner and in terms of a geometric 
language$,$ 
a comparison can be made to the line-element as it would be without the 
gravitational field. It will be argued in article 3$,$ that $dS^2$ is actually 
a physical invariant only due to a special property associated with 
infinitesimal modeling. \par
The geometric approach actually contradicts the creator of the mathematics 
employed. The universe is NOT controlled by geometry. {Geometry is a human 
construct} that is used to model other intuitive concepts. Newton tells us 
that the intuitive concepts of mechanics come first. After experiencing  
behavior of forces$,$ velocities$,$ accelerations$,$ resistance to motion$,$ and the 
hundreds of other purely physical concepts$,$ then the mathematical model is
devised that will aid in logical argument and prediction. The mathematical 
model itself is not reality. The same holds for a Riemannian generalization
of the geometry. Physical intuition should come first. Then certain rules 
for infinitesimal modeling are applied. The resulting 
constructions should aid in comprehending and predicting physical behavior 
as it is compared to the  ``simpler'' and original physical intuition.
The rules for such modeling are completely controlled by what is perceived 
of as {simple behavior  ``in the small''} that leads to complex behavior 
``in the large.''  We  ``understand'' the more complex behavior by comparing 
it to the simpler behavior. Now Riemannian geometry can still be used$,$ if it 
is properly interpreted in terms of the actual physical entities. But$,$ if the 
foundation for your ``gravitational forces'' is stated in terms of a 
geometric language rather than in terms of intuitive physical qualities$,$ then 
you would lose the great power of the infinitesimal calculus as a predictor of 
intuitive physical behavior. \par
For physical reasons$,$ {Patton and Wheeler} (1975) also reject geometry as the ultimate 
foundations for our physical universe. Wheeler coined a new term 
 {``pregeometry''} for the actual foundations.\parm
{\leftskip 0.5in \rightskip 0.5in \noindent Riemannian geometry likewise 
provides a beautiful vision of reality; but it will be useful as anything we 
can do to see in what ways geometry is inadequate to serve as primordial 
building material. . . .  ``geometry'' is as far from giving an understanding 
of space as ``elasticity'' is from giving an understanding of a solid. . . .
(1975$,$ p. 544$,$ p. 557-558 \par}\parm
In what follows$,$ from the viewpoint of infinitesimal light-clocks$,$ acceptable 
coordinate changes are first changes in the orientation of the measuring light 
clocks. After this is done$,$ certain physical processes associated with the 
intuitive idea of a physical {``potential velocity''}  
are postulated. These are modeled by a simple linear correspondence. 
This leads to a general 
 invariant line-element. As will be demonstrated in article 3$,$ in order for 
this line-element 
to remain invariant$,$ the differently oriented light-clocks would need to 
have one and only one aspect altered when substituted into this line-element. 
If one then associates the ``potential velocities'' with those that would be 
produced by a specific gravitational field$,$ then solutions are obtained that 
are the same
as those obtain from the Einstein equations and the Riemannian 
geometry approach$,$ except that the solutions are stated in terms of actual 
physical processes that are being altered by these potentials. 
Other potential velocity substitutions lead to
vreified predictions for the Special Theory. The predictions obtained are in terms of 
infinitesimal light-clocks and how the altered light-clock behavior   
compares to that of the unaltered light-clock behavior.  
\parm    
\leftline{\bf 7. Why Different Derivations?}\parm
There is a definite need for different derivations 
for theories that predict local events since all such theories are based upon 
philosophic foundations that impinge upon personal belief-systems. When 
the basic hypotheses for theory construction are identified$,$ then these 
hypotheses 
always have a broader meaning called their {{\it descriptive content.}} 
The content of a 
collection of written statements$,$ diagrams and other symbolic forms 
is defined as all of the mental impressions 
that the collection evokes within the mind of the reader. These impressions 
are$,$ at the least$,$ based upon an individual's experiences. \par
A personal {belief-system} also has content and this content can be contradicted 
by the content generated by the hypotheses or predictions of a scientific 
theory. Suppose that 
you have theories $T_1$ and $T_2$ each based upon different foundational 
concepts but the verified logical consequences of these two theories are the 
same. Further$,$ the content of the foundations for $T_1$ does not contradict 
the content of your belief-system$,$ while the content of the foundations for 
$T_2$ does. You now have
 two choices. You can accept theory $T_1$ as 
reasonable and reject $T_2;$ or you can change your personal belief-system$,$  
accept theory $T_2$ and reject $T_1.$  Whichever theory you select can be 
analyzed relative to its humanly comprehensible technical merits. Such an 
analysis does not necessarily imply that the selected theory is the correct 
theory. \par
 There is absolute evidence that much of 
what passes for scientific theory is designed to force a rejection of various 
philosophic concepts$,$ a rejection of well-founded belief-systems. For this reason$,$ 
if for none other$,$ it is important to identify the philosophic foundations 
of all scientific theories and to allow individuals a free-choice as to 
which they wish to accept.\par 

With respect to both the Special and General Theories of Relativity$,$ it is 
now possible to use the correct infinitesimal concepts$,$ a NSPPM that does not reveal all of its properties,
and obtain verified consequences of both theories without 
logical error. The error is eliminated by predicting alterations in clock behavior rather than by the error of inappropriate generalization. This new 
mathematical model alters the basic philosophic 
assumption of no privilege observer.\parm 

\leftline{\bf 8. A Corrected Derivation} \parm
Notwithstanding Einstein's inability$,$ using the philosophy (S)$,$ to describe 
a natural world \ae{ther} that will lead to verified 
Special Theory 
effects$,$ there does exist  
a description that includes the necessary infinitesimals.  
 A new derivation for the Lorentz
transformation based upon absolute time with infinitesimals 
has been published
(Herrmann$,$ 1993). 
  In this derivation$,$ there is a {{\it privileged 
observer using an inertial Cartesian coordinate system}} as well as  
an additional$,$ almost trivial$,$ Galilean 
infinitesimal effect based upon natural world laboratory observations.    
The coordinate system lies in a 
portion of an NSP-world the  
 {{\it nonstandard photon-particle medium}} -- the NSPPM - which is contained in an NS-substratum. This ``medium'' is a portion of the subparticle field and it yields relativistic alterations in natural-system behavior through an interaction with natural world entities associate with simple electromagnet propagation properties. 
The term {{\it inertial}} refers only to the weakest aspect of 
Galilean-Newtonian mechanics where$,$ with respect to this NSPPM$,$ a state of 
rest or uniform motion can be altered only by {(force-like) interactions.} \par
Few things can be known about the NSPPM. What is known is that certain 
basic expressions for 
Newtonian mechanics must be altered and even the general descriptions for such 
laws are slightly different. This NSPPM can be considered as part of the 
natural world if one wishes but it would be a very distinct part. 
The basic derivation is obtained using an 
 {absolute Newtonian time concept} within NSPPM and a {Galilean 
photon propagation theory} that includes an infinitesimal statement 
which assumes that the velocity of light IS dependent 
upon an additional NSP-world velocity which can be the source of the photon. The apparent inability of 
measuring the velocity of the Earth through this medium using certain 
electromagnetic 
techniques and$,$ hence$,$ only being able to determine by such natural world 
techniques relative velocities is incorporated into this basic derivation.
Further$,$ the {constancy of the to-and-fro natural world measurements  
of the velocity of light} is included where$,$ however$,$ this velocity 
need not be  constant with respect to the absolute time. Simply 
stated$,$ this new derivation adjoins to certain electromagnetic behavior a 
simple additional infinitesimal property related to observed electromagnetic 
propagation. \par
The basic properties for the NSPPM time as measured by a to-and-fro
electromagnetic propagation experiment$,$ once obtained from this derivation$,$  
are now applied to the  natural world.  
In order to retain the electromagnetic character of the Special Theory 
effects$,$ it is necessary that the transformation be derived via hyperbolic geometry.   
This derivation 
method retains in its time and distance measures  the electromagnetic 
propagation language. Such  measures are termed as {{\it Einstein measures}}
 and are 
denoted by the subscript $E$. The 
Prokhovnik (1967) interpretation of the results relative to the Hubble textural 
expansion is$,$ however$,$ totally rejected since the $\omega$ is not related in 
any manner to such an assumed expansion. I point out that the use of the NSPPM
eliminates not only the incomprehensible physical world length contraction and 
time dilation relations but even the difficulties associated with the 
reciprocals of 
these relations (i.e. the twin paradox).\par
A new refinement of the concept of Einstein measure is mentioned in this soon 
to be published 
1993 paper.  This refinement retains completely the electromagnetic 
propagation 
language by introducing the analogue model of the light-clock. 
The analogue  light-clock model
is composed of a fix rigid  arm, of various lengths, that has a beginning light source attached 
to one end  
with a counter and returning mirror$,$ and simply a returning mirror at the 
other end. The  ``counter'' counts the number of to-and-fro paths an 
electromagnetic pulse ``makes'' from some fixed beginning count setting. Notice that from 
field properties an electromagnetic pulse's speed (at least one photon) is measured as constant. This concept is passed to the infinitesimal 
world where the arm's ``length'' can be considered a positive infinitesimal. What this yields 
is the {{\it  infinitesimal light-clock.}} This does not preclude the possibility that under various conditions the speed of photons in the infinitesimal light-clocks is altered.
\par
In the natural world$,$ the light-clock concept can only 
 {approximate the continuum} associated Einstein measures where 
the approximation is improved as the 
length of the arm is reduced. As will be discussed in article 2$,$ 
this approximation may be made exact within the infinitesimal world. 
Identical infinitesimal light-clocks are used to 
measure the 
``time'' by corresponding this concept to the counter number$,$ the counter 
``ticks.''  Twice the arm's length multiplied by the counter number gives a 
measurement of an apparent   
 distance$,$ in terms of linear units$,$  the electromagnetic pulse has traversed for 
a specific count number. Tracing the path of electromagnetic radiation  
leads to the basic interpretation that these light-clock counts can be used to
measure$,$ within the NSP-world$,$ the to-and-fro electromagnetic path length 
within one moving light-clock.   
It is a measure in 
terms of linear units for nonlinear behavior. \par 
Technically$,$ in the NSP-world$,$ twice the ruler-like 
measurement$,$ in terms of private units$,$ of the arm's length can be considered a 
positive infinitesimal 
$L.$  It is known (Herrmann$,$ 1991b$,$ p. 108)  that for every positive real 
number $r$ 
there exists an ``infinite'' 
count number $\Pi,$ where $\Pi$ is a Robinson infinite number$,$ such that $r$ 
is {{\it infinitesimally (or infinitely) close to}} $L\Pi.$ Notationally this is 
written as $r \approx L\Pi,$ where  $\approx$ is at the least an equivalence 
relation. What this means 
is that $r - L\Pi$ is an infinitesimal. There also exist infinite numbers
$\Lambda$ such that $L\Lambda$ is infinitesimally near to zero.  \par
 If $L\Pi$ is known to be 
infinitesimally near to real number$,$ then $L\Pi$ is said to be {{\it finite}}. 
There is a process that can be used to 
capture the real number $r$ when $L\Pi$ is known to be finite. This process 
uses the {{\it standard part operator}} that is denoted by ``st.'' 
Many  properties 
of the operator ``st'' are obtained from results in abstract algebra and 
these properties include the same formal properties as the ``limit'' viewed as an 
operator. The infinitesimal light-clocks yield an exact analogue for Einstein 
measures when the standard part operator is applied.\par  
The discussion in 
Herrmann (1992) shows how the use of infinitesimal light-clocks
allows for a {return to the concept of absolute Newtonian time.}   
Natural world observations lead to infinitesimal properties for the NSPPM 
expressed in terms of  absolute time. 
Then$,$ counter to the Einstein 
claim$,$ these light-clock approximations are used to define infinitesimal 
light-clocks$,$ which in turn lead to unique NSPPM times.  
Special light propagation properties lead to an  Einstein time definition. 
But Einstein time 
can also be  successively approximated within  the natural world in terms of 
light-clock measures and only such measures. In terms of 
infinitesimal light-clocks$,$ Einstein time can be exactly obtained. 
This discussion  also shows that 
known  Special Theory effects associated with uniform relative velocity 
(i.e. not incorporating possible gravitational effects) 
can be interpreted as manifestations of the electromagnetic character 
of natural world 
entities and how they interact with the NSPPM. What has not been 
investigated$,$
is what relation gravity has with respect to the NSPPM  
and 
how such a gravity relation might influence the 
physical behavior within our natural world.\par 
In Article 3$,$ based upon a privileged  
observer 
located within the NSPPM$,$ the infinitesimal chronotopic line-element 
is derived from light-clock properties and shown to 
be related to the propagation of electromagnetic radiation.  A general 
expression is derived$,$
without the tensor calculus$,$ 
from basic infinitesimal theory 
applied to obvious Galilean measures for  distances traversed by an
electromagnetic pulse. Various line-elements are obtained from this general 
expression. These include the Schwarzschild (and modified) line-element$,$ which 
is obtained by merely substituting a Newtonian  gravitational velocity into 
this expression; the de Sitter and the Robertson-Walker which are obtained by 
substituting a velocity 
associated with the cosmological constant or an expansion (contraction)
process. The relativistic (i.e. 
transverse Doppler)$,$ gravitational and cosmological
redshifts$,$ and alterations of the radioactive decay rate are derived 
from a general behavioral model associated with atomic systems$,$ 
 and it is predicted that similar types of shifts will take 
place for other specific cases. Further$,$ the mass alteration expression is 
derived in a similar manner. 
  From these derivations$,$ locally verified 
predictions of the Einstein Special and General Theories of Relativity can 
be obtained. A process is also given that minimizes the problem of the
``infinities'' associated with such concepts as the Schwarzschild radius.
These ideas are applied to the formation of black holes and pseudo-white holes.\parm

\centerline{\bf References for Article 1.}\parm
\id{B}arnes$,$ T. G. 1986. Space medium theory. Geo/Space Research Foundation.
\id{B}enton$,$ D. J. 1988. The Special Theory of Relativity: its 
assumptions and implications. {\it C. R. S. Quarterly} 
25:88--90. 
\id{B}ergmann$,$ P. 1976. Introduction to the theory of relativity. Dover$,$ 
New York. 
\id{C}auchy$,$ A. 1821. Cours d'Analyse (Alalyse Alb{\accent 19 
e}brique).
\id{E}instein$,$ A. 1905. Zur Elektrodynamik bewegter K{\accent "7F 
o}rper. Annalen der Physik$,$
 17:891--921. 
\id{}{\vbox{\hrule width 0.5in}} 1907. \"Uber das Relativit\"atsprinzip und die aus demselben gezogenen Folgerungen, {\it Jahrbuch der Radioaktivit\"at und Electronik}, 4:411-462. [In this paper, Einstein uses the Minkowski method to obtain the transformation and looses the special character, even required for this approach, for the measurements.]  
\id{}{\vbox{\hrule width 0.5in}} 1945. The Meaning of Relativity. Princeton University 
Press$,$ Princeton$,$ NJ. 
\id{}{\vbox{\hrule width 0.5in}}  1979. Albert Einstein autobiographical notes. Transl. 
and ed. by Paul
Schilpp$,$ Open Court$,$ La Salle$,$ IL.
\id{}{\vbox{\hrule width 0.5in}}  1982. How I created the theory of relativity. Transl. 
Yoshimasa 
A. Ono. {\it Physics Today} 35:45--47. 
\id{E}instein$,$ A. and L. Infeld. 1938. The Evolution of Physics. Simon 
and Schuster$,$ New York. 
\id{F}ock$,$ V. 1959. The Theory of Space Time and Gravity. Pergamon Press. New 
York.
\id{F}okker$,$ A. D. 1955. Albert Einstein$,$ inventor of chronogeometry. {\it 
Synth\'e se}$,$ 9:442-444.
\id{H}errmann$,$ R. A. 1986. D-world evidence. {\it C. R. S.
Quarterly}$,$ 23:47-53.
\id{}{\vbox{\hrule width 0.75in}} 1989. Fractals and ultrasmooth microeffects. {\it 
J. Math. Physics}. 30:805--808. 
\id{}{\vbox{\hrule width 0.75in}} 1991a. Some applications of nonstandard analysis to 
undergraduate mathematics: infinitesimal modeling and elementary physics. {\it 
Instructional Development Project$,$ Mathematics Department$,$ U. S. Naval 
Academy$,$ Annapolis$,$ MD 21402-5002.}\hfil\break http://arxiv.org/abs/math/0312432 
\id{}{\vbox{\hrule width 0.75in}} 1991b. The Theory of Ultralogics. \hfil\break
http://www.arxiv.org/abs/math.GM/9903081\hfil\break
http://www.arxiv.org/abs/math.GM/9903082\par
 \id{}{\vbox{\hrule width 0.75in}} 1992. A corrected derivation for the Special 
Theory of 
Relativity. Presented before the Mathematical Association of America$,$ Nov. 14.
at Coppin State College$,$ Baltimore$,$ MD.
Preprints from IMP$,$ P. O. Box 3268$,$ Annapolis$,$ MD 21403.                         
\id{}{\vbox{\hrule width 0.75in}} 1994$,$ The Special Theory and a nonstandard substratum. 
{\it Speculations in Science and Technology.} 17(1):2-10. 
\id{L}awden$,$ D. F. 1982. An introduction to tensor calculus$,$ relativity 
and cosmology. John Wiley \& Sons$,$ New York. 
\id{M}axwell$,$ J. C. 1965a. The scientific papers of James Clark 
Maxwell. Ed. W. D. Niven$,$ Vol. 1$,$ Dover$,$ New York. 
\id{}{\vbox{\hrule width 0.6in}} 1965b. The scientific papers of James Clark 
Maxwell. Ed. W. D. Niven$,$ Vol. 2$,$ Dover$,$ New York. 
\id{N}ewton$,$ I. 1934. Mathematical principles of natural philosophy. 
Transl. by Cajori$,$ University of Cal. Press$,$ Berkeley. 
\id{P}atton$,$ C. M. and J. A. Wheeler. 1975. Is Physics Legislated by a 
Cosmogony? In {\it Quantum Gravity.} ed. Isham$,$ Penrose$,$ Sciama. Clarendon 
Press$,$ Oxford. pp. 538--605. 
\id{P}lanck$,$ M. 1932. The mechanics of deformable bodies. Vol II$,$ 
Introduction to theoretical physics. Macmillan$,$ New York. 
\id{P}rokhovnik$,$ S. J. 1967. The logic of Special Relativity. Cambridge 
University Press$,$ Cambridge.
\id{R}obinson$,$ A. 1966. Non-standard analysis. North-Holland$,$ Amsterdam.\par\bigskip

\vfil\eject
\centerline{{\bf A Corrected 
Derivation for the Special Theory of Relativity, Article 2}\footnote*{This is an 
expanded version of the paper presented before 
the Mathematical Association of America$,$ Nov. 14$,$ 1992$,$ Coppin State College$,$ 
Baltimore Md}} \medskip 
\centerline{Robert A. Herrmann} 
\medskip 
Abstract: Using properties of the 
nonstandard physical world$,$ a new fundamental derivation for the 
effects of the Special Theory of Relativity is given. This fundamental 
derivation removes all the contradictions and logical errors in the original 
derivation and leads to the fundamental expressions for the Special Theory Lorentz transformation. Necessarily, these are obtained by means of hyperbolic geometry.  
It is shown that the Special Theory effects are manifestations of the 
interaction between our natural world and a nonstandard medium, the NSPPM.  This 
derivation  eliminates the controversy associated with any physically  
unexplained absolute time dilation and length contraction. It is shown that 
there is no such thing as a absolute time dilation and length contraction 
but$,$ rather$,$ alterations in pure numerical quantities associated with an 
electromagnetic-type interaction with an NSP-world NSPPM. 
\medskip
\leftline{\bf 1. The Fundamental Postulates.} 
\medskip
There are various Principles of Relativity. The most general and least 
justified is the one as stated by Dingle ``{\it There is no meaning in 
absolute motion.} By saying that such motion  has {\it no meaning}$,$
 we assert that 
there is no observable effect by which we can determine whether an object is 
absolutely at rest or in motion$,$ or whether it is moving with one velocity or 
another.''[1:1] Then we have Einstein's statements that  ``I. The laws of 
motion are equally valid for all inertial frames of reference. II. The velocity 
of light is invariant for all inertial systems$,$ being independent of the 
velocity of its source; more exactly$,$ the measure of this velocity (of light) 
is constant$,$ $c$$,$ for all observers.''[7:6--7] I point out that Einstein's 
original derivation in his 1905 paper ({\it Ann. der Phys.} {\bf 17}: 
891) uses certain well-known processes related to partial differential 
calculus.\pars 
In 1981 [5] and 1991 [2]$,$  
it was discovered that 
the intuitive concepts associated with the  Newtonian laws of motion were 
inconsistent with respect to the mathematical theory of infinitesimals when 
applied to a theory for light propagation. The apparent nonballistic nature for light propagation when transferred to infinitesimal world would 
also yield a nonballistic  behavior. 
Consequently$,$ {\bf there is an absolute contradiction between Einstein's 
postulate II and the derivation employed.} This contradiction would not have 
occurred if it had not been assumed that the \ae{ther} followed the 
principles of Newtonian physics with respect to electromagnetic propagation. 
[Note: On Nov. 14$,$ 1992$,$ when the information in this article was formally presented$,$ I listed various predicates that Einstein used and showed the specific places within the derivations where the predicate's domain was altered without any additional argument. Thus$,$ I gave specific examples of  
the model theoretic error of 
generalization. See page 49.] 
\pars
 I mention that Lorentz speculated that \ae{ther} 
theory need not correspond directly to the mathematical structure but could 
not show what the correct correspondence would be. Indeed$,$ if one assumes that the 
NS-substratum satisfies the most basic concept associated with an 
inertial system that a body can be considered in a state of rest or 
uniform motion 
unless acted upon by a force$,$ then the expression $F=ma,$ among others$,$ 
may be altered for infinitesimal NS-substratum behavior. Further$,$ the 
NSPPM portion of the NS-substratum$,$ when light propagation is discussed$,$ does not follow the
Galilean rules for velocity composition. The additive rules are followed but 
no negative real velocities exterior to 
the Euclidean monads are used since we are only interested in the propagation properties for electromagnetic radiation. 
The derivation in 
section 3 removes all contradictions by applying the most simplistic Galilean 
properties of motion$,$ including the ballistic property$,$ but only to behavior 
within a Euclidean monad. \pars
As discussed in section 3$,$ the use of an NSP-world (i.e. {\it nonstandard physical world}) 
NSPPM allows for the  elimination of the well-known 
Special Theory ``interpretation'' contradiction that the 
mathematical model 
uses the concepts
of Newtonian absolute time and space$,$ and$,$
yet$,$ one of the major interpretations is 
that there is no such thing as absolute time or absolute space.\pars 
Certain  general principles for NS-substratum light propagation will be  
specifically stated  
in section 3.  These principles  
can be gathered together as follows: (1) {\sl There is a porton of the 
nonstandard photon-particle medium - the NSPPM - that sustains 
N-world  (i.e. natural = physical world) photon propagation. Such propagation 
follows the infinitesimally presented laws of Galilean dynamics$,$ when 
restricted to monadic clusters$,$  
and the 
monadic clusters follow an additive and an actual metric property 
for linear relative motion when considered collectively.} [The term 
``nonstandard electromagnetic field'' should only be construed as a NSPPM notion, where the propagation of electromagnetic radiation follows 
slightly different principles than within the natural world.] 
(2) {\sl The motion of light-clocks
 within the N-world (natural world) is associated with one single 
effect. This effect 
is an alteration in an appropriate light-clock mechanism.}  [The 
light-clock concept will be explicitly defined at the end of section 3.] 
It will be shown 
later that an actual physical cause may be associated with Special Theory 
verified physical alterations. 
 {Thus the 
Principle 
of Relativity$,$ in its general form$,$ and the inconsistent
portions of the Einstein principles are
eliminated from consideration and$,$ as will be shown$,$ 
the existence of a special type of medium, the NSPPM, can be assumed without 
contradicting  experimental evidence.} \pars
In modern Special Theory interpretations [6]$,$ it is  claimed 
that the effect of ``length contraction''  has no physical meaning$,$ whereas 
time dilation does. 
This is probably true if$,$ indeed$,$ the Special Theory is actually based upon 
the intrinsic N-world  concepts of length and time. What follows will further
demonstrate that 
 the  Special Theory is  a light propagation theory$,$ 
as has been previously 
argued by others$,$ and that the so-called ``length contraction'' and time 
dilation can both be 
interpreted as  physically real effects when they are described in terms of the NSPPM. The effects are only relative to a theory of light 
propagation.\pars   
\medskip
\leftline{\bf2. Pre-derivation Comments.} 
\medskip
Recently [2]--[4]$,$ nonstandard analysis [8] has proved to be 
a very significant tool in investigating the mathematical foundations for various 
physical theories. In 1988 [4]$,$ we discussed how the methods of 
nonstandard analysis$,$ when applied to the symbols that appear in statements 
from a physical theory$,$ lead formally to a pregeometry and the entities termed 
as subparticles. One of the goals of  
 NSP-world research is the 
re-examination of the foundations for various controversial 
N-world theories and the eventual elimination of such 
controversies by viewing such theories as but restrictions of more simplistic 
NSP-world concepts. This also leads to indirect evidence for the actual 
existence of the NSP-world. \par
The Special Theory of Relativity still remains a very controversial theory due 
to its philosophical implications. Prokhovnik [7] produced a derivation that 
yields all of the appropriate transformation formulas based upon a light 
propagation theory$,$ but unnecessarily includes an interpretation of 
the so-called Hubble textural expansion  of our 
universe as an additional ingredient. 
The new derivation we give in this article shows that  
properties of a NSPPM also lead to 
Prokhovnik's expression  
(6.3.2) in reference [7] and from which all of the appropriate equations 
can 
be derived. However$,$ rather than considering the Hubble expansion as directly 
related to Special Relativity$,$ it is shown that one only needs to consider 
simplistic 
NSP-world behavior for light propagation and the measurement of 
time by means of N-world light-clocks. This leads to the 
conclusion that Special Theory effects may be produced by a dense NSPPM 
within the  NSP-world. Such a NSPPM -- an \ae{ther} -- 
yields 
N-world 
Special Theory effects. \pars   
\medskip
\leftline{\bf 3. The derivation}
\medskip
The major natural system in which we exist locally is a space-time 
system. {``Empty'' space-time} has only a few characterizations when viewed 
from an Euclidean perspective. We investigate$,$ from the NSP-world viewpoint$,$ 
 {electromagnetic propagation through a Euclidean neighborhood} of space-time. 
Further$,$ we assume that light is such a propagation. One of the basic precepts 
of infinitesimal modeling is the experimentally verified {{\it simplicity}} 
for such a local system. For actual time intervals$,$ certain physical processes 
take on simplistic descriptions. These NSP-world descriptions are represented 
by the exact same description restricted to infinitesimal intervals. Let 
$[a,b],\ a \not= b,\ a > 0,$ be an objectively real conceptional time interval and let $t 
\in (a,b).$\par
The term  ``time'' as used above is very misunderstood. There are 
various viewpoints relative to its use within mathematics. Often$,$ it is but 
a term used in mathematical modeling$,$ especially within the calculus. It is a 
catalyst so to speak. It is a modeling technique used due to the necessity for 
infinitesimalizing physical measures. The idealized concept for the 
``smoothed out'' model for distance measure appears acceptable. Such an 
acceptance comes from the  use of the calculus in such areas as quantum 
electrodynamics where it has great predictive power. In the subatomic region$,$
the assumption that geometric measures have physical meaning$,$ even without 
the ability to measure by external means, is justified as an appropriate 
modeling technique. Mathematical procedures applied to 
regions ``smaller than'' those dictated by the uncertainty principle are 
accepted although the reality of the infinitesimals themselves need not be 
assumed. On the other hand$,$ for this modeling technique to be applied$,$ the 
rules for ideal infinitesimalizing should be followed.\par
 The infinitesimalizing 
of ideal geometric measures is allowed. But$,$ with respect to the time concept 
this is not the case. Defining measurements of time as represented by the 
measurements of some physical periodic process is not the definition upon 
which the calculus is built. Indeed$,$ such processes cannot be 
infinitesimalized. To infinitesimalize a physical measurement using physical 
entities$,$ the entities being  
observed must be capable of being smoothed out in an ideal sense. 
This means that only the macroscopic is considered$,$ the atomic or microscopic 
is ignored. Under this condition$,$ you must be able to subdivide the device  
into ``smaller and smaller'' pieces. The behavior of these pieces can then be 
transferred to the world of the infinitesimals. Newton based the calculus not 
upon geometric abstractions but upon observable mechanical behavior. It was 
this mechanical behavior that Newton used to define physical quantities that 
could be infinitesimalized. This includes the definition of ``time.''\par
All of Newton's ideas are based upon velocities as the defining concept. 
The notation that uniform (constant) velocity exists for an object when that object is 
not affected by anything$,$ is the foundation for his mechanical observations.
This is an ideal velocity$,$ a universal velocity concept. The modern approach 
would be to add the term ``measured'' to this mechanical concept. This will 
not change the concept$,$ but it will make it more relative to natural world 
processes and a required theory of measure. This velocity concept is coupled 
with a smoothed out scale$,$ a ruler$,$ for measurement of distance. Such a 
ruler can be infinitesimalized. From observation$,$ Newton then 
infinitesimalized his uniform velocity concept. This produces the theory of 
fluxions. \par
Where does observer time come into this picture? 
It is simply a defined quantity based upon the length and velocity concept. 
Observationally$,$ it is the ``thing'' we call time that has passed when a test 
particle with uniform velocity first crosses a point marked on a scale and 
then crosses a second point marked on the same scale. This is in the absence of any 
physical process that will alter either the constant velocity or the scale.
Again this definition would need to be refined by inserting the word 
``measured.'' Absolute time is the concept that is being measured and cannot be altered as aconcept. \par
Now with Einstein relativity$,$ we are told that measured quantities are effected 
by various physical processes. All theories must be operational in that the 
concept of measure must be included. But$,$ the calculus is used. Indeed, used by Einstein in his original derivation. Thus$,$ 
unless there is an actual physical entity that can be substituted for the 
Newton's ideal velocity$,$ then any infinitesimalizing process would 
contradict the actual rules of application of the calculus to the most basic 
of physical measures. But$,$ the calculus is used to calculate the measured 
quantities. Hence$,$ we are in a quandary. Either there is no physical basis for 
mathematical models based upon the calculus$,$ and hence only selected 
portions can be realized while other selected portions are simply parameters 
not related to reality in any manner$,$ or the calculus is the incorrect 
mathematical structure for the calculations. Fortunately$,$ nature has provided 
us with the answer as to why the calculus$,$ when properly interpreted$,$ remains 
such a powerful tool to calculate the measures that describe observed physical 
behavior.\par   
In the 1930s$,$ it was realized that the measured uniform velocity of the 
to-and-fro velocity of electromagnetic radiation$,$ (i.e. light) is the only 
known natural entity that will satisfy the Newtonian requirements for an ideal 
velocity and the concepts of space-time and from which the concept of time 
itself can be defined.  The first to utilize this in relativity theory was 
Milne. This fact I learned after the first draughts of this paper were
written and gives historical 
verification of this paper's conclusions. Although$,$ it might be assumed that 
such a uniform velocity (i.e. velocity) concept as the velocity of light or light paths {\it in vacuo} 
cannot be infinitesimalized$,$ this is not the case. Such infinitesimalizing 
occurs for light-clocks and from the simple process of ``scale changing'' for a smoothed out ruler. What this means is that$,$ at its most basic {\it physical} level$,$ {\it conceptually} absolute or universal Newton time can have operational meaning as a physical foundation for a restricted form of ``time'' that can be used within the calculus. \par

As H. Dingle states it$,$ ``The second point is that the conformability of light to 
Newton mechanics . . . makes it possible to define corresponding units of 
space and time in terms of light instead of Newton's hypothetical `uniformly 
moving body.' '' [The Relativity of Time$,$ {\it Nature$,$} 144(1939): 888--890.]
It was Milne who first (1933) attempted$,$ for the Special Theory$,$ to use this definition for a ``Kinematic 
Relativity'' [{\it Kinematic Relativity}$,$ Oxford University Press$,$ Oxford$,$ 
1948] but failed to extend it successfully to the space-time environment. 
In what follows such an operational 
 time concept is being used and infinitesimalized. 
It will be seen$,$ however$,$ that based 
upon this absolute time concept another time notion is defined$,$ and this is the actual 
time notion that must be used to account for the physical changes that 
seem to occur due to relativistic processes. In practice$,$ the absolute time  is 
eliminated from the calculations and is replaced by defined ``Einstein time.'' 
It is shown that Einstein time can be infinitesimalized through the use of the 
definable  ``infinitesimal light-clocks'' and gives an exact measurement.\par

 Our first assumption is based entirely upon the logic of 
infinitesimal analysis$,$ reasoning$,$ modeling and subparticle theory. \parm 
{\leftskip 0.5in \rightskip 0.5in \noindent (i) { ``Empty'' space within 
our universe$,$ from the NSP-world viewpoint$,$ is composed of a dense-like 
NSPPM that sustains$,$ comprises and yields  N-world Special Theory  
effects. These NSPPM effects are electromagnetic in character.} \par}\parm 
\noindent This medium through which 
the effects appear to propagate comprise the objects that yield 
these effects. The next assumption is convincingly obtained from a simple and 
literal translation of the concept of infinitesimal reasoning. \parm
 {\leftskip 0.5in \rightskip 0.5in \noindent  (ii) {Any 
N-world position from or through which an electromagnetic effect appears to 
propagate$,$ when viewed from the NSP-world$,$ is embedded into a disjoint 
``monadic cluster'' of the NSPPM, where this monadic cluster 
mirrors the same unusual order properties$,$ with respect to propagation$,$ as the 
nonstandard ordering of the nonarchimedian field of hyperreal numbers 
$\hyperreal.$ [2] A monadic cluster may be 
a set of NS-substratum subparticles located within a monad of the standard 
N-world position. The propagation properties within each such monad are 
identical.}\par}\parm 
In what follows$,$ consider two (local) fundamental pairs of N-world 
positions $F_1,\ F_2$ 
that are in nonzero uniform (constant) NSP-world linear and 
relative motion. Our interest is in what effect such nonzero velocity might 
have upon electromahnetic propagation. Within the NSP-world$,$ this uniform and linear 
motion is measured by the number $w$ that is near to a standard number 
$\omega$ and this velocity is measured with respect to conceptional NSP-world time 
and a stationary subparticle field. [Note that field expansion can be additionally incorporated.] The same NSP-world linear ruler is used in 
both the NSP-world and the N-world. The only difference is that the ruler is 
restricted to the N-world when such measurements are made. 
N-world  time is  
measured by only one type of machine -- the light-clock. The concept of the 
light-clock is to be considered as any clock-like apparatus that utilizes 
either  
directly or indirectly an equivalent process. As it will be 
detailed$,$ due to the   
different propagation effects of electromagnetic radiation within the 
two  ``worlds$,$'' measured N-world light-clock time need not be the same 
as the NSP-world 
time. Further$,$ the NSP-world ruler is the measure used to define the 
N-world light-clock. \pars
Experiments show that for small time intervals $[a,b]$  
the {Galilean theory} of average velocities suffices to give 
accurate information relative to the compositions of such velocities. Let there 
be an internal function $q \colon \Hyper [a,b] \to \hyperreal,$ where  $q$ 
represents in the NSP-world a {distance function}. Also$,$ let nonnegative  and 
internal $\ell\colon \Hyper [a,b] \to \hyperreal$  be a function 
that yields the NSP-world 
velocity of the electromagnetic propagation at any  $t \in \Hyper [a,b].$ As usual 
$\monad t$ denotes the monad of standard $t,$ where ``$t$'' is an absolute NSP-world ``time'' parameter. \par

The general and correct methods of infinitesimal modeling state that$,$ within the internal 
portion of the NSP-world$,$ two measures $m_1$ and $m_2$ are {{\it 
indistinguishable for dt}} (i.e. infinitely close of order one) (notation $m_1 \sim m_2$) if and only if $0 \not= dt \in \monad 0,$ ($\monad 0$ the set of infinitesimals)
$${{m_1}\over {dt}} - {{m_2}\over{dt}} \in \monad 0. \eqno(3.1)$$
Intuitively$,$ indistinguishable in this sense means that$,$ although within the 
NSP-world the two measures are only equivalent and not necessarily equal$,$ 
the {{\it first 
level}} (or first-order) effects  these measures represent over $dt$ are 
indistinguishable within 
the N-world (i.e. they appear to be equal.) \par
In the following discussion$,$ we continue to use {{\it photon}}
 terminology. Within the N-world our photons need not be conceived of 
as particles in the sense that there is a nonzero finite N-world distance 
between individual photons. Our photons {\it may be} finite combinations of 
intermediate subparticles that exhibit$,$ when the standard part operator is 
applied$,$ basic electromagnetic field properties. Subparticle need not be discrete 
objects when viewed from the N-world$,$ but rather they could just as well give 
the {\it appearance} of a dense NS-substratum. Of course$,$ this dense NSPPM portion 
is not the usual notion of an  ``\ae{ther}'' (i.e. ether) 
for it is not a subset of 
the N-world. This dense-like portion of the NS-substratum contains the {\it nonstandard 
particle medium -- the NSPPM}. Again ``photon'' can be considered as 
but a convenient term used to discuss electromagnetic propagation. 
Now for another of our simplistic physical assumptions. \parm
{\leftskip 0.5in \rightskip 0.5in \noindent (iii) {In an N-world convex 
space neighborhood $I$ traced out over the time interval $[a,b],$  the NSPPM  
disturbances (or photons) appear to propagate linearly. }\par}\parm 
\noindent As we proceed through this derivation$,$ other such assumptions will be 
identified.\par
The functions $q,\ \ell$ need to
 satisfy some simple mathematical 
characteristic. The best known within nonstandard analysis is the concept of 
S-continuity [8]. So$,$ where defined$,$ let $q(x)/x$ (a velocity type expression) and $\ell$
be S-continuous$,$ and $\ell$ limited (i.e. finite) at each $p \in [a,b],\ (a\! +\ {\rm at}\ a,\ b\! - \ {\rm at}\ b).$ 
From compactness$,$ $q(x)/x$ and $\ell$ are S-continuous, and 
$\ell$ is limited on $\Hyper [a,b].$  
Obviously$,$ both $q$ and $\ell$ may  have 
infinitely many totally different NSP-world characteristics of which we could 
have no knowledge. But the function $q$ represents, within the NSP-world, the 
distance traveled with linear units by an identifiable NSPPM disturbance. The notion of ``lapsed-time'' is used. The $x \not= 0$ is the lapsed-time between two events.  
 It 
follows from all of this that  
for each $ t \in [a,b]$ and $t^\prime \in \monad t 
\cap \Hyper [a,b],$
$$ {{q(t^\prime)}\over{t^\prime}} - {{q(t)}\over{t}} \in \monad 0;\ \ell 
(t^\prime) - \ell (t) \in \monad 0.\eqno (3.2)$$
Expressions (3.2) give relations between nonstandard $t^\prime \in \monad 
t$ and the standard $t.$  Recall that 
if $x,\ y \in \hyperreal,$ then  $x \approx y$ iff $ x - y\in \monad 0.$ \par
 From (3.2)$,$ it follows that for each $dt \in \monad 0$ such that  $t + dt 
\in \monad t \cap \Hyper [a,b]$
$${{q(t+dt)}\over{t+dt}} \approx {{q(t)}\over {t}}, \eqno (3.3)$$ 
$$\ell( t+dt) + {{q(t+dt)}\over{t+dt}} \approx \ell (t) +{{q(t)}\over {t}}. 
\eqno (3.4)$$
\indent  
From (3.4)$,$ we have that  
$$\left(\ell( t+dt) + {{q(t+dt)}\over{t+dt}}\right)dt \sim \left(\ell (t) 
+{{q(t)}\over {t}}\right)dt. 
\eqno (3.5)$$
\indent It is now that we begin our application of the concepts of classical 
Galilean composition of velocities but restrict these ideas to the NSP-world 
monadic clusters and the notion of indistinguishable effects. You will notice 
that within the NSP-world the transfer of the classical concept of equality  
of constant or average quantities is replaced by the idea of
 indistinguishable. At the moment 
$t\in [a,b]$ that 
the standard part operator is applied$,$ an effect is transmitted through the  NSPPM as follows: 
\parm
{\leftskip 0.5in \rightskip 0.5in \noindent (iv) {For each $dt 
\in \monad 0$ and $t \in [a,b]$ such that $t +dt \in \Hyper [a,b],$ the 
NSP-world distance $q(t+dt) - q(t)$ (relative to $dt$) traveled by the 
NSPPM-effect 
within a monadic cluster is indistinguishable for $dt$ 
from the distance produced by the Galilean composition of 
velocities.}\par}\parm
\noindent From (iv)$,$ it follows that 
$$ q(t+dt) - q(t) \sim \left(\ell (t+dt)+ {{q(t+dt)}\over{t+dt}}\right)dt.
\eqno (3.6)$$                                                                 
And from (3.5)$,$ 
$$ q(t+dt) - q(t) \sim  \left(\ell (t)+ {{q(t)}\over{t}}\right)dt.\eqno
 (3.7)$$\par
Expression (3.7) is the basic result that will lead to conclusions relative 
to the Special Theory of Relativity. In order to find out exactly what 
standard functions will satisfy (3.7)$,$ let arbitrary $t_1 \in [a,b]$ be the standard 
time at which electromagnetic propagation begins from position $F_1.$ Next$,$ let 
$q= \hyper s$ be an extended standard function 
and $s$ is continuously differentiable on $[a,b].$ Applying the definition of $\sim,$ yields 
$${{\hyper s(t + dt) -s(t)}\over{dt}} \approx \ell (t) + 
{{s(t)}\over {t}}. \eqno(3.8)$$
Note that $\ell$ is microcontinuous on $\Hyper [a,b].$  For each $t \in 
[a,b],$ the value $\ell (t)$ is limited. Hence$,$ let $\st {\ell(t)} = v(t)\in 
\real.$ From Theorem 1.1 in [3] or 7.6 in [10]$,$ $v$ is continuous on $[a,b].$ 
[See note 1 part a.] Now (3.8) may be rewritten as 
$$\Hyper {\left({{d(s(t)/t)}\over{dt}}\right)} = {{\Hyper {v(t)}}\over{t}},
\eqno (3.9)$$
where all functions in (3.9) are *-continuous on $\Hyper [a,b].$
Consequently$,$ we may apply the *-integral to both sides of (3.9). [See note 
1 part b.] Now 
(3.9) implies 
that for $t \in [a,b]$
$${{s(t)}\over{t}} =\hyper {\int_{t_1}^t}{{\Hyper v (x)}\over{x}}dx, \eqno
(3.10)$$
where$,$ for $t_1 \in [a,b],$ $s(t_1)$ has been initialized to be zero.\par
 Expression (3.10) is of 
interest in 
that it shows that although (iv) is a simplistic requirement for monadic 
clusters and the requirement that $q(x)/x$ be S-continuous is a customary 
property$,$ they do not lead to 
a simplistic NSP-world function$,$ even when view at standard NSP-world times. 
It also shows 
that the light-clock assumption was necessary in that the time represented by 
(3.10) is related to the distance traveled and unknown velocity  
of an identifiable NSPPM disturbance. It is also obvious that for pure 
NSP-world 
times the actual path of motion of such propagation  
effects is highly nonlinear in 
character$,$ although within a monadic cluster the distance $\hyper s(t + dt) - 
s(t)$ is indistinguishable from that produced by the linear-like Galilean 
composition of velocities. \par
Further$,$ it is the standard function in (3.10) that allows us to cross over 
to other monadic clusters. Thus$,$ substituting into (3.7) yields$,$ since the 
propagation behavior in all monadic clusters is identical$,$\pars  
\line{\hfil $ \hyper s(t+dt) - s(t) \sim  \left(\Hyper v (t)+  
{\left(\hyper {\int_{t_1}^t}{{\Hyper v 
(x)}/{x}}dx\right)}\right)dt,$\hfil (3.11)}\smallskip
\noindent for every $t \in [a,b],\  t+dt \in \monad t \cap \Hyper [a,b]$  \par
Consider a second standard position $F_2$ at which electromagnetic reflection 
occurs at $t_2 \in [a,b],\ t_2 > t_1,\ t_2 + dt \in \monad {t_2} \cap \Hyper [a,b].$  
Then (3.11) becomes\pars
\line{\hfil $\hyper s(t_2+dt) - s(t_2) \sim \left(\Hyper v (t_2)+  
\left(\hyper {\int_{t_1}^{t_2}}{{\Hyper v 
(x)}/{x}}dx\right)\right)dt.$\hfil (3.12)}\pars
\indent Our final assumption for monadic cluster behavior is that the classical 
ballistic property holds with respect to electromagnetic propagation.\parm 
{\leftskip 0.5in \rightskip 0.5in \noindent (v) {From the exterior 
NSP-world viewpoint$,$ at standard time $t \in [a,b],$ the velocity $\Hyper v 
(t)$ acquires an added finite velocity $w$.}\par}\parm
Applying the classical statement (v), with the indistinguishable concept,  
means that the distance traveled $\hyper s(t_2 + dt) - s(t_2)$ is 
indistinguishable from $(\Hyper v (t_2) + w)dt.$ Hence$,$
$$(\Hyper v (t_2) + w)dt\sim \hyper s(t_2+dt) - s(t_2) \sim  \left(\Hyper v (t_2)+  
{\left(\hyper {\int_{t_1}^{t_2}}{{\Hyper v 
(x)}\over{x}}dx\right)}\right)dt.\eqno 
(3.13)$$                                                        
Expression (3.13) implies that 
$$\Hyper v (t_2) + w\approx \Hyper v (t_2)+  
{\left(\hyper {\int_{t_1}^{t_2}}{{\Hyper v (x)}\over{x}}dx\right)}.\eqno 
(3.14)$$                                                        
Since $\st w$ is a standard number$,$ (3.14) becomes after taking the standard 
part operator$,$                                                                 
$$ \st w=  
\St {\left(\hyper {\int_{t_1}^{t_2}}{{\Hyper v (x)}\over{x}}dx\right)}.
\eqno(3.15)$$\pars  
After reflection$,$ a NSPPM disturbance returns  to the first position $F_1$ arriving at 
$t_3  \in [a,b],\ t_1 < t_2 < t_3.$ Notice that the function $s$ does not 
appear in equation (3.15). Using the nonfavored position concept$,$ 
a reciprocal argument entails that 
$${{s_1(t_3)}\over{t_3}} = \St {\left(\hyper {\int_{t_2}^{t_3}}{{\Hyper v_1          
(x)}\over{x}}dx\right)}, \eqno(3.16)$$
$$ \st w=  
\St {\left(\hyper {\int_{t_2}^{t_3}}{{\Hyper v_1(x)}\over{x}}dx\right)},
\eqno(3.17)$$ 
where $s_1(t_2)$ is initialized to be zero. It is not assumed that $\Hyper v_1 = \Hyper v.$ \par
We now combine (3.10)$,$ (3.15)$,$ (3.16)$,$ (3.17) and obtain an 
interesting  
nonmonadic view of the relationship between distance traveled by a NSPPM disturbance and relative velocity. 
$$s_1(t_3) - s(t_2) = \st w (t_3 -t_2). \eqno(3.18)$$
Although reflection has been used to determine relation (3.18) and a 
linear-like interpretation involving reflection seems difficult to express$,$ there 
is a simple nonreflection analogue model for this behavior. \par

Suppose that a NSPPM disturbance is transmitted from a position $F_1,$ to 
a position $F_2.$ Let $F_1$ and $F_2$ have no NSP-world relative motion. 
Suppose that a NSPPM disturbance is transmitted from $F_1$ to $F_2$ with 
a constant velocity $v$ with the 
duration of the transmission $t^{\prime\prime} - t^\prime,$ where the path 
of motion is considered as linear. The disturbance continues linearly after it 
passes point $F_2$ but has increased during its travel through the monadic 
cluster at $F_2$ to the velocity $v+\st w.$ The disturbance then travels 
linearly for the same duration  $t^{\prime\prime} - t^\prime.$ The linear 
difference in the two distances traveled is $w(t^{\prime\prime} - t^\prime).$
Such results in the NSP-world should be construed only as behavior mimicked by the analogue NSPPM model.  
\par

Equations (3.10) and (3.15) show that in the NSP-world NSPPM 
disturbances propagate. Except
 for the effects of material 
objects$,$ it is assumed that in the N-world the path of motion displayed by a 
NSPPM disturbance is linear. This includes the path of motion within an 
N-world light-clock.  Let $s$ be the distance from $F_1$ of the disturbance and $w$ the NSPPM relative $F_1,\ F_2$ separation velocity. 
We continue this derivation based upon what$,$ at present$,$ appears to be additional parameters$,$ a private NSP-world time and an NSP-world rule. 
Of course$,$  
the idea of the N-world light-clock is being used as a fixed means of identifying the different 
effects the NSPPM is having upon these two distinct worlds. A question 
yet to be answered is how can we 
compensate for differences in these two time measurements$,$ the NSP-world 
private time measurement of which we can have no knowledge and N-world 
light-clocks.\par
The weighted mean value theorem for 
integrals in nonstandard form$,$ when applied to equations  (3.15) and (3.17)$,$ states that 
there are two NSP-world times $t_a,\ t_b \in \Hyper [a,b]$ such that $t_1 \leq t_a \leq t_2 
\leq t_b \leq t_3$ and  
$$ \st w= \st {\Hyper v (t_a)} \int_{t_1}^{t_2}{{1}\over{x}}dx = 
\st {\Hyper v_1(t_b})\int_{t_2}^{t_3}{{1}\over{x}}dx.\eqno (3.19)$$
[See note 1 part c.]  
Now suppose that within the local N-world an $F_1 \to F_2,$ $F_2 \to F_1$ light-clock styled
measurement for the velocity of light using a fixed instrumentation  
yields equal quantities. (Why this is the case is established in Section 6.) Model this by (*) $\st {\Hyper v (t_a)} = \st {\Hyper v_1(t_b)}=c$ a NSPPM constant quantity.  
I point out that there are many  nonconstant *-continuous functions 
that satisfy property (*). For example$,$ certain standard nonconstant linear 
functions and nonlinear modifications of them.   
Property (*) yields  
$$ \int_{t_1}^{t_2}{{1}\over{x}}dx = 
\int_{t_2}^{t_3}{{1}\over{x}}dx.\eqno (3.20)$$ 
And solving (3.20) yields
$$\ln \left({{t_2}\over{t_1}}\right) =  \ln 
\left({{t_3}\over{t_2}}\right).\eqno (3.21)$$ 
From this one has
$$ t_2= \sqrt{t_1t_3}.  \eqno (3.22)$$ 
\par
Expression (3.22) is Prokhovnik's equation (6.3.3) in reference [7].
However$,$ the interpretation of this result and the others that follow 
cannot$,$ for the NSP-world$,$ be those as 
proposed by Prokhovnik. The times $t_1,\ t_2,\ t_3,$ are standard NSPPM 
times. Further$,$ it is not logically 
acceptable when considering how to measure such time in the NSP-world or 
N-world to 
consider just any mode of measurement. The mode of light velocity measurement 
must be carried out within the confines of the language used to obtain this 
derivation. Using this language$,$ a method for time calculation 
that is permissible in the N-world is the light-clock method. Any other 
described method for 
time calculation should not include significant terms from other sources. 
Time as 
expressed in this derivation is not a mystical {\it absolute} something or other. 
It is a measured quantity based entirely upon some mode of measurement. \pars

They are two major difficulties with most derivations for expressions  
used in the Special Theory. One is the above mentioned absolute time concept. 
The other is the ad hoc nonderived N-world relative velocity. 
In this case$,$ no consideration is given as to how such a relative velocity is 
to be measured so that from both $F_1$ and $F_2$ the same result would be 
obtained. It is possible to achieve such a measurement method because of the 
logical existence of the NSPPM.\pars

In a physical-like sense, the ``times'' can be considered as the numerical values recorded by single device stationary in the NSPPM. It is conceptual time in that, when events occur, then such numerical event-times ``exist.'' It is the not yet identified NSPPM properties that yield the unusual behavior indicated by (3.22). One can use light-clocks and a counter that 
indicates$,$ from some starting count$,$ the number of times the light pulse 
has traversed back and forth between the mirror and source of our light-clock. 
 Suppose that $F_1$ and  $F_2$ can coincide. When they do coincide$,$ the $F_2$ light-clock counter number that appears conceptually first after that moment can be considered to coincide with the counter number for the $F_1$ light-clock. \pars

After $F_2$ is perceived to no longer coincide with $F_1,$ a light pulse is 
transmitted from $F_1$ towards $F_2$ in an assumed linear manner.  
The ``next'' $F_1$ counter number after this event is $\tau_{11}.$  {We  
assume that the relative velocity of $F_2$ with respect to $F_1$ may have  
altered the light-clock counter numbers, compared to the count at $F_1,$ for a light-clock riding with $F_2$. The length $L$ used to define a 
light-clock is measured by the NSP-world ruler and would not be altered. Maybe the light velocity $c$, as produced by the standard part 
operator$,$ is altered by N-world relative velocity.}  Further$,$ these two  
N-world light-clocks are only located at the two positions $F_1,\ F_2,$ and  this light pulse is represented by a NSPPM disturbance. 
The light pulse is reflected back to $F_1$ by a mirror similar to the light-clock itself. The 
first counter number on the $F_2$ light-clock to appear$,$ intuitively$,$  
``after'' this reflection is approximated by $\tau_{21}.$ The $F_1$ counter number first 
perceived after the arrival of the returning light pulse is $\tau_{31}.$ \pars

From a linear viewpoint$,$ 
at the moment of reflection$,$ denoted by $\tau_{21},$ the pulse has traveled 
an 
operational linear light-clock distance of $(\tau_{21} - \tau_{11})L.$ After 
reflection$,$ under our assumptions and nonfavored position concept$,$ 
a NSPPM disturbance would trace out the same operational linear 
light-clock distance
measured by $(\tau_{31} - \tau_{21})L.$  Thus the operational 
light-clock distance from 
$F_1$ to $F_2$ would be at the moment of operational reflection$,$ 
under our linear assumptions$,$  1/2 the sum 
of these two distances or $S_1 = (1/2) (\tau_{31} - \tau_{11})L.$ 
Now we can also determine the appropriate operational relation between these 
light-clock counter numbers for $S_1 = (\tau_{21} - \tau_{11})L.$  Hence$,$
$\tau_{31} = 2\tau_{21} -\tau_{11},$ and $\tau_{21}$ operationally behaves like an Einstein measure.\pars
After$,$ measured by light-clock counts$,$ the pulse has been received back 
to $F_1,$ a second light pulse (denoted by a second subscript of 2) is 
immediately sent to $F_2.$ Although $\tau_{31}\leq \tau_{12}$$,$ it is  
assumed that $\tau_{31}= \tau_{12}$ [See note 2.5]. 
The same analysis with new 
light-clock count numbers yields a different operational distance $S_2 = (1/2) 
(\tau_{32} -\tau_{12})L$  and $\tau_{32} = 2\tau_{22} -\tau_{12}.$  
One
 can determine the operational light-clock time intervals by considering
$\tau_{22}-\tau_{21} = (1/2)((\tau_{32} - \tau_{31}) + (\tau_{12}-\tau_{11}))$ 
and the operational linear light-clock distance difference $S_2 -S_1 = 
(1/2)((\tau_{32} - \tau_{31}) - (\tau _{12}-\tau_{11}))L.$  Since we can only 
actually measure numerical quantities as discrete or terminating numbers$,$ it 
would be empirically sound to write the N-world time intervals for these
scenarios as $t_1 = \tau_{12}-\tau_{11},t_3=(\tau_{32} - \tau_{31}).$ 
This yields the operational Einstein measure expressions in (6.3.4) of [7] as
$\tau_{22}-\tau_{21} =t_E$ and operational light-length $r_E = S_2 - S_1,$  
using our specific 
light-clock approach. This allows us to define$,$ operationally$,$ the N-world 
relative velocity as $v_E= r_E/t_E.$ [In this section$,$ the $t_1,\ t_3$ are not the same Einstein measures$,$ in form$,$ as described in [7]. But$,$ in section 4$,$ 5$,$ 6 these operational measures are used along with infinitesimal light-clock counts to obtain the exact Einstein measure forms for the time measure. This is: the $t_1$ is a specific starting count and the $t_3$ is $t_1$ plus an appropriate lapsed time.] 
\pars
Can we theoretically turn the above approximate operational approach for 
discrete N-world light-clock time  into a time continuum? Light-clocks can be 
considered from the NSP-world viewpoint. In such a case$,$ the actual NSP-world
length used to form the light-clock might be considered as a nonzero 
infinitesimal. Thus$,$ 
at least$,$ 
the numbers $\tau_{32}, \ \tau_{21}, \ \tau_{31}, \ \tau_{22}$ are 
infinite hyperreal numbers$,$ various differences would be finite and$,$ 
after taking the standard part operator$,$ all of the N-world times and lengths 
such as $t_E,\ r_E,\ S_1,\ S_2$
 should be exact and not approximate in character.
These concepts will be fully analyzed in section 6.  
Indeed$,$ as previously indicated$,$  for all of this to hold the 
velocity $c$ cannot be measured by any means. As indicated in section 6$,$ the 
actual numerical quantity $c$ as it appears in (3.22) is the standard part of 
pure NSP-world quantities. Within the N-world$,$ one obtains an ``apparent''  
constancy for the velocity of light since$,$ for this derivation$,$ 
it must be measured by means of a to-and-fro light-clock styled procedure 
with a 
fixed instrumentation.  \pars
As yet$,$ we have not discussed relations between N-world light-clock measurements 
and N-world physical laws. It should be self-evident that the assumed linearity of the light 
paths in the N-world can be modeled by the concept of projective geometry. 
Relative to the paths of motion of a light path in the NSP-world$,$ the NSPPM disturbances$,$ the
N-world path behaves as if it were a projection upon a 
plane. Prokhovnik analyzes such projective behavior and comes to 
the conclusions that in two or more dimensions the N-world light paths would 
follow the rules of hyperbolic geometry. In Prokhovnik$,$  
the equations (3.22) and the statements establishing the relations between 
the operational or exact Einstein measures $t_E,\ r_E$ and $v_E$ lead to 
the Einstein expression 
relating the light-clock determined 
 relative velocities for three 
linear positions 
having three NSP-world relative and uniform velocities $w_1,\ w_2,\ w_3.$ \pars 

In the appendix$,$ {in terms of light-clock determined Einstein measures}
and based upon the projection idea$,$ the basic Special Theory coordinate 
transformation is correctly obtained. Thus$,$ all of the NSP-world times have 
been removed from the results and even the propagation differences with 
respect to light-clock measurements. Just use light-clocks in the N-world to 
measure all these quantities in the required manner and the entire Special 
Theory is forthcoming. \par
I mention that it can be shown that  $w$ and $c$ may be measured by probes that 
are not N-world 
electromagnetic in character. Thus $w$ need not be obtained 
in the 
same manner as is $v_E$ except that N-world light-clocks would be used for 
N-world time measurements. For this reason$,$ $\st w=\omega$ is not directly related to the 
so-called textual expansion of the space within our universe. The NSPPM is not to be taken as a nonstandard translation of the Maxwell EMF equations. \par
\medskip
\leftline{\bf 4. The Time Continuum.}\par 
\medskip
With respect to models that use the classical continuum approach (i.e. 
variables are assumed to vary over such things as an interval of real 
numbers) does the mathematics perfectly measure quantities within 
nature -- quantities that cannot be perfectly measured by a human being? Or is the 
mathematics only approximate in some sense? Many would believe that if  
``nature'' is no better than the human being$,$ then classical mathematics 
is incorrect as a perfect measure of natural system behavior. However$,$ this 
is often contradicted in the limit. That is when individuals refine their 
measurements$,$ as best as it can done at the present epoch$,$ then 
the discrete human 
measurements seem to approach the classical as a limit.   
Continued exploration of this question is a 
philosophical problem that will not be discussed in this paper$,$ but it is 
interesting to model those finite things that can$,$ apparently$,$ be accomplished 
by the human being$,$ transfer these processes to the NSP-world and see what 
happens. 
For what follows$,$ when the term  ``finite''  (i.e. limited) hyperreal 
number is used$,$ since it is usually near to a nonzero real number$,$ it will 
usually refer to the ordinary  
nonstandard notion of finite except that the infinitesimals have been 
removed. This allows for the existence of finite multiplicative inverses. 
\pars  
First$,$ suppose that $t_E= \st {t_{Ea}},\ r_E= \st {r_{Ea}},\ 
S_1= \st 
{S_{1a}},\ S_2 = \st {S_{2a}}$ and each is a nonnegative real number. 
Thus $t_{Ea},\ r_{Ea},\ {S_{1a}},\  S_{2a}$ are all nonnegative  
finite hyperreal numbers. Let $L = 1/10^\omega > 0, 
\omega\in \nat^+_\infty.$ By transfer
and the result that $S_{1a},\ S_{2a},$ are  considered finite (i.e. 
near standard)$,$ then 
$S_{1a} \approx (1/2)L(\tau_{31} - \tau_{11}) \approx L
(\tau_{21} - 
\tau_{11})\Rightarrow (1/2)(\tau_{31} - \tau_{11}),\ (\tau_{21} - \tau_{11})$ 
cannot be finite. Thus$,$ by Theorem 11.1.1 [9]$,$ it can be assumed
that there exist $\eta, \gamma \in 
\nat^+_\infty$ such that
$(1/2)(\tau_{31} - \tau_{11}) =\eta,\  (\tau_{21} - \tau_{11})= \gamma.$ This 
implies that each $\tau$ corresponds to an infinite light-clock count and that 
$$\tau_{31} = 2\eta + \tau_{11},\ \tau_{21} = \gamma + 
\tau_{11}.\eqno (4.1)$$\par
In like manner$,$ it follows that
$$\tau_{32} = 2\lambda + \tau_{12},\ \tau_{22} = \delta + 
\tau_{12},\ \lambda, \delta \in \nat^+_\infty.\eqno (4.2)$$
\noindent Observe that the second of the double subscripts being 2 indicates 
the light-clock counts for the second light transmission. \pars
Now for $t_{Ea}$ to be finite requires that 
the corresponding  nonnegative $t_{1a},\ t_{3a}$ be finite. Since a different mode of conceptual time might be used in the NSP-world$,$ then there is a need for 
a number $u=L/c$ that adjusts NSP-world conceptual time to the 
light-clock count numbers. [See note 18.] By transfer of the case where  
these are real number counts$,$ this yields that $t_{3a} \approx u(\tau_{32} - 
\tau_{31}) = 2u(\lambda -\eta) +u(\tau_{12} - \tau_{11}) 
\approx 2u(\lambda -\eta) + t_{1a} $  and $t_{Ea} \approx u(\tau_{22} - 
\tau_{21}) \approx u(\delta- \gamma) + t_{1a}.$
 Hence for all of this to hold in 
the NSP-world $u(\delta -\gamma)$ must be finite or that there 
exists some $r\in \real^+$ such that $u(\delta -\gamma) \in \monad 
r.$  Let $\tau_{12} = \alpha,\ \tau_{11} = \beta.$ 
Then $t_{Ea} \approx u(\delta - \gamma) + u(\alpha - \beta)$ implies 
that $u(\alpha - \beta)$ is also finite.  \pars

The requirement that these infinite 
numbers exist in such a manner that the standard part of their products with
$L$ [resp. $u$] exists and satisfies the continuum requirements of classical mathematics 
is satisfied by Theorem 11.1.1 [9]$,$ where in that theorem $10^\omega = 1/L$ 
[resp. $1/u$]. [See note 2.]
It is obvious that the 
nonnegative numbers needed to
satisfy this theorem are nonnegative infinite 
numbers since the results are to be nonnegative and finite. 
Theorem 11.1.1 [9] allows for the appropriate $\lambda,\ \eta,\ \delta,\ \gamma$ to 
satisfy a bounding property in that we know two such numbers exist such that
$\lambda, \ \eta < 1/L^2,\ \delta,\ \gamma < 1/u^2.$ [Note: It is important to 
realize that due to this correspondence to a continuum of real numbers that 
the entire analysis as it appears in section 3 is now consistent with a mode 
of measurement. Also the time concept is replaced in this analysis with a 
``count'' concept. This count concept will be interpreted in section 8 as a 
count per some unit of time measure.] \pars
Also note that the concepts are somewhat simplified if it is assumed that
$\tau_{12} =\tau_{31}.$ In this case$,$ substitution into 4.1 yields that 
$t_{1a} \approx 2u\eta$ and $t_{3a} \approx 2u\lambda.$ Consequently$,$ $t_{Ea} = (1/2)(t_{1a} + 
t_{3a})
\approx u(\lambda+\eta).$ This predicts what is to be expected$,$ that, in this 
case$,$ the value of $t_E$ from the NSP-world viewpoint is not related to the first 
``synchronizing'' 
light pulse sent.\parm
\leftline{\bf 5. Standard Light-clocks and c.}
\medskip
I mention that the use of subparticles or the concept of the NSPPM is not 
necessary for the derivation in section 3 to hold. One can substitute for 
the NSPPM the term   ``NS-substratum'' or the like and for the  term
``monadic cluster'' of possible subparticles  just the concept of a 
``monadic neighborhood.''  It is not necessary that one assume that the 
NS-substratum contains subparticles or any identifiable entity$,$ only  
that NSPPM transmission of such radiation behaves in the 
simplistic manner stated.\pars
It is illustrative to show by a diagram of simple light-clock counts how this 
analysis actually demonstrates the two different modes of  propagation$,$ the NSP-world
 mode and the different mode when viewed from 
the N-world. In general$,$ $L$ is always fixed and for the following analysis
and$,$ for this particular scenario$,$ inf. light-clock $c$ 
may change. This process of using N-world light-clocks to approximate
the relative velocity should only be done once due to the necessity of 
``indexing'' the light-clocks when $F_1$ and $F_2$ coincide. 
In the following diagram$,$ the numbers represent actual 
light-clock count numbers as perceived in the N-world.
 The first column are those recorded at $F_1,$ the 
second column those required at $F_2.$  The arrows 
and the numbers above them represent our $F_1$ comprehension of 
what happens when the transmission is considered to take place in the 
N-world. The Einstein measures are only for the $F_1$ position.  
 \parm
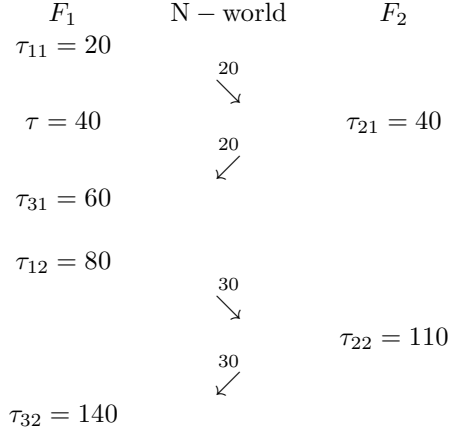
$$ \matrix{F_1&&{\rm N-world}&&F_2\cr
\tau_{11} = 20&&&&&\cr
&&{\buildrel 20 \over \searrow}&&\cr
\tau = 40 &&&&\tau_{21} = 40\cr
&&{\buildrel 20 \over \swarrow}&&\cr
\tau_{31} = 60&&&&&\cr
&&&&&&\cr
\tau_{12} = 80&&&&&\cr
&&{\buildrel 30 \over \searrow}&&\cr
&&&&\tau_{22} = 110\cr
&&{\buildrel 30 \over \swarrow}&&\cr
\tau_{32} =140&&&&&\cr}$$   
\medskip 
Certainly$,$ the above diagram satisfies the required light-clock count 
equations.  
The only light-clock counts that actually are perceivable are those at $F_1.$
And$,$ for the  transformation equations$,$ the scenario is altered. 
When the Special Theory transformation equations are obtained$,$ two distinct N-world observers are used and a third 
N-world distinct fundamental position. All light-clock counts made at each of 
these three positions are entered  into the  appropriate expressions 
for the Einstein measures {\bf as obtained for each individual position.}\par 
\bigskip
\leftline{\bf 6. Infinitesimal Light-clock Analysis.}
\medskip
In the originally presented Einstein derivation$,$ time and length are 
taken as absolute time and length. It was previously 
pointed out that this 
assumpt yields logical error. The scientific
  community extrapolated 
the language used in the derivation$,$ a language stated only in terms of 
light propagation behavior$,$ without logical reason, to the ``concept'' of 
Newtonian absolute time and length. Can the actual meaning of the 
``time'' and ``length'' expressed in the Lorentz transformation be 
determined?\parm

In what follows$,$ a measure by light-clock counts is used to analyze 
the classical transformation as derived in the Appendix-A and, essentially, such ``counts'' will replace conceptional time. [See note 1.5] The superscripts 
indicate the counts associated with the light-clocks$,$ the Einstein measures$,$ 
and the like$,$  
at the positions $F_1, \ F_2.$ The $1$ being the light-clock measures at $F_1$ 
for a light pulse event from $P,$ the $2$ for the light-clock measures at the $F_2$ 
for the same light pulse event from 
$P,$ and the 3 for the light-clock measures and its corresponding Einstein 
measures at $F_1$ for the  velocity of $F_2$ relative to $F_1.$ The 
NSP-world measured angle$,$ assuming linear projection due to the constancy of 
the velocities$,$ from $F_1$ to the light pulse event from  
$P$ is $\theta,$ and that from $F_2$ to $P$ is an exterior angle $\phi.$  \par

The expressions for our proposes are $x^{(1)}_E = v^{(1)}_Et^{(1)}_E\cos 
\theta,\  x^{(2)}_E = -v^{(2)}_Et^{(2)}_E\cos 
\phi.$ [Note: The negative is required since $\pi/2 \leq \phi \leq \pi$ and use of the customary coordinate systems.] In all that follows$,$ $i$ 
varies from 1 to 3. We investigate what happens when the standard model is 
now embedded back again into the {\it non-infinitesimal finite} NSP-world. All of the  
``coordinate'' transformation equations are in the Appendix and they  actually only involve 
$\omega_i/c.$ These equations are interpreted in the NSP-world. 
But  
as far as the light-clock counts are concerned$,$  their appropriate 
differences   
are only infinitely near to a standard number. The appropriate expressions are 
altered to take this into account. For simplicity in notation$,$ it is again 
assumed that   ``immediate'' in the light-clock count process means $\tau^{(i)}_{12} 
= \tau^{(i)}_{31}.$ Consequently$,$  
$t^{(i)}_{1a} 
\approx 2u\eta^{(i)},\ t^{(i)}_{3a} 
\approx 2u\lambda^{(i)},\ \eta^{(i)},\lambda^{(i)}\in \nat^+_\infty.$ Then
$$t^{(i)}_{Ea}\approx  u(\lambda^{(i)}+\eta^{(i)}),\ \lambda^{(i)},\eta^{(i)}\in \nat^+_\infty. \eqno (6.1)$$  \par
Now from our definition $r^{(i)}_E \approx L(\lambda^{(i)} - 
\eta^{(i)}),\ (\lambda^{(i)} - 
\eta^{(i)})\in \nat^+_\infty.$  Hence, since all of the numbers to which {\tt st} is applied 
are nonnegative and finite and $\st {v^{(i)}_{Ea}}\ \st {t^{(i)}_{Ea}} = \st 
{r^{(i)}_{Ea}},$ it follows that
$$v^{(i)}_{Ea}\approx L{{(\lambda^{(i)} - 
\eta^{(i)})}\over{u(\lambda^{(i)}}+\eta^{(i)})}.\eqno (6.2)$$ 
Now consider a set of two
4-tuples $$(\st {x^{(1)}_{Ea}},\st {y^{(1)}_{Ea}},\st {z^{(1)}_{Ea}},
\st {t^{(1}_{Ea})},$$
$$(\st {x^{(2)}_{Ea}},\st {y^{(2)}_{Ea}},\st {z^{(2)}_{Ea}},\st {t^{(2)}_{Ea}}),$$ where 
they are viewed as Cartesian coordinates in the NSP-world. First$,$ we have
$\st {x^{(1)}_{Ea}} = \st {v^{(1)}_{Ea}}\st {(t^{(1)}_{Ea}}\st {\hyper \cos \theta},\ 
\st {x^{(2)}_{Ea}} = \st {v^{(2)}_{Ea}}\st {t^{(2)}_{Ea}}\st {\hyper \cos \phi}.$
Now suppose the local constancy of  $c$. The N-world Lorentz transformation 
expressions are
$$\st {t^{(1)}_{Ea}} =  {\beta_3}(\st {t^{(2)}_{Ea}} + \st {v^{(3)}_{Ea}}\st 
{x^{(2)}_{Ea}}/c^2),$$ 
$$\st {x^{(1)}_{Ea}}= {\beta_3}(\st {x^{(2)}_{Ea}} + \st {v^{(3)}_{Ea}}\st 
{t^{(2)}_{Ea}}),$$
where $\beta_3 = \st {(1 - (v^{(3)}_{Ea})^2/c^2)^{-1/2}}.$ 
Since $L(\lambda^{(i)} - \eta^{(i)}) \approx cu(\lambda^{(i)} - \eta^{(i)})$$,$ 
the finite character of $L(\lambda^{(i)} - \eta^{(i)}),\ u(\lambda^{(i)} - 
\eta^{(i)})$ yields that $c= \st {L/u}$ [See note 8].  When transferred to the NSP-world 
with light-clock  counts$,$ substitution yields
$$t^{(1)}_{Ea}\approx u(\lambda^{(1)}+\eta^{(1)})
\approx \beta[u(\lambda^{(2)}+\eta^{(2)}) -u(\lambda^{(2)}+\eta^{(2)})K^{(3)}K^{(2)}\hyper \cos 
\phi],\eqno (6.3)$$
where $K^{(i)} = (\lambda^{(i)} - \eta^{(i)})/(\lambda^{(i)}+\eta^{(i)}),\  \beta 
= (1 - (K^{(3)})^2)^{-1/2}.$ \par
For the  ``distance''
               transformation$,$ we have
$$x^{(1)}_{Ea} \approx L(\lambda^{(1)} - \eta^{(1)})\hyper 
\cos \theta 
\approx$$ 
$$\beta(-L(\lambda^{(2)} - \eta^{(2)})
\hyper \cos \phi +{{L(\lambda^{(3)} - 
\eta^{(3)})}\over{u(\lambda^{(3)}+\eta^{(3)})}}u(\lambda^{(2)}+\eta^{(2)})).\eqno (6.4)$$
Assume 
in the NSP-world 
that $\theta \approx \pi/2,\ \phi \approx \pi.$ Consequently$,$ 
substituting into 6.4 yields
$$-L(\lambda^{(2)} - \eta^{(2)})\approx
{{L(\lambda^{(3)} - 
\eta^{(3)})}\over{u(\lambda^{(3)}+\eta^{(3)})}}u(\lambda^{(2)}+\eta^{(2)}).\eqno (6.5)$$\par 
Applying the finite property for these numbers$,$ and, for this scenario, taking into account the different modes of the corresponding light-clock measures, yields
$$ {{L(\lambda^{(3)} - 
\eta^{(3)})}\over{u(\lambda^{(3)}+\eta^{(3)})}}\approx {{-L(\eta^{(2)} - 
\lambda^{(2)})}\over{u(\lambda^{(2)}+\eta^{(2)})}} \Rightarrow v^{(3)}_{Ea} \approx - v^{(2)}_{Ea}. 
\eqno (6.6)$$
Hence$,$ $\st {v^{(3)}_{Ea}} =-\st {v^{(2)}_{Ea}}.$ [Due to the coordinate-system  selected, these are directed velocities.] This predicts that$,$ in the 
N-world$,$ the light-clock determined relative velocity of $F_2$
as measured from the $F_1$ and $F_1$ as measured from the $F_2$ positions would be 
the same if these 
special infinitesimal light-clocks are used. If noninfinitesimal 
N-world light-clocks are used$,$ then the values will be 
approximately the same and  equal in the limit. \par 
Expression 6.4 relates the light-clock counts relative to the measure of the 
to-and-fro paths of light transmission. By not  substituting for 
$x^{(2)}_{Ea},$ 
it is easily seen that $x^{(2)}_{Ea}\approx LG,$ where $G$ is an expression 
written entirely in terms of various light-clock count numbers. This implies 
that 
the so-called 4-tuples $(\st {x^{(1)}_{Ea}}, \st {y^{(1)}_{Ea}},$  
$\st {z^{(1)}_{Ea}},$ $\st {t^{(1}_{Ea})},$
$(\st {x^{(2)}_{Ea}},$$\st {y^{(2)}_{Ea}},$ $\st {z^{(2)}_{Ea}},
\st {t^{(2)}_{Ea}})$ 
are not the absolute Cartesian 
type coordinates determined by  Euclidean geometry and used to model 
Galilean dynamics. These coordinates are dynamically 
determined by the behavior of electromagnetic radiation within the 
N-world. Indeed$,$ in [7]$,$ the analysis within the  (outside of the 
monadic clusters) that leads to Prokhovnik's conclusions is only relative to 
electromagnetic propagation and is done by pure number Galilean dynamics. 
Recall that the 
monadic cluster analysis is also done by Galilean dynamics.  
\par    
In general$,$ when it is claimed that  ``length contracts'' with respect to 
relative velocity the  ``proof'' is stated as follows: $x^\prime =
\st {\beta}(x + vt); \ {\overline{x}}^\prime = \st {\beta}({\overline{x}} + 
\overline{v}\overline{t}).$ Then these two 
expressions are subtracted. Supposedly$,$ this yields $                            
{\overline{x}}^\prime -x^\prime= \st {\beta}(\overline{x} - x)$ since its assumed 
that $\overline{v}\overline{t} = vt.$
 A more complete expression would be 
$${\overline{x}}^{(1)}_E -x^{(1)}_E= \st {\beta}((\overline{x}^{(2)}_E - x^{(2)}_E)+ 
(\overline{v}^{(3)}_E\overline{t}^{(2)}_E - v^{(3)}_Et^{(2)}_E)).\eqno (6.7)$$\par
In this particular analysis$,$ it has been assumed that all NSP-world relative 
velocities $\omega_{i}, \overline{\omega_i} \geq 0.$ To obtain the classical 
length contraction expression$,$ let $\omega_{i}= \overline{\omega_i},\ 
i=1,2,3.$  Now this  implies that 
$\overline{\theta} = \theta,\ \overline{\phi} = \phi$ as they appear in the 
velocity figure on page 52 and that  
$${\overline{x}}^{(1)}_E -x^{(1)}_E= \st {\beta}(\overline{x}^{(2)}_E - 
x^{(2)}_E).\eqno (6.8)$$\par
The difficulty with this expression has been its interpretation.
Many modern treatments of Special Relativity [6] argue that (6.8) has no 
physical meaning. But in these arguments it is assumed that 
${\overline{x}}^{(1)}_E -x^{(1)}_E$ means ``length'' in the Cartesian 
coordinate sense as related to Galilean dynamics. As pointed out$,$ such a   
physical meaning is not  
the case. Expression (6.8) is a relationship between light-clock counts and$,$ 
in general$,$ displays properties of electromagnetic propagation within the 
N-world. Is there a difference between the right and left-hand sides of 6.8 
when viewed entirely from the NSP-world. First$,$ express 6.8 as
$\overline{x}^{(1)}_E -x^{(1)}_E= \st \beta\overline{x}^{(2)}_E - 
\st \beta x^{(2)}_E.$ In terms of operational light-clock counts$,$ this expression becomes
$$L(\overline{\lambda}^{(1)}\hyper {\cos {\theta}} - 
\overline{\eta}^{(1)}
\hyper {\cos {\theta}}) - L(\lambda^{(1)}\hyper {\cos \theta} - 
\eta^{(1)}\hyper {\cos \theta})
\approx\eqno (6.9)$$ 
$$L(\overline{\lambda}^{(2)}\beta{\vert \hyper {\cos {{\phi}}}\vert} - 
\overline{\eta}^{(2)}\beta
{\vert \hyper {\cos {\phi}\vert })} -  L(\lambda^{(2)}\beta {\vert\hyper {\cos 
\phi\vert}} 
- \eta^{(2)}\beta
{\vert \hyper {\cos \phi\vert}}),$$ 
where finite $\beta = (1 - (K^{(3)})^2)^{-1/2}$ and $\vert \cdot \vert$ is 
used so that the Einstein velocities are not directed numbers and the 
Einstein distances  are comparable. Also as  long as $\theta,\ \phi$ satisfy 
the velocity figure on page 45$,$ then (6.9) is independent of the specific angles
chosen in the N-world since in the N-world expression (6.8) no angles 
appear relating the relative 
velocities. That is, the velocities are not vector quantities in the N-world$,$ 
but scalars.\par
\vfil\eject
Assuming the nontrivial case that $\theta \not\approx 
\pi/2, \ \phi \not\approx \pi,$ we have from Theorem 11.1.1 [9]
that there exist $\overline{\Lambda}^{(i)},\ \overline{N}^{(i)},\ 
\Lambda^{(i)},\ N^{(i)} \in \nat_\infty,\ i=1,2$ such that
$\hyper  {\cos {\theta}}\approx \overline{\Lambda}^{(1)}/
\overline{\lambda}^{(1)} \approx 
\overline{N}^{(1)}/\overline{\eta}^{(1)}\approx$
${\Lambda}^{(1)}/
{\lambda}^{(1)} \approx {N}^{(1)}/{\eta}^{(1)},$
$\beta\vert \hyper {{\cos {\phi}}}\vert\approx \overline{\Lambda}^{(2)}/
\overline{\lambda}^{(2)} \approx \overline{N}^{(2)}/\overline{\eta}^{(2)}$
$\approx {\Lambda}^{(2)}/
{\lambda}^{(2)} \approx {N}^{(1)}/{\eta}^{(2)}.$  Consequently$,$ using the 
finite character of these quotients and the finite character of 
$L(\overline{\lambda}^{(i)}),\ L(\overline{\eta}^{(i)}),\ 
L({\lambda}^{(i)}),\ L({\eta}^{(i)}), \  i = 1,2,$ the general 
three body NSP-world
view 6.9 is
$$L(\overline{\Lambda}^{(1)} -\overline{N}^{(1)})-L(\Lambda^{(1)}-N^{(1)}) 
= L\Gamma^{(1)}\approx$$ $$ L\Gamma^{(2)}_1 =L(\overline{\Lambda}^{(2)} -\overline{N}^{(2)})-L(\Lambda^{(2)}-
N^{(2)}). \eqno (6.10)$$\par
The obvious interpretation of 6.10 from the simple NSP-world 
light propagation viewpoint is displayed by taking the standard part of 
expression 6.10.\pars

$$\st {L(\overline{\Lambda}^{(1)} -\overline{N}^{(1)})}-\st {L(\Lambda^{(1)}-
N^{(1)})}=\st {L\Gamma^{(1)}} 
=$$ $$\st {L\Gamma^{(2)}_1}= \st {L(\overline{\Lambda}^{(1)} -\overline{N}^{(1)})}-\st {L(\Lambda^{(1)}-
N^{(1)})}. \eqno (6.11)$$\par
This is the general view as to the equality of the standard NSP-world 
distance traveled by a light pulse moving to-and-fro within a light-clock as 
used to measure at $F_1$ and 
$F_2,$ as viewed from the NSPPM only$,$ the occurrence of the light pulse event from $P$. In order to interpret 6.9 for the N-world and a single NSP-world 
relative velocity$,$ you consider additionally that $\omega_1 
=\omega_2=\omega_3.$  Hence$,$
$\theta = \pi/3$ and 
correspondingly ${\phi} = 2\pi/3.$ In this case$,$  
 $\beta$ is unaltered and since $\cos \pi/3,\ \cos 2\pi/3$ are nonzero 
and finite$,$  
 6.9 now yields
$$\st {L(\overline{\lambda}^{(1)} -\overline{\eta}^{(1)})}-\st {L(\lambda^{(1)}-
\eta^{(1)})} 
=$$ $$ \st {\beta}(\st {L(\overline{\lambda}_1^{(2)} -\overline{\eta}_1^{(2)})}-
\st {L(\lambda_1^{(2)}-
\eta_1^{(2)})})\Rightarrow$$
$$(\st{L\overline{\lambda}^{(1)}} -\st {L\overline{\eta}^{(1)}})-(\st 
{L\lambda^{(1)}}-
\st {L\eta^{(1)}}) 
=$$ $$ \st {\beta}((\st {L\overline{\lambda}_1^{(2)}} -\st 
{L\overline{\eta}_1^{(2)}})-
(\st {L\lambda_1^{(2)}}-
\st {L\eta_1^{(2)}})). \eqno (6.12)$$
Or
$$\st { L(\overline{\lambda}^{(1)} -\overline{\eta}^{(1)})- L(\lambda^{(1)}-
\eta^{(1)})} =$$                                               
$$\st { L[(\overline{\lambda}^{(1)} -\overline{\eta}^{(1)})- (\lambda^{(1)}-
\eta^{(1)})] }=$$
$$\st { L{\Pi}^{(1)} }=\st {\beta}\st {L{\Pi}_1^{(2)} }=\st 
{\beta L{\Pi}_1^{(2)} }=\eqno (6.13)$$ 
$$\st { L[(\overline{\lambda}^{(1)} -\overline{\eta}^{(1)})- (\lambda^{(1)}-
\eta^{(1)})] }=$$
$$
\st {\beta L[(\overline{\lambda}_1^{(2)} -\overline{\eta}_1^{(2)})-
(\lambda_1^{(2)}-
\eta_1^{(2)})]}. $$
In order to obtain the so-called  ``time dilation'' expressions$,$ follow the 
same procedure as above. Notice$,$ however$,$ that (6.3) leads to a contradiction unless 
$$u((\overline{\lambda}^{(1)} + \overline{\eta}^{(1)})- (\lambda^{(1)}+\eta^{(1)})) \approx\beta 
u((\overline{\lambda}^{(2)}+\overline{\eta}^{(2)}) - (\lambda^{(2)}+\eta^{(2)})).\eqno (6.14)$$ 
It is interesting$,$ but not surprising$,$ that this procedure yields (6.14) without hypothesizing a relation between the $\omega_i,\ i= 1,2,3$ and implies that the timing infinitesimal light-clocks are the fundamental constitutes for the analysis.
In the NSP-world$,$ 6.14 can be re-expressed as
$$u((\overline{\lambda}^{(1)}+\overline{\eta}^{(1)}) - (\lambda^{(1)}+\eta^{(1)})) \approx 
u(\overline{\lambda}_2^{(2)} - \lambda_2^{(2)}).\eqno (6.15)$$  
\noindent Or
$$\st {u((\overline{\lambda}^{(1)}+\overline{\eta}^{(1)})} =\st 
{u\Pi_2^{(1)}}=$$ $$ 
\st{u\Pi_3^{(2)}}= \st {u(\overline{\lambda}_2^{(2)} - 
\lambda_2^{(2)})}.\eqno (6.16)$$ 
[See note 4.] 
From the N-world$,$ the expression becomes$,$ taking the standard part operator$,$
$$\st {u(\overline{\lambda}^{(1)}+\overline{\eta}^{(1)})} - \st {u(\lambda^{(1)}+\eta^{(1)})}=
$$ $$
\st{\beta}(\st{u(\overline{\lambda}^{(2)}+\overline{\eta}^{(2)})} - \st {u(\lambda^{(2)}+\eta^{(2)}))}.
\eqno (6.17)$$\par
\noindent Or
$$\st {u\Pi_2^{(1)}}=\st{\beta}\st {u\Pi_4^{(2)}}= \st 
{ \beta u\Pi_4^{(2)}}=$$
$$\st {u((\overline{\lambda}^{(1)}+\overline{\eta}^{(1)}) - (\lambda^{(1)}+\eta^{(1)}))} =    
\st {\beta u[(\overline{\lambda}^{(2)}+\overline{\eta}^{(2)}) - (\lambda^{(2)}+\eta^{(2)})]}.\eqno (6.18)$$
Note that using the standard part operator in the above expressions$,$ yields  
continuum time and space coordinates to which the calculus can now be applied. 
However$,$ the time and space measurements are not to be made with respect to an 
 universal (absolute) clock or ruler. The measurements are relative to electromagnetic 
propagation. The Einstein time and length are not the NSPPM time and length$,$ 
but rather they are concepts that incorporate a mode of measurement into 
electromagnetic field theory. This mode of measurement follows from the one wave property used for Special Theory scenarios, the property that, in the N-world, the propagation of a photon do not take on the velocity of its source. It is this that helps 
clarify properties of the NSPPM. 
 Expressions such as (6.13)$,$ (6.18) will be interpreted in the  next 
sections of this paper.\par
\medskip 
 \leftline{\bf 7. An Interpretation.}
\medskip
In each of the expressions  $(6.i),\ i=10,\ldots, 18$ 
the infinitesimal numbers $L,\ u$ are unaltered. If this is the case, then the light-clock counts would appear to be altered. As shown in Note [2], alteration of $c$ can be represented as alterations that yield infinite counts. Thus, in one case, you have a specific infinitesimal $L$ and for the other infinitesimal light-clocks a different light-clock $c$ is used. But, $L/u = c.$  Consequently the only 
alteration that takes place in N-world expressions $(6.i),
\ i=12,13,17,18$ is the infiniteimal light-clocks that need to be employed. This is exactly what (6.13) 
and (6.18) state if you consider it written as say, $(\beta L)\  \cdot$ rather than $L (\beta \ \cdot)$. Although these are external 
expressions and cannot be  ``formally'' transferred back to the N-world$,$ the 
methods of infinitesimal modeling require the concepts of  ``constant'' and 
``not constant'' to be preserved. \par 

These N-world expressions can be 
re-described in terms of N-world approximations. Simply substitute $\doteq$ for 
$=,$ a nonzero real $d$ [resp. $\mu$] for $L$ [resp. $u$] and 
real natural numbers for each 
light-clock count in equations $(6.i)$,$
\ i=12,17.$ Then for a particular $d$ [resp. $\mu$] any change 
in the light-clock
measured relative velocity $v_E$ would dictate a change in the 
the light-clocks used. 
Hence$,$ the N-world need not be concerned with the idea that  ``length''
contracts but rather it is the required light-clocks change. It is the required change in infiniteimal light-clocks that lead to real physical changes in behavior as such behavior is compared to a standard behavior. {\it But$,$ in many cases$,$  the use of light-clocks is 
not intended to 
be a literal use of such instruments.} For certain scenarios$,$ light-clocks are 
to be considered as  {\it analog 
models} that incorporate  electromagnetic energy properties. [See note 18, first paragraph.]  
\par 
The analysis given in the section 3 is done to discover 
a general property for the transmission of electromagnetic radiation. 
It is clear that property (*) does not require that the  measured velocity of 
light be a universal constant. All that is needed is that for the two 
NSP-world  times $t_a,\ t_b$ that $\st {\ell (t_a)} = \st {\ell_1 (t_b)}.$
 This means that all that is required for the most basic aspects of the 
 Special Theory 
to hold is that at two NSP-world times in the $F_1 \to F_2,\ F_2 \to F_1$  reflection 
process $\st {\ell (t_a)} = \st {\ell_1 (t_b)}$$,$ $t_a$ a time 
during the transmission prior to reflection and $t_b$ after 
reflection.  If $\ell,\ \ell_1$ are nonstandard extensions of standard functions 
$v,\ v_1$ continuous on $[a,b]$$,$ then given any $\eps \in \real^+$ 
there is a $\delta$ 
such that
for each $t,\ t^\prime \in [a,b]$ such that $\vert t - t^\prime \vert < 
\delta$ it follows  that 
$\vert v(t)-v(t^\prime) \vert < \eps/3 $ and 
$\vert v_1(t)-v_1(t^\prime) \vert < \eps/3.$ Letting $t_3 -t_1 < 
\delta,$ then $\vert t_a - t_b \vert < \delta.$ Since $\Hyper v(t_a) = 
\ell(t_a) \approx \Hyper v_1(t_b) = \ell_1(t_b),$ *-transfer implies
$\vert \Hyper v(t_2) - \Hyper v_1(t_2) \vert < \eps.$ [ See note 5.] 
Since $t_2$ is a 
standard number$,$ $\vert  v(t_2) - v_1(t_2) \vert < \eps$  implies that 
$v(t_2) = v_1(t_2).$ Hence$,$ in this case$,$ the two functions $\ell,\ \ell_1$ do 
not differentiate between the velocity $c$ at $t_2.$ But $t_2$ can be considered 
an arbitrary (i.e. NSPPM) time such that $t_1 < t_2 < t_3.$ {\bf This does not require  
$c$ to be the same 
for all cosmic times} only that $v(t) = v_1(t),\ t_1 < t < t_3.$ \par 
{\it The restriction that $\ell,\ \ell_1$ are extended standard functions 
appears necessary for our derivation.} Also$,$ this analysis is not 
related to what $\ell$ may be for a stationary laboratory. In the case of 
stationary $F_1,\ F_2,$ then the integrals are zero in equation (19) of 
section 3. The 
easiest thing to do is to simply postulate that $\st {\Hyper v (t_a)} $ is a 
universal constant. This does not make such an assumption correct. 
\par
One of the properties that will allow the Einstein velocity transformation 
expression to be derived is the {\it equilinear} property. {\bf This property 
is 
weaker than the $c$ = constant property for light propagation.} 
Suppose that you have within the NSP-world three observers 
$F_1,\ F_2,\ F_3$ that are linearly related. Further$,$ suppose that
$w_1$ is the NSP-world velocity of $F_2$ relative to $F_1$ and 
$w_2$ is the NSP-world velocity of $F_3$ relative to $F_2.$ It is assumed 
that for this nonmonadic cluster situation$,$ that Galilean dynamics also 
apply and that $\st {w_1} + \st {w_2} = \st {w_3}.$ Using the description  
for light propagation as given in section 3$,$ let $t_1$ be the 
cosmic time when a light pulse leaves $F_1$$,$ $t_2$ when it  ``passes''  
$F_2,$  and $t_3$ the cosmic time when it arrives at $F_3.$\par
From equation (3.15)$,$ it follows that
$$\st {w_1} = \st {\Hyper v_1(t_{1a})} 
\St {\left(\hyper {\int_{t_1}^{t_2}}{{1}\over{x}}dx\right)}+$$ 
$$[\st {w_2} =] \st {\Hyper v_2(t_{2a})} 
\St {\left(\hyper {\int_{t_2}^{t_3}}{{1}\over{x}}dx\right)}=$$ 
$$\st {w_3} = \st {\Hyper v_3(t_{3a})} 
\St {\left(\hyper {\int_{t_1}^{t_3}}{{1}\over{x}}dx\right)}.\eqno 
(7.1)$$ 
If $\st {\Hyper v_1(t_{1a})}=\st {\Hyper v_2(t_{2a})}=\st {\Hyper v_3(t_{3a})},$ then we 
say that the velocity functions $\Hyper v_1, \Hyper v_2, \Hyper v_3$ are {\it equilinear.} 
The constancy of $c$ implies equilinear$,$ but not conversely.  
In either case$,$ functions such as $\Hyper v_1$ and $\Hyper v_2$ need not be the same 
within a stationary laboratory after interaction. \par
Experimentation indicates that electromagnetic propagation does  
``appear'' to behave in the N-world in such a way that it does not accquire the velocity of the source. 

The light-clock analysis is consistent with the following speculation.
{\bf Depending upon the scenario$,$
the uniform velocity 
yields an effect via interactions with the subparticle field (the \underbar{NSPPM}) that uses a photon particle behavioral model.  
This is termed the \underbar{(emis)} effect.}
Recall that a ``light-clock''  can 
be considered as an analog model for the most basic of the 
electromagnetic properties. On the other hand$,$ only those experimental  methods 
that replicate or are equivalent to the methods of Einstein measure
would be relative to the Special Theory. This is one of the basic logical 
errors in theory application. The experimental language must be related to 
the language of the derivation. The concept of the light-clock$,$ linear paths 
and the like are all intended to imply substratum interactions. Any 
explanation for experimentally verified Special Theory effects should be stated in 
such a language and none other. 
I also point out that there are no paradoxes in this derivation for you cannot 
simply  ``change your mind'' with respect to the NSPPM. For example$,$ 
an observer is either in motion or 
not in motion$,$ and  not both {\bf with respect to the NSPPM.} \pars   
\leftline{\bf 8. A Speculation and Ambiguous Interpretations}\pars
Suppose that the correct principles of infinitesimal modeling were 
known prior to the M-M (i.e. Michelson-Morley) experiment. Scientists would 
know that the (mathematical) NSPPM is not an N-world entity. They would 
know that they could have very little knowledge as to the refined workings of 
this NSP-world NSPPM since $\approx$ is not an $=.$ They would have been 
forced to accept the statement of Max Planck that ``Nature does not allow 
herself to be exhaustively expressed in human thought.''[{\it The Mechanics of 
Deformable Bodies$,$ Vol. II$,$ Introduction to Theoretical Physics}$,$ Macmillian$,$ 
N.Y. (1932)$,$p. 2.]\par
Further suppose$,$ that human comprehension was advanced enough so that all 
scientific experimentation always included a theory of measurement.  
The M-M experiment would then have been performed to 
learn$,$ if possible$,$ more about this NSP-world NSPPM. When a null finding 
was obtained then a derivation such as that in section 3 might have been 
forthcoming. Then the following two expressions would have emerged from the derivation.\par

The Einstein method for measurement - the ``radar'' method - is used (see A3, p. 52) to determining the relative velocity of the moving light-clock. Using Appendix-A equations (A14), let $P$ correspond to $F_2$. Then $\theta = 0, \ \phi = \pi/2.$ Since, $x^{(2)} = 0$ from page 54, then $F_2$ is the $s$-point Hence, $t^2_E =t^{(2)}$. The superscript and subscript $s$ represents local measurements about the $s$-point, using various devices, for laboratory standards (i.e. standard behavior) and using infinitesimal light-clocks or approximating devices such as atomic-clocks. [Due to their construction atomic clocks are effected by relativistic motion and gravitational fields approximately as the infinitesimal light-clock's counts are effected.] Superscript or subscript $m$ indicates local measurements, using the same devices, for an entity considered at the $m$-point in motion relative to the $s$-point, where Einstein time and distance via the radar method as registered at $s$ are used to investigate $m$-point behavior. For example, $m$-point time is measured at the $s$-point via infinitesimal  light clock and the radar method and this represents time at the $m$-point. To determine how physical behavior is being altered, the $m$ and $s$-measurements are compared. Many claim that you can replace each $s$ with $m$, and $m$ with $s$ in what follows. This leads to various controversies which are elimianted in part 3.  A specific interpretation of
$$\st {\beta}^{-1}({\overline{t}}^{(s)} -t^{(s)})= 
\overline{t}^{(m)}_E - 
t^{(m)}_E \eqno (8.1)$$
or the corresponding  
 $$\st {\beta}^{-1}({\overline{x}}^{(s)} -x^{(s)})= 
\overline{x}^{(m)}_E - 
x^{(m)}_E \eqno (8.2)$$  
seems necessary. However, (8.2) is unnecessary since $v_E(\st {\beta}^{-1}({\overline{t}}^{(s)} -t^{(s)}))= v_E(
\overline{t}^{(m)}_E - 
t^{(m)}_E)$ yields (8.2), which can be used when convienient. Thus, only the infinitesimal light-clock ``time'' alterations are significant. Actual length as measured via the radar method is not altered. It is the clock counts that are altered.  \par

\underbar{If}$,$ in (8.2)$,$ which is employed for convenience, ${\overline{x}}^{(s)} - x^{(s)}= U^s$ (note that $x^{(s)} = v_Et^{(s)}$ etc.) is 
interpreted as ``any'' standard unit  
for length measurement at the $s$-point and 
$\overline{x}^{(m)}_E - 
x^{(m)}_E)= U^m$ the same ``standard'' unit for length measurement in a system 
moving with respect to the NSPPM (without regard to direction)$,$ 
 then for equality to take place the unit of measure $U^m$ may seem to be altered 
in the moving system. Of course$,$ it would have been immediately realized 
that the error in this last statement is that 
$U^s$ is ``any'' unit of measure. Once again$,$ the error in these two statements is the term
``any.'' (This problem is removed by application of $(14)_a$ or $(14)_b$ p. 60.)\par

 Consider 
experiements such as the M-M$,$ Kennedy-Thorndike and many others.  
When viewed from the wave state$,$ the  
interferometer measurement  technique is  determined completely by a 
light-clock type process -- the \underbar{number} of light waves in the 
linear path.  We need to use $L_{sc}^m,$ a scenario associated 
light unit$,$ for  $U^m$ and use a $L_{sc}^s$ for $U^s.$ 
It appears for this particular scenario$,$ that $L_{sc}^s$ may be considered 
the private unit of length in the NSP-world, such as $L,$  used to measure NSP-world 
light-path length. 
The  ``wavelength'' $\lambda$ 
of any 
light source must also 
be measured in the same light units. Let $\lambda = N^sL_{sc}^s.$  
Taking into consideration a unit 
conversion 
factor $k$ between the unknown NSP-world private units$,$ such that $\st {k
L_{sc}^s }= U^s$, the number of light waves
in $s$-laboratory would be $A^s \st {k L_{sc}^s}/N^s\st {kL_{sc}^s} = 
A^s/N^s,$ 
where $A^s$ is a pure number such that $A^s\st {kL_{sc}^s}$ is the ``path-length'' using the 
units in the $s$-system. In the moving system$,$ assuming that this 
simple aspect of light propagation holds in the NSP-world and the 
N-world which we did to obtain the derivation in section 3$,$ 
it is claimed that substitution yields  
$\st {\beta^{-1}A^skL_{sc}^s}/\st{\beta^{-1}N^s kL_{sc}^s}= A^s\st {\beta^{-1}kL_{sc}^s}/N^s\st{{\beta}^{-1}kL_{sc}^s}= A^s/N^s.$ 
Thus there would be 
no difference in the number of light waves in any case where the 
experimental set up 
involved the sum of light paths each of which corresponds to the to-and-fro 
process [1: 24]. Further$,$ the same conclusions would be reached using (8.2). 
not 
relevant to a Sagnac type of experiment. However$,$ this does not mean that a
similar derivation involving a polygonal propagation path cannot be 
obtained. [Indeed$,$ this may be a consequence of a result to be derived in 
article 3. However$,$ see note 8 part 4$,$ p. 80.] \par  
Where is the logical error in the above argument?
The error is the object upon which the $\st {\beta}^{-1}$ operates. Specifically
(6.13) states that \pars

$$ \st {\beta}^{-1}(A^skL_{sc}^s)\ \ \ {\buildrel ({\rm emis}) \over \longleftrightarrow}
\ \ \ \beta^{-1}(L\Pi^{(s)}) = (\beta^{-1}L)\Pi^{(s)} \ {\rm and}\eqno (8.3)$$ 
$$\st {\beta}^{-1}(N^skL_{sc}^s)\ \ \ {\buildrel ({\rm emis}) \over \longleftrightarrow}
\ \ \ \beta^{-1}(L\Pi_1^{(s)})= (\beta^{-1}L)\Pi_1^{(s)}.\eqno (8.4)$$ 
\noindent It is now rather obvious that the two (emis) aspects of the M-M experiment 
nullify each other.  
 Also for no finite  
$w$ can $\beta \approx 0.$ There is a great difference between the 
propagation properties in the NSP-world and the N-world. For example$,$ the 
classical Doppler effect is an N-world effect relative to linear propagation. 
{\bf Rather than indicating that the NSPPM is not 
present$,$ the M-M results indicate indirectly that the NSP-world NSPPM 
exists.}\par 

Apparently$,$ the well-known Ives-Stillwell$,$ and all similar$,$ experiments 
used in an attempt to verify such things as the relativistic redshift are of 
such a 
nature that they 
eliminate other effects that motion is assumed to have  upon the scenario 
associated electromagnetic 
{\it propagation.}  What was shown is that 
the frequency $\nu$ of the canal rays vary with respect to a representation for 
$v_E$ measured from electromagnetic theory in the form 
$\nu_m=\st {\beta}^{-1} \nu_s.$ First$,$ we must investigate what the so-called time 
dilation statement (8.2) means. What it means is exemplified by (6.14) and 
how the human mind comprehends the measure of  ``time.'' 
In the scenario associated (8.2) expression$,$ 
for the right and left-sides to be comprehensible$,$
the expression  should be 
conceived of as  a measure that originates with infinitesimal light-clock behavior.
It is the 
experience with a specific unit and the number of them that ``passes'' that 
yields the intuitive concept of ``observer time.'' On the other hand$,$ for some purposes 
or as some authors assume$,$ (8.2) might be 
viewed as a change in a time unit $T^s$ rather than in an infinitesimal light-clock. 
Both of these interpretations can be incorporated into a frequency statement.
First$,$ relative to the frequency of light-clock counts$,$
for a fixed stationary unit of time 
$T^s,$ (8.2) reads 
$$\st {\beta}^{-1}C_{sc}^s/T^s\doteq C^m_{sc}/T^s \Rightarrow\st {\beta}^{-1}
C_{sc}^s\doteq C^m_{sc}.\eqno (8.5)$$ \par

\underbar{But} according to (6.18)$,$ the $C_{sc}^s$ and $C^m_{sc}$ correspond to 
infinitesimal light-clocks measures and nothing more 
than that. Indeed$,$ (8.5) has nothing to do with  the concept of absolute ``time''
only with the different infinitesimal light-clocks that need to be used due to relative motion. 
This requirement may be due to (emis). Indeed$,$ the ``length contraction'' 
expression (8.2) and the  ``time dilation'' expression (8.1) have nothing 
to do with either absolute length or absolute time. These two expressions are both saying the 
same thing from two different viewpoints. There is an alteration due to the 
(emis). [Note that the  second $\doteq$ in (8.5) depends upon the $T^s$ 
chosen.] \par  
On the other hand$,$ for a relativistic redshift type experiment$,$ the usual interpretation is that 
$\nu_s\doteq p/T^s$ 
and $\nu_m\doteq  p/T^m.$ This leads to $p/T^m \doteq \st {\beta}^{-1} 
p/T^s\Rightarrow T^m  
\doteq \st {\beta}T^s.$   
Assuming that all frequency alterations due 
to (emis) have been eliminated then this is interpreted to mean that  
``time'' is slower in the moving excited hydrogen atom 
 than in the ``stationary'' laboratory. When compared to (8.5)$,$ 
there is the ambiguous interpretation in that the $p$ is considered the same for both sides
(i.e. the 
concept of the frequency is not altered by NSPPM motion). It is 
consistent with all that has come before that the Ives-Stillwell result be 
written as $\nu_s\doteq p/T^s$ 
and that $\nu_m\doteq q/T^s,$ where ``time'' as a general notion is not altered. This leads to the expression 
$$\st {\beta}^{-1}p \doteq q\  [=\ {\rm in\ the\ limit}].\eqno (8.6)$$\par 
Expression (8.6) does not correspond to a concept of  ``time'' but rather to the 
concept of alterations in emitted frequency due to (emis).   
One$,$ therefore$,$ has an ambiguous  
interpretation that in an Ives-Stillwell scenario
the number that represents the frequency of light emitted from an atomic unit moving with velocity 
$\omega$ with respect to the NSPPM is altered due to (emis).  
This (emis) alteration depends upon $K^{(3)}.$ It is critical that the two 
different infinitesimal light-clock interpretations be understood. One 
interpretation is relative to electromagnetic {\it propagation} theory. In this case$,$ the light-clock concept is taken in its most literal form. The 
second interpretation is relative to an infinitesimal light-clock as an {\it analogue} 
model. 
This means that the cause need not be related to propagation but is more 
probably due to how individual constituents interact with the NSPPM. 
The 
exact nature of this interaction and a non-ambiguous approach needs further investigation based upon 
constituent models since the analogue model specifically denies that there is 
some type of 
 {\it absolute time} dilation but$,$ rather$,$ signifies the existences of other 
possible causes. [In article 3$,$ the  $\nu_m=\st {\beta}^{-1} \nu_s$ 
is formally and non-ambiguously derived from a special line-element$,$ a universal 
functional requirement and Schr\"odinger's equation.]  \par  
In our analysis it has been 
assumed that $F_1$ is stationary in the NSP-world NSPPM. It is clear$,$ 
however$,$ that under our assumption that the scalar velocities in the NSP-world
are additive with respect to linear motion$,$ then if $F_1$ has a velocity
$\omega$ with respect to the NSPPM and  $F_2$ has the velocity 
$\omega^\prime,$ then 
it follows that the light-clock counts for $F_1$ require the use of a different light-clock with respect to a stationary $F_0$ due to the (emis) and the light-clocks for 
$F_2$ have been similarly changed with respect to a stationary $F_0$ 
due to (emis). Consequently$,$ a light-clock related expressed by $K^{(3)}$ is 
the result of the combination$,$ so to speak$,$ of these two (emis) influences. The relative NSPPM velocity $\omega_2$ of $F_1$ with respect to $F_2$ which yields the difference between these influences is that which would satisfies the additive rule for three linear positions.  
\par
 As previously stated$,$ within the NSP-world
relative to electromagnetic propagation, observer scalar velocities are either 
additive or related as discussed above. Within the N-world$,$ this last 
statement need not be so. Velocities of individual entities are modeled 
by either 
vectors or$,$ at the least$,$ by signed numbers. Once the N-world expression is 
developed$,$ then it can be modified in accordance with the usual (emis) 
alterations$,$ in which case the velocity statements are N-world Einstein 
measures. For example$,$ deriving the so-called relativistic Dopplertarian 
effect$,$
the combination of the classical and the relativistic redshift$,$
by means of a NSPPM argument such as appears in [7] where it is 
assumed that the light propagation laws with respect to the photon concept
in the NSP-world
are the same as those in the N-world$,$ is in logical error. Deriving the 
classical Doppler effect expression then$,$ when physically justified$,$ 
making the wave number alteration
in accordance with the (emis) would be the correct logic needed to obtain 
the relativistic Dopplertarian effect. [See note 6.]\par
Although I will not$,$ as yet$,$  re-interpreted all of the Special Relativity results with 
respect to this purely electromagnetic interpretation$,$ it is interesting to 
note the following two re-interpretations. The so-called variation of  
``mass'' was$,$ in truth$,$  originally derived for imponderable matter 
(i.e. elementary 
matter.) This would lead one to believe that the so-called rest mass and its 
alteration$,$ if experimentally verified$,$ is really 
a manifestation of the electromagnetic nature of such elementary matter. 
Once again the 
so-called mass alteration can be associated with an (emis)
concept. 
The $\mu$-meson decay rate may also show the same type of alteration as 
appears to be the case in an Ives-Stillwell experiment.   It does not take 
a great 
stretch of the imagination to again attribute the apparent
alteration in this rate to an (emis) process. This would lead to 
the possibility that such decay is controlled by electromagnetic properties. 
Indeed$,$ in order to conserve various things$,$  $\mu$-meson decay
is said to lead to the 
generation of the neutrino and antineutrino. 
[After this paper was completed$,$ a method was discovered that establishes that 
predicted mass and decay time alterations are (emis) effects. The derivations 
are found in article 3.]\par
I note that such things as neutrinos and antineutrinos need not exist. Indeed$,$ 
the nonconservation of certain quantities for such a scenario leads to the
conclusion that subparticles exist within the NSP-world and carry off 
the ``missing'' quantities. Thus the invention of such objects may 
definitely be considered as only a bookkeeping technique.\par
As pointed out$,$ all such experimental verification of the properly interpreted 
transformation equations can be considered as indirect evidence that the 
NSP-world NSPPM exists. But none of these results should be extended 
beyond the experimental scenarios concerned. Furthermore$,$ I conjecture that 
no matter how the human mind attempts to explain the (emis) in terms of a 
human language$,$ it will always be necessary to postulate some interaction
process with the NSPPM without being able to specifically describe this 
interaction in terms of more fundamental concepts. Finally$,$ the MA-model
specifically states that the Special Theory is a local theory and should not
be extended, without careful consideration, beyond a local time interval $[a,b].$\pars
\leftline{\bf 9. Reciprocal Relations}\pars
As is common to many mathematical models$,$ not all relations generated by 
the mathematics need to correspond to physical reality. This is the modern
approach to the length contradiction controversy [6]. Since this is a 
mathematical model$,$ there is a theory of correspondence between the physical 
language and the mathematical structure. This correspondence should be 
retained throughout any derivation. This is a NSPPM theory and what 
is stationary or what is not stationary with 
respect to the NSPPM must be maintained throughout any correspondence.
This applies to such reciprocal relations as 
$$\st {\beta}^{-1}({\overline{t}}^{(m)}_E -t_E^{(m)})= 
\overline{t}^{(s)} - 
t^{(s)} \eqno (9.1)$$
and 
$$\st {\beta}^{-1}({\overline{x}}^{(m)}_E -x^{(m)}_E)= 
\overline{x}^{(s)} - 
x^{(s)} \eqno (9.2)$$       
Statement (8.1) and (9.1) [resp. (8.2) and (9.2)] both hold from the 
NSPPM viewpoint only when 
$v_E = 0$ since it is not the question of the N-world viewpoint of relative 
velocity but rather the viewpoint that $F_1$ is fixed and $F_2$ is not fixed 
in the 
NSPPM or $\omega \leq \omega^\prime.$ The physical concept of the $(s)$ 
and $(m)$ must be maintained throughout the physical correspondence. 
 Which expression would hold for a particular scenario depends upon 
laboratory confirmation. This is a scenario associated theory.  
All of the laboratory scenarios discussed in this paper use infinitesimalized (9.1) and (9.2) as derived from line-elements and the ``view'' or comparison is always made relative to the $(s)$. Other authors$,$ such as Dingle [1] and Builder [7]$,$ have, in a absolute sense, excepted one of these sets of equations, without derivation, rather the other set. I have not taken this  
stance in this paper.\pars

One of the basic controversies associated with the Special Theory is whether 
(8.2) or (8.1) [resp. (9.1) or (9.2)] actually have physical meaning. The 
notion is that either ``length'' is a fundamental concept and ``time'' is 
defined in terms of it$,$ or ``time'' is a fundamental concept and length is 
defined in terms of it. Ives$,$ and many others assumed that the 
fundamental notion is length contraction and not time dilation. The modern 
approach is the opposite of this. Length contraction in the N-world has no physical meaning$,$ but time dilation does [6].  
We know that time is often defined in terms of 
length and velocities. But$,$ the length or time being considered here is 
Einstein length or Einstein time. This is never mentioned when this problem is 
being considered. As discussed at the end of section 3$,$ Einstein length is 
actually defined in terms of infinitesimal light-clocks or in terms of 
the Einstein velocity and Einstein time. As shown after equation (8.2) is considered, it is only infinitesimal light-clock ``time'' that is altered and length altertions is but a technical artefact.   The changes in the infinitesimal light-clock  counts yields an analogue model for physical changes that cause Special 
Theory effects.  [See note 7.] \parm

$\{$Remark: Karl Popper notwithstanding$,$ it is not the sole purpose of 
mathematical models to predict natural system behavior. The major purpose is to 
maintain logical rigor and$,$ hopefully$,$ \underbar{when applicable} to discover 
new properties for natural systems. I have used in this speculation a 
correspondence theory that takes the stance that any verifiable Special Theory
effect is electromagnetic in character rather than a problem in measure. 
However$,$ whether such effects are simply effects relative to the
propagation of electromagnetic information or whether they are effects
relative to the constituents  involved cannot be directly obtain from the 
Special Theory. All mathematically stated effects involve
the Einstein measure of relative velocity$,$ $v_E$ -- a propagation related 
measure. The measure of an effect should also be done in accordance with 
electromagnetic theory.   
As demonstrated$,$ the Special Theory should not 
be unnecessarily applied to the behavior of all nature systems since it is 
related to electromagnetic interaction; unless$,$ of course$,$ all natural systems 
are electromagnetic in character. 
Without strong justification$,$ the assumption that one theory does 
apply to all scenarios is one of the 
greatest errors in mathematically modeling.    
But$,$ if laboratory experiments verify 
that alterations are  taking place in measured quantities and these 
variations are 
\underbar{approximated} in accordance with the Special Theory$,$
then this would indicate that either the alterations are related to 
electromagnetic propagation properties or the constituents  have an 
appropriate electromagnetic character.$\}$
\parm\vfil\eject

\centerline{NOTES}
\medskip
[1] (a) Equation (3.9) is obtained as follows: since $t \in [a,b],$  
$t$ finite and not infinitesimal. Thus  division by $t$ preserves $\approx.$
Hence$,$
$$\left[t\left({{\hyper s(t + dt) - s(t)}\over {dt}}\right)- s(t)\right]/{t^2} 
\approx {{\ell(t)}\over{t}}.\eqno (1)$$ Since $t$ is an arbitrary standard number 
and $dt$ is 
assume to be an arbitrary and appropriate nonzero infinitesimal and 
the function $s(t)/t$ is differentiable$,$ the standard part of the 
left-side equals 
the standard part of the right-side. [For the end-points, the left and right derivatives are used.] Thus 
$${{d(s(t)/t)}\over {dt}} = {{v(t)}\over {t}}, \eqno (2)$$ 
for each $t \in [a,b].$ 
By *-transfer$,$ equation (3.9) holds for each $t \in \Hyper [a,b].$\par
(b) Equation 
(3.10) is then obtained by use of the *-integral and the fundamental theorem 
of integral calculus *-transferred to the NSP-world. It is useful to 
view the definite integral over a standard interval say $[t_1,t]$ as an 
operator$,$ at least$,$ defined on the set $C([t_1,t], \real )$ of all continuous real valued 
functions 
defined on $[t_1,t].$ Thus$,$ in general$,$ the fundamental theorem of integral 
calculus can be viewed as the statement that $(f^\prime , f(t) - f(t_1)) \in 
\int_{t_1}^t.$ Hence $\Hyper (f^\prime , f(t) - f(t_1)) \in 
\Hyper \int_{t_1}^t \Rightarrow (\hyper f^\prime , \Hyper (f(t) -
f(t_1))) \in \hyper \int_{t_1}^t \Rightarrow   (\hyper f^\prime , f(t) -
f(t_1)) \in \hyper \int_{t_1}^t.$\pars
(c)  To obtain the expressions in (3.19)$,$
consider $f(x) = 1/x.$ Then $\hyper f$  is limited and S-continuous on 
$\Hyper [a,b].$ Hence 
$(\hyper f, \ln t_2 - \ln t_1) \in \hyper \int_{t_1}^{t_2}.$  Hence
$\st {(\hyper f, \ln t_2 - \ln t_1)} = (f, \ln t_2 - \ln t_1) \in 
\int_{t_1}^{t_2}.$ Further (3.19) can be interpreted as an interaction property. \pars
[1.5] Infinitesimal light-clocks are based upon the QED model as to how electrons are kept in a range of distances in a hydrogen atom proton. The back-and-forth exchanges of photons between a proton and electron replaces ``reflection'' and the average distance between the proton and electron is infinitesimalized to the $L$. In this case, the proton and electron are also infinitesimalized. The large number of such interchanges over a second, in the model, is motivation for the use of the members of $\nat^+_\infty$ as count numbers.\pars

[2] The basic theorem that allows for the entire concept of infinitesimal 
light-clocks and the analysis that appears in this monograph has not been 
stated. As taken from ``The Theory of Ultralogics,''the theorem, for this application, is:\par 
{\bf Theorem 11.1.1} {\sl Let $10^\omega \in \nat_\infty.$ Then for each $ r \in \real$ there exists an $x \in \{2m/10^\omega \mid(2m \in \Hyper {\b Z}) \land
(\vert 2m\vert < \lambda 10^\omega)\},$ for any $\lambda \in \nat_\infty,$ such that $x \approx r$ (i.e. $ x \in \monad r.)$}\par
Theorem 11.1.1 holds for other members of $\nat_\infty.$ Let $L =1/10^\omega$ where $\omega$ is any hyperreal infinite natural number 
(i.e. $\omega \in \nat_\infty).$ Hence$,$ by this theorem$,$ for any positive real 
number $r$ there exists some $m \in \nat_\infty$ such that $2\st {m/10^\omega} = r.$
I point out that for this nonzero case it is necessary that $m \in 
\nat_\infty$ for if  $m \in \nat,$ then $\st{m/10^\omega} = 0.$ Since $c =\st {L/u},$ then $2\st {um} = 2\st {(L/c)m} = t = r/c$ as required. Thus, the infinitesimal light-clock determined length $r$ and interval of time $t$ are determined by the difference in infinitesimal light-clock counts $2m = (\lambda -\eta)$. Note that our approach allows the calculus to model this behavior by simply assuming that the standard functions are differentiable etc. \par
What occcurs in the infinitesimal light-clock to alter the counts? Whithin the infinitesimal light-clock, of linearity is assumed, the velocity of light can be considered as altered. The constant $c$ always denotes the N-world measured invariant velocity of light. But, using (9.2), then for zeroed light-clocks $u(\beta^{-1})(\lambda_m -\eta_m) = u(\lambda_s-\eta_s).$ From light-clock construction one reason for the changes in light-clock counts may be the velocity in the $s$-point light-clock be less than that in the $m$-point by the factor $\beta^{-1}.$ If (8.2) is employed, then the velocity in the 
$m$-point light-clock is less than that in the $s$-point by the factor $\beta.$ \pars

[2.5] (4 JUN 2000) Equating these counts here and elsewhere is done so that the ``light pulse'' is considered to have a ``single instantaneous effect'' from a global viewpoint and as such is not a signal in that globally it contains no information. 
Thus additional analysis is needed before one can state that the Special Theory applies to informational transmissions. It's obvious from section 7 that the actual value for $c$ may depend upon the physical application of this theory.\par
[3] At this point and on$,$ the subscripts 
on the $\tau$ have a different meaning than previously indicated. The 
subscripts denote process numbers while the superscript denotes the position 
numbers. For example$,$ $\tau^2_{12}$ means the light-clock count number when the second 
light pulse leaves
$F_2$ and $\tau^2_{31}$ would mean the light-clock count number when the first 
light pulse returns to position $F_2.$\par
  The additional piece of each subscript
denoted by the $a$ on this and the following pages indicates$,$ what I thought 
was obvious from the lines that follow their introduction$,$ that these are 
approximating numbers that are 
infinitesimally near to standard NSP-world number obtained by taking the 
standard part. \par
[4] Note that such infinite hyperreal numbers as $\Pi^{(2)}_3$ (here and 
elsewhere) denote the 
difference between two infinitesimal light-clock counts and since we are 
excluding the finite number infinitesimally near to 0$,$ these numbers must be 
infinite hyperreal. Infinitesimal light-clocks can be assumed to measure this 
number by use of a differential counter. BUT it is always to be conceived of 
as an infinitesimal light-clock  ``interval'' (increment$,$ difference$,$ etc.) It 
is important to recall this when the various line-elements in the next 
article are considered. \par
[5] This result is 
obtained as follows: since $t_a \leq t_2 \leq t_b,$  it follows that 
$\vert t_a - t_2 \vert < \delta,\ \vert t_b - t_2 \vert 
< \delta.$ Hence by *-transfer, $\vert \Hyper v(t_2) - \Hyper v(t_a) \vert < 
\eps/3,\ \vert \Hyper v_1(t_b) - \Hyper v_1(t_2) \vert < 
\eps/3.$ Since we assume arbitrary $\eps/3$ is a standard positive number$,$ 
then $\Hyper v(t_a) = 
\ell(t_a) \approx \Hyper v_1(t_b) = \ell_1(t_b) \Rightarrow 
\vert \Hyper v(t_a) - \Hyper v_1(t_b) \vert< \eps/3.$ Hence $\vert \Hyper 
v(t_2) - \Hyper v_1(t_2)\vert < \eps.$                                                           
 \par

[6] In this article$,$ I mention that all previous derivations for 
the complete Dopplertarian effect (the N-world and the transverse) are in 
logical error. Although there are various reasons for a redshift not just
the Dopplertarian$,$ the electromagnetic redshift based solely upon properties 
of the NSPPM can be derived as follows: \par
 (i) let $\nu^s$ denote the  ``standard'' laboratory frequency for radiation 
emitted from an atomic system. This is usually determined by the observer. 
The  
NSP-world alteration in emitted frequency at an atomic structure 
due to (emis) is  $\gamma\nu^s = \nu^{\rm radiation},$ where $\gamma = 
\sqrt {1-v_E^2/c^2}$ and $v_E$ is the Einstein measure of the relative velocity 
using light-clocks only. \par
(ii) Assuming that an observer is observing this emitted radiation in a 
direct line with the propagation and the atomic structure is receding with 
velocity $v$ from the 
observer$,$ the frequwncy of the electromagnetic 
propagation$,$ within the 
N-world$,$ is altered compared to the observers standards. This alteration is $\nu^{\rm radiation}
(1/(1+v/c)) = \nu^{\rm received}.$ Consequently$,$ this yields the total 
alteration as $\gamma\nu^s(1/(1+v/c)) = \nu^{\rm received}.$ Note that $v$ 
is measured in the N-world and can be considered a directed velocity. 
Usually$,$ if due to 
the fact that we are dealing with electromagnetic radiation$,$  
we 
consider $v$ the Einstein measure of linear velocity (i.e. $v = v_E$)$,$ 
then the total Dopplertarian effect for $v \geq 0$ can be written as
$$\nu^s\left({{1-v_E/c}\over{1 + v_E/c}}\right)^{1/2} = \nu^{\rm 
received}.\eqno (3)$$\par
It should always be remembered that there are other reasons$,$ such as the 
gravitational redshift and others yet to be analyzed$,$ that can mask 
this total Dopplertarian redshift. \par

[7]  A question that has been asked relative to the new 
derivation that yields verified Special Theory results is why in the 
N-world do we have the apparent nonballistic effects associated with 
electromagnetic radiation? In the derivation$,$ the opposite was assumed for the 
NSP-world
monadic clusters. The constancy of the {\it measure}$,$ by light-clocks and the 
like$,$ of the $F_1 \to F_2,\ F_2 \to F_1$ velocity of electromagnetic radiation was modeled by 
letting $\st {t_a} = \st {t_b}.$ As mentioned in the section on the Special 
Theory$,$  the Einstein velocity measure transformation expression can be 
obtained prior to embedding the world into a hyperbolic velocity space. It is 
obtained by considering three in-line standard positions $F_1,\ F_2, \ F_3$ 
that have the NSP-world velocities $w_1$ for $F_2$ relative to $F_1$$,$
$w_2$ for $F_3$ relative to $F_2$ and the simple composition $w_3= w_1 + w_2$
for $F_3$ relative to $F_1$. Then simple substitution in this expression
yields 
$$v_E^{(3)}= (v_E^{(1)}+v_E^{(2)})/\left(1 + 
{{v_E^{(1)}v_E^{(2)}}\over{c^2}}\right).\eqno (4)$$
This relation is telling us 
something about the required behavior in the N-world of electromagnetic 
radiation. To see that within the N-world we need to assume   for 
electromagnetic radiation effects the nonballistic property$,$  simply  
let  $v_E^{(2)}\doteq c,$ where always $v_E^{(2)}<c.$ Then $v_E^{(3)} 
\doteq c.$ Of course$,$ the reason we do not have a contradiction is that 
we have two distinctly different views of the behavior of electromagnetic 
radiation$,$ the NSP-world view and the N-world view. Further$,$ note how$,$ for 
consistency$,$ the velocity of electromagnetic radiation is to be measured. It 
is measured by the Einstein method$,$ or equivalent$,$ relative to a 
to-and-fro path and measures of ``time'' and ``distance'' by means of a 
(infinitesimal) light-clock counts. Since one has the NSPPM, then letting $F_1$ be fixed in that medium, assuming that ``absolute'' physical standards are measured at $F_1,$ equation (4) indicates why, in comparison, physical behavior varies at $F_2$ and $F_3$. The hyperbolic velocity space properties are the cause for such behavior differences.  \pars\par
  
I am convinced that the dual character of the Special 
theory derivation requires individual reflection in order to be understood fully. In the NSP-world$,$ 
electromagnetic radiation behaves in one respect$,$ at least$,$ like a particle in 
that it satisfies the ballistic nature of particle motion.  The reason that 
equation (3) is derivable is due to the definition of Einstein time. But 
{\it Einstein time$,$ as measured by electromagnetic pulses$,$ models the 
nonballistic 
or one and only one wave-like property in that a wave front does not partake of the 
velocity of the source.} This is the reason why  
I wrote that {\it a NSPPM disturbance would trace the 
same operational linear light-clock distance.} The measuring light-clocks are 
in the N-world in this case. $F_1$ is modeled as fixed in the NSPPM and 
$F_2$ has an NSP-world relative velocity. The instant the light pulse is 
reflected 
back to $F_1$ it does not$,$ from the N-world viewpoint$,$ partake of the N-world
relative velocity and therefore traces out the exact same apparent N-world 
linear path. The position $F_2$ acts like a virtual position having no other 
N-world effect upon the light pulse except a reversal of direction. \par
[8] This expression implies that the ``$c$'' that appears here and elsewhere is to be measured by infinitesimal light-clocks.  As noted $u \approx L/c,$ but infinitesimal light-clock construction yields that $u =L/c.$ For a fixed $L$, from the NSPPM  viewpoint, $u$ is fixed. Notice that $t^{(i)}\approx u(2\eta^{(i)}) = u(\gamma^{(i)}),\ \gamma^{(i)}\in \nat^+_\infty.$\par
[9] In this monograph, conceptual time is used and NSPPM and gravitational field processes yield non-classical relations between these times. For example, $t_2 = \sqrt{t_1t_3}.$ For the Special Theory, there is only one aspect of physical-world behavior that corresponds to the infinitesimal-world behavior. This is the sudden photon interaction with other particles. Hence, such interactions are particle-like, which predicts the QED assumption. An NSPPM velocity for the source is always necessary for photon emission due to a photon's momentum. In the actual derivation, the wave-property, where classical wave-mechanics can be applied, is not a property within the infinitesimal-world. Classical wave-mechanics model photon paths of motion within our physical-world. Wave-behavior emerges after the "st" operator is applied. The particle behavior takes place only for the interactions. Hence, the probability interpretation that comes from a photon's wave-property can be used to predict the number of photon interactions. Consequently, there is neither a contradiction between these two interpretations nor the particle assumption.\par 
[10] Modern derivations attempt to remove the mode of measurement, but by so doing the twin anomaly occurs that cannot by removed even by using GR [7]. It can be removed by using the method presented here and Einstein measures. The solution of the twin anomaly [7, pp:108-111] clearly implies that, with certain exceptions, all Special Theory ``time'' dependent effects are relative. The effects have physical meanings only when compared. The major and maybe only exception within our universe are electromagnetic radiation  (photon) related effects, which satisfy the inertia frame of reference $\rm I_S$ requirements. When such effects are separated from their source and prorogate with respect to these requirements, then they can affect the behavior of other physical objects. When compared, such effects can become significant. \parm

\centerline{\bf REFERENCES}
\parm
\noindent {\bf 1} Dingle$,$ H.$,$ {\it The Special Theory of Relativity}$,$ Methuen's 
Monographs on Physical Subjects$,$ Methuen \& Co. LTD$,$ London$,$ 1950.\pars
\noindent {\bf 2} Herrmann$,$ R. A.$,$ {\it Some Applications of Nonstandard 
Analysis to Undergraduate Mathematics -- Infinitesimal Modeling}$,$ 
Instructional Development Project$,$ Mathematics Dept. U. S. Naval Academy$,$ 
Annapolis$,$ MD 21402-5002$,$ 1991. http://arxiv.org/abs/math/0312432 \pars
\noindent {\bf 3} Herrmann$,$ R. A.$,$ Fractals and ultrasmooth microeffects. {\it 
J. Math. Physics.}$,$ {\bf 30}(1989)$,$ 805--808.\pars
\noindent {\bf 4} Herrmann$,$ R. A.$,$  Physics is legislated by a cosmogony. {\it
Speculat. Sci. Technol.}$,$ {\bf 11}(1988)$,$ 17--24.\par
\noindent {\bf 5} Herrmann$,$ R. A.$,$ Rigorous infinitesimal modelling$,$ {\it Math. 
Japonica} 26(4)(1981)$,$ 461--465.\pars
\noindent {\bf 6} Lawden$,$ D. F.$,$ {\it An Introduction to Tensor Calculus$,$ 
Relativity and Cosmology}$,$ John Wiley \& Son$,$ New York$,$ 1982.\par
\noindent {\bf 7} Prokhovnik$,$ S. J.$,$ {\it The Logic of Special Relativity}$,$ 
Cambridge University Press$,$ Cambridge$,$ 1967.\pars
\noindent {\bf 8} Stroyan$,$ K. D. and Luxemburg$,$ W. A. J. 1976. {\it Introduction to 
the Theory of Infinitesimals}$,$ Academic Press$,$ New York.\pars
\noindent {\bf 9} Herrmann$,$ R. A.$,$ {\it The Theory of Ultralogics}$,$ (1992)
 \hfil\break
http://www.arxiv.org/abs/math.GM/9903081\hfil\break
http://www.arxiv.org/abs/math.GM/9903082\pars

\noindent {\bf 10} Davis M.$,$ {\it Applied Nonstandard Analysis$,$} John Wiley 
\& Sons$,$ New York$,$ 1977.
\medskip
\noindent NOTE: Portions of this monograph have been published in various journals. \pars
\parm

 {Overhead material relative to the paper
``A corrected derivation for the Special Theory of relativity'' 
as presented at the above mentioned MAA meeting of Nov. 14$,$ 1992.}\pars
 \centerline{\bf Relativity and Logical Error}\pars 
 {In a 1922 lecture Einstein stated}
{the bases of his Special Theory.}
{``Time cannot be absolutely defined$,$}
{and there is an inseparable relation}
{between time and signal velocity.''}
{In the paper you're going to receive$,$}
{the first phrase is shown to be {\bf false}.}
{Thus with respect to natural models$,$}
{the stated hypotheses yield an}
{inconsistent theory.}\par
{Originally$,$ Einstein did NOT reject}
{an \ae{ther} or medium concept.}
{In the same lecture$,$ he said}  
{``Since then [1905] I have come to}
{believe that the motion of the}
{Earth [through the \ae{ther}] cannot be}
{detected by any optical experiment}
{though the Earth is revolving about}
{the Sun.''}
{Einstein also stated that he simply}
{couldn't describe the properties}
{of such an \ae{ther}.}\par  
{In his derivation$,$ he first uses the}
{term  ``clock'' as meaning  {\bf any}}
{measure of time within the natural}
{world without further defining}
{the apparatus.}
{But then he restricts the}
{characterization of such clocks}
{by adding light propagation}
{terminology relative to their}
{synchronization and$,$ hence$,$ creates}
{{\bf a new predicate model}.}\par
{Einstein now uses these restricted}
{clocks to measure a new time$,$} 
{\bf the proper time$,$}
{in terms of}
{additional light propagation language.}
{This is a third predicate model.}
{After this$,$ he \underbar{assumes} that the second}
{predicate model is Newtonian infinitesimal}
{time$,$ another type of absolute time}
{which is a different fourth predicate model.}
{Thus substituting one predicate model for}
{another$,$ as if they are the same$,$}
{he obtains the Lorentz transformation.}
{Of course$,$ this substitution is a}
{logical error.}\par {Now$,$ Einstein's}
{form of the Lorentz transformation}
{has proper time on the left-side and}
{Newtonian absolute time on the right.}
{Then to apply this transformation$,$}
{the predicate model for proper time with its}
{light propagation language is extended to include}
{an absolute {\bf any} time concept.}
{The logical error of substituting one}
{predicate model for another predicate is}
{compounded by the error of}
{{\bf model theoretic generalization}.}
{The statement that what holds for one}
{domain (time restricted by the language}
{of light propagation) cannot be extended}
{ad hoc to a larger domain.}\par
\bigskip

\centerline{\bf Appendix-A}\bigskip

\noindent {\bf 1. The Need for Hyperbolic Geometry}\parm
In this appendix$,$ it is shown that from equations (3.21) and (3.22) the Lorentz transformation are derivable. All of the properties for the Special Theory are based upon ``light'' propagation. In Article 2$,$ the concern is with two positions $F_1,\ F_2$ in the NSPPM within the NSP-world and how the proposed NSPPM influences such behavior. Prior to applications to the N-world$,$ with the necessity for the N-world Einstein measures$,$ the NSPPM exhibits infinitesimal behavior and special NSPPM non-classical global behavior. The behavior at specific moments of NSPPM time for global positions and classical uniform velocities are investigated. \par 

The following is a classical description for photon behavior. Only NSPPM relative velocities (speeds) are being considered. 
Below is a global diagram for four points that began as the corners of a square, where $u$ and $\omega$ denote uniform relative velocities between point locations and no other point velocities are considered. The meanings for the symbolized entities are discussed below. \par\bigskip

\line{\hskip 1.35in$\bullet F_1\ t$ $\sim\kern -0.25em\sim\kern -0.25em\longrightarrow $\hskip 1.25in$\omega \longrightarrow \bullet F_2(t(p_1))\sim\kern -0.25em\sim\kern -0.25em\longrightarrow$\hfil}
\vskip 0.35in 
\line{\hskip 1.35in $\bullet F_1'\ t'\sim\kern -0.25em\sim\kern -0.25em$\hfil}
\line{\hskip 1.20in $u \ \downarrow$\hfil}
\line{\hskip 1.35in $\bullet F_1'\ t(p_2)$\hskip 2.00in $\omega \longrightarrow \bullet F_2'(t(p_2))$\hfil} 

\line{\hskip 1.20in $u \ \downarrow$ \hskip 2.76in$\ \downarrow\ u$\hfil}\par\bigskip

Consider the following sequence of (conceptual) NSPPM time-ordered events. First, the N-world position points $F_1,\ F_2,\ F_1',\ F_2'$ are stationary with respect to each other and form the corners of a very small rhombus, say the side-length is the average distance $d$ between the electron and proton within an hydrogen atom. The sides are $\overline{F_1,F_2}, \ \overline{F_2F_2'},\ \overline{F_2',F_1'},\ \overline{F_1'F_1}$. At the NSPPM time $t_g$, the almost coinciding $F_2,F_2'$ uniformally recede from the almost coinciding $F_1,F_1'$ with constant velocity $\omega$.  At a time $t > t_g,$ where the distance between the two groups is significantly greater than $d$, one process occurs simultaneously. The point $F_1'$ separates from $F_1$ with relative velocity $u$ and $F_2'$ separates from $F_2$ with a relative velocity $u$. [Using NSP-world processes, such simultaneity is possible relative to a non-photon transmission of information (Herrmann, 1999).] At any time $\geq t,$ the elongating line segments $\overline{F_1F_1'}$ and $\overline{F_2F_2'}$ are parallel and they are not parallel to the parallel elongating line segments $\overline{F_1F_2}$ and $\overline{F_1'F_2'}.$ \par

At NSPPM time $t,$ a photon $p_1$ is emitted from $F_1$ towards $F_2$ and passes through $F_2$ and continues on. As $F_1'$ recedes from $F_1$, at $t' >t,$ a photon $p_2$ is emitted from $F_1'$ towards $F_2'$. The original classical photon-particle property that within a monadic cluster photons prorogate with velocity $\omega +c$ is extended to this global environment. [Again there are NSP-world processes that can ensure that the emitted photons acquire this prorogation velocity (Herrmann, 1999).] Also, this classical photon-particle property is applied to $u$. Thus, photon $p_2$ is assumed to take on an additional velocity component $u$. Photon, $p_1$, passes through $F_2$ at the NSPPM time $t(p_1)$. Then $p_2$ is received at point $F_2'$ at time $t(p_2)$.  \par

Classically, $t(p_1') > t(p_1).$ From a viewpoint relative to elongating $\overline {F_1F_1'}$, the distance between the two photon-paths of motion measured parallel to elongating $\overline {F_1F_1'}$  is $u(t(p_2) -t).$ On the other hand, from the viewpoint of elongating $\overline{F_1F_1'}$, the distance between photon-paths, if they were parallel, is $u(t' - t).$ 
By the relativity principle, from the viewpoint of $F_1'$, the first equation in (3.19) should apply. Integrating, where $\st {\hyper {v(t_a))}} = c$, one obtains $u(t(p_2)- t)) = ue^{\omega/c}(t' - t).$ [Note: No reflection is required for this restricted application of (3.19).] This result is not the classical expression $u(t(p_2)- t)).$  For better comprehension, use infinitesimal light-clocks to measure NSPPM time. Then using the same NSPPM process that yields information instantaneously throughout the standard portion of the NSPPM, all clocks used to determine these times can be set at zero when they indicate the time $t$. This yields that the two expressions for the distance are $ut(p_2)$ and $ue^{\omega/c}t'.$ However, the classical expression $ut(p_2)$ has the time $t(p_2)$ dependent upon both $\omega$ and, after the  $t'$ moment, upon $u.$ But, for the relativistic expression, the $t'$ is neither dependent upon the $u$ velocity after $t'$ nor the $\omega$ and the factor $e^{\omega/c}$ has only one variable $\omega.$ What property does this NSPPM behavior have that differentiates it from the classical?\par

Consider the two velocities $u$ and $ue^{\omega/c}.$ These two velocities only correspond when  $\omega =0.$ Hence, if we draw a velocity diagram, one would conclude that, in this case, the velocities are trivially ``parallel.''  Using Lobatchewskian's horocycle construction$,$ Kulczycki (1961) shows that for ``parallel  geometric'' lines in hyperbolic space$,$ the distance between each pair of such lines increases (or decreases) by a factor $e^{x/k},$ as one moves an ordinary distance $x$ along the lines and $k$ is some constant related to the $x$ unit of measurement. Phrasing this in terms of velocities, where $x = \omega$ and $k = c$, then, for this case, the velocities, as represented in the NSPPM by standard real numbers, appear to satisfy the properties for an hyperbolic velocity-space. Such velocity behavior would lead to this non-classical NSPPM behavior.\par

When simple classical physics is applied to this simple Euclidian configuration within the NSPPM, then there is a transformation $\rm \Phi\colon NSPPM \to$ N-world, which is characterized by hyperbolic velocity-space properties. This is also the case for relative velocity and collinear points, which are exponentially related to the Einstein measure of relative velocity in the N-world. In what follows, this same example is used but generalized slightly by letting $F_1$ and $F_2$ coincide. \parm 
\noindent {\bf 2. The Lorentz Transformations}\parm
Previously, we obtained the expression that $t_2 = \sqrt {t_1t_3}.$ The Einstein measures are defined formally as 
$$\cases{t_E = (1/2)(t_3 + t_1)&\cr
         r_E = (1/2)c(t_3 - t_1)&\cr
         v_E = r_E/t_E,\ {\rm where\ defined.}&\cr}\eqno (A1)$$
Notice that when $r_E = 0,$ then $v_E = 0$ and $t_E = t_3 = t_1= t_2$ is not Einstein measure.\par 

\noindent The Einstein time $t_E$ is obtained by considering the ``flight-time'' that would result from using one and only one wave-like property not part of the NSPPM but within the N-world. This property is that the $c$ is not altered by the velocity of the source. This Einstein approach assumes that the light pulse path-length from $F_1$ to $F_2$ equals that from $F_2$ back to $F_1.$ Thus, the Einstein flight-time used for the distance $r_E$ is $(t_3 - t_1)/2$. The $t_E,$ the Einstein time corresponding to an infinitesimal light-clock at $F_2,$ satisfies $t_3 - t_E = t_E - t_1.$ From $(A1),$ we have that 
$$t_3 = (1+v_E/c)t_E\ {\rm and } \ t_1 = (1 - v_E/c)t_E, \eqno (A2)$$ and$,$ hence$,$ $t_2 = (\sqrt {1 - v_E^2/c^2})t_E.$ Since $e^{\omega/c} = \sqrt {t_3/t_1},$ this yields
$$e^{\omega/c} = \left({{1 + v_E/c}\over{1 - v_E/c}}\right)^{(1/2)}. \eqno (A3)$$ \par

Although it would not be difficult to present all that comes next in terms of the nonstandard notions$,$ it is not necessary since all of the functions being consider are continuous and standard functions. The effect the NSPPM has upon the N-world are standard effects produced by application of the standard part operator ``st.''\par

From the previous diagram, let $F_1$ and $F_2$ coincide and not separate. Call this location $P$ and consider the diagram below. This is a three position classical NSPPM light-path and relative velocity diagram used for the infinitesimal light-clock analysis in section 6 of Article 2. This diagram is not a vector composition diagram but rather represents linear light-paths with respect to medium measures for relative velocities. It is also a relative velocity diagram to which hyperbolic ``geometry'' is applied. \par\bigskip
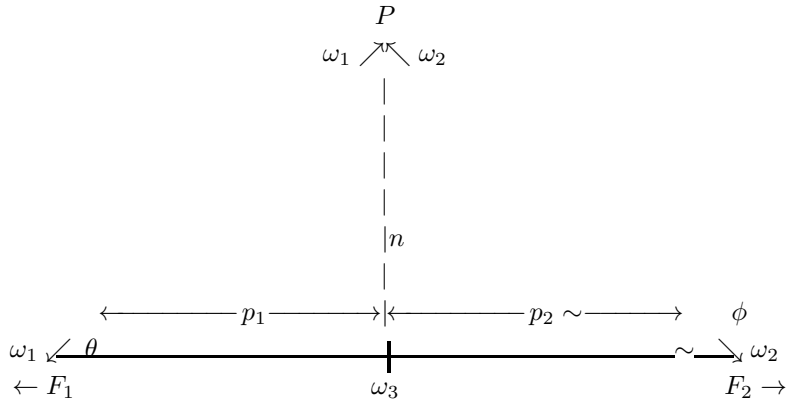
\line{\hfil$P$\hfil} 
\line{\hfil $\omega_1\nearrow\nwarrow \omega_2$\hfil}
\line{\hskip 0.04in\hfil $\vert$ \hfil}
\line{\hskip 0.04in\hfil $\vert$ \hfil}

\line{\hskip 0.04in\hfil $\vert$ \hfil}

\line{\hskip 0.04in\hfil $\vert$ \hfil}

\line{\hskip 0.13in\hfil $\vert n$ \hfil}

\line{\hskip 0.04in\hfil $\vert$ \hfil}

\line{\hskip 1.50in {\hbox to 0.75in{\leftarrowfill}}$p_1${\hbox to 0.60in{\rightarrowfill}}$\vert$\kern -0.1em{\hbox to 0.75in{\leftarrowfill}}$p_2\sim${\hbox to 0.53in{\rightarrowfill}}\hskip 0.25in$\phi$\hfil}
\line{\hfil\hskip 0.10in$\omega_1 \swarrow$\kern -.6em\vbox{\hrule width .15in}$\theta$\kern -.6em\vbox{\hrule width 1.60in}{\vrule height6pt width1pt depth0pt}\kern -0.1em{\vrule height0pt width1pt depth6pt}\kern-.1em\vbox{\hrule width 1.5in}$\sim$\vbox{\hrule width 0.20in}\kern -.6em $\searrow \omega_2$\hfil}
\line{\hskip 1.05in$\leftarrow F_1 $\hskip 1.55in $\omega_3$\hskip 1.70in  $ F_2\rightarrow$\hfil}\par \bigskip

Since Einstein measures are to be associated with this diagram$,$ then this diagram should be obtained relative to infinitesimal light-clock counts and processes in the NSPPM. The three locations $F_1,\ F_2,\ P$ are assumed$,$ at first$,$ to coincide. When this occurs$,$ the infinitesimal light-clock counts coincide. The positions $F_2, \ F_1$ recede from each other with velocity $\omega_3$. The object denoted by location $P$ recedes from the $F_1,\ F_2$ locations with uniform NSPPM velocities$,$ in standard form$,$ of $\omega_1,\ \omega_2$, respectively.  Further, consider the special case where both are observing the pulse sent from $P$ at the exact some $P$-time. This produces the internal angle $\theta$ and exterior angle $\phi$ for this velocity triangle. 
The segments marked $p_1$ and $p_2$ are the projections of the velocity representations (not vectors) $F_1P$ and $F_2P$ onto the velocity representation $F_1F_2$. The $n$ is the usual normal for this projection. We note that $p_1 + p_2 = \omega_3.$  We apply hyperbolic 
trigonometry in accordance with [2]$,$ where we need to consider a particular $k$. We do this by scaling the velocities in terms of light units and let $k = c$. From [2, p. 143]\par
$$\cases{\tanh(p_1/c) = (\tanh (\omega_1/c)) \cos \theta&\cr
          \tanh (p_2/c) = - (\tanh (\omega_2/c))\cos \phi&\cr},\eqno (A4)$$
and also 
$$\sinh (n/c) = (\sinh (\omega_1/c))\sin \theta = (\sinh (\omega_2/c))\sin \phi. \eqno (A5)$$
Now$,$ eliminating $\theta$ from $(A4)$ and $(A5)$ yields [1, p. 146]
$$\cosh (\omega_1/c) = (\cosh (p_1/c))\cosh (n/c). \eqno (A6)$$
Combining (A4), (A5) and (A6) leads to the hyperbolic cosine law [2, p. 167].  
$$\cosh (\omega_1/c)=(\cosh (\omega_2/c))\cosh (\omega_3/c) + (\sinh(\omega_2/c))(\sinh(\omega_3/c))\cos \phi. \eqno (A7)$$
From $(A3),$ where each $v_i$ is the Einstein relative velocity$,$ we have that 
$$e^{\omega_i/c} = \left({{1 + v_i/c}\over{1 - v_i/c}}\right)^{(1/2)}, i = 1,2,3.\eqno (A3)'$$
From the basic hyperbolic definitions$,$ we obtain from $(A3)'$
$$\cases{\tanh (\omega_i/c) = v_i/c&\cr
          \cosh (\omega_i/c) = (1 - v_i^2/c^2)^{-1/2} = \beta_i&\cr
          \sinh (\omega_i/c) = \beta_iv_i/c&\cr}. \eqno (A8)$$
Our final hyperbolic requirement is to use 
$$\tanh (\omega_3/c) = \tanh(p_1/c + p_2/c)= {{\tanh (p_1/c) + \tanh (p_2/c)}\over{1 + (\tanh (p_1/c))\tanh (p_2/c)}}. \eqno (A9)$$\par
Now into $(A9),$ substitute $(A4)$ and then substitute the first case from $(A8).$ One obtains
$$v_1\cos \theta = {{v_3 - v_2\cos \phi}\over{1 - \alpha}},\ \alpha = {{v_3v_2\cos \phi}\over{c^2}}. \eqno (A10)$$
Substituting into $(A7)$ the second and third cases from $(A8)$ yields
$$\beta_1=\beta_2\beta_3(1-\alpha), \ \beta_i = (1-v_i^2/c^2)^{-1/2}.\eqno (A11)$$
From equations $(A11), \ (A5)$ and the last case in $(A8)$ is obtained
$$v_1\sin \theta ={{v_2\sin\phi}\over{\beta_3(1-\alpha)}}. \eqno (A12)$$\par

 For the specific physical behavior being displayed, the photons received from $P$ at $F_1$ and $F_2$ are ``reflected back'' at the  
NSPPM $P$-time $t^r.$ We then apply to this three point scenario our previous results. [Note: For comprehension, it may be necessary to apply certain relative velocity viewpoints such as from $F_1$ the point $P$ is receding from $F_1$ and $F_2$ is receding from $P$.  In this case, the NSPPM times when the photons are sent from $F_1$ and $F_2$ are related. Of course, as usual there is assumed to be no time delay between the receiving and the sending of a ``reflected''  photon.] In this case$,$ let $t^{(1)},\ r^{(1)},\ v_1$ be the Einstein measures at $F_1$ for this $P$-event$,$ and $t^{(2)},\ r^{(2)},\ v_2$ be the Einstein measures at $F_2$. Since $t^r = \beta_1^{-1}t^{(1)},\ t^r = \beta_2^{-1}t^{(2)}$ (p. 52), then  
$${{t^{(1)}}\over {\beta_1}} = {{t^{(2)}}\over {\beta_2}}\ {\rm and}\ r^{(1)}=v_1t^{(1)},\ r^{(2)}=v_2t^{(2)}.\eqno (A13)$$\par

Suppose that we have the four coordinates$,$ three rectangular$,$ for this $P$ event as measured from $F_1 = (x^{(1)}, y^{(1)}, z^{(1)}, t^{(1)})$ and from $F_2 = (x^{(2)}, y^{(2)}, z^{(2)}, t^{(2)})$ in a three point plane. It is important to recall that the $x,y,z$ are related to Einstein measures of distance. Further$,$ we take the $x$-axis as that of 
$F_1F_2$. The $v_3$ is the Einstein measure of the $F_2$ velocity as measured by an inf. light-clock at $F_1$. To correspond to the customary coordinate system employed [1, p. 32], this gives
$$\cases{x^{(1)}=v_1t^{(1)}\cos\theta,\ y^{(1)}=v_1t^{(1)}\sin\theta, \ z^{(1)} = 0&\cr
         x^{(2)}=-v_2t^{(2)}\cos\phi,\ y^{(2)}=v_2t^{(2)}\sin\phi, \ z^{(2)} = 0&\cr}.\eqno (A14)$$\par
It follows from $(A10), \cdots, (A14)$ that 
$$t^{(1)} = \beta_3(t^{(2)} -v_3x^{(2)}/c^2),\ x^{(1)} = \beta_3(x^{(2)} -v_3t^{(2)}),\ y^{(1)} = y^{(2)},\ z^{(1)} = z^{(2)}. \eqno (A15)$$\par

Hence, for this special case, $\omega_1,\ \omega_2, \ \theta, \ \phi$ are eliminated and the Lorentz Transformations are established. If $P \not= F_1,  P\not= F_2$, then the fact that $x^{(1)},\ x^{(2)}$ are not the measures for a physical ruler but are measures for a distance related to Einstein measures$,$ which are defined by the properties of the propagation of electromagnetic radiation and infinitesimal light-clock counts$,$ shows that the notion of actual \underbar{natural} world ``length'' contraction is false. For logical consistency$,$ Einstein measures as determined by the light-clock counts are necessary. This analysis is relative to a ``second'' pulse when light-clock counts are considered. The positions $F_1$ and $F_2$ continue to coincide during the first pulse light-clock count determinations. \par

Infinitesimal light-clock counts allow us to consider a real interval as an interval for ``time'' measure as well as to apply infinitesimal analysis. This is significant when the line-element method in Article 3 is applied to determine alterations in physical behavior. All of the coordinates being considered must be as they would be understood from the Einstein measure viewpoint. The interpretations must always be considered from this viewpoint as well. Finally, the model theoretic error of generalization is eliminated by predicting alterations in clock behavior rather than by the error of inappropriate generalization. \parm 
\centerline{\bf REFERENCES} 
\noindent {\bf Dingle, H.} {\it The Special Theory of Relativity,} Methuem \& Co., London, 1950.\pars
\noindent {\bf Herrmann, R. A.} ``The NSP-world and action-at-a-distance,'' In ``Instantaneous Action at a Distance in Modern Physics: `Pro' and `Contra,' '' Edited by Chubykalo, A.,  N. V. Pope and R. Smirnov-Rueda, (In CONTEMPORARY FUNDAMENTAL PHYSICS - V. V. Dvoeglazov (Ed.)) Nova Science Books and Journals, New York, 1999, pp. 223-235.  Also see http://arxiv.org/abs/math/9903082 pp. 118-119.\pars

\noindent {\bf Kulczycki, S.} {\it Non-Euclidean Geometry,} Pergamon Press, New York, 1961.
\vfil\eject
%The break for web book 2. 
 
\centerline{\bf Foundations and Corrections to Einstein's}
\centerline{{\bf Special and General Theories of Relativity, Article 3.}\footnote*{Partially funded by a grant from the United States Naval Academy 
Research Council.}} 
\parm
Abstract: In Article 3 of this paper, based upon a privileged  
observer 
located within a nonstandard substratum, the infinitesimal chronotopic line-element 
is derived from light-clock properties and shown to 
be related to the propagation of electromagnetic radiation.  A general 
expression is derived,
without the tensor calculus, 
from basic infinitesimal theory 
applied to obvious Galilean measures for  distances traversed by an
electromagnetic pulse. Various line-elements (i.e. ``physical metrics,'' not obtained via tensor analysis) are obtained from this general 
expression. These include the  Schwarzschild (and modified) line-element, which 
is obtained by merely substituting a Newtonian gravitational velocity into 
this expression; the de Sitter and the Robertson-Walker which are obtained by 
substituting a velocity 
associated with the cosmological constant or an expansion (contraction)
process. The relativistic (i.e. 
transverse Doppler), gravitational and cosmological
redshifts, and alterations of the radioactive decay rate are derived 
from a general behavioral model associated with atomic systems, 
 and it is predicted that similar types of shifts will take 
place for other specific cases. Further, the mass alteration expression is 
derived in a similar manner. 
  From these and similar derivations, the locally verified 
predictions of the Special and General Theories of Relativity should be  
obtainable. A process is also given that minimizes the problem of the
``infinities'' associated with such concepts as the Schwarzschild radius.
These ideas are applied to black holes and pseudo-white holes. \parm
\leftline{\bf 1. Some Special Theory Effects}\parm 
 Recall that it does not appear possible to give a detailed 
description for the 
behavior of the NSPPM. For example, Maxwell's equations are based upon 
infinitesimals. Deriving these equations using only the NSP-world language gives but approximate NSP-world information about infinitesimals for they would be expressed 
in terms of $\approx$ and not in terms of $=.$ These facts require that a new approach be used in order change $\approx$ into $=$ within the NSP-world and to determine other properties of the NSPPM,
properties that 
are originally approximate in character and gleaned from observations within the 
natural world$.$  This is the 
view implied by the Patton and Wheeler statements and taken within this research. The view is that space-time geometry is but 
a {convenient language} and actually tells us nothing about the true 
fundamental causes for such behavior. As mentioned and as will be 
demonstrated, this ``geometric'' language description is but an 
analogue model for 
properties associated with electromagnetic radiation. In most cases, {Riemannian geometry} 
will not 
be used for what follows. However, certain Riemannian concepts can still be 
utilized \underbar{if} they are properly 
interpreted in terms of light-clock behavior. 
\par
The so-called {``Minkowski-type line-element''} is usually defined. 
But, using this new approach,
it is derived relative to light propagation behavior. 
As shown in Herrmann (1992) (i.e. article 2), there is in the NSP-world a unit conversion 
$u$ that relates the private fixed absolute time measurements to 
corresponding light-clock measurements. Further, it is shown that $c= \st 
{L/u}= L/u,$  where c is 
 a  {{\it local} measure of the velocity} of electromagnetic radiation {\it in vacuo}.  
  \par
Suppose a timing infinitesimal light-clock  
is at a standard point in the NSPPM and $\Pi_s$ counts have occurred, where a subscript or superscript $s$ denotes ``standard laboratory measurements'' (measurements not considered as affected by the physical processes being considered). A {{\it potential velocity}} is a  
velocity that {\it may be} produced by a physical process. Whether or not 
motion actually occurs depends upon the physical scenario. If the same 
light-clock is moved with a standard relative (potential) velocity $v_E$ 
as infinitesimal light-clock measured in the natural world, then the (8.1) and (8.2) scenario dictates that the infinitesimal light-clock \underbar{intervals}, $u\Pi_s,\ u\Pi_m,$ are related as follows:
 $$\st {u\Pi_s}\gamma = \st {(\gamma u)\Pi_s} = \st {u\Pi_m},\ \gamma = \sqrt {1- v_E^2/c^2}.\eqno (1)$$
[See note [12] before proceeding.]
The numbers $u\Pi_s, \ u\Pi_m, \ \gamma, \ L/u$ being finite in character
allows the standard part operator to be dropped and $=$ to be replaced by 
$\approx.$ This yields the following NSP-world statement. 
$$u\Pi_s\approx u\Pi_m\gamma^{-1} \approx \eqno (2)$$
$${{L\Pi_m}\over{(L/u)\gamma}}\approx {{L\Pi_m}\over{c\gamma}}, \eqno (3)$$
where a subscript or superscript $m$ denotes potential
``velocity'' relative to the $s$-point. It is assumed that $v_E^s$ is generalized to a type of velocity $v+d$ (i.e. has velocity units) and satisfies, in the NSPPM, the Galilean definition for uniform velocity. Then 
$(v+d)^2 = (\Delta r^s)^2)/(\Delta t^{s})^2=((\Delta x^{s})^2 +(\Delta y^{s})^2 + (\Delta z^{s})^2)/(\Delta t^{s})^2$ 
and $\st {u\Pi_s} = \Delta t^{s}$. Combining (2) and (3) 
yields $L\Pi_m \approx c\,\gamma \,u\Pi_s.$ This can be re-written as 
$$(\st {L\Pi_m})^2 = (1 -((v+d)^2/c^2)(\Delta t^{s})^2 c^2 =((\Delta t^{s})^2 c^2 - ((\Delta x^{s})^2 +(\Delta y^{s})^2 +
 \Delta z^{s})^2). \eqno (4)_a$$\par 

$$(\st {u\Pi_m})^2 =(\Delta t^{s})^2  - (1/c^2)((\Delta x^{s})^2 
+(\Delta y^{s})^2 + (\Delta z^{s})^2). \eqno (4)_b$$\par
 The left side of equations (4)$_a$ and (4)$_b$ are only relative to 
electromagnetic properties as being analogue modeled by Einstein measures and 
equivalent infinitesimal light-clocks. The well-known right hand side of (4)$_a$ has been termed the {{\it chronotopic interval}}, a term that indicates its relationship to electromagnetic propagation.  It is important to always keep in mind, that statements such as 
$(4)_a,\ (4)_b$ and the forthcoming statements $(5)_a,\ (5)_b$ refer to the 
use of infinitesimal light-clocks, or an approximating device, to measure time. \par

Although (4)$_a$ is similar to expression (21) in Ives 
(1939), this interpretation is completely distinct from the Ives' assumption 
that $\ell$ is altered within  the N-world by relative motion. 
To measure the velocity of light by means of infinitesimal light-clocks 
and the 
Euclidean length expression, simply consider (4)$_a$ written as 
$0 =(\Delta t^{s})^2 c^2 - 
((\Delta x^{s})^2 +(\Delta y^{s})^2 + (\Delta z^{s})^2).$
Expressions such as (4)$_a$ and (4)$_b$ always incorporate both the length 
and time 
infinitesimal light-clocks  due to the definition of (scalar) velocity. 
\par 
For nonuniform motion and its local effects, one passes (4)$_a$, (4)$_b$ 
to the 
infinitesimal world, where $v+d$  
is considered not a constant but a differentiable function that 
behaves as if it is a constant in the infinitesimal world. Such a re-statement of $(4)_a,\ (4)_b$ does not come from the more formal process of ``infinitesimalizing.''  It is a physical infinitesimal light-clock hypothesis. 
Further, this implies that  $L\Pi_m$  is an infinitesimal 
although $\Pi_m$ can still be an infinite number
and that all other similar finite 
quantities in 
(3) are nonzero infinitesimal numbers representing infinitesimal light-clock 
measures. 
 Writing these infinitesimals 
in the customary form, yields
$$dS^2 =(dt^{s})^2 c^2 - 
((dx^{s})^2 +(dy^{s})^2 + (dz^{s})^2), \eqno (5)_a$$
$$d\tau^2 =(dt^{s})^2  - (1/c^2)
((dx^{s})^2 +(dy^{s})^2 + (dz^{s})^2), \eqno (5)_b$$
where $dZ^s, dZ^m, \ Z \in \{t,x,y,z\}$ are infinitesimal light-clock measures. Notice that the {quasi-time-like (5)$_a$ and 
proper-time-like (5)$_b$} 
are not metrics as these terms are generally understood within mathematics. 
[The standard approach (Bergmann, 1976, p. 44) is to consider real numbers as represented by variables without the $s$ or $m$ superscripts. But, for any real number $r\not= 0$, there exists an infinite integer $\Gamma_L,$ [resp. $\Gamma_u$] such that $L\Gamma_L \approx r,$ [resp. $L\Gamma_u \approx r$]. Obviously, $ds\approx dS.$ For hyperrational numbers, $ds = dS$ and the two $\tau$ are equal. [Also see note 22b.]]  The basic goal is to determine to what the left-hand sides of $(5)_a,\ (5)_b$ correspond. This is done by examining two NSP-world infinitesimal views that compare infinitesimal physical world behavior with NSPPM altered infinitesimal behavior as both are viewed from the NSP-world. \par
Equations $(4)_a, \ (4)_b$ use unaltered infinitesimal light-clocks. The $L, \ u,$ for such clocks are infinitesimalized.   
 For this reason, the right-hand 
side of equations $(4)_a,\ (4)_b$ can be expressed in terms of infinitesimal
concepts. However, the left-hand side can only be considered as near to the 
right hand side. Does this matter?    
The differentials that appear in 
$(5)_a$ and $(5)_b$ represent infinitesimal Einstein time intervals
 and associated distance measures. 
As shown in Herrmann (1985, p. 175), classical differential calculus 
cannot differentiate between types of differentials if either the concept of infinitesimal  
``indistinguishable affects'' or Riemann integration is considered. 
Under these conditions, this gives an additional 
freedom in differential selection. Consequently, as shown in Herrmann (1992, Article 2), each of these  
differentials can be assumed to represent {\it exactly}, rather 
than approximately, the infinitesimal light-clock behavior. 
Expressions 
(5)$_a,$ (5)$_b$ do not refer to the geometry of the universe in which we 
dwell. They and $dS$  refer totally to the restricted concept of 
electromagnetic propagation within an infinitesimal light-clock. These 
expressions  do not reveal what natural world relation might be operative 
unless other considerations are introduced. For this reason,  certain basic 
properties are imposed upon the NSPPM. \par 
As done in Herrmann (1991a p. 170, arxiv p. 162; 1992, 
1994a), the concept of NSPPM {{\it infinitely close of order one}} effects (i.e. indistinguishable for $dt$ effects) is the simplest  and most successful modeling condition to impose 
upon the NSPPM.  In order to investigate what affect $(5)_a$ produces in the 
infinitesimal world over $dt$ of NSPPM infinitesimal changing $t,$  this 
concept says, for $dS,\ dS^m,$ that for each infinitesimal $dt$ there is an infinitesimal $\eps$ such that $dS^m = dS + \eps\, dt.$ This infinitesimal approach allows infinitesimal changes about a point to be extended to a local environment. In this case, the to-and-fro property is observationally 
subdivided into a 
finite collection of ``to''s followed by a finite collection of ``fro''s.
If one subdivides a NSPPM t-interval $\Hyper {[a,b]}$ into 
infinitesimal pieces of  ``size'' $dt$ and considers the $(dS^m)^i$ and $dS^i,$ where the superscript, not exponent, $i$ varies 
over the   
number $\omega$ of these subdivisions, it is not difficult to show using the notion of indistinguishable
effects (Herrmann, 
1991a, p. 87, arxiv p. 60) that\par
 $$\sum_{i=1}^\omega (dS^m)^i = 
\sum_{i=1}^\omega dS^i + \lambda, \eqno (6)$$ 
\noindent where $\lambda$ is an 
infinitesimal. \par 
Suppose that the entities 
$t_E,\ x_E,\ y_E,\ z_E,$ that appear in (4)$_a$, are functions 
expressed in parameter $t.$  For nonuniform relative velocity, these 
functions might be but restrictions of differentiable functions that yield an 
integrable $dS= f(t)dt$ over standard $[a,b].$ Assuming this, the standard part 
operator yields $$\st {\sum_{i=1}^\omega (dS^m)^i} = \st {\sum_{i=1}^\omega dS^i} = \int_a^b f(t)\, 
dt. \eqno (7)$$ 
Equation (7) can be interpreted in the exact manner as is 
(4)$_a$ assuming the conditions imposed upon its derivation. [See important note 21.] Expressions such as (6) and (7) and the properties of the standard part operator imply that when standard methods are used to derive an expression, then the expression can be re-expressed in terms of infinitesimal light-clock counts. 
\parm
\leftline{\bf 2. General Effects}\parm
{Surdin} (1962)  states that it was {Gerber} in 1898 who first 
attempted to adjoin to Newton's theory of gravity a time-varying potential so 
as to explain the additional advance of Mercury's perihelion.              
Is it possible that an 
infinitesimal effect such as $\gamma dt^{m}$ or $\gamma dt^{s}$ 
is a significant 
part of a general NSPPM effect for all natural system behavior? \par
Various investigators (Barnes and Upham, 1976),  assume 
that ``clocks'' 
and ``rods'' are altered in a specific way and select the appropriate $\gamma$ 
expression that will transform (5)$_b$ into a proper time-like Schwarzschild line-element. 
However, from a more general  relation derived solely from some  
simple NSP-world 
assumptions, the Schwarzschild 
relation will be obtained.  
To 
investigate a possible and simple  effect, we improve significantly upon a 
suggestion of Phillips (1922) 
and use a simple 
monadic (i.e. infinitesimal) world behavioral concept. 
Recall that the collection of NSP-world entities infinitely close to a 
(standard) natural world position is called {{\it a monad, monadic 
neighborhood}}, or a {{\it monadic cluster}} when considered as composed of 
various types of {subparticles}.\par
Further, recall that subparticles should never be visualized as  ``particles'' 
(Herrmann, 1986b, p. 50). Indeed, they are often simply called  ``things.''
The {ultimate subparticles}, those with all but two coordinates denoted by 
$\pm 1/10^\delta,\ \delta \ {\rm an\ infinite\ number},$ are combined  
into intermediate subparticles as 
modeled by a well-defined process that mirrors finite 
combinations and is expressed by 
equational 
system 
(2) in Herrmann (1986b, p. 50). Various relations between 
subparticle coordinates determine the types of matter or fields that such 
combinations produce. In what follows, we investigate a simple theoretical
 relation 
between the electromagnetic and velocity coordinates of the subparticles that comprise the NSPPM. \par 
Although for certain behavior an expression that models natural system 
behavior within the N-world (i.e. natural (physical) world) may be considered as an invariant form, the 
General Principle of Relativity is not assumed for our basic line-element 
derivations 
and its 
model tensor analysis is not applied. Without some reasonable physical 
basis, not all smooth curvilinear coordinate transformations need be allowed. 
However, for applications of these line 
elements to 
specific physical problems, certain {invariant forms} and {solution methods} 
will 
be assumed. 
Although infinitesimal world alterations are allowed they represent 
various physical effects and do not correspond to alterations in the geometry 
of space-time, but only to alterations in the behavior of natural world 
entities.\par

Suppose that timing infinitesimal light-clocks are used as 
an analogue model to  
investigate how the NSPPM behavior is related to 
a  
physical or physical-like process denoted by $P.$ The process $P$ influences various infinitesimal 
light-clocks as they are specially oriented. In the first case, we consider the 
``distance'' measuring light-clocks as oriented in a radial and rotational 
direction as compared to  a Euclidean (Cartesian) system. 
The {``timing'' infinitesimal 
light clocks have two non-coordinate  ``increment'' orientations}, nonnegative or 
nonpositive. Refined meanings for the superscripts or subscripts $s$  
and $m$ are discussed in note 12 and the Appendix B page 93. 
\par 
Using the established methods of infinitesimal modeling as in Herrmann 
(1994a),  
suppose that  $v$ and $d$ behave within a monadic neighborhood as if they 
are constant with respect to $P.$ [See note 14.] {\bf Moreover, as a physical principle, since behavior in a monadic neighborhood is a proposed simple behavior, the simple Galilean velocity-distance law and (4)$_a$, (4)$_b$ not just for Einstein measures but for other potential velocities hold.} 
Hence, for photon behavior and a proposed potential velocity  
$$((v+d) + c)dt^s = (v+d)dt^s + cdt^s = $$ $$dR^s + dT^s,\ dR^{s} = (v+d)dt^{s},\ dT^{s} = cdt^{s},\ {\rm and}\ {{dR^s}\over{dT^s}}={{v+d}\over{c}},\ dt^s\not= 0.\eqno (8)$$ 
Suppose that for physical effects, not just for photons,   
a microeffect (Herrmann, 1989) alters the $dR^s$ and $dT^s$ in (8). 
This alteration is 
characterized by a linear transformation
(A): $dR^{s} = (1-\alpha \beta)dR^m -  \alpha dT^{m}$ and (B): $dT^s = \beta dR^m + dT^m.$ This is 
conceived of within a monadic cluster as determining a
subparticle coordinate  relation relative to the $P$-process. 
The $\alpha, \ \beta$ are to be determined. Since the effects 
are to be observed in the natural world, $\alpha, \ \beta$ have 
standard values. [Equations (A) and (B) represent an ``infinitesimal'' linear 
transformation of the infinitesimal light-clocks measurements.] 
\par  
Assuming simple 
NSP-world  behavior with respect to radial  motion,    
transform only 
the length portion of (5)$_a$ into 
spherical coordinates. This yields, not in terms of Einstein  
measures, but infinitesimal light-clock counts at points that is form invariant and infinitely close to $dS^2$, which by choice, is equated to $dS^2$. 
$$dS^2 = (dT^{s})^2 - (dR^{s})^2 - 
(R^{s})^2(\sin^2\theta^{s}(d\phi^{s})^2 +
(d\theta^{s})^2).  \eqno (9)$$ 
In (9), the infinitesimals $d\phi^s,\ d\theta^s$ are assumed to be 
infinitesimal light-clocks for the two
rotational aspects. Such ``clock'' behavior can be viewed as spherical transformed values for Cartesian coordinate infinitesimal light-clock values. \par 
 Consider the radial portion of (9) and let
$k = (dT^{s})^2 - (dR^{s})^2.$ Substituting (A) and (B) into $k$ yields 
$$k = (1-\alpha^2)(dT^{m})^2 + 2(\alpha +\beta(1-
\alpha^2))dR^{m}dT^{m}+$$
$$(\beta^2 -(1-\alpha\beta)^2)(dR^{m})^2. \eqno (10)$$\par
For real world time interval measurements, it is assumed that 
timing counts can be added or subtracted. This is 
transferred to a monadic neighborhood and requires  
$dT^{m}$ to take on two increment orientations represented by nonpositive or 
nonnegative infinitesimal values. As done for space-time, suppose that the 
$P$-process is symmetric with respect to the past and future sense of a time variable. This implies that $dS^2$ is unaltered when $dt^m$ is replaced by 
$-dt^m$ (Lawden, 1982, p. 143). Hence, $k,$ is not altered in infinitesimal value when 
$dT^{m}= cdt^m$ is positive or negative. This implies a transformation restriction that $2(\alpha + \beta(1-
\alpha^2)) = 0.$  
For simplicity of calculation,  let $\alpha = -\sqrt{1 - \eta}.$ Hence, 
$\beta = \sqrt {1 - \eta}/\eta.$ Substituting these expressions into (A) and (B) yields\par
$$dR^{s} = {{1}\over{\eta}}dR^{m} + \sqrt{1-\eta}\, dT^{m}$$
$$dT^{s}= {{\sqrt{1-\eta}}\over{\eta}}dR^{m} + dT^{m}. \eqno (11)$$
\noindent Combining both equations in (11) produces 
$${{dR^{s}}\over{dT^{s}}}=
\left({{1}\over{\eta}}{{dR^{m}}\over{dT^{m}}} +\sqrt{1-\eta}\right)\div
\left({{\sqrt{1-\eta}}\over{\eta}}{{dR^{m}}\over{dT^{m}}}+ 1\right).\eqno (12)$$ 
For this derivation, a static condition is assumed. This is modeled by letting $dR^{m}/dT^m= 0$ over some time interval.  [This time interval can be finite or potentially infinite.] Hence, (12) yields $dR^{s}/dT^{s} = (v+d)/c 
= \sqrt {1- \eta}$ or $\eta = 1-(v+d)^2/c^2=\lambda\not=0$ for a standard neighborhood. 
We note that using $\alpha = \sqrt {1 - \eta}$ yields the contradiction $(v+d)/c < 0$ for the only case considered in this article 
that $0\leq v+d.$ [See note 1.] By substituting $\eta$ into (11) and then (11) into (9), we 
have, where $dT^m = cdt^m,$\vfil\eject
$$dS^2 = \lambda(cdt^{m})^2 - (1/\lambda)(dR^{m})^2 -$$
$$(R^{s})^2(\sin^2\theta^{s}(d\phi^{s})^2 +
(d\theta^{s})^2),\ {\rm or}  \eqno (13)_a$$\par 
$$(d\tau^s)^2 = \lambda(dt^{m})^2 - (1/(c^2\lambda))(dR^{m})^2 -$$ 
$$(R^{s}/c)^2(\sin^2\theta^{s}(d\phi^{s})^2 +
(d\theta^{s})^2).  \eqno (13)_b$$\par 
Under our assumption that instantaneous radial behavior is being investigated, 
then $\theta^{s},\ \phi^{s}$  are not affected directly by the radial properties of 
$P$ and the superscript $(s)$ in the parts of (13)$_a,$ and (13)$_b$ 
containing these angles can be replaced by a superscript $(m)$. Further, 
since the monadic neighborhood is only relative to radial behavior, then an infinitesimal rotation of the monadic neighborhood by $d\theta^s$ should have no effect upon the standard radial measure. Hence, considering all other infinitesimal light-clock counts to be zero, this implies that $R^sd\theta^s = R^md\theta^m = R^md\theta^s.$ Hence, $R^s = R^m$ and this substitution is considered throughout all that follows in this article. [Also see note 12.] [However, if $R$ is not considered as an independent spacetime parameter, then, as will be discussed if it is assumed to vary in ``time,'' it will be considered as a ``universal'' function.] 
I again mention that  
(13)$_a,$  (13)$_b$
do not determine metrics under 
the general mathematical  
meaning for this term. What they do represent is 
restricted to an instantaneous 
 effect relating radial effects of $P$ to electromagnetic propagation 
where a subparticle 
coordinate relation that yields this  is partially 
identified. [See note 2 and for the important case where we consider pure complex $(v+d)i$. Also see note 15.]\par
There are, of course, many conceivable $P$-processes. Suppose 
that there exists such a process that is only related, 
in general, to an objects relative velocity with respect to an observer 
and possibly an objects distance from such an observer.  Not transforming to 
spherical coordinates, but using the  same 
argument used to obtain (13)$_a$ and (13)$_b$ yields the {{\it linear effect} 
line-element} 
$$dS^2 = \lambda(cdt^{m})^2 - (1/\lambda)(dr^{m})^2, \eqno (14)_a$$
$$(d\tau^s)^2 = \lambda(dt^{m})^2 - (1/(c^2\lambda))(dr^{m})^2, \eqno (14)_b$$
where $(dr^m)^2= (dx^m)^2 + (dy^m)^2 +(dz^m)^2,$ where $x^m,\ y^m,\ z^m$ are not functions in $t^m$ (Herrmann, 1995). This linear effect line 
element can be used, among other applications, to determine 
linear effects solely attributable to the Special Theory. For that use, one usually considers $v = v_E, \ d = 0$ and Einstein measures are used. The introduction of two aspects for Special Theory behavior 
will lead to a simplification of $(14)_a,\ (14)_b$. [See note 4,5,6.]\parm

\leftline{\bf 3. Relativistic Alterations}\parm

There are many possible coordinate effects produced by the Lorentz 
transformation. 
In Herrmann (1992)(i.e. Article 2), it is shown that since a NSPPM is used for this new 
derivation that there is no contradictory Einstein Special Theory reciprocal effects.
For example, consider the so-called {time dilation effect}. This effect is a 
``light-clock'' (electromagnetic) effect and has no relation to the concept 
of natural world time dilation. An infinitesimal light-clock is used as a measure of the concept of time. 
This implies that such time is measured by $u\Pi$ with a possible alteration in $c$ and $\Pi$ is an infinite Robinson number. 
\par

Within the NSPPM, relative motion is 
measured with respect to two entities $F_1$ and $F_2.$ With respect to 
relativistic effects,  
it is always the 
case that $F_1$ is selected as denoting the entity that has NSPPM scalar
velocity $\omega_1$ less than or equal to $\omega_2$ for $F_2,$ and 
$\omega_1,\omega_2 \geq 0$. 
NSPPM scalar velocities are {not modeled by directed numbers} and, hence, 
do not follow the same arithmetic as natural world 
scalar velocities. Such velocities are additive but not subtractive. 
Subtraction 
is replace by the (true) metric $d(\omega_1,\omega_2) = \vert \omega_1 
- \omega_2 \vert,$ which represents NSPPM relative velocity. 
 Such effects are {instantaneous  ``snapshot''  effects}. \par
For the non-infinitesimal change point of view, a timing infinitesimal light-clock 
corresponding to $F_1$ is denoted by $u\Pi^{s}.$ The timing infinitesimal light-clock 
that corresponds to motion is relative to 
$F_2$ and is denoted by $u\Pi^{m}.$ If only the Einstein time 
coordinate is affected by some NSPPM process, then there are two possible 
alterations in these infinitesimal light-clock counts due to motion as it can be viewed 
from the NSPPM. Either,
(I) $\gamma u\Pi^{s} = (\gamma u)\Pi^{s} \approx u\Pi^m$ or
(II) $\gamma u\Pi_1^{m} = (\gamma u)\Pi_1^{m}) \approx u\Pi_1^{s},$ but both cannot hold 
unless $v_E = 0.$ Since (I) has lead to (5)$_a$ and (5)$_b$ and relativity notions must be maintained, it is conjectured that when $m$-points are considered that derivations will lead to line elements related to (9.1) and (9.2).   \par 

Once a derivation is obtained, it requires interpretation, 
although the basic electromagnetic properties that lead to a consequence are fixed.  
In all cases, the nature of the NSPPM and 
its associated 
effects are characterized by the  ``time'' measuring infinitesimal 
light-clocks, 
while the ``length'' measuring infinitesimal light-clocks characterize the 
to-and-fro path traversed concept within the NSPPM and its associated effects. As this 
investigation progresses, it will become more evident that {Builder} (1960) 
may be correct in that all physical phenomena are  associated with properties 
of an electromagnetic  field; in this case, the NSPPM. \par 
Laboratory determined                 
 selection of one or more of these possibilities is the major 
modeling technique used by numerous investigators, including {Barnes and 
Upham, Builder, Lorentz,  Dingle,
 but, especially, Herbert Ives,} in their attempts to understand 
relativistic properties, where the velocities considered are the Einstein 
relative velocities $v_E.$  
This same accepted method 
is used in Herrmann (1992) to 
explain the  Michelson-Morley and Kennedy-Thorndike type 
experiments. However, for the Ives-Stillwell type experiment (Ives and 
Stillwell, 1938), it is argued from 
empirical evidence that 
an alteration of the unit of ``time'' is not the correct interpretation, but, 
rather, that such experiments display an alteration of the frequency of 
radiation associated with an 
atomic system, as it corresponds to the timing 
infinitesimal 
light-clock, 
 due to an {{\it electromagnetic interaction with the 
NSPPM}} (emis). For special relativity, (emis) refers to hyperbolic velocity behavior within the NSPPM (p. 50, (A3)), and it  
implies that this intimate relation exists. 
 Thus  
 certain details as to how the NSPPM's behavior influences these and similar 
types of experimental scenarios is indirectly known. However, rather than 
simply postulating the correct infinitesimal relation, the correct relation 
will be predicted from well-established 
atomic system behavior.\par

Suppose that certain aspects of a natural system's behavior are 
governed by a function
$T(x_1,x_2,\ldots,x_n,t)$ 
that satisfies an 
expression $D(T) = k(\partial T/\partial t),$ where $D$ is a {(functional) 
t-separating operator and $k$ is a universal constant.} 
Such an expression is actually saying something about the 
infinitesimal world. With respect to electromagnetic effects, 
the stated variables are replaced by variables with superscripts $s$ if 
referred to $F_1$ behavior where $s$ now indicates that no (emis) 
modifications occur. In what follows, all measures are Einstein 
measures.\par
In solving such expressions, 
the function $T$ is often considered as 
separable and $D$ is not related to the coordinate $t.$ In this case,
let $T(x_1,x_2,\ldots,x_n,t) = 
h(x_1,x_2,\ldots,x_n)f(t).$ 
Then  $D(T(x_1,x_2,\ldots,x_n,t)) = D(h(x_1,x_2,\ldots,x_n))f(t) = 
(kh(x_1,x_2,\ldots,x_n))(df/dt)$ and is an invariant separated 
form.  
\par
Let $(x_1^{s},x_2^{s},\ldots,x_n^{s},t^{s})$ correspond to measurements 
taken of the behavior of a natural system that is influenced by (9) [resp. 
(5)$_a$] and using 
identical modes of measurement let 
$(x_1^{m},x_2^{m},\ldots,x_n^{m},t^{m})$ correspond to measurements 
taken of the behavior of a natural system that is influenced by (13)$_a$ 
[resp. (14)$_a$]. [Notice that this uses the language of ``measurements'' and 
not that of transformations. The term ``modes'' means that identically 
constructed devices are used.] Now suppose 
that   $T(x_1^{s},x_2^{s},\ldots,x_n^{s},t^{s}) = 
h(x_1^{s},x_2^{s},\ldots,x_n^{s})f(t^{s}).$ We assume that $T$ is a {universal 
function} and that separation is an invariant procedure. What this means is 
that the same solution method holds throughout the universe and any 
alterations in
the
measured quantities preserves the functional form; in this 
case, preserves the separated functions. [See note 14.]
Let the values
$h(x_1^{s},x_2^{s},\ldots,x_n^{s})= H(x_1^{m},x_2^{m},\ldots,x_n^{m})$ and the values $f(t^s) = F(t^m)$ and $T(x_1^{m},x_2^{m},\ldots,x_n^{m},t^m) =
H(x_1^{m},x_2^{m},\ldots,x_n^{m})F(t^m).$ 
One differentiates with respect to $t^s$ and obtains by use of the chain rule
$$ 
\delta^{s}=\left({{D_s(h(x_1^{s},x_2^{s},\ldots,x_n^{s}))}\over{h(x_1^{s},x_2^{s},\ldots,x_n^{s})}}\right) 
= k{{1}\over{f(t^{s})}}{{df}\over{dt^{s}}}=
k{{1}\over{F(t^{m})}}{{dF}\over{dt^{m}}}{{dt^{m}}\over{dt^{s}}}.\eqno (15)$$  
With respect to $m,$ 
$$\left({{D_m(H(x_1^{m},x_2^{m},\ldots,x_n^{m}))}\over{H(x_1^{m},x_2^{m},
\ldots,x_n^{m})}}\right) 
= k{{1}\over{F(t^{m})}}{{dF}\over{dt^{m}}} 
=\delta^{m}.\eqno (16)$$ \par 
First, consider physical structure. With 
respect to the NSPPM, consider the (emis) effects of this 
$P$-process
caused by the NSPPM velocity $\omega.$  This is a physical process of the nonsigned 
relative velocity. Further, this (emis) effect is 
considered as occurring within the physical structure itself and, due to the use of infinitesimals, it is the 
general 
practice to assume the modeling concept that when such physical alterations occur 
the structure is
momentarily at rest with respect to both the observer 
and its immediate environment. This is modeled with respect to the 
linear effect line-element by letting $dr^m=0$ in (14)$_a,$  
$v =v_E,\ d = 0,$ and  $dr^s = 0,$ in (5)$_a.$  Comparing the resulting invariant $dS^2$
yields for this 
case, that this Special Theory  $P$-process
requires via application of (15) and (16) that 
$$\sqrt \lambda\, dt^m =  \gamma\, dt^m =
dt^s,\eqno (**)$$ 
which is all that is needed for such relativistic effects. [See note 22c.] The assumption that $d = 0$ is taken to mean that we are either interested only in local 
effects where the $d$ effect would be removed from the problem or effects 
where the $d$ is exceeding small in character. [See note 3.]
Notice that this is  one of the many possible Special Theory coordinate alterations 
and  establishes that the eigenvalues are related to  
physical NSPPM 
properties as they are measured by infinitesimal light-clocks. Also notice that from the infinitesimal viewpoint (**) is similar to (9.2) [Note: 21b]. Finally, 
assume that for a physical structure that a $P$-process is modeled by the above 
$D$ operator equation. Then
 $$\delta^s=\delta^{m}/\gamma. \eqno (17)$$ 
Consequently,  
$\gamma\, \delta^s =  \delta^{m}.$ \par
Suppose that $T = \Psi$ is the total wave function, $D$ for (17) is the 
operator $\nabla^2 - p,$ where $n = 3$ and the constant $k,$ function 
$p$ are those 
associated with the classical {time-dependent Schr\"odinger equation} 
for an atomic system.  It is not assumed, as yet, that such 
a Schr\"odinger type equation  predicts any other behavior except that it 
reasonably approximates the energy associated with electromagnetic radiation 
and that the frequency of such radiation may   
be obtained, at least approximately, from this predicted energy variation.
The eigenvalues (Pohl, 1967, p. 31)  
for this separable solution correlate to energies $E^s$ 
and $E^{m}$ for such a radiating atomic system. Hence, 
$$\gamma \Delta E^s =\Delta E^{m}. \eqno (18)$$
Since comparisons as viewed from one location are used, divide (18) by Planck's constant  and the comparative relativistic {redshift 
(transverse Doppler, where $\alpha = \pi/2$) result}
$\gamma\nu^s=  \nu^{m}, d = 0,$ is predicted. [Note that, as mentioned, 
pure Special 
relativistic effects would, usually, have the effects of any $d$ removed from 
consideration. However, this may not be the case, in general, for all such 
effects.]  This is the same expression, where $\nu^s$ means a stationary (laboratory determined value), first 
verified by Ives and Stillwell and which is attributed to observer time 
dilation. But 
it has been established that it can be interpreted as an (emis) effect and 
not an effect produced by absolute time dilation. Observer time is dilated via alterations in the machines that ``measure'' time.\par

The above Schr\"odinger equation approach does not just apply to 
atomic and molecular physical processes that exhibit uniform 
frequencies associated with such energy changes. Indeed, for consistency, all time related behavior must undergo similar alterations. This implies that the Schr\"odinger equation approach is universal. The same alterations are produced by gravitational fields.
If needed to verify this conjecture, as demonstrated shortly, the (15), (16) method may need to be modified. Again these would be an electromagnetic or (emis) effects. [See note 22.]\par

We next apply the linear effect line-element to the problem of {radioactive 
and similar decay rates.} The usual arguments for the alteration of such rates are in 
logical error. 
Let $N(t^s)$ denote a measure for the number of active entities at the 
light-clock count time $t^s$ and $\tau_s$ be the 
(mean) lifetime. These measures are taken within a laboratory and are used as 
the standard measures. This is equivalent to saying that they are, from the 
laboratory viewpoint, not affected by relativistic alterations.    
The basic statement is that there exists some $\tau \in (0,B]$ such that (*) $(-\tau)dN/dt = N.$ Even though the number of active entities is a natural number, this 
expression can only have meaning if $N$ is differentiable 
on some time interval. But, since the $\tau$ are averages and the number of entities is usually vary large, then such a differential function is a satisfactory approximation.  
Recall that the required operator expression is $$D(T) = k (\partial/\partial t)(T).\eqno (19)$$ \par 
Let $k = 1$ and $h(r) = 0\cdot r^2 + 1 = 1.$ Then define $T(r,t) = h(r)N(t) = (0\cdot r^2 + 1){N}(t),$ where $r^2 = x^2 + y^2 + z^2,$ and let $D$ be the identity map $I$ on $T(r,t).$ Then $D(T(r,t))= D(h(r))N(t) = D(0\cdot r^2 +1)N(t) = 1\cdot N(t)$ and, in this form, $D$ is considered as only applying to $h$ and it has no effect on $N(t).$ In this required form, first let $r= r^s$ and $t = t^s.$  Then, consider $T(r^m, t^m) = H(r^m)\overline{{N}}(t^m), \ H(r^m) = 0\cdot (r^m)^2 + 1 = 1.$ (Notice that it is not necessary to explicitly define $h$ and  $H$ when one assumes the such a $T$ is a universal function, since the $h$ and $H$ are factored from the final result. One simply assumes that there are  functions $h$ and  $H$ such that $h(r^s) = H(r^m).$) In order to determine whether there is a change in the $\tau_s$, one considers the value $N(t^s) = \overline{N}(t^m).$ This yields the final requirement for $T$. Notice that $t^m$ is Einstein time as measured from the $s$-point. This is necessary in that the $v = v_E$, which is a necessary requirement in order to maintain the hyperbolic-velocity space behavior of $v$. \par 
Applying (19) to $T$ and considering a corresponding differentiable equation (*) and the chain rule, one obtains  that there exist a real number $\tau_s$ such that  
$${N}(t^s) = (-\tau_s)(d/dt^s){N}(t^s)=$$ $$(-\tau_s) (d/dt^m)\overline{{N}}(t^m) (d/dt^s)(t^m)=$$ $$ (-\tau_s/\gamma)(d/dt^m)\overline{{N}}(t^m).\eqno (20)$$
\noindent And, with respect to $m$, and for $\tau_m$ 
$$\overline{{N}}(t^m)= (-\tau_m)(d/dt^m)\overline{{N}}(t^m). \eqno 21)$$\par
Using (20) and (21) one obtains that ${\tau_m} = \tau_s/\gamma.$ (In the linear effect line element $d = 0.$) This is one of the well-known expressions for the prediction for the alteration of the decay rates due to relative velocity (that is $v_E$).  The $\tau_s$ can always be taken as measured at rest in the laboratory since the relative velocity of the active entities is determined by experimental equipment that is at rest in the laboratory. 
\par

As another example of the previous procedures, we consider the so-called
Special Theory {mass alteration expression.}
Consider two perfectly elastic objects of mass M and moving in opposite directions with the same velocity and colliding. Let this occur at both the s-position and m-positions. Then at the moment they collide they are momentary at rest. Thus, at that moment, $dr^s = 0 = dr^m$ expression (**) holds. Consider one of these colliding objects. For a Hamilton characteristic function
$S^\prime,$ the classical Hamilton-Jacobi equation becomes $(\partial 
S^\prime/\partial r)^2 = -
2M(\partial S^\prime/\partial t).$  Suppose that $S^\prime(r,t)
= h(r)f(t).$ Again consider a $P$-process that yields this 
isotropic behavior and that  $S^\prime$ is universal in character and the 
solution method holds throughout our universe. This yields 
that  
$h(r^s) =H(r^m),\ f(t^s) = F(t^m), \ S^\prime(r^s,t^s) = 
S^\prime(r^m,t^m).$ Let $D= (\partial(\cdot)/\partial r)^2.$ 
The same procedure used previously yields
$$ 
\left({{\partial h(r^s)}\over{\partial 
r^s}}\right)^2\left({{1}\over{h(r^s)}}\right) 
= -2{{M^s}\over{f^2(t^{s})}}{{df}\over{dt^{s}}}= $$ $$
-
2{{M^s}\over{F^2(t^{m})}}{{dF}\over{dt^{m}}}{{dt^{m}}\over{dt^{s}}}
=M^s\lambda^m/\gamma.\eqno (22)$$  
With respect to $m,$ 
$$\left({{\partial H(r^m)}\over{\partial 
r^m}}\right)^2\left({{1}\over{H(r^m)}}\right) 
= -2{{M^m}\over{F^2(t^{m})}}{{dF}\over{dt^{m}}}=M^m\lambda^m.\eqno (23)$$ 
In (22) and (23), the quantities $M^s$ and $M^m$ are obtained by means of  
identical modes of measurement that 
characterizes   
``mass.'' Assuming that the two separated forms on the left of (22) and (23) 
are invariant, leads to the Special Theory mass expression 
$M^m = (1/\gamma)M^s.$ This result is postulated to hold in general for  identical objects stationary at the s-point and m-point. This result indirectly demonstrates an actual cause 
for the so-called rest mass alteration. It indicates 
the possible existence of a $P$-process that is produced by  
the NSPPM and yields an alteration in the mass effect which is either 
electromagnetic in nature or, at the least, an (emis) effect.  
 \parm
\leftline{\bf 4. Gravitational Alterations}\parm
The results derived in the previous section are all relative to Special  
Theory alterations produced by  (emis) effects and these are not associated 
with a {Newtonian 
gravitational potential.} 
Suppose, however,  that a varying Newtonian 
gravitational potential additionally influences electromagnetic behavior and 
that this potential is determined by a $P$-process. 
An entity of mass $M,$ as analogue modeled by a homogeneous spherical 
object of radius $R_0 \leq R^m,$ has an 
instantaneous Newtonian potential 
at a distance $R\geq R^m,$ where we are not concerned with the question of 
whether this is an action-at-a-distance or a field propagation effect. 
Let the 
general potential difference expression for a distance $R$ and a mass $m_0$ be 
$$U(R) = {{GMm_0}\over{R^m}} - {{GMm_0}\over{R}}, \eqno (24)$$
where the minute potential due to all other matter in the universe is 
omitted.\par
By the usual techniques of 
nonstandard analysis, (24) has 
meaning in the NSP-world at points where $R= \Lambda$ is an infinite number. 
In this case, 
$$U(\Lambda) \approx {{GMm_0}\over{R^m}}.\eqno (25)$$ \par
\noindent Of course, depending upon the observer's viewpoint, you can consider
$U(\Lambda) \approx 0$ and $U(R^m) = -GMm_0/R^m.$\par 
Viewed as escape velocity, $\st {U(\Lambda)} = GMm_0/R^m$ indicates the 
potential energy associated with escaping totally from the natural world to 
specific points within the NSP-world. Since the NSP-world   ``size'' of the 
our universe at the present epoch, whether finite or not, is not known, such a radius $\Lambda$ 
from the massive body 
is 
necessary in order to characterize a total escape. On the other hand,  
$\st {U(\Lambda)} = GMm_0/R^m$ can be viewed as the potential energy associated 
with a ``potential velocity'' that is attained at $R^m$ 
when finite mass  $m_0$ is    ``moved''
from such a specific point within the NSP-world.
 Suppose that this potential energy is characterized 
as kinetic energy. Then the potential velocity is $v_p = 
\sqrt{2GM/R^m}.$ In order to incorporate an additional gravitational (emis) effect, 
substitute  
$v_p =v_E$ into (13)$_a$ rather than the possibility that $v_p = \omega,$ where 
$\omega$ is 
the NSPPM 
velocity.  
 This expression holds even if we assume infinitesimal masses. 
The use of 
infinitesimal light-clocks  allows us to apply the calculus. Further, although it need only be considered as a 
modeling technique, the concept of viewing the behavior of our cosmos from 
a single external 
position within the NSP-world appears relevant to various cosmologies. 
\par
As is customary, for this single object 
derivation, (9) applies at an ``infinite'' distance from the homogeneous
object. But when astronomical and atomic distances are compared, then (9)
 can 
be assumed to apply approximately to many observers within the universe. 
This is especially the case if an observer is affected by a second much weaker 
Newtonian potential, in which case  (9) is used as a local 
line-element relating  measures of laboratory standards.\par
Following the usual practice for radiation purposes, 
the representative  atomic 
system, as well as other physical systems, is considered as momentarily at rest with respect to the spherical object and 
the observer. Hence
$dR^m= d\phi^m= d\theta^m = dR^s= d\phi^s= d\theta^s =0.$ Equation    
(13)$_a$ yields  
that $dS^2 = \lambda (cdt^m)^2$ and (9) yields $dS^2= 
(cdt^s)^2.$ Thus, for this atomic system case, the differentials $dt^m$ 
and $dt^s$ are again related by the expression  
$\gamma dt^m=dt^s,$                   
where $\gamma=\sqrt \lambda.$\par 
Using the same operator as in the relativistic redshift case, yields
the basic gravitational redshift expression $\gamma \nu^s = \nu^m$ (Bergmann, 
1976, p. 222, where $d = 0$).  
 Using the General Theory of Relativity, this same 
expression is obtained for what is termed a  ``weak'' gravitational field 
that can be approximated by a Newtonian potential although it is often 
applied to strong fields and is derived using time dilation. 
The Schr\"odinger type equation   
derivation is a different approach and holds for Newtonian 
potentials in general. Indeed, even if radiation is not the immediate 
product of atomic emission, it may be assumed to be controlled by the 
Schr\"odinger equation if there is any energy alteration. Hence, 
this approach may be applied anywhere within a gravitational field.\par
 Also repeating the derivations in the previous section, we have 
the predicted {gravitational alterations in the radioactive decay} 
rates, and atomic clocks, etc. This will yield laboratory verified
variations that are attributed to time dilation, but in this theory they are 
all attributed to electromagnetic or (emis) effects.
\par
 For the case where
$2GM/c^2 < R^m,$ when one substitutes $v=v_p = 
\sqrt{2GM/R^{m}}$ into 
(13)$_b$ 
[resp. (13)$_a$], one obtains, for $d = 0,$   
the so-called proper-time-like {Schwarzschild line-element} (13)$_{bp}$ [resp. 
(13)$_{ap}$].     
Using the Schwarzschild relation (13)$_{bp}$, many physical predictions have 
been made. These include
the advance of the perihelion of Mercury 
and the deflection of a light ray by a massive body are predicted (Lawden, 
1982, pp. 147--152),  where the language of geodesics does not refer to 
space-time curvature but rather to a $P$ effect. That is Riemannian geometry 
is but an analogue model for the behavior of infinitesimal light-clocks.  
A relation such as 
(13)$_{bp}$ should only be 
applied to 
a Newtonian potential for which other such potentials can be neglected. Thus 
such predictions would be restricted to special physical scenarios. To see 
that this identified (emis) is an intimate subparticle coordinate relation 
and not a cause and effect concept, in one 
scenario the Newtonian potential appears to alter NSPPM behavior, but for 
another scenario properties of the NSPPM appear to alter 
the Newtonian potential. {\bf It is important to note that other 
(gravitational) line-elements can be obtained for non-homogeneous bodies if 
one is able to find an appropriate  ``potential'' velocity  $v_p$ 
for such bodies.} For nonzero values of $d$, one obtains 
what I term the {{\it quasi-Schwarzschild}} line-element.\par 
In general when the Newtonian potential velocity is used and $d$ is constant, 
this is not 
exactly the same as the modified Schwarzschild line-element that contains the 
Einstein cosmological constant $\Lambda$ (Rindler, 1977, p. 184). There is one
additional term in the radial coefficients of the quasi-Schwarzschild line 
element. 
Also, $d = 0$ if and only if $\Lambda =0$ and, more 
generally, $d \geq 0$ if and only if $\Lambda \geq 0.$ However, there is a 
$d=f(R^m)$ that will yield the {{\it modified Schwarzschild}} line-element 
$(13)_a$ where $\lambda = 1 - 2GM/(R^mc^2) - (1/3)\Lambda (R^m)^2/c^2$ (Rindler, 
1977, p. 184, Eq. 8.151), (the $\Lambda$ unit is (time)$^{-2}$). Further note that when one puts $M=0$ and uses this 
modified Schwarzschild, the {de Sitter line-element} is obtained (Rindler, 1977, 
p. 184, eq. 8.155). Thus, major line-elements are obtained by this 
approach.\par
One of approximations discovered for the Einstein gravitational field 
equations is when gravitational potentials are used and $v_p^2 \ll c^2.$  
The coefficients $g_{ij}$ are determined by considering the system of bodies 
that produce the field to be at a great distance from the point being 
considered. From the viewpoint of the Newtonian potential, this would yield an 
approximating parallel force field relative to the center of mass for the 
system. Hence, consider substituting into the 
linear effect line-element $(14)_a$, $v^2 = 2GM/r^m,\ d = 0.$
 From this, one obtains
$$dS^2= (1- 2GM/(r^mc^2))(cdt^m)^2 - \left[{{1}\over{(1- 2GM/(r^mc^2)}}
\right](dr^m)^2.\eqno 
(26)$$
But, letting $(dr^m)^2= dx_m^2 + dy_m^2 + dz_m^2,$ then (26) would be the proper 
line-element for this approximation.\par    
What happens when the coefficient $[1/((1- 2GM/(r^mc^2))]$ is 
approximated? Expanding, we have that $1/((1- 2GM/(r^mc^2))= 1 + 
2GM/(r^mc^2) + \cdots.$ Substituting this approximation into (26) yields
$$dS^2\approx (1- 2GM/(r^mc^2))(cdt^m)^2 - (1+ 2GM/(r^mc^2))(dr^m)^2.\eqno 
(27)$$ 
Equation (27) is the customary and well-known Newtonian first approximation 
associated with the Einstein law of gravity. It is used by many authors for 
various purposes.\par       
Notwithstanding the methods used to model 
Einstein's General Principle of Relativity that allow 
for numerous coordinate transformations to be incorporated into the theory 
and for certain so-called ``singularities'' to be considered as but 
coordinate system anomalies, there would be a natural world Newtonian bound 
$0\leq (2GM)/c^2 \leq R^m$ when the 
operator expression  $D(T) = k (\partial T/\partial t)$ is used to investigate
such things as the frequency change $\gamma \nu^s =\nu^m,$ and gravitational timing device alterations.   
The bound would be retained, at this epoch, for entities  
 to which the Schr\"odinger equation applies. 
On the other hand, unless for other reasons such a bound 
is shown to be necessary, 
then an alteration in atomic system behavior can 
certainly 
be 
incorporated into theorized mechanisms if one assumes that it is 
possible in the natural world for $(2GM)/c^2 > R^m.$  
\par
It is also interesting to note that $D$ can be replaced by the 
{Laplacian.} The same analysis would yield that if only NSPPM behavior is considered, then gravitational effects would alter 
the nature of internal heat transfer. General Relativity also predicts 
certain temperature shifts (Misner, 1973, p. 568).\par 
In this section, the use of {ad hoc coordinate transformations is not 
allowed.} 
Thus it is necessary that $\lambda \not=0$ in (13)$_a$ and (13)$_b.$ From 
continuity considerations this, at present,  forces 
upon us for Newtonian potentials the slightly move restrictive bound 
$0\leq 2GM/c^2 < R^m.$  
 Using the proper initial condition $U(\Lambda) \approx 0$ and 
the light-clock interpretation does not lead to a theoretical impasse as 
suggested  by Lawden (1982, p. 156).      
%See end of paper for stuff on black hole transformation which, at the %present 
%time I don't believe is the proper way to go. 
Relative to action-at-a-distance, where natural world alterations  of 
Newtonian potentials are to take place instantaneously with respect 
to a natural world time 
frame, there may be  NSP-world informational entities, the 
hyperfast subparticles, 
that mediate these instantaneous changes (Herrmann, 1986b, p. 51).  Today, the concept
of action-at-a-distance is still accepted, partially or wholly, as the only 
known description that accurately predicts certain natural system behavior
(Graneau, 1990. Herrmann, 1999). 
On the 
other hand, since Maxwell's field equations are stated in the language of 
infinitesimals and reveal basic behavior within a monadic cluster, then it 
is not ad hoc to speculate that gravitational alterations, that in the 
gravostatic case yield the Newtonian potential,  
are controlled by 
a similar set of field expressions and that gravitational effects are 
propagated with a finite velocity. 
This speculation has been advanced by various researchers as far back as 1893
(Heaviside, 1922).  There are, as well,  more recent advocates of this approach 
(Barnes and Upham, 1976; Barnes, 1983; Jefferson, 1986; Brietner, 1986). \parm
\leftline{\bf 5. NSPPM Analysis}\parm
It may be assumed that numerous distinct processes contribute to the 
present epoch behavior of the cosmos. 
 One  possibility 
is the {metamorphic (i.e. sudden) structured appearance} of our universe or a 
natural system contained within it at various levels of complexity. 
(Schneider, 1984; Herrmann, 1990, p. 13.) The rational existence of  
 {ultimate ultrawords} and a specific NSP-world process, an {ultralogic,}
that combines these processes together rationally and yields
the present epoch behavior has been established (Herrmann, 1991b, p. 103).
[See (Herrmann, 1986a, p. 194) where such entities are described using older 
terminology.]  It might be argued that one aspect of present day cosmological
behavior superimposed upon all of the other aspects  is a smoothed out galactic gas aspect.  
A {cosmological redshift} verification for an ``expanding'' 
universe model might not actually be attainable by means of a redshift  
measurement (Ettari 1988, 1989). However, the superimposition of a NSPPM  
physical effect can model such an expansion as it occurs throughout the development of our universe.  
 Of course, this need not be extrapolated backwards to the beginning of,  say, the ``Big Bang.'' \par 
We solve this problem by the exact same method used in section 3 to obtain all 
of the previous line-elements. {\it The assumption is that there is a type of 
expansion of the subparticle field taking place that directly influences the general 
behavior of our universe and this affect is  superimposed upon all other 
behavior.}  For the derivation of this new line-element, consider $v + d 
 = (R^s/a), \ v=0,\ \ a\not=0, \ dR^s = (R^s/a)dt^s.$ The exact same derivation yields the line-element
$$ dS^2_e = (1-(R^m)^2/c^2a^2)(cdt^m)^2- {{(dR^m)^2}\over{1- (R^m)^2/(ca)^2}} -$$
$$ (R^m)^2(\sin^2 \theta^m(d\phi^m)^2 + (d\theta^m)^2). \eqno (28)$$   
What is needed is to relate measures within the NSPPM to physical 
alterations within the natural world. From all of our previous derivations for 
the natural world alterations in physical behavior produced by the NSPPM, 
it was argued in each case that $(dt^s)^2 = \lambda (dt^m)^2.$ Thus consider 
substituting into the line-element that relation and obtaining 
$$dS^2_e = (cdt^s)^2 -  {{(dR^m)^2}\over{1- (R^m)^2/(ca)^2}} - $$
$$
(R^m)^2(\sin^2 \theta^m(d\phi^m)^2 + (d\theta^m)^2)= (cdt^s)^2 - d\ell_e^2.
 \eqno (29)$$   
Equation (29) is one form of the famous {Robertson-Walker line-element,} when it 
is assumed that $R^m$ is only ``time'' dependent so that it can be used for 
the concept of an isotropic and homogeneous universe (Ohanian, 1976, p. 392). 
It can also be used to 
satisfy the {Copernican principle (i.e. the Cosmological principle),} in this 
case, if one 
wishes to apply it at \underbar{any}  position. Of course, this general application is \underbar{not} 
necessary. For that reason, it is not assumed, in general, that such
NSPPM behavior as it affects the natural world only depends upon $t^s.$   
Functions such as $R^m$ are considered as possibly dependent upon a fixed 
frame of reference within the NSPPM as it is reflected in the 
non-gravitational and non-expanding natural world by the coordinates 
$t^s, x^s, y^s, z^s.$ To correlate to {Ohanian}, the units are changed to  
to light units (i.e. c = 1). In comparing (29) to the Einsteinian approach, it is the 
line-element for the {positive curvature Friedmann model} with the Riemannian curvature  
 $K = 1/a^2$ (Ohanian, 1976, p. 391-2),where $a$ behaves like the ``radius'' of a four-dimensional hypersphere. In the analysis as it appears in Lawden (1982), not using light units, $a = S/c.$ The definition of $d$ clearly indicates that the expanding NSPPM ``velocity'' forces this hypersphere behavior upon our universe. (I mention that {Fock} (1959, p. 173, 371) and 
others who require 
$\Lambda = 0$,  consider only the 
Friedmann models as viable scenarios.)  
\par
The selection of the above $d$ is not ad hoc but is consistent with astronomical 
observation.  The observation is the apparent expansion redshift as 
modeled by the {Hubble Law.} This law states that the {velocity of the 
expansion (or contraction)} is $v(t^s) = 
H(t^s)\ell_e(t^s),$ where $\ell_e(t^s)$ is a distance to 
the center of a spherical mass. 
The influence this 
expansion has within the natural world is a direct reduction (or 
addition) to the Newtonian gravitational potential. Assuming Newtonian 
potentials and Second 
Law of motion, Ohanian (1976, pp. 
370-373)  derives easily the usual relation between the {{\it deceleration 
parameter}} $q = -(1 +(1/H^2)dH/dt^s)$ and its relation to a uniform, but 
changing in $t^s,$ density $\rho.$  This relation is $-qH^2 = -(4\pi G 
\rho(t^s))/3.$ Moreover, by 
letting $dt^s = ad\eta,$ one obtains the usual differential equation for the 
 {Friedmann closed universe model,} $3((da/d\eta)^2 + a^2)=(4GM/\pi)a,$ 
where the mass of the universe 
$M= 2\pi^2a^3\rho(t^s).$ If we let $d= (R^s/a)i,$ then the {Friedmann open 
universe model} is obtained. These Friedmann models are, therefore, based
entirely upon the Newtonian potential concept and not upon Einstein's theory. 
Further, the Hubble function $H(t^s) = 
(1/a)(da/dt^s).$ These results are all consistent with the $d$ selected. 
[See note 7.] The concept of 
an expanding NSPPM will be utilized explicitly for other purposes in the 
last section of this article.   \par
Relative to gravitational collapse, {\bf unless specific physical processes are 
first introduced}, the methods discussed in this section {preclude collapse 
through  
the Schwarzschild surface} (i.e. radius). This yields but a {restricted 
collapse.} However, if such a restricted 
collapse occurs, 
then the {optical appearance} to an exterior observer would be as it is usually 
described (Misner 1973, p. 847), but no actual black hole would be formed. Any 
observational black hole effects, would mostly come from the gravitational 
effects 
such an object might have upon neighboring objects. Black hole internal 
effects, such as the retention of photons and material particles generated 
within the body as well as those that pass over the Schwarzschild surface,  
depend upon a strong gravitational field. In this research, this 
would depend upon {\it physical concepts} that lead to a coordinate 
transformation, not, in general, conversely. For this reason,
there would not exist, as yet, the concept of the 
 {Einstein--Rosen bridge (wormholes)} (Misner 1973, p. 837) and the like which 
tend to be associated with regular coordinate transformations (Rindler, 1977, 
p. 320.)\par 
In the next section, a method will be investigated that will allow certain 
physical transformations to occur in such a way that the so-called 
Schwarzschild radius or singularity $(2GM)/c^2=R^m$ can essentially be 
bypassed.  \parm
\leftline{\bf 6. Minimizing Singularities}\parm
In section 2, the quasi-Schwarzschild line-element $dS^2$ is derived  
by considering Newtonian gravity as being a $P$-process 
considered as emanating from the center of a 
homogeneous spherical configuration    
and its interaction with the 
NSPPM. Adjoined to this gravitational effect was a possible NSPPM 
expansion (contraction) effect. The (emis) interaction is modeled by taking 
the Special Theory 
chronotopic interval and modifying its spherical coordinate transformation
by a type of damping of the basic  infinitesimal 
light-clocks. This damping is characterized by an infinitesimal
nonsingular linear transformation of the infinitesimal light-clock mechanism.
[Note that the length determining infinitesimal 
light-clock is the same as the ``time'' determining light-clock with the 
exception of a different unit of measure. Thus, in all of our cases, there is a 
close relationship between these measuring devices.]\par 
In General Relativity, appropriate {differentiable coordinate 
transformations with nonvanishing Jocabians} may be applied to a solution of 
the Einstein law of gravity and, with respect to the same physical 
constraints, this would yield, at least, a different view of the gravitational field. In this and the next section, this notion is assumed.\par
 The {Eddington-Finkelstein} transformation was one of the first   
of the many purposed transformations. 
But, in (Lawden, 1982),
the derivation and argument for using the this simple 
transformation 
(**) $dU^m = dt^m + f_M(R^m)dR^m$ to obtain a black hole line-element appears 
to be  flawed. 
The apparent flaw is caused by the usual ad hoc logical 
errors in  ``removing infinities.'' Equation 57.11 in (Lawden (1982), p. 157), 
specifically requires that $R^m> 2GM/c^2.$ However, in arguing for the use of 
the transformed Schwarzschild line-element, Lawden assumes that it is 
possible for 
$R^m\leq 2GM/c^2.$ But the assumed real valued function defined by equation 
57.11 is not 
defined for $R^m$ such that $R^m\leq 2GM/c^2.$  
 Hence, new and  
rigorously correct procedure for such transformations might be useful. 
Such a procedure is accomplished by showing 
that (**) can be considered as a hypercontinuous and hypersmooth 
transformation associated with a new $P$-process that 
yields an alteration to the gravitational  field in the vicinity of the 
Schwarzschild surface  during the process of gravitational collapse. \par This 
 {\it speculation} 
is modeled by the expression (**) which is conceived of as an alteration in the 
time measuring light-clock. [See note 13.] Further, this alteration is conceptually the same 
as the ultrasmooth microeffects model for fractal behavior (Herrmann, 1989) 
and thus has a similar physical bases. This 
transformation takes the Schwarzschild line-element, which  applies only 
to the case where $R^m> 
2GM/c^2,$ and yields an NSP-world 
 black hole line-element that only applies for the case 
where $R^m\leq 2GM/c^2.$ Like ultrasmooth microeffects, the nonstandard 
transformation process is considered as an ideal model of behavior that 
approximates the actual natural world process.  
Thus we have two district line-elements connected by such a transformation 
and each applies to a specific $R^m$ domain. \par
To establish that an internal function $f_M(R^m)$ exists with the 
appropriate properties proceed as follows: let $\cal I$ be the set of all 
nonsingleton intervals in $\power{\real},$ where $\real$ denotes the real 
numbers. Let $\cal F\subset \power {\real \times 
\real}$  be the set of all nonempty functional sets of ordered pairs. For each 
$I \in \cal I,$ let $C(I,\real) \subset \cal F$ be the set of all real valued 
continuous functions (end points included as necessary) defined on $I.$ For 
each $k > 0,$  $\exists f_k\in C((-\infty,0],\real),\ (-\infty,0]\in
\cal I,$  such that $\forall x \in (-\infty,0],\ f_k(x) = 1/(x-k).$  Further,
$\exists g_k\in C((0,2k],\real),\ (0,2k]\in
\cal I,$  such that $\forall x \in (0,2k],\ g_k(x) = 
-x^3/(2k^4)+7x^2/(4k^3)-x/k^2-1/k.$ Then 
$\exists h_k\in C((2k,+\infty),\real),\  (2k,+\infty) \in
\cal I,$  such that $\forall x \in (2k,+\infty),\ h_k(x) = 0.$ Finally, it 
follows that $\lim_{x \to 0^-} f_k(x) = \lim_{x \to 0^+} g_k(x),\ 
\lim_{x \to 2k^-} g_k(x) =\lim_{x \to 2k^+} h_k(x).$ Hence 
$$H_k(x) = \cases{f_k(x);&$x\in (-\infty,0]$\cr
                g_k(x);&$x\in (0,2k]$\cr
                h_k(x);&$x\in (2k,+\infty)$\cr}$$ 
is continuous for each $x \in \real$ and has the indicated properties.\par
Now $H^\prime_k(x)$ exists and is continuous for all $x \in \real$ and 
$$H^\prime_k(x) = \cases{f_k^\prime(x);&$x\in (-\infty,0]$\cr
                g_k^\prime(x);&$x\in (0,2k]$\cr
                h_k^\prime(x);&$x\in (2k,+\infty)$\cr}$$ \par
All of the above can be easily expressed in a first-order language 
and all the statements hold in our superstructure enlargement (Herrmann, 
1991b). Let $0<\eps \in 
\monad 0.$ Then there exists an internal hypercontinuous hypersmooth $H_\eps \colon 
\hyperreal \to \hyperreal$ such that $\forall x \in 
\Hyper {(-\infty,0]},\ 
H_\eps (x)= 1/(x - \eps)$ and  $\forall x \in  \Hyper {(-\infty,0)}\cap    
\real,\ \st {H_\eps (x)} =\st {1/(x - \eps )} = 1/x;$  and for $x = 0,\ 
H_\eps(0)$ exists, although $\st {H_\eps(0)}$ does not exist as a real 
number. Further, 
$\forall x \in  {(2\eps ,+\infty )}\cap \real=(0,+\infty),   
\ \st {H_\eps(x)} =0.$  To obtain the hypercontinuous hypersmooth $f_M,$  
simply let $cf_M 
= H_\eps,\ x = \lambda,\ R^m \in \hyperreal.$ \par
In order to motivate the selection of these functions, first   
recall that a function $f$ defined on interval $I$ is 
standardizable (to $F$) on $I$ if 
$\forall x \in I \cap \real,\  F(x) = \st {f(x)}\in \real.$ 
Now, consider the transformation (**) in the nonstandard form 
$dU^m = dt^m +f_M(R^m)dR^m$ where internal $f_M(R^m)$ is a function defined 
on $A \subset \hyperreal,$ and $\lambda=\lambda(R^m).$ There are infinitely 
many nonstandard functions that can be standardized to produce the line 
element 
$dS^2.$ In this line-element, consider substituting  for the function $\lambda = 
\lambda (R^m,)$ the function $\hyper \lambda -\eps.$ The transformed 
line-element then  
becomes, prior to standardizing the coefficient functions (i.e. restricting 
them the natural world),\parm 
$$T= (\hyper\lambda -\eps)c^2((dU^m)^2 -2f_MdU^mdR^m + f^2_M(dR^m)^2) - $$
$$(1/(\hyper\lambda -\eps))(dR^m)^2-$$
$$(R^{m})^2(\sin^2\theta^{m}(d\phi^{m})^2 +
(d\theta^{m})^2)=$$ 
$$(\hyper\lambda -\eps)c^2(dU^m)^2 -2(\hyper \lambda -\eps)c^2f_MdU^mdR^m  + 
$$ $$
\overbrace{((\hyper\lambda -\eps)c^2 f^2_M- 
(1/(\hyper\lambda-\eps)))dR^m}^b dR^m -$$
$$(R^{m})^2(\sin^2\theta^{m}(d\phi^{m})^2 +
(d\theta^{m})^2).\eqno (30)$$ \par
Following the  procedure outlined in Herrmann (1989), 
first consider the partition 
$\real= (-\infty, 0] \cup (0,2\eps] \cup (2\eps, +\infty),$  where
$\eps$ is a positive infinitesimal. Consider the required constraints. 
(2) As required, for specific real intervals, 
all coefficients of the terms of the transformed
line-element are to be standardized and, hence, the simplest are extended 
standard functions. 
(3) Since any  line-element 
transformation, prior to standardization, should retain its infinitesimal 
character with respect to an 
appropriate interval $I$, then for any 
infinitesimal $dR^m$ and for each value $R^m\in I$ terms such as  
$G(R^m)dR^m,$ where $G(R^m)$ is a coefficient function, 
must be of  infinitesimal value.\par
For the important constraint (3), 
Definition 4.1.1, and theorems 4.1.1, 4.1.2 in Herrmann (1991a) 
imply that for a fixed
infinitesimal $dR^m$  in order to have expression $b$ infinitesimal as $R^m$ varies, the 
coefficient $h(R^m)= (\hyper\lambda -\eps)c^2 f^2_M- 
1/(\hyper\lambda-\eps)$ must be infinitesimal on a subset $A$
of an appropriate interval $I$ such 
that  $0 \in A.$ The simplest case would be to assume that $A = \Hyper (-
\infty,0].$ Let 
standard $r \in A \cap \real.$ Then it follows that $h(r) \subset \monad 
0.$ Thus $\st {h(r)} = 0.$ Indeed, let $x \in (\cup \{\monad{r}\mid r <0,\ r 
\in \real \})\cup (\monad 0 \cap A).$ Then $\st {h(x)} = 0.$ Since we are seeking 
a transformation 
process that is  hypercontinuous, at least on $\Hyper (-\infty, 0],$ this last 
statement suggests the simplest to consider would be that  
 on $\Hyper (-\infty,0],\  h = 0.$ Thus the basic constraint yields the 
basic requirement that on $\Hyper (-\infty, 0]$ the simplest function to choose 
is $cf_M(x) =1/(x - \eps).$ Since standardizing is required 
on $\Hyper (-
\infty, 0)\cap \real,$ we have for each $x \in \Hyper (-
\infty, 0)\cap \real,$ that $\st {cf_M(x)} = c\st {f_M(x)} =\st {1/(x - 
\eps)}= 1/x.$
This leads to the assumption that 
on 
$(-\infty, 0]$ the function $f_k(x) = 1/(x - k),\ k >0,$ 
should be considered. 
After *-transferring and prior to standardizing, this selection
would satisfy (3) for both of the coefficients in which
$f_M$ appears and for the interval $I=\Hyper (-\infty,0].$ 
The function $g_k$ is arbitrarily selected to satisfy the hypercontinuous 
and hypersmooth property and, obviously, $h_k$ is selected to preserve 
the original  line-element  for the interval $(2\eps, +\infty).$ Finally, it 
is necessary that the resulting new coefficient functions, prior to 
standardizing, all satisfy (3) at least for a fixed $dR^m$ and a varying 
$R^m\in \Hyper (-\infty, 0]$ for the expression (**).   
It is not difficult to show 
that $\vert H_\eps (x) \vert \leq 2/\eps$ for all $ x \in \hyperreal.$ 
Consequently,  for 
$\eps = (dR^m)^{1/3}$ expression (**) is an infinitesimal for all $R^m \in 
\hyperreal.$ \par
Let $1-v^2/c^2=\lambda.$ For the collapse, scenario $R^m =2GM/c^2.$ If
$2GM/(R^mc^2) < 1,$ substituting $2GM/R^m = v^2,$ into $(13)_k,$ 
yields  
the so-called Schwarzschild line-element. With respect to the transformation, 
(A) if $R^m <2GM/c^2,$ then for $ 
\st {f_M(R^m)} = 1/(c\lambda), 
\ \lambda = 1 - 2GM/(R^mc^2);$ 
for (B) $R^m > 2GM/c^2;\ \st {f_M(R^m)} = 0,$ and for the case that 
(C) $R^m = 2GM/c^2,$ 
the 
function $f_M$ is defined and equal to a
NSP-world value $f_M(R^m).$ [It is a Robinson infinite number. Such numbers have very interesting algebraic properties and are not equivalent to the behavior of the concept being considered when one states that something or other ``approaches infinity.'' (Also see Article 2.)] 
Now, for case  (C), $\st {f_M(R^m)}$ does not exist as a real number. 
 Hence, (C) 
has no direct effect within the natural world when $R^m = 2GM/c^2,$ although 
the fact that 
$f_M(R^m)dR^m$ 
is an infinitesimal implies that $\st {f_M(R^m)dR^m} = 0.$ Using 
these NSP-world functions and (30),
 cases (A) and (C) yield  
$$dS_1^2 = \lambda(cdU^{m})^2 - 2cdU^m\,dR^m -$$ 
$$(R^{m})^2(\sin^2\theta^{m}(d\phi^{m})^2 +
(d\theta^{m})^2). \eqno (31)$$
But case (B), leads to the Schwarzschild line-element. 
The two constraints are met by 
$f_M(R^m),$ and indeed the standardized (A) form for $f_M(R^m)$ 
is unique if (3) is to be 
satisfied for a specific interval.\par
Since this is an ideal 
approximating model,  
in order to apply this ideal model to the natural world, 
one must select an 
appropriate real $k$ for the real valued function $H_k.$ Finally, it is not 
assumed that the function $g_k$ is unique. For solutions using line-element
(31), the $dU^m$ [resp. $dR^m$] refers to the timing [resp. 
length] infinitesimal light-clocks. \par
If one is concerned with how this theory of gravity compares with 
the Einsteinian theory, then for a given positive non-infinitesimal 
 $k$ 
the transformation in 
terms of such a $k$ do not satisfy the (restricted) Einstein gravitational 
law. The new transformation required in the {transitional 
zone} is not a proper transformation at but the one point $k/3$  
where it does not have a
local inverse. Continuity considerations can be 
applied on either side of this point.  
Since within this zone the $R^m$ is assumed to change 
very little, then to analyze effects in this zone it seems reasonable to 
assume that $a$ behaves as a constant.\par 
The above method for transforming one line-element into another to  
minimize singularities, so to speak, can be done for other transformations. 
Our use 
of the Eddington-Finkelstein transformation is but an illustration,  
although in the next section, based upon the concept of simplicity, 
it will be used in a brief investigation of the non-rotating black hole 
concept. It will be shown how a black hole can exhibit interesting white hole 
characteristics when NSPPM expansion is included. \parm
\leftline{\bf 7. Applications}\parm 
How should distance be measured when such line-elements as the Schwarzschild, 
(26), (27) or (31) are investigated? For measurements relative to a Special Theory scenario, such as an essentially constant field, Einstein distances, 
velocities or times are appropriate for the measurements.  {\bf The 
infinitesimal light-clocks yield how the physical universe is being 
altered by the (emis) effects.} This is all that this entire theory represents.  
As pointed out by {Fock (1959, p. 177, 351), such  concepts only have 
meaning when they 
are compared to similar  concepts more closely related with human intuition.}
Basic intuition at the most local level is obtained through our use of a 
Euclidean system. But we also live in an  approximate  
exterior  (gravitational) Schwarzschild field. To investigate how new 
gravitational 
effects alter behavior, these new effects should be compared with  
behavior essentially not affected by the specific gravitational field effect being considered. Such comparisons are made with 
respect to the measurements made by local infinitesimal 
light-clocks used as a standard. Depending upon the problem, these may be taken  
as the ``s.'' For comparison purposes, the measurement of distance between 
two fixed 
points, where the distance measure is considered as being physically altered 
by (emis) effects, the {``radar'' 
or ``reflected light pulse'' method} would be appropriate. The same would hold 
for the linear 
effect line-element. Of course, this is actually the behavior of 
an electromagnetic pulse (a photon) as it is being affected by (emis). Thus 
to obtain such measurements, let $dS^2$ or $dS^2_1 = 0$ 
and consider only the radial part for the Schwarzschild and (31), or 
the entire line-element for the linear effect (26) and (27). For $(13)_a$, 
this gives $(dR^m)^2 = (\lambda)^2 (cdt^m)^2.$ Taking, say $(1/\lambda)dR^m =  
cdt^m.$  Let $R^m_1,\ R^m_2$ be light clock measured distances from the 
center, assuming both are exterior to the Schwarzschild radius, taken at 
light-clock times $t^m_1,\ t^m_2.$ These light-clocks have been altered 
by the exterior Schwarzschild field. For the Schwarzschild field, 
this gives
$$c(t^m_2-t^m_1) = R^m_2-R^m_1 + r_0\ln {\left[{{R^m_2 - r_0}\over{R^m_1 - 
r_0}}\right]},\eqno (32)$$
where $r_0 = 2GM/c^2.$ 
Of course, for the quasi-Schwarzschild case, (32) has a slightly more complex 
form. The reason we have such an expression is that, compared to 
no gravitational field, the path of motion of such a pulse is not 
linear. \par
One of the more interesting coordinate transformations used in the Einsteinian 
theory is the {Kerr transformation}  for a rotating body. Analysis of the pure 
vacuum Kerr solution yields some physical science-fiction type conclusions 
(Ohanian, 1976, p. 334). Ohanian states that these conclusions are
physically meaningless unless ``. . . a white hole was -- somehow --
created when the universe began'' (1976, p. 341). 
He also mentions that such 
objects would probably not survive for a long ``time'' due to the instability of 
white holes, and the like, within the pure Kerr geometry.  The possibility that a
non-rotating black hole can be altered into a {pseudo-white hole} will be 
investigated using the notions from the previous section. This will not be 
done using another science-fiction type transformation, the {Kruskal space.}
As pointed out by Rindler (1977, p. 164), this space also suffers from 
numerous {``insurmountable difficulties''} and the full Kruskal space probable 
does not exist in nature. \par
Referring back to the previous section, consider a {\it 
transition zone} for a given $k.$ This zone is for all $\lambda = 1 - 
(v_p+d)^2/c^2 \in (0,2k].$  The transition zone motion of photons and 
material  ``particles'' can be investigated. As in the case of expressions 
such as (32), the mathematical analysis must be properly interpreted. To do 
this, physical entities must be characterized as behaving like a ``photon'' or
like a {``material particle.''} You will almost never find a definition of the 
notion of the ``material particle'' within the literature. Since this is the 
infinitesimal calculus, then, as has been shown (Herrmann, 1991a), the 
``material'' particle is an (not unique) infinitesimal entity that carries
hyperreal material characteristics. However, such a particle does not affect 
the standard field (i.e. the field exterior to the Schwarzschild surface
still behaves like a gravitational vacuum). 
 When one  
investigates the paths of motion or other physical concepts for 
such material particles as they are measured by 
infinitesimal light-clocks, the proper interpretation is to state that the 
physical property is affected by the (emis) as it is altered by the NSPPM  
\underbar{and} as it is compared with the non-altered (emis) effects. \par 
To analyze properly particle behavior, we have {three partial line-elements,} 
that represent radial behavior.
These partial line-elements are
$$d\ell^2 = (\lambda -k)(cdt^m)^2 - (1/(\lambda -k))(dR^m)^2,$$
$$d\ell_1^2= (\lambda - k)(cdU^m)^2 -2cdU^mdR^m,$$
$$d\ell_t=(\lambda - k)c^2(dU^m -(1/c)g_kdR^m)^2 - (1/(\lambda -
k))(dR^m)^2,\eqno (33)$$ 
\noindent where $d\ell_1$ is applied when $\lambda <0,$ $d\ell_t$ is applied when 
$0 <\lambda <2k,$ and $d\ell$ is applied when $\lambda > 2k.$ For most 
analysis, the conclusions are extended to $\lambda = 0, 2k$ by continuity 
considerations. \par
Rather than give a complete analysis of the above line-elements in this article, we 
illustrate how it would be done by considering the general behavioral 
properties for  
electromagnetic radiation and material particles. If there was no 
altering gravitational field present,
then we have that, for radiation, $dS^2 = 0 = c^2(dt^m) - d(R^m), s = m.$ 
Hence, we have that 
$dR^m/dt^m = \pm c.$ The selection is made that $dR^m/dt^m = c>0$ represents 
``outgoing'' radiation and that $dR^m/dt^m = -c<0$ represents ``incoming''
radiation. For material particles, we have the concept of incoming and 
outgoing particles as they are modeled by the general statement that 
$\pm dR^m/dt^m < c.$ Lawden (1982, pp. 155--157) correctly analyzes the cases 
for $d\ell$ and $d\ell_1.$ His analysis holds for the case of the 
quasi-Schwarzschild field as well. Both incoming radiation and particles can 
pass 
from the exterior quasi-Schwarzschild field into the  
transition zone and if they can leave the transition zone, they be 
forced to continue towards the center of attraction, where, for this analysis,
one may assume that the mass is concentrated. Further, no material particle or 
radiation pulse can leave the region controlled by $d\ell_1.$ What happens 
in the transition zone?\par
If one analyzes what effects occur in the transition zone in terms of 
$t^m$ and $R^m$, one finds that there is a general 
photon turbulence depending upon the different photon ``families'' predicted 
by the statement that $dR^m/dt^m = \pm c(\lambda -k),$ where $0 < \lambda 
<2k.$ A Lawden type analysis for material particle behavior also leads to a 
particle turbulence within the transition zone.\par
 Now as  suggested by Ohanian, 
what happens if the universe is formed and a homogeneous sphere 
appears such that its physical radius $R$ implies $\lambda \leq 2k$? 
For simplicity, it is assumed that for the following speculative scenario 
that, with respect to these pure gravitational line-elements, $d = 0$  \par
(1) Since the entity being considered, under a special process, could appear when the universe is formed, 
then the actual material need not be the same as a neutron star and the like.
Of course, this does not necessarily assume the Big Bang cosmology but 
does assume that 
the Laws of Nature also appear at the moment of formation. [This possibility 
comes from the {MA-model scenario} (Herrmann, 1994b).]  
Suppose that the material is malleable with respect 
to the gravitational force. If there is 
any ``space'' outside of the radius $R^m,$ then the Schwarzschild field would 
apply in this vacuum. [Recently using a different argument and 
assuming 
the entire Einsteinian theory, {Humphreys}  (1994) has speculated that the material 
could be ordinary water!] \par
(2) During the {transitional phase, the apparent  
turbulent physical behavior} for the material particles would tend  
to keep a certain amount of material within this zone as the gravitational 
collapse occurs with respect to the remaining material.  
This has the effect of producing a {``halo'' or a spherical
shell} that might or might not 
collapse along with the material that has been separated from the 
material that has remained within the transition zone. Collapse of this shell 
would tend to depend on its density, the actual continuation of the collapse
and numerous other factors. \par
 Unless something were to intervene after collapse began, then the usual 
black hole scenario would occur. In the next and last section of this article, 
we speculate upon such an intervention. \parm
\noindent {\bf 8. Prior to Expansion, Expansion and Pseudo-White Hole Effects.}\parm
The concept of the white hole is almost always developed by considering 
 {``time reversal.''} Although the  timing infinitesimal light-clock
can be oriented with a non-positive orientation, this process may not be 
necessary for certain cosmological events to occur. If one assumes that the investigated
black hole is formed, then there is a vast number of 
 {diverse black hole scenarios} highly dependent upon many parameters that have 
been held fixed for our previous investigations. The notions of spherical 
symmetry, uniform density, uniform expansion, gravitational vacuum states, and 
the like, will in reality be violated. In this theory, as well as the Einstein 
theory, the assumptions made as to the simplicity of the Newtonian potentials 
will certainly be violated in a real expanding universe. Although it might be 
possible to substitute a new potential velocity $v_p$ for the simple one 
presented here, say in the case of a bounded universe, in reality other 
competing potentials will preclude such a simplistic approach.
 Thus any speculative 
description for how a black hole could actually change in such a manner 
to produce white hole 
effects, and even produce some effects associated with a Big Bang cosmology, 
must be 
very general in character. Further, these speculations probably cannot 
be arrived at by 
strict mathematical analysis. This leaves these speculations to 
descriptions generalized by means of human intuition. This further weakens 
their application since it assumes that the human mind is capable of 
comprehending, even on the most general level, the actual processes that yield 
the development of the universe in which we dwell. I caution the reader
to consider these facts while examining the next speculative descriptions.\par
First, the real values for the basic expansion function $a$ (i.e. $\st {\hyper 
a}=a$) are 
controlled from the NSP-world. 
Indeed, it is a connection 
between the 
NSPPM and the natural world. Since observation from the NSPPM is the 
preferred observation, then, as mentioned, one might assume that
 there is a preferred coordinate system within the NSPPM and that the 
values of $d$ can depend 
upon specific locations within the natural world. This, of course, includes the 
special case where $d$ depends only upon the $t^s$ and the Copernican 
principle is applied, as well as the restrictive 
case that the universe actually has a center and expansion is with respect 
to that center. Even if one assumes 
the Copernican principle,  
it need not be assumed that the radial expansion follows this simplistic pattern 
 {\it when viewed 
from every position.}\par
In what follows, certain simplistic assumptions are 
made relative to the behavior of the NSPPM and such assumptions are
only assumed to hold in a neighborhood of the black hole. Five possible 
scenarios are presented. It might not be possible to differentiate 
scientifically these 
scenarios one from another by any observational means from our present epoch.
For this reason, an individual's selection of one of the  numerously many 
scenarios will probably involve philosophical considerations rather than 
scientific data. \par
(1) Suppose that at formation scenario (2) of section 7 occurs 
and the halo is separated from the remaining material and the malleable 
material undergoes gravitational collapse. No superimposed expansion occurs. 
That is $\hyper a$ is an infinite hyperreal number. Suppose that at a 
particular moment based upon collapse factors, an {extreme rate of spatial 
expansion occurs.} This would radically {alter the density of the collapsed malleable 
material.} The material captured in the transition zone, since the zones size is fixed 
and the gravitational effects in this zone are based upon the gravitational 
field produced by the collapsing material, would expand beyond that zone and 
continue to ``move outward'' with respect to the center of the collapsing 
material (i.e. the center of attraction). Under exceptional critical values 
for the parameters involved, it is possible that the expansion rate is so 
great and the density is reduced to such an extent that the assumed 
symmetric collapsing material actual increases its radius to a point that 
would rapidly go beyond the original Schwarzschild radius, and beyond the transition 
zone. This {reduces the original Schwarzschild radius} due to the loss of mass 
and would alter greatly the predicted behavior since the simple potential 
velocity 
and spherical symmetry would certainly be altered. \par
As previously derived, the alterations in many physical processes due to the 
gravitation field 
would 
probably still occur  ``near'' to the new Schwarzschild radius. These 
alterations would tend to {``slow'' down certain processes} when compared to a 
standard. On the other hand, if the standard is considered as having values near 
to the new Schwarzschild surface, then the processes would seem to 
 {``increase''  
within the material that expanded to great distances beyond the new 
Schwarzschild surface.} Depending upon critical values associated with the 
expansion, the composition of the material, and the {``explosive'' effects,} 
the new Schwarzschild surface could slowly continue to shrink. \par
The composition of the material that has expanded beyond the 
original Schwarzschild 
surface depends upon whether or not the collapsing material follows the known 
laws of quantum mechanics. These laws need not apply since we have no 
laboratory verified knowledge of 
exactly how extreme {gravitational potentials affect quantum behavior.} 
What can be assumed is that a great deal of ``cooling'' might take place along 
with the great increase in expansion. I mention the theoretically derived 
result that a {subparticle composed NSPPM is capable of drawing off 
energy,} 
even in vast amounts. Depending upon many factors, there could be a great deal 
of collapsed material forced beyond the previously located Schwarzschild 
surface. This process might, of course, stop the gravitational collapse. To an 
outside observer there might appear to be a spherical material shell expanding 
from the previous black hole with a great deal of additional material between 
this shell and the remnants. This ``explosion'' of material from a specific 
position in space corresponds to one white hole aspect. This is a pseudo-white 
hole effect. The effect would also be similar to the appearance of a 
supernova. I point out that extreme and varying expansion rates have been 
postulated previously (Guth, 1981).\par
(2) Consider the scenario described in (1) with the exception that the 
transition zone is empty and collapse begins at the edge of the transition 
zone at formation. \par
(3) Consider the scenario described in (1) with material in the 
transition zone, but at formation there is a separation between the transition 
zone and the collapsing material at the moment of formation. \par
(4) Consider the scenario described in (1) with the exception that there is no 
material in the transition zone and there is a separation between the 
transition zone and collapsing material at the moment of formation.\par
(5) Modify the previous scenarios by considering the infinitely many  
variations brought about by diversions from the ideal uniform and homogeneous 
behavior.\par
Of course, if there was no black hole formed, then the above 
speculations are vacuous.  Further, I am convinced 
that this corrected theory is capable of predicting all of the actual  
effects that are associated with the Einstein theory. It specifically shows 
that such effects are caused by an (emis).  
Also, this corrected theory shows, once again,  
that there is not one theory that predicts such behavior. {Hence, this article 
establishes that the selection of any theory for the development of our 
universe must be based upon considerations exterior to theoretical science.}\par
As previously mention, recently {Humphreys} (1994) has proposed a theory for the 
formation of the entire universe by assuming that the entire universe was 
produced by such a black hole. His theory uses concepts from Einstein's 
theory, which I reject, and assumes that the rapid expansion is caused by a 
sudden 
change in the cosmological ``constant.'' His scenario is very similar to (4) 
above. He assumes that the {cosmos is bounded} and that formation occurs at a 
specific point, say at the center of the black hole, which, at the moment of formation, is composed of ordinary water.  
The above corrected 
theory also applies to a bounded cosmos with its gravitational center of 
mass. Humphreys includes various descriptions for physical changes that 
might occur within the black hole while undergoing gravitational collapse.
His selection of one of the above general scenarios is based entire upon a 
philosophical stance. Indeed, it appears that scenario (3) might have been a slightly better choice. \par
[Note added 1 June 1998, corrected 28 March 1999. It appears that Humphreys' model as stated may fail to achieve the goals claimed in a few instances. First, the present day cosmological constant $\Lambda$ is estimated to be no larger than $10^{-56}{\rm cm^{-2}}.$ Humphreys uses the Schwarzschild configuration, the vacuum solution and the classical Schwarzschild surface (i.e. event horizon) throughout his discussions, especially relative to the geometry \underbar{exterior} to such a surface. Due to the dust-like properties of matter interior to this surface and due to a comparatively large cosmological constant, the collapse scenario for the dust-like material would be overcome and the material would escape through the event horizon and give a white hole effect. (However, this scenario does not appear to have all of the actual white hole properties.) He states, ``I suggest that the event horizon reached earth early in the morning of the fourth day.'' (Humphreys, 1994, p. 126) The earth here is a type of ``water-world'' that has stayed ``coherently together.'' (Humphreys, 1994, p. 124) The event horizon also remaines approximately in that position the entire ``fourth day.'' Humphreys discusses the Klein line-element and the Schwarzschild line-element without the cosmological constant. However, it is necessary that the cosmological constant be considered and other line-elements that include the cosmological constant should be  investigated prior to acceptance of this model. For example, consider the modified Schwarzschild solution where the significant expression is $1-2GM/(c^2r)-(1/3)\Lambda r^2 = 0$. A simple calculation, using $6.67 \times 10^8$cm as the radius of earth and $.889$cm as the value for $2GM/c^2$ yields $\Lambda = 7.39\times 10^{-18}/{\rm cm}^2.$ Humphreys states that the cosmological constant is to be set at a large value on day two of his creation model in order to produce a ``rapid, inflationary expansion of space.'' (Humphreys, 1994, p. 124) This does not appear to be the large value of the cosmological constant that Humphreys is considering in order to obtain the necessary rapid inflationary expansion as shown in the Moles paper cited by Humphreys. More importantly, if this value is inserted into this modified Schwarzschild expression with the mass of the universe, then the event horizon that was at $450 \times 10^6$ lyr suddenly vanishes, indeed, no event horizon exist. This yields a direct contradiction. (Including a term for ``charge,'' in the above, will not significantly affect these results.)]\par
I would like to thank T. G. Barnes and R. J. Upham who supplied certain 
important references that led to many improvements in this paper. [See note
8, parts 1, 2, 3, 5.]\parm
 \centerline{\bf NOTES}\parm
[1] I have been asked to more fully justify my statement in  
Article 3 of the  ``Foundations . . . .'' paper, just below equation (12).  
The basic characterizing equation in line 14 is 
$\alpha + \beta(1-\alpha^2) = 0.$ The notation used assumes that we are working 
in real and not complex numbers. The expression 
$\alpha = \pm \sqrt {1 - \eta}$ characterizes $\eta$ such that $\sqrt {1-n} 
\geq 0.$ The  combined velocity $v+d$ has the property that $0\leq v+d.$ 
Consider the case that $0 < v+d$ and suppose that $\alpha = \sqrt{1 - \eta}.$ 
Then from the line 14 expression we have that $\eta\beta = -\sqrt {1-\eta}.$ 
This also tells us that $\eta \not= 0.$ Hence we may write 
$\beta = -\sqrt {1-\eta}/\eta.$  Now substituting into (B) 
($dT^s = \beta dR^m + dT^m)$ yields (i) $dT^s 
=-(\sqrt {1-\eta}/\eta)dR^m+dT^m$ and into (C) ($ dR^s = (1- \alpha\beta)dR^m 
- \alpha dT^m)$ 
yields (ii) 
$dR^s = dR^m/\eta -\sqrt {1-
\eta}dT^m.$ Combining the two differentials yields the requirement that  
$${{dR^{s}}\over{dT^{s}}}=
\left[{{{{1}\over{\eta}}{{dR^{m}}\over{dT^{m}}} -\sqrt{1-\eta}}}\right]\div     
\left[-{{\sqrt{1-\eta}}\over{\eta}}{{dR^{m}}\over{dT^{m}}}+ 1\right].\eqno 
(1)$$ 
Now consider $dR^m = 0$ in (i) and (ii). Then we have that for $v+d \not= 0,\ (v+d)/c = -\sqrt {1-\eta} <0.$ This 
contradicts the original requirement for this case that $0 < v+d.$ \par
Note that in modern infinitesimal analysis the actual combining process to 
obtain (1) is not division, although the result is the same but the function 
interpretation is important. (i) states that $dT^s/dT^m 
=-(\sqrt {1-\eta}/\eta)(dR^m/dT^m)+1,$ (and nonzero since $t^m$ is 
dependent upon $v$) $\Rightarrow
dT^m/dT^s =(-(\sqrt {1-\eta}/\eta)(dR^m/dT^m)+1)^{-1}$ and from (ii) 
$dR^s/dT^m =(1/\eta)dR^m/dT^m -\sqrt {1-\eta}.$ 
Application of the chain rule leads to (1).\par

[2] One of the more difficult aspects of this research is to disregard   
the classical  interpretations of Einstein's theories and to  
interpret the mathematical statements in a simple and consistent manner. 
Notice that the chronotopic 
interval expressions (4)$_a$, (4)$_b$, (5)$_a$,  (5)$_b$ are stated in terms 
of subscripts and superscripts. The infinitesimal light-clocks being 
modeled by these expressions have two distinct interpretations. One is that 
they represent the ``actual'' $u$ or $L$ associated with an infinitesimal light-clock, 
the other is that they are but analogue models for an {\it electromagnetic 
interaction with the NSPPM} or (emis) effect. \par

Usually, the left-hand sides represent ``actual'' 
infinitesimal light-clock and how the light-clock is altered with 
respect to motion 
by the NSPPM. For these four expressions,
the right hand side represents unaffected infinitesimal light-clocks 
located in the proper coordinate positions and how they would 
``measure,'' the ``actual'' changes that take place 
within the 
the two ``actual'' NSP-world measuring infinitesimal light-clocks.\par

The expressions (4)$_a$, (4)$_b$, (5)$_a$,  (5)$_b$ are what is gleaned from 
the most basic laboratory observations within the N-world relative to  
electromagnetic propagation and are not, as yet, related to  possible 
(NSP-world) NSPPM physical-like properties that produce these effects. Possible 
NSPPM physical properties are investigated by using expressions such as (8) and 
(A), (B).  Coordinate transformations, if made, are relative to 
possible NSPPM physical properties.
Suppose that these measuring infinitesimal light-clocks 
are affected by 
this $P$-process. We seek a relationship $\phi(dR^m,d\theta^m,d\phi^m,dt^m)$  
between altered behavior of these infinitesimal light-clocks as they would appear for the measuring infinitesimal light-clocks 
so that $L\Pi_m = \phi(dR^m,d\theta^m,d\phi^m,dt^m).$\par

 The results of this approach indicate that 
the quantities now denoted by the superscript $m$ represent N-world 
measurements with respect to coordinate infinitesimal light-clocks 
that incorporate the (emis) effects. These measurements are taken 
relative to the coordinate transformation. It is clear from our letting 
$R^m = R^s, \ \theta^s = \theta^m,\ \phi^s = \phi^m$ and other considerations
that the measuring devices used to measure these qualities are not, at the moment, being considered themselves 
as being altered by the $P$-process. However, the reason for now using such a variable in the form $R^m$ is that relative to the hyperfinite approach to integrals, in general, the values of $R^m$ or any function in $R^m$ are not constant, in general, but depend upon each subdivision. Thus, such variables as $R^m$ represent a type of cumulative alteration from the standard $R^s.$ 
The method used to find the value of $\eta$ after equation (12),
where $dR^m/dT^m$ has been set to zero, should be considered as an initial condition for the static behavior of $P$ or other behavior that yields $dR^m/dT^m = 0.$  \par

After the line-element is obtained, the next step is to 
interpret the element relative to its effects upon other entities. There are 
arguments that show that such a {\it re-interpretation} is viable 
by means of the 
test particle concept which itself does not essentially alter the field. The behavior to be 
calculated is N-world behavior. Thus, a test particle can be conceived of as 
a physical infinitesimal light-clock. \par

When this derivation of the  Schwarzschild line-element is used to investigate 
planetary motion, $\theta^m =$ constant. Then the $dt^m$ and $dS$ are 
eliminated to obtain an expression in terms of $R^m$ and $\phi^m.$ Rather than 
use the concept of geometric geodesics, Fermat's general principle of least 
``time'' or ``action'' can be 
applied. This leads to the same partial differential equations as would the 
geodesic approach. The resulting variables $\phi$ and $R$ can be measured with 
respect to any standard by the observer. Exactly the same variables are used 
for the measurement of the possible deflection of electromagnetic radiation by 
massive bodies. \par

[3] There is a slight confusion as to my use of the term ``invariant'' with 
respect to $dS^2.$ The term``invariant'' only refers to its ``value.''
How does one relate this inf. light-clock approach to the classical one? For clocks, the classical approach claims to present expressions for time measurements for all clocks. But, for many circumstances, time needs to vary continuously. There is no such clock that has this property. So, the classical approach is but an approximation. But, it uses the notion of infinitesimalizing. The only physical behavior that most closely approximates infinitesimal behavior is subatomic photon behavior. Hence, I chose the light-clock for both time and distance measures. The inf. light-clock is, of course, a conceptual model. It also allows us to show how photons locally behavior for various $v_p$ . \par
Of course, one can immediately return to the classical approach by symbol substitution but must explain them. For example, for a gravitational field if the $m$ and $s$ superscripts are dropped in (13)$_a$, let $dS = ds,$ and, state the variables that now appear represent measures where there is no gravitational field, you can derive the classically  expressed Schwarzschild line element.\par 

These articles are intended as a mere beginning and as an indication of an 
appropriate method. 
Not all reasonable and superimposed physical processes have been considered. 
Further research, by individuals other then myself, relative to other 
reasonable physical NSPPM processes should provide additional confirmation that 
this approach is viable. \par    
[4] The concept of the linear effect line-element is significant in that 
it yields a NSPPM physical reason for Special Theory effects. Notice that so 
as not to confuse this with the spherically transformed line-element, $dr^s$ replaces $dR^s$ and $dr^m$ replaces $dR^m.$ 
Further, the 
method used to obtain the transverse Doppler and the mass variation 
predictions replaces the $r^m$ and the $t^m$  (emis) effect quantities with 
measurable N-world quantities. (See note [19]). \par

[5] If for various investigations variables such as $\tau$ and $t^m$ cannot 
be replaced by other measurable quantities, then the  concept of the proper 
time is replaced  
by an ``actual'' NSP-world measure of the infinitesimal light-clock 
time as registered by the test particle itself. The ``coordinate'' 
time $t^m$  is conceived as related to an   
N-world light-clock that is affected by the field and that can only
approximate an  infinitesimal light-clock. Further, the infinitesimal 
light-clock used to measure this ``coordinate time'' is measuring just one aspect of 
the (emis) effect. Under certain circumstances, such as investigating the 
Schwarzschild line-element for test particles approaching the Schwarzschild 
surface, the NSP-world time measure ${u\Pi}$ has different properties
than the N-world time measure. Thus these two differences in  
``time'' measurements does not lead to a contradictory. This is similar to the  
Special Theory statement that the relative velocity $w$ is an 
unbounded NSP-world measure while the N-world Einstein measure $v_E$ is 
bounded. \par
[6] For the case of electromagnetic radiation, the linear effect line-element refers only to the 
Ives' interferometer case of two fixed positions $F_1$ and $F_2$ having $v_E = 
0$ (i.e. a to-and-fro ``linear'' light path scenario).
In which case, the Fermat's principle (or ``geodesic'') approach leads to the 
two 
differential equations $dt^m/d\tau = K$ and $K^2(dr^m/dt^m)^2 -c^2K^2 = 0.$
Since, in general, $K\not= 0,$ this yields that $(dr^m/dt^m)= c.$ This 
corresponds to the laboratory measurement of the to-and-fro velocity of light. 
$\rm (14)_a,\ (14)_b$ must be carefully considered when applied to the physical world. (See note [19].) Usually, they applied only to infinitesimal light-clocks where there is an alteration in counts since $c$ is probably altered. 
This eliminates the need to interpret (8.2), (9.2) as ``unit'' alterations.
\par 
[7]  The actual derivation 
uses the ballistic property within a monadic cluster. This is modeled by a 
moving point. In the case where the point's velocity is the velocity of the 
NSPPM itself, then using the velocity   
$R^s/a$ result (29) is obtained. However, it is just as possible that the 
point is not moving
and the P-process is affecting the velocity of $c$ in the manner 
indicated. This can be produced by various motions of the 
NSPPM material.  
Unless there 
is some other confirmation that $d$ is motion of the point, then this cannot 
be assumed. Thus this                                               
need only be an 
\underbar{apparent} general textual expansion due to this special effect associated with 
electromagnetic radiation. What this means that such general expansion need 
not occur in reality, but rather it would indicate a superimposed property 
produced by an interaction with the underlying NSPPM. If this were the 
case, it would invalidate some of these speculations. \par 
[8] It should be obvious that in these articles I have not attempted to 
duplicate the 
approximately 80 years of General and Special Theory work. What has been done 
is to point out the absolute logical errors in these theories. Then:\par 
(1) a method 
is given that retains the concepts associated with electromagnetic 
propagation. \par
(2) Based upon laboratory observations within our local environment and a 
privileged observer within the NSPPM, certain conclusions are developed
by a strict interpretation of a mathematical structure. \par
(3) Specific physical descriptions for behavior within the NSPPM (i.e. the 
P-process) that alter natural world behavior
are introduced. These physical processes are then modeled by means of a line-elements that predict the behavior of physical entities as this behavior is 
measured by infinitesimal light-clocks. This yields intimate relations 
between such behavior and the properties of electromagnetic propagation. 
These physical processes and the paths along which they operate have taken the 
place of the ad hoc coordinate transformations of the Einstein theory.\par
(4) It is very important to realize that all of this analysis relative to the 
Special Theory (ST) is local. 
Indeed, it is relative only to effects within an empty universe. 
The effects can only be properly measured over local regions where the 
gravitational potential is considered to be constant and can, thus, be
``factored,'' so to speak, from the measurements. 
Measurements that might be taken by what could be considered as ``large 
light-clocks'' are not analyzed and could give different results. Predictions 
associated with ST for all of the alterations are based entirely upon very 
local measurements and uniform relative velocity. These predictions become 
less accurate when these conditions are altered. Further, as I have shown 
previously (Herrmann, 1989), so-called constant quantities can 
actually be considered as but piecewise constant. This piecewise constancy 
may be 
considered as a natural world restriction of hypercontinuous and hypersmooth 
NSP-world processes.  We infinitesimalize according to our views of natural 
world behavior. But as shown by our ST derivation, physical processes within 
the natural world and the NSP-world might be considered as contradictory if 
they 
occurred solely within one of these worlds. Consequently, physical observation 
and theoretical constructs are needed prior to infinitesimalizing for physical 
theories. It is very possible that the ST results cannot themselves 
be infinitesimalized in the mathematical sense that a continuous curve is the 
standard part of an infinitesimal polygonal path. Only experimental evidence 
would imply such a process. This means that there can be a considerable 
difference in (emis) effects for rotation. These effects need not be 
predictable by application of locally linear ST effects. 
\par
(5) I have not combined together the pure ST relative velocity as a $d$ and 
the gravitational potential velocity $v_p$ due to what I feel is a logical 
difficulty with the intuitive difference between a ``directed'' and  
``non-directed'' effect. The $d$ expression used to obtain our 
Robertson-Walker type line-element assumes, in order for it to have any 
meaning, that the universe is \underbar{not} empty. Further, 
always remember that the predictions derived for the Newtonian potentials 
would not be correct if the potentials took on a more complex character. 
Considering the alterations in the Newtonian potentials associated  
with a rotating body could certainly be considered and would lead to another 
line-element that might compare favorably with the Einsteinian theory.
[These effects need not be gravitational but may be centrifugal or the like. 
It is possible that there are no actual local (emis) gravitational field 
rotational effects except those countering effects produced by pure rotation 
itself.] 
However, in general, their character would be so complex that only very 
general conclusions could be 
predicted. The analysis in this article can be extended to include 
the effects of a Newtonian potential propagated with some 
specific velocity such as $c.$ This could lead to a Newtonian theory of 
gravitational propagation. For example, see equation (1) in Surdin (1962, p. 
551).\par
Finally, the idea of potentials and trajectories as modeled by the complex plane has been well established. Hence, it may be profitable to consider, within various line-elements, complex valued $v$, $d$ or $v+d$ when they model potentials or tensions. Indeed, for pure complex $d$ and $v=0$ this would lead, in many cases, to a reverse in the alterations of various measures. If this is a NSPPM effect, then a P-process within the NSPPM need not be considered except that such processes should be predicted within the natural world prior to determining which 
measures will be altered. \par
[9] The derivation that appears here is the correct simplified derivation and replaces the 
the overly complex derivation that appears in the published version Herrmann, R. A. An operator equation and relativistic alteration in radioactive decay, Internat. J. Math. and Math. Sci., 19(2)(1996):397-402. This published version also contains some notational errors. \par
[10] Relative to equation (28), the $d$ is selected for the obvious 
purpose. Although the $d$ is used, with $v=0,$ this, of course, can be reversed.
Further, there are two superimposed aspects within the universe. One is where
$v$ is considered as an actual velocity and the other is where $v$ is 
considered as a potential velocity. For consistency, it might be that the same 
two aspects can be used with equation (28). This is the line-element that 
leads to the expansion redshift. But, it is possible that a superimposed second 
effect can occur. A natural property of some fields is a tension property.
Thus there might be a superimposed second aspect relative to such a tension 
property. Then the $d$ would represent such a tension. In this case, the tension 
could be both NSP-world position and time dependent. Following the exact same 
procedures as previously presented, such a tension could also alter physical 
processes. For example, it could  ``slow them down'' so to speak with respect 
to a standard. There is also one other interesting aspect to all of these
line-elements. If either $v$ or d are say pure nonzero complex numbers, then 
the exact same methods used to determine alterations in physical measures can 
be applied and will yield the exact reverse of such alterations. If an 
alteration for a real $v$ or $d$ is a decrease in some measure, then 
the alteration in a measure using a $v$i or $d$i would be an 
increase in the measure.\par
[11] (29 MAY 1997) The most absurd statement ever made by intelligent 
individuals is that physical processes or entities actually alter the behavior 
of the concept called ``absolute time'' and such an alteration is reflected by 
alterations in natural processes or in the characteristics of various 
physical entities. The same can be said for those that rejected such a 
``time'' alteration and replaced it with an alteration in the concept of 
``absolute length.'' The most basic assumption within natural science is that 
all natural-system behavior is altered either by an alteration in 
characteristics of the entities themselves or such behavior is altered by an 
interaction with natural entities. Unless ``time'' is a natural entity, a 
particle or a field or whatever, then such a statement is absurd. It has never 
been demonstrated that ``time'' corresponds to one of these types of physical entities. What is altered is observer time. 
\par
Although I have mentioned it previously, the results in this book \underbar{do not} 
overthrow the Einsteinian General Theory of Relativity. What has been altered 
is the very basic interpretation and foundations of the Einsteinian theory. It 
is well-known that the Hilbert-Einstein gravitational field equation has 
solutions that do not correspond to the universe that many members of the 
scientific community accept. Full Kerr or Kruskal geometries 
contradict the standard cosmological model. 
What is being established is that the language of Riemannian geometry is not 
the language of reality. Riemannian geometry is but an analogue model for 
behavior within a gravitational field. The Patton and Wheeler remarks 
state strongly this same conclusion. ``Riemannian geometry likewise 
provides a beautiful vision of reality; but it will be useful as anything we 
can do to see in what ways geometry is inadequate to serve as primordial 
building material . . . `geometry' is as far from giving an understanding of 
space as `elasticity' is from giving an understanding of a solid.'' Such 
phrases as ``curved space-time'' are but technical phrases that are not to be 
associated with real physical entities according to the basic ideas associated 
with the concept termed a ``pre-geometry.'' As demonstrated this NSPPM model may be the true objective reality.\par
 After the basic concepts that 
associate gravitational fields with geometric terminology are introduced, 
Marzke and Wheeler, using such concepts, describe the construction of a type 
of absolute clock they call the ``geometrodynamic clock.'' This clock is 
similar to my infinitesimal light-clock. However, the infinitesimal 
light-clock is a fundamental entity within the theory presented here and is not 
introduced after the model is constructed. The generation of certain line 
elements from physical considerations, {\bf not} related to the language of 
Riemannian geometry, is significant in that it tends to indicate that all 
physically meaningful line-elements that satisfy the Hilbert-Einstein 
gravitational equation could also be so generated. From these considerations, 
one could conclude that all of the consequences of the Einsteinian theory that 
apply to an actual physical universe should be interpreted in terms of 
infinitesimal light-clocks. That is; these conclusions alter such 
light-clock behavior. How does this influence our view of the actual relationship 
between gravitational fields and those physical entities the fields are 
predicted to affect?\par
In 1904 at the International Congress of Arts and 
Sciences held in St. Louis, Poincar\'e gave a talk entitled ``The present and 
future of mathematical physics.'' In his researches, he could not eliminate a 
certain constant c, the velocity of light in a vacuum, from any of his 
conclusions. He explained this fact as follows: ``(1) Either there exists 
nothing in the universe that is not of electromagnetic origin; (2) or, this, 
which is common to all physical phenomena, appears only because it relates our 
methods of measurements.'' As discussed above, the material in this book shows 
that Poincar\'e may have been partially correct. For, relative to this
re-interpretation and the NSPPM pre-geometry, it is probable that all 
natural-system behavior is related to the properties of electromagnetic 
radiation and this relation is developed by considering one and only one mode 
of measurement; the light-clock. \par
%baselineskip=12pt
What appears in this little book is that 
there is a measure that can be used as an model for alterations in certain 
behavior of the natural process called the propagation of electromagnetic 
radiation. This is what leads to alterations in natural-system behavior. The 
alterations can best be comprehended by introducing the nonstandard physical 
world (NSP-world) model. Further, verification of each prediction made by the 
theory presented here or the re-interpreted Einsteinian theory gives 
strong, albeit indirect, evidence that something like the NSP-world might exist 
in objective reality. However, rejection of the reality of the NSP-world does 
not preclude the use of such a concept as a model. \par
(18 JULY 1999) An external approach has been used to obtain the general line-elements discussed in this book. This means that no attempt has been made to associate changes of infinitesimal light-clocks with any modern field or particle theory other than the basic requirement that it all be subparticle controlled. There has been an approach that postulates a natural world vacuum electromagnetic zero-point field (ZPF) (Puthoff, 1989) and attributes 
the potential used above to obtain the Schwarzschild line-element to that induced by {\it Zitterbewegung} motion on a charged particle, where it is assumed that all matter is composed of charged particles and the {\it Zitterbewegung} motion is induced by the ZPF. This potential is associated directly with the kinetic energy associated with {\it Zitterbewegung} motion 
and this kinetic energy is considered as the gravitational mass effect. Of course, infinitesimal light-clocks are ideal approximators for natural world light-clocks. Haisch, Rueda and Puthoff (1997), state that the ZPF theory needs to be correlated to curved spacetime. If the results using ZPF can be related to light-clocks, then such a correlation would be, at least partially, achieved.\par    
[12] (18 JULY 1999)   The superscript and subscript $s$ represents local measurements about the $s$-point, using various devices, for laboratory standards and using infinitesimal light-clocks or approximating devices such as atomic-clocks. [Due to their construction atomic clocks are effected by relativistic motion and gravitational fields approximately as the infinitesimal light-clock's counts are effected.] All measures, rate of changes and the like, should be viewed via comparison. [See note 18.] The superscript or subscript $m,$ for the Special Theory, indicates how, with respect to the measures $s$, the motion of the $m$-point with a specific relative velocity yields physical behavior that differs from that at the $s$-point. For Special Theory, $\Pi_m$ for expressions (2) and (3) are obtained as follows: $t^{(m)}_E = (t_1 + t_3)/2.$ Then consider $ \Pi_m =(\overline{\Pi}^{(1)}_s + \overline{\Pi}^{(3)}_s)/2 - ({\Pi}^{(1)}_s + {\Pi}^{(3)}_s)/2).$ Now $\Pi_m$ need not be a member of $\nat^+_\infty.$ If any of the $\Pi$ is an odd hyperinteger, then replace it with $\Pi +1$. Then for nonzero $\st {u\Pi_m},$ each term when divided by 2 is a member of $\nat^+_\infty.$ Hence, in this case, $\Gamma_m$ can be used in place of $\Pi_m$ since 
$\st {u\Pi_m} = \st {u\Gamma_m}.$ The $\Gamma_m$ is an equivalent infinitesimal light-clock count. Indeed, the inf. light-clocks used could have a counter that incorporates this process. \par

Recall, that due to the presence of the NSPPM, there are absolute physical standards. These are what would be measured from a point $f$-point fixed in the NSPPM, where the NSPPM relative velocities when viewed in our physical world follow the unusual behavior indicated by equation (4) on page 48. (Additional material 22 NOV 2007 and extended on 24 JUL 2009, corrected on 8/25/09.) \par 

For the interpretation of the ``s'' and ``m'' for the General Theory see the appendix starting at page 93 on ``Gravitational Time-dilation.'' Previously, infinite Robinson numbers are used to model Einstein measures. Let the infinite number $\Pi$ denote the infinitesimal light-clock counts that yield, as an example, a time measure. The difference between two such measures is denoted by $\Pi'$. These differences can be a natural number or even, in some cases, an infinite number and in both cases $u\Pi' \in \mu (0).$ Then the symbol $dt^s$, generally, means $dt^s \approx u\Pi'$ which includes the interpretation $dt^s = u\Pi'$ in the NSP-world. But for this behavior to be realized in the physical world, then point behavior expressed in terms of such notation needs to be related to standard intervals as is done (6) and (7). Depending upon the usual function requirements that lead to the notions of the definite integral, when passing from the monadic NSP-world to the standard world, the behavior of $dt^s$ and the other differentials in a metric are considered as being infinitely close of order one to the term containing the corresponding altered differential. Further, since all functions are considered as continuous, a gravitational field is closely approximated by a constant field over ``small'' neighborhoods. This is why comparisons are made relative to the chronotopic or a transformed chronotopic line element.\par
[13] (9 AUG 1999) In my view, one of the possible causes for ambiguity or contradiction within the General Theory is that, although mathematically such transformations are allowed, physically they need to be interpreted via the theory of correspondence. It is this interpretation that leads to such difficulties. This does not occur in this analysis.\par
[14] (10 AUG 1999) It is assumed that it is the operator equation that will reveal the alterations in physical behavior. In order to determine these alterations, something needs to be fixed. The technique used is to investigate what physical behavior modifications would be needed so that the equational statements, that follow this note in the article, hold. For example, it can be argued that if the $R$ varies in ``time,'' then $R$ may be a factor of a universal function. In this case, $R$ may be represented by two possibly distinct functions $r^s$ and $R^m$ that have the property that  $r^s(t^s)= R^m(t^m).$ This technique is based upon the acceptance that a concept can be represented by a universal function. Notice, however, that this requirement is not very restrictive. It is only a ``form'' restriction. There may be confusion when the ``technique'' is applied to determine a ``relation'' between the standard infinitesimal light-clock measured quantities and the quantities that are measured by other altered infinitesimal light-clocks. It is this technique that allows for the equational statements $h(x_1^s,\ldots,x_n^s) = H(x_1^m,\ldots,x_n^m),\ f(t^s)= F(t^m)$ and, hence, the technique statement that $T(x_1^s,\ldots,t^s) = T(x_1^m,\ldots,t^m).$ It is usually a simple matter to argue for the universal function restriction. It's implicit in this work that the Special and General Theories need not be universal in application. The scenario and physical property investigated must satisfy all of the constraints of this technique; constraints that formally yield alterations in physical behavior. These theories may not apply to certain informational transmissions. On the other hand, since the selection of a line-element is often scenario dependent, the universal function or other aspects of this technique could be postulated. Such postulation would need to be verified by experimentation before actual acceptance could be considered. If verified, this would imply, at the least, one new property for the physical entity under investigation. \par
 Notice that the method used to derive such expressions as (13) and (14) are obtained by starting with the Minkowski Special Theory line-element. The potential velocity statement $(v + d)$ assumes the usual infinitesimal statement that $(v + d)$ behaves as if it is a constant over the monadic neighborhood. Using the fact that most functions, and especially continuous ones defined on a compact domain, can be approximated to any degree of approximation by step (constant) functions, for the monadic neighborhood $\mu(t_0),$ the same thing holds where a standard continuous function, say $v(t),$ is considered as behaving like the constant value  $v(t_0)$ over $\mu(t_0)$ (i.e. $\hyper v(t) \approx v(t_0).$) (Of course, if $v$ is differentiable, then $v(t)$ can also be considered as behaving like a linear function over $\mu(t_0).$) Using this result, the extension of the ``constant'' potential-type velocity statement to a piecewise continuous function is an acceptable modeling technique. \par
(3 FEB 2004) Rather than use $v + d$, for certain line-elements, a more complex monadic behavior function $f(v,d)$ may be necessary. For example, considering the vector $\langle v,d\rangle$ and the Euclidean norm $\Vert \langle v,d\rangle \Vert$ in place of $v+d$, the modified Schwarzschild line-element is obtained using one form for the cosmological constant $\Lambda.$\par 
[15] (16 AUG 1999) I thought that the following was obvious. But, it appears that a formal presentation is necessary. If the physical effect is the reverse of the P-process used to obtain a particular line-element, then this can be modeled by considering a complex velocity vector $((v + d)i +c)dt^s$ in expression (8) and within a monadic neighborhood, as it is done in many two dimensional fluid flow problems. All of the analysis that leads from equation (8) through and including equations (13) and (14) holds with only one alteration in the 
$\lambda.$ In this case, $\lambda = 1 + (v+d)^2/c^2.$ \par

The basic derivation method uses the notion of physical-like photon behavior.  I do not accept  unconstrained line element transformations. Transformations need to be constrained due to physical conditions, one of which is a reasonable photon position requirement. \par 
[16] (23 AUG 1999) From the  derivations for the line-elements ${\rm(13)_a,\ (13)_b,\  (14)_a,\  (14)_b}$, it is obvious that application of these line-elements must be carefully considered. Please note the behavior being investigated does not alter $\phi^s$ nor $\theta^s.$ Hence, $\phi^s(t^s) = \phi^m(t^m)$ and $\theta^s(t^s) = \theta^m(t^m)$ although these are not considered as universal functions.
These line-elements are only applicable to alterations in radial or linear behavior. Further, one must be convinced that they hold for a particular scenario. Notice that they can be written entirely in terms of the standard measures and then compared with the Minkowski form. For example, this yields for (13)$_{\rm a}$ the general identity $(cdt^s)^2 - (dR^s)^2 = \lambda (cdt^m)^2 - (1/\lambda)(dR^m)^2,$ which cannot be solved unless other conditions are met. The arguments, for the variations in measures, are
based upon the necessity to show that the line-elements apply and then to solve this expression for the variations in infinitesimal light-clock measures. \par
[17] (17 SEPT 1999) [Ref. note 12.] Notice that the expression $\gamma dt^m = dt^s$ (i.e. NSPPM proper time) used along with the universal function concept to obtain the alterations in physical behavior is a unique representation for one of the ambiguous forms discussed in Article 2 section 8. Although it represents an alteration in infinitesimal light-clock, by considering $u\,dt^m = (u/\gamma)dt^s,$ one \underbar{might comprehend} this alteration in terms of what \underbar{appears} to be an N-world alteration in the time unit.\par
[18] (31 DEC 1999)  In article 1, the unit adjusting number $u$ is used. To obtain this number, where it is assumed that no relativistic physical alterations occur, simply note that $u = L/c$ from the NSPPM viewpoint. Thus the $u$ units being considered are entirely related to the units used to express the $c$ velocity. Hence, it is not the ``units'' of measure that is being altered, but as shown in Special Theory note [2], it is a NSPPM hyperbolic velocity-space behavior that yields infinitesimal light-clock count alterations and alterations in the physical world.\par
Almost no mathematical structure involving real or complex numbers is a perfect model for the measured natural system behavior, macroscopic, large scale or otherwise. One must always constrain the predicted results by physical conditions and arguments. This often has to do with such ideas as negative length, negative mass, negative energy and the like. If you substitute the derived $\alpha$ and $\beta$ into equations (A) and (B) you 
get statements relative only to infinitesimals, in general, infinitesimals that, in some cases, would not lead to N-world effects. For example, the expression $dR^m/dT^m$ is the derivative of an assumed differentiable standard function $R^m(T^m)$ in terms of an altered distance and time measuring infinitesimal light-clocks. Technically, one needs to consider both positive and negative infinitesimals as well as all members of $\monad {0}.$  We have not used a physical argument for all infinitesimals, but have extended our results to all members of $\monad{0}$ by considering physical behavior for those infinitesimals obtained from infinitesimal light-clock behavior. When $\Hyper {dR^m/dT^m}$ is evaluated over a monadic neighborhood, we have set the standard part equal to zero. Thus, there is some standard real neighborhood where the rate of change of $R^m$ with respect to $T^m$ is zero. This, however, from the most basic aspect of relativity is what one would expect since alterations in the two infinitesimal light-clocks caused by the $P$-process would need to counter each other with respect to this rate of change. When it comes to relativistic alterations in natural system behavior, one is interested in comparing non-influenced effects with influenced effects.\par
[19] (7 JAN 2000) [Cross reference note [4]]. Concerning the existence of a subparticle (i.e. sub-quantum) region, Special Theory alterations are also distinct from those of gravitational alterations since any such sub-quantum region that might produce these alterations does so in such a manner that measuring devices only yield the same relative velocity measures and also, due to the countering of the alterations, no device using electromagnetic properties will reveal any constant linear motion through this sub-quantum region. This should not be the case with respect to acceleration, however.\par

The derivation of the linear effect line-element is different from that of the Schwarzschild and this fact is profound in its consequences. We argued for the Schwarzschild that $R^s = R^m.$ The same argument does not hold for the linear effect line-element. However, what does hold is the fact that this line-element deals entirely with but two measurements within the natural world ``time'' and ``length'' and how such measurements are altered by properties of the sub-quantum monadic cluster. The facts are that to measure a distance alteration relative to a moving object, the alteration would only appear if measurements where made simultaneously. Also the ``timing'' and ``distance'' measuring infinitesimal light-clocks are the same ``clock'' undergoing alterations in $c$. To incorporate this ``simultaneous'' measurement requirement into this line-element, we need to consider that when $dt^m=0,$ then $dt^s = 0$. In this case, the invariance leads to the requirement that $(dr^m)^2 = \lambda (dr^s)^2.$ Equating the two line-element forms, would require that we replace $dr^s$ in the derivation with a $dr^m$ to obtain (14)$_a$ and (14)$_b$. (A second justification for this requirement is that $dr^s/dt^s = dr^m/dt^m = v.$) \par

[20] (2/9/2005) This is relative to article 2.  Expression (6.9) (p. 37) is the standard two coordinate location ``length'' contraction expression (6.8) viewed using operational infinitesimal light-clocks relative to the corrected NSPPM hyperbolic diagram (p. 52). The basic scenario, as mentioned, is that when $F_1$ and $F_2$ coincide a single light pulse is sent to $P$ that does not coincide with $F_1,\ F_2$ (p. 52) In the diagram, since it is a position and velocity diagram, the ``lengths'' of the sides of the triangle correspond to the NSPPM velocities. To obtain (6.5), consider the scenario where $P$ coincides with $F_1$. This is equivalent to $\omega_1 = 0.$ This gives the stated values for $\theta$ and $\phi,$ where if standard values are considered the $\approx$ is replaced by $=$, and yields $L(\lambda^{(1)}-\eta^{(1)}) \approx 0.$\par
  The expression (6.8) is obtained by considering two coordinate locations for $F_1,\ F_2$ while all other aspects of the diagram remove fixed. That is, while the velocities $\omega_i = \overline{\omega}_i$ and $\overline{\theta} = \theta,\ \overline{\phi} = \phi.$ In general, the values of the NSPPM velocities $\omega_i, \ i=1,2,3,$ are not related and it is these values that determine the infinitesimal light-clock counts. Within the diagram, the $\omega_i,\ i=1,2,3,$ are related to the $\theta,\ \phi$ in that $\omega_1\cos \theta + \omega_2\vert \cos \phi \vert = \omega_3$. The Einstein measures of the relative velocities are, in general, not related. That is, they are, at least, relationally independent. If, in (6.9), $\hyper {\cos \theta}$ is not infinitesimal, then the proper definition of length also depends upon the diagrammed circumstances.\par\smallskip
{\leftskip=0.5in \rightskip=0.5in \noindent (a) Alterations in the values for $\theta,\ \phi$ simply imply that the values for the standard Einstein coordinates, such as $\st {x^{(i)}_{Ea}},\ i = 1,2,$ as well as length, must be altered under altered circumstances. \par}\par\smallskip
  \noindent The over-long analysis that leads to (6.12) or (6.13) (p. 38) shows how to obtain, for the infinitesimal light-clocks, the same ``form'' within the NSPPM as that of expression (6.8). \par

I have not given an argument that yields (6.14) on page 38. That is, $$u((\overline{\lambda}^{(1)}+\overline{\eta}^{(1)}) - (\lambda^{(1)}+\eta^{(1)})) \approx\beta 
u((\overline{\lambda}^{(2)}+\overline{\eta}^{(2)}) - (\lambda^{(2)}+\eta^{(2)})).\eqno (6.14)$$ 
\par
To derive (6.14), first consider the general expression
$${\overline{t}}^{(1)}_E -t^{(1)}_E= \st {\beta}(\overline{t}^{(2)}_E - 
t^{(2)}_E), \eqno [20.1]$$
and make the substitution from (6.3). This yields $$u[(\overline{\lambda}^{(1)}+\overline{\eta}^{(1)}) - (\lambda^{(1)}+\eta^{(1)})]
\approx (\beta u)[(\overline{\lambda}^{(2)}+\overline{\eta}^{(2)})(1- \overline{\alpha}) -(\lambda^{(2)}+\eta^{(2)})(1 -\alpha)],\eqno [20.2]$$
where finite $\alpha = K^{(3)}K^{(2)}\hyper {\cos\phi},\ \overline{\alpha} = \overline{K}^{(3)}\overline{K}^{(2)}\hyper {\cos\overline{\phi}}.$\par
Let $\omega_1,\ \overline{\omega_1},$ be fixed and consider the standard part of [20.2]. Then we have the contradiction that the infinitesimal light-clock measurement for the event at $P$ as determined by $F_1$ depends upon $\st{(\lambda^{(2)}+ \eta^{(2)})\alpha -(\overline{\lambda}^{(2)} + \overline{\eta}^{(2)})\overline{\alpha}}$. However, for this scenario, $\phi$ and $\overline{\phi},$ are arbitrary, where it is only required that $\pi/2 \leq \phi,\overline{\phi}\leq \pi.$ 
Thus  
$$[(\lambda^{(2)}+ \eta^{(2)})\alpha -(\overline{\lambda}^{(2)} + \overline{\eta}^{(2)})\overline{\alpha}]\approx 0. \eqno [20.3]$$ 
This yields (6.14). \par\smallskip

{\leftskip=0.5in \rightskip=0.5in \noindent (b) The infinitesimal light-clock expression (6.14) is universal in that it is not altered by the diagrammed circumstances.\par}\par\smallskip
It is claimed that expression (6.14) (i.e. [20.1]) written as 
(\dag) $\Delta t^m = \st {\beta}\Delta t^s$ indicates that ``rates of change,'' via a change in the ``time'' unit, are altered with respect to relative motion. However, more than rates of change might be altered and so as to eliminate the model theoretic error of generalization all individual alterations in behavior need to be derived. Technically, such an expression as (\dag) \underbar{cannot be} physically infinitesimalized for \underbar{any} clock. Expressions $(5)_a$ and $(5)_b$ (p. 56) are relative to light propagation via the chronotopic interval and along with the linear-effect line-element leads to (**) (p. 62). It is the light propagation determined linear effect line-element that yields the infinitesimalized version of (\dag), where $\gamma = \beta^{-1}, \ d=0.$ It is the line-element method that displays the correct alterations, alterations that are employed to derived altered physical behavior, where the basic alteration is the alteration of $c$ in the inf. light-clocks. This method leads to alterations in measures of mass, energy, etc. \par

As previously mentioned, Einstein originally used partial derivatives in his derivation and, regardless of the derivation method used, the differential calculus is the major classical tool. Einstein had to synchronize his clocks using ``light'' signals via the ``radar' method. This and other requirements signify that this theory and any theory that must reduce to the Special Theory, under certain circumstances, is but a light propagation theory. \par

In order to apply the differential calculus to any physical measure in a reasonably accurate many, especially for macroscopic behavior, one uses an ideal physical entity and measures that can be``infinitesimalized.'' Almost no ``clock'' corresponds to an ``ideal'' clock that can be infinitesimalized. Almost no clock mechanism that displays standard time can be ``smoothed-out" and infinitesimalized in a reasonable manner. But, the basic notion of ``length'' is smoothed-out and infinitesimalized. Further, since the theory uses the language of light (i.e. electromagnet) propagation, then one should not generalize beyond this language. There is one physical mechanism that uses light propagation as its basic mechanism and can be infinitesimalized relative to lengths of light-paths. The clock is the light-clock.\par     
One needs only consider two viewpoints relative to the implications of this theory to see that the use of the standard methods and incorrectly generalizing leads to contradictory yet viable philosophic views.\par
Dingle (1950) considers the contraction of ``length'' as physical nonsense. He states relative to the usual expression for length-contraction $(\overline{x}'-x') =(\overline{x}-x)(1-v^2/c^2)^{-(1/2)}$, ``The implication of this choice is often expressed by the statement that a body contracts on moving, but the expression is unfortunate: it suggests that something happens to the body, whereas the `movement' may be given it merely by our mental change of the standard of rest, and we can hardly suppose that the body shrinks on becoming aware of it'' (p. 30). He also rejects the notion that``space'' changes and contends that ``. . . our province is simply that of physical measurements, and our object is simply to relate them with one another accurately and consistently . . . . this is completely achieved by a re-definition of length . . .'' (p. 31). All this comes about since it claimed that ``length is not an intrinsic property of the body'' (p. 30). But these remarks assume that there is no \ae{ther}, no privileged observer with privileged frame of reference. \par

As to time-dilation, Dingle claims that this comes about only due to the way science defines velocity, a defining method that need not be used. A magnetic form of speedometer on a specific vehicle can simply be marked off in speedometer units as a measure of velocity is one of his examples. One then uses the length-contraction statement and obtains the necessary time-dilation expression. Indeed, he claims for this time expression the following: ``A very familiar expression of this result is that statement that the rate of a clock is changed by motion, and by this we are intended to understand that some physical change occurs in the clock. How false this is can be seen, just as the falsity of the corresponding statement for space-measuring rods, by remembering that we can change the velocity of the clock merely by changing our minds'' (p. 39-40).  He does not mention and explain in this book, using his definition notion, the Ives-Stillwell 1938 experiment (frequency changes in emitted ``light'' from hydrogen canal rays) nor changes in decay rates that are attributed to time-dilation. \par 
Then we have Lawden's (1982) statements about length contraction. ``The contraction is not to be thought of as the physical reaction of the rod to its motion and as belonging to the same category of physical effects as the contraction of a metal rod when it is cooled. It is due to a changed relationship between the rode and the instruments measuring its length. . . . It is now understood that length, like every other physical quantity, is defined by the procedure employed for its measure and it possesses no meaning apart from being the result of this procedure. [Notice the use of the term ``quantity.'' Physical properties exist in reality. But, the properties are distinct from the methods used to measure the properties.] . . .  [I]t is not surprising that, when the procedures must be altered to suit the circumstances, the result will also be changed. It may assist the reader to adopt the modern view of the Fitzgerald contraction if we remark that the length of the rod considered above can be alerted at any instant simply by changing our minds and commencing to employ the $S$ frame rather than the $S'$ frame. Clearly, a change of mathematical description can have no physical consequences'' (p. 12).\par

Lawden takes the observer time-dilation expression as the one that has physical significance ``. . . all physical processes will evolve more slowly when observed from a frame relative to which they are moving'' (p.13). However, no further explanation is possible for this effect since once again no \ae{ther} is assumed.  But, are these time-dilation alterations actual mean physical changes in the measuring machines or are they but observational illusions? \par

It appears significant that there was no rigorous basis for the philosophic stances of Dingle and Lawden. This all changes with the use of infinitesimal light-clock theory and this note. When infinitesimal light-clocks are used to measure the length of a rod with respect to linear N-world relative motion, (a) shows that there is no fixed \underbar{N-world} expression for such length ``contraction'' independent of the parameters.  The idea that the expressions simply imply altered ``definitions'' is viable, but only for the N-world.\par

On the other hand, with respect to linear relative motion, (b) shows that the alterations in the infinitesimal light-clocks used to measure ``time'' are independent from the circumstances (i.e. the two parameters) and they are a universal requirement in that the time-dilation model determines what physical changes occur. Hence, (b) rejects the Dingle notion that (infinitesimal light-clock measured) time-dilation is simply a problem of measure and definition and verifies that it must represent actual physical or observed alterations in behavior if such a notion is actually applicable to the physical world. This ``time'' alteration is not in absolute time but rather in observer time and is related to the NSPPM \ae{ther} via $L/u =c.$ The linear-effect line-element is used to derive Special Theory alterations in behavior. But, using the line-elements associated with the General Theory, the exact same infinitesimal time-light clock differential expression is obtained as used for the Special Theory. The exact same derivations yield the gravitational relativistic alterations in behavior where the ``velocities'' are but potential. Few doubt that the gravitational alterations are physically real. This yields additional very strong evidence that similarly derived Special Theory alterations are physically real, at the least, with respect to the NSPPM, while infinitesimal light-clock alterations lead to physical manifestations within the natural-world. \par

[21] (25 NOV 2007) (a) It is assumed that this approach is a strong classical approach. As such, solutions to the line-elements presented are related to, at least, the Riemann integral and the additional requirement that they be integrable over a closed ``interval'' and continuous at the point within the interval under consideration. As shown by Theorem 5.1.1 in Herrmann (1994a), such integrals are independent from the differential chosen. Thus, restricting  each $d\xi$ to specific $L\Pi^\prime$ or $u\Pi^\prime$ infinitesimals, where $\Pi^prime \in \hypernat,$ is sufficient. \par
(b) It is clear from statement (*) as it relates to (9.2) and the actual observed behavior that (*) predicts that the idea of the infinitesimal light-clock is the correct approach to relativistic physics and the only one that includes the required light propagation language. This includes the General Theory.\par

[22] (18 JUL 2009) (a) Since we are comparing results at different locations in a gravitational field or flat-space and the medium, the $L$ and $u$ do not change, then Planck's constant is not altered from what it would be at an s-point or its assumed value at the m-position (or other ``viewing'' positions). \par
(b) It should not be forgotten that the times used for the Special Theory line elements is Einstein time. The alterations in inf. light-clock measurements are required to maintain a NSPPM hyperbolic velocity-space. Notice that local (or flat-space) measurements or equivalent will always yield the constant $c$. This is modeled by (5)$_a$, where $dS = 0$.\par
(c) For gravitational effects, each inf. light-clock count alteration follows from $\gamma dt^m = \gamma u \Pi_m = dt^s = u \Pi_s,\ \Pi_s,\Pi_m \in \hypernat.$ Since $u$ does not vary, then, from light-clock construction transferred to inf. light-clocks, the only physical way that this can happen is that, for the s-clock, the actual velocity of light within the clock is $\gamma c.$ Since inf. light-clocks are used, this implies that such alterations in behavior probably correspond to subatomic photon behavior. Notice that for the Robinson-Walker line element, the alteration follows from $c\gamma(1/\lambda.)$  \parm   

\centerline{\bf References for article 3.}\parm
\id{B}arnes, T. G. 1983. Physics of the future. ICR El Cajon, CA.
\id{B}arnes, T. G. and R. J. Upham. 1976. Another theory of gravitation: an 
alternative to Einstein's General Theory of Relativity. {\it C. R. S. 
Quarterly} 12:194-197. 
\id{B}ergmann, P. 1976. Introduction to the theory of relativity. Dover, 
New York. 
\id{B}reitner, T. C. 1986. Electromagnetics of gravitation. 
{\it Speculations in Science and Technology} 9:335--354.
\id{B}uilder, G. 1960. Lecture notes for a colloquium held at the 
University of New 
South Wales.
\id{E}ddington, A. S. 1924. A comparison of Whitehead's and Einstein's 
formulas. {\it Nature} 113:192.
\id{E}ttari. 1988. Critical thoughts and conjectures concerning the Doppler 
effect and the concept of an expanding universe -- part I. {\it C. R. S. Quarterly} 
25:140--146.       
\id{} {\vbox{\hrule width 0.75in}} Critical thoughts and conjectures concerning the Doppler 
effect and the concept of an expanding universe -- part II. {\it C. R. S. Quarterly} 
26:102--109.       
\id{F}inkelstein, D. 1958. Past-future asymmetry of the gravitational field of a 
point particle. {\it Phys. Rev.} 110:965--967.
\id{F}ock, V. 1959. The Theory of Space Time and Gravity. Pergamon Press. New 
York.
\id{G}reneau, P. 1990. Far-action versus contact action. {\it 
Speculations in Science and Technology} 13:191--201.
\id{G}uth, A. H. 1981. Inflationary universe: a possible solution to the 
horizon and flatness problem. {\it Physical Review D} 23(2):347-356.
\id{H}aisch, B. H, A. Rueda and H. E. Puthoff. 1997. Physics of the zero-point field: implications for inertia, gravitation and mass. {\it 
Speculations in Science and Technology} 20:99--114.
\id{H}eaviside, O. 1922. Electromagnetic theory. Vol(1). Benn Brothers, 
Limited, London.
\id{H}errmann, R. A. 1985. Supernear functions. {\it Mathematica Japonica} 
30:169--185. 
\id{} {\vbox{\hrule width 0.75in}} 1986a. Developmental paradigms. {\it C. R. S. Quarterly} 22:189--198. 
\id{}{\vbox{\hrule width 0.75in}} 1986b. D-world evidence. {\it C. R. S. Quarterly} 23:47--54. 
\id{}{\vbox{\hrule width 0.75in}}
 1988. Physics is legislated by a cosmogony. {\it 
Speculations in Science and Technology} 11:17--24.
\id{}{\vbox{\hrule width 0.75in}}
 1989. Fractals and ultrasmooth microeffects. {\it 
J. Math. Physics}. 30:805--808. \par
\id{}{\vbox{\hrule width 0.75in}}
 1990. World views and the metamorphic model: their 
relation
 to the concept of variable constants. {\it C. R. S. Quarterly} 27:10--15.
\id{}{\vbox{\hrule width 0.75in}} 1991a {Some} {application} 
{of} {nonstandard} {analysis} {to} 
{undergraduate} {mathematics:} {infinitesimal} {modeling} {and} 
{elementary} {physics}, {Instructional} {Development} 
{Project}, {Mathematics} {Department}, U. S. Naval 
Academy, Annapolis, MD, 21402-5002. http://arxiv.org/abs/math/0312432 
\id{}{\vbox{\hrule width 0.75in}}
 1991b. The Theory of Ultralogics\hfil\break http://www.arxiv.org/abs/math.GM/9903081\hfil\break 
http://www.arxiv.org/abs/math.GM/9903082.
\id{}{\vbox{\hrule width 0.75in}} 1992. 
 A corrected derivation for the Special Theory of 
relativity. Presented before the Mathematical Association of America, Nov. 14, 
1992 
at Coppin State College, Baltimore, MD.
\id{}{\vbox{\hrule width 0.75in}} 
 1994(a). The Special Theory and a nonstandard substratum. {\it 
Speculations in Science and Technology} 17(1):2-10.
\id{}{\vbox{\hrule width 0.75in}} 
1994(b). A Solution to the ``General Grand Unification Problem'' or How to Create 
a Universe. Presented before the Mathematical Association of America, Nov. 12
1994 at Western Maryland College, Westminister MD. 
http://www.arxiv.org/abs/astro-ph/9903110
\id{}{\vbox{\hrule width 0.75in}} (1995) 1996. An operator equation and relativistic 
alterations in the time for radioactive decay. {\it Intern. J. Math. \& Math. 
Sci.} 19(2):397-402
\id{}{\vbox{\hrule width 0.75in}} 1995. Operator equations, separation of variables and relativistic alterations. {\it Intern. J. Math. \& Math. 
Sci.} 18(1):59-62.
\id{}{\vbox{\hrule width 0.75in}} 1997. A hypercontinuous, hypersmooth Schwarzchild line-element transformation. {\it Intern. J. Math. \& Math. 
Sci.} 20(1):201-204.
\id{}{\vbox{\hrule width 0.75in}} 1999. 
The NSP-world and action-at-a-distance, In ``Instantaneous Action at a Distance in Modern Physics: `Pro' and `Contra,' ''
Edited by Chubykalo, A.,  N. V. Pope and R. Smirnov-Rueda,
(In CONTEMPORARY FUNDAMENTAL PHYSICS - V. V. Dvoeglazov (Ed.))
Nova Science Books and Journals, New York: 223-235.
\id{H}umphreys, R. 1994. Starlight and time: solving the puzzle of distant 
starlight in a young universe. Master Books, Colorado Springs, Colorado.
\id{I}ves, H. 1932. Derivation and  significance of the so-called 
``Chronotopic interval'', {\it J. Opt. Soc. Am.} 29:294-301. 
\id{I}ves, H. and G. Stillwell. 1938. Experimental study of the rate of a 
moving atomic clock, {\it J. Opt. Soc. Am.} 28:215--226.
\id{J}efferson, S. N. 1986. Is the source of gravity field an 
electromagnetic mechanism?  {\it Speculations in Science and Technology} 
9:334.
\id{L}andau, L. D. and E. M. Lifshitz. 1962. The classical theory of fields. 
Transl. by M. Hamermesh, Pergamon Press, Oxford and Addison-Wesley 
Publishing Co. Reading MA.
\id{L}awden, D. F. 1982. An introduction to tensor calculus, relativity 
and cosmology. John Wiley \& Sons, New York. 
\id{M}isner, C. W., K. S. Thorne and J. A. Wheeler. 1973. Gravitation.
W. H. Freeman and Co., San Francisco. 
\id{P}hillps, H. B. 1922. Note on Einstein's theory of gravitation. {\it 
{\rm(}MIT{\rm )} 
Journal of Mathematics and Physics}  I:177--190. 
\id{P}ohl, H. 1967. Quantum mechanics for science and engineering, Prentice-Hall, Englewood Cliffs, NJ.
\id{P}uthoff, H. E. 1989. Gravity as a zero-point-fluctuation force.
{\it Phys. Rev. A} 39(5):2333--2342. 
\id{O}hanian, H. C. 1976. Gravitation and Space-Time. W. W. Norton \& Company, New 
York.
\id{P}rokhovnik, S. J. 1967. The logic of Special Relativity. Cambridge 
University Press, Cambridge.
\id{R}indler, W. 1977. Essential relativity. Springer-Verlag, New York.
\id{S}chneider, H. 1984. Did the universe start out structured. {\it C. R. S. Quarterly}
21:119--123.
\id{S}urdin, M. 1962. A note on time-varying gravitational potentials. 
{\it Proc. Cambridge Philosophical Society} 58:550-553.\par\vfil\eject
\centerline {\bf NOTES}
\vfil\eject

\centerline{\bf Appendix-B}
\centerline{\bf Gravitational Time-dilation.}
\bigskip 
\noindent{\bf 1. Medium Time-dilation Effects.}\parm

\noindent {\it Within a gravitational field, the superscript $s$ represents local measurements taken at a spatial point $Q$ using specific devices. Using identically constructed devices, the superscript $m$ represents local measurements taken at an m-point $Q_1$ or at $+\infty$, where there is no gravitational field. These local m-point measures are compared with the local measures made at point $Q$.}\par 
Equation (**) on page 62, when expressed for the general physical metric for a fixed spatial point is
$$\sqrt{g_1} dt^m= dt^s. \eqno (B1)$$
(Note that many authors denote $g_1$ as $g_4$ and the term ``clock'' is specifically defined.) Consider another spatial point $R$ within the gravitational field, where for the two points $P,\ R$ the expression $g_1$ is written as $g_1(P),\ g_1(R),$ respectively. Considering the point effect at each pont and applying the relativity principle, this gives, in medium $t_3$ time, that
$${{\Delta t^s_P}\over{\sqrt{g_1(P)}}}=\Delta t^m = {{\Delta t^s_R}\over{\sqrt{g_1(R)}}} \eqno (B2)$$
Equations like (B2) are comparative statements. This means that identical laboratories are at $P$ and $R$ and they employ identical instrumentation, definitions, and methods that lead to the values of any physical constants. Since infinitesimal light-clocks are being used, standard ``clock'' values can take on any non-negative real number value. The $\Delta t^s_P,\ \Delta t^s_R$ represent the comparative view of the gravitationally affected ``clock'' behavior as observed from the medium where there are no gravitational effects. (The $1/\sqrt{g_1}$ removes the effects.) 
  Assume a case like the Schwarzschild metric where real $\sqrt {g_1} <1$. Consider two different locations $P, R$ along the radius from the ``center of mass.'' (A ``tick'' is a one digit change in a light-clock counter. As discussed below, there can be ``portions'' of a tick.) Then there is a constant $r_s$ such that 
$$\sqrt {1- {{r_s}\over{r_P}}}\ \Delta t_{R}^s\ ({\rm in\ R\!\!-\!ticks}) = \sqrt {1- {{r_s}\over{r_R}}}\ \Delta t_{P}^s \ ({\rm in \ P\!\!-\!ticks}). \eqno (B3)$$
where $r_s \leq r_P,\ r_R.$
[The cosmological ``constant'' $\Lambda$ modification ($\Lambda$ is not assumed constant) is 
$$\sqrt {1 - {{r_s}\over{r_P}} - (1/3)\Lambda{{r_P^2}\over{c^2}}}\ \Delta t_{R}^s =\sqrt {1 - {{r_s}\over{r_R}} - (1/3)\Lambda_1{{r_R^2}\over{c^2}}}\ \Delta t_{P}^s.]\eqno (B3)'$$\par

Of course, these equations [(B3), (B3)'] are comparisons that must be done with the same type of ``clocks.'' As an example for (B3), suppose that $r_s/r_P = 0.99999$ and $r_R = 100,000r_P.$ Then $r_s/r_R = .000009999$. This gives 
$0.003162278 \Delta t_{R}^s = 0.999995 \Delta t_{P}^s$. Hence, $\Delta t_{R}^s = 316.2262 \Delta t_P^s.$ Thus, depending upon which ``change'' is known, this predicts that ``a change in the number of R-ticks'' equals ``316.2262 times a change in the number of P-ticks.'' Suppose that at $P$ undistorted information is received. Observations of both the ``P-clock'' and the ``R-clock'' digit changes are made. (The fact that it takes ``time'' for the information to be transmitted is not relevant since our interest is in how the ticks on the ``clocks'' are changing.) Hence, if the ``clock'' at $P$ changes by 1-tick (the 1-tick), then the change in the R-ticks is 316.2262. (What it means to have a ``portion'' of a tick is discussed later in this appendix.)\par

The careful interpretation of such equations and how their ``units'' are related is an important aspect of such equations since (B2) represents a transformation. Using a special `` clock'' property, if the ``R-clock'' changes its reading by 1, then at $P$ the ``P-clock'' shows that only 0.003162 ``P-clock'' time has passed. If you let $R = \infty,$ then $\Delta t_{R}^s = 316.2278 \Delta t_{P}^s$ and, in a change in the reading of 1 at $P$, the R-reading at $\infty$ is 316.2278.  Is this an incomprehensible mysterious results? No, since it is shown previously, the gravitational field is equivalent to a type of change in the infinitesimal light-clock itself that leads to this result. But, for our direct physical world, thus far, the answer is yes if there is no physical reason why our clocks would change in such a manner. The equation (B3) [$(B3)'$] must be related to physical clocks within the physical universe in which we dwell.\par

\noindent{\bf 2. The Behavior of Physical Clocks.}\parm

\noindent Einstein did not accept general time-dilation for the gravitation redshift but conjectured that such behavior, like the gravitational redshift, is caused by changes within atomic structures rather than changes in photon behavior during propagation. This was empirically verified via atomic-clocks. To verify Einstein's conjecture theoretically and to locate the origin of this atomic-clock behavior, the comparative statement that $dt^s = \sqrt {g_1} dt^m$ is employed. It has been shown using time-dependent 
Schr\"{o}dinger equation, that certain significant energy changes within atomic structures are altered by gravitational potentials.  Once again, consider identical laboratories, with identical physical definitions, physical laws, construction methods etc. at two points $P$ and $R$ and within the medium. When devices such as atomic-clocks are used in an attempt to verify a statement such as (B3), the observational methods to ``read'' the clocks are chosen in such a manner  that any known gravitational effects that might influence the observational methods and give method-altered readings is eliminated. For point $P$, let $E_P^s,$ denote measured energy. In all that follows, comparisons are made.  Using the principle of relativity, the following equation (B4) (A) holds, in general,  and if $g_1$ is not time dependent, then (B) holds. 
$${\rm (A)}\ \sqrt {g_1(P)}dE_P^s = dE^m= \sqrt {g_1(R)}dE_{R}^s , \ {\rm (B)}\ \sqrt {g_1(P)}\Delta E_P^s = $$ $$\Delta E^m =\sqrt {g_1(R)}\Delta E_{R}^s.\eqno (B4)$$ 
This is certainly what one would intuitively expect. It is not strange behavior. Hence, in the case that $g_1$ is not time dependent, then   
$$\sqrt {g_1(P)}\Delta E_P^s = \sqrt {g_1(R)}\Delta E_{R}^s,\eqno (B5)$$\par

For this application, equation (B4) corresponds to the transition between energy levels relative to the ground state for the specific atoms used in atomic-clocks. But, {\it for this immediate approach, the atomic structures must closely approximate spatial points.} Further, at the moment that such radiation is emitted the electron is considered at rest in the medium and, hence, relative to both $P$ and $R$. The actual aspect of the time-dependent Schr\"{o}dinger equation that leads to this energy relation is not the spatial ``wave-function'' part of a solution, but rather is developed from the ``time-function part.''  \par

{\it The phrase ``measurably-local'' means, that for the measuring laboratory the gravitational potentials are considered as constants.} Diving each side of (B) in (29) by Planck's (measurably-local) constant in terms of the appropriate units, yields for two observed spatial point locations $P, \ R$ that
 $$\sqrt {g_1(P)}\nu_P^s = \sqrt {g_1(R)}\nu^s_{R}, \eqno (B6)$$
 Equation (B6) is one of the expressions found in the literature  for the gravitational redshift [6, p. 154] but (B6) is relative to medium ``clocks.'' Originally, (B6) was verified for the case where $\sqrt {g_1(R)} \ll \sqrt {g_1(P)}$ using a physical clock. Note that since the $P$ and $R$ laboratories are identical, then the numerical values for $\nu^s_P$ and $\nu^s_R$ as measured using the altered medium ``clocks'' and, under the measurably-local requirement, are identical. Moreover, (B6) is an identity that is based upon photon behavior as ``clock'' measured. \par

What is necessary is that a comparison be made as to how equation (B6) affects the measures take at $R$ compared to $P,$ or at $P$ compared to $R$.  Suppose that $\vert \nu^s_A \vert$ indicates the numerical value for $\nu^s_A$ at any point $A$.  To compare the alterations that occur at $P$ with those at $R$, $\vert \nu^s_R \vert_P$ is symbolically substituted for the $\nu^s_P$ and the expression $\sqrt {g_1(P)}\vert \nu^s_R \vert_P = \sqrt {g_1(R)}\nu^s_{R}$ now determines the frequency alterations expressed in $R$  ``clock'' units.  As will be shown for specific devices, this is a real effect not just some type of illusion. {\it This substitution method is the general method used for the forthcoming ``general rate of change'' equation.} As an example, suppose that for the Schwarzschild metric $R = \infty$ and let $\vert \nu^s_R \vert = \nu_0.$ Then $\nu^s_\infty = \nu = \sqrt {1 - r_s/r_P}\,\nu_0.$ This result is the exact one that appears in [1, p. 222]. 
However, these results are all in terms of the behavior of the ``clocks'' and how their behavior ``forces'' a corresponding alteration in physical world behavior and not the clocks used in our physical world. These results need to be related to physical clocks. \par

Consider atomic-clocks. At $P$, the unit of time used is related to an emission frequency $f$ of a specific atom. Note that one atomic-clock can be on the first floor of an office building and the second clock on the second-floor or even closer than that. Suppose that the identically constructed atomic-clocks use the emission frequency $f$ and the same decimal approximations are used for all measures and $f$ satisfies the measurably-local requirement. The notion of the ``cycle'' is equivalent to ``one complete rotation.'' For point-like particles, the rotational effects are not equivalent to gravitational effects [8, p. 419] and, hence, gravitational potentials do not alter the ``cycle'' unit C. Using the notation ``sec.'' to indicate a defined atomic-clock second of time, the behavior of the $f$ frequency relative to the ``clocks'' requires, using equation (B6), that 
$$\sqrt {g_1(P)}{{1\r C}\over{{\rm P\!\!-\!sec.}}}=\sqrt {g_1(R)}{{1\r C}\over{{\rm R\!\!-\!sec.}}}.\eqno (B7)$$
$$\sqrt {g_1(P)}{{1}\over{{\rm P\!\!-\!sec.}}}= \sqrt {g_1(R)}{{1}\over{{\rm R\!\!-\!sec.}}}.\eqno (B7)'$$
For measurably-local behavior, this unit relation yields that  
$$\sqrt {g_1(P)}(\overline{t_R} - t_R)({\rm R\!\!-\!sec.})= \sqrt {g_1(R)}(\overline{t_P} - t_P)({\rm P\!\!-\!sec.}).\eqno (B8)$$
Hence, in terms of the atomic-clock seconds of measure 
$$\sqrt {g_1(P)}\Delta t_R= \sqrt {g_1(R)}\Delta t_P.\eqno (B9)$$\par
Equation (B9) is identical with (B3), for the specific $g_1$, and yields a needed correspondence between the ``clock'' measures and the atomic-clock unit of time.   Corresponding ``small'' atomic structures to spatial points, if the gravitational field is not static, then, assuming that the clocks decimal notion is but a consistent approximation,  (B9) is replaced by a (B3) styled expression 
$$ \sqrt {g_1(P,t_P)}dt^s_R = \sqrt {g_1(R,t_{R})}dt^s_{P}\eqno (B10)$$
and when solved for a specific interval correlates directly to atomic-clock measurements. Also, the Mean Value Theorem for Integrals yields 
$$\sqrt {g_1(R,t'_P)}(\overline{t_R} - t_R) = \sqrt {g_1(R,t'_{R})}(\overline{t_{P}} - t_{P}), \eqno (B11)$$ 
for some $t'_P \in [t_P,\overline{t_P}],\ t'_{R} \in [t_{R},\overline{t_{R}}].$ 
Equations (B9), (B10) and (B11) replicate, via atomic-clock behavior, the exact ``clock'' variations obtained using the medium time, but they do this by requiring, relative to the medium, an actual alteration in physical world photon behavior. The major interpretative confusion for such equations is that the ``time unit,'' as defined by a specific machine, needs to be considered in order for them to have any true meaning. As mentioned, the ``unit'' notion is often couched in terms of ``clock or observer'' language. The section 1 illustration now applies to the actual atomic-clocks used at each location. 
\par

For quantum physically behavior, how any such alteration in photon behavior is  possible depends upon which theory for electron behavior one choices and some accepted process(es) by which gravitational fields interaction with photons. Are the alterations discrete or continuous in character? From a quantum gravity viewpoint, within the physical world, they would be discrete if one accepts that viewpoint. This theoretically establishes the view that such changes are real and are due to ``the spacings of energy levels, both atomic and nuclear, [that] will be different proportionally to their total energy'' [3, pp. 163-164]. Further, ``[W]e can rule out the possibility of a simple frequency loss during propagation of the light wave. . . .'' [8, p. 184]. \par

Although the time-dependent Schr\"{o}dinger equation applies to macroscopic and large scale structures via the de Broglie ``guiding-wave'' notion, the equation has not been directly applied, in this same manner, to such structures since they are not spatial points. However, it does apply to all such point-approximating atomic structures since it is the total energy that is being altered. One might conclude that for macroscopic and large scale structures there would be a cumulative effect for a collection of point locations. (From the ``integral'' point of view relative to material objects, such a cumulative procedure may only apply, in our physical world, to finitely many such objects.) Clearly, depending upon the objects structure, the total effect for such objects, under this assumption, might differ somewhat at different spatial points.  However, the above derivation that leads to (B9) is for the emission of a photon ``from'' an electron and to simply extend this result to all other clock mechanisms would be an example of the model theoretic error in generalization unless some physical reason leads to this conclusion.\par

As mentioned, equation (B9) is based upon emission of photons. Throughout all of the atomic and subatomic physical world the use of photon behavior is a major requirement in predicting physical behavior, where the behavior is not simply emission of the type used above. This tends to give more credence to accepting that, under the measurably-local requirement, each material time rate of change, where a physically defined unit U that measures a $Q$ quality has not been affected by the gravitational field, satisfies 
$$\sqrt {g_1(P)}\Delta Q_P {\rm\ in\ a \ P\!\!-\!sec.} = \sqrt {g_1(R)}\Delta Q_R\ {\rm\ in\ an \ R\!\!-\!sec.},$$
 $$\Delta Q_R\ {\rm\ in\ an \ R\!\!-\!sec.}={{\sqrt {g_1(P)}}\over{\sqrt {g_1(R)}}}\vert \Delta Q_R\vert_P, \eqno (B12)$$
Equations (B12) give comparative statements as to how gravity alters such atomic-clock time rates of change including rates for other types of clocks.\par

 Prior to 1900, it was assumed that a time unit could be defined by machines that are not altered by the earth's gravitational field. However, this is now known not to be fact and as previously indicated, such alterations in machine behavior is  probably due to an alteration in photon behavior associated with a NSPPM source that undergoes two types of physical motion, uniform or accelerative. There is a NSPPM process that occurs and that alters photon behavior as it relates to the physical world. These alterations in how photons physically interact with atomic structures and gravitational fields is modeled (mimicked) by the defining machines that represent the physical unit of time, when the mathematical expressions are interpreted. The observed accelerative and relative velocity behavior is a direct consequence of this process. As viewed from the NSPPM, every a nonzero uniform velocity obtained from a zero velocity first requires acceleration. This is why the General Theory and the Special Theory are infinitesimally close at a standard point. \par
It is claimed by some authors that regular coordinate transformations for the Schwarzschild solution do not represent a new gravitational field but rather allows one to investigate other properties of the same field using different modes of observation. When such transformations are discussed in the literature another type of interpretation appears necessary [6, p. 155-159]. Indeed, what occurs is that the original Schwarzschild solution is rejected based upon additional physical hypotheses for our specific universe that are adjoined to the General Theory. For example, it is required that certain regions not contain physical singularities under the hypotheses that physical particles can only appear or disappear at chosen physical ``singularities.'' Indeed, if these transformations simply lead to a more refined view of an actual gravitational field, then the conclusions could not be rejected. They would need to represent actual behavior. One author, at least, specifically states this relative to the Kruskal-Szekeres transformation. In [10, p. 164], Rindler rejects the refined behavior conclusions that would need to actually occur within ``nature.''  ``Kruskal space would have to be {\it created in toto}: . . . . There is no evidence that full Kruskal spaces exist in nature.'' \par
One way to interpret the coordinate transformation that allows for a description of ``refined'' behavior is to assume that such described behavior is but a ``possibility'' for a specific gravitational field and that such behavior need not actually occur. This is what Rindler appears to be stating. Of course, such properly applied coordinate transformations also satisfy the Einstein-Hilbert gravitational field equations. For the medium view of time-dilation, this leads to different alterations in the atomic-clocks for any collection of such ``possibilities.'' \par 

These results, as generalized to the behavior exhibited by appropriate physical devices, imply that no measures using these devices can directly determine the existence of the medium. Although Newton believed that infinitesimal values did apply to ``real'' entities and, hence, such measures exist without direct evidence, there is a vast amount of indirect evidence for existence of such a medium.\parm

\noindent {\bf 3. Infinitesimal Light-Clocks, in General.}\parm

For applications, everything in the infinitesimal-world is composed of ``simple'' Euclidean or physical notions. There are no curves or curved surfaces and the like only objects that are ``linear'' in character. This follows from the nonstandard version of the Fundamental Theorem for Differentials in terms of a linear operator. For physical problems, when one attempts to find the appropriate nonstandard approach the closer the approach approximates physical behavior the more likely it will yield acceptable results. (Sometimes the approach used today by those that do not use the formal infinitesimals leads to statements that may appear to be mathematical but they are not. A foremost example of this is the Feynman integral.)  \par

Using a photon language, we further analyzed light-clocks. In the case being considered, light-clocks have a light source and two reflecting surfaces (A) and (B). At one ``mirror'' (A) a very short pulse (of photons) is emitted. The pulse is reflected at (B) and returns to (A). A detector at (A) registers its return. Then immediately the pulse returns to (B), etc. Of course, the pulse may need to be replenished with a new pulse after a while. Light-clocks are used since they mimic the behavior of a type of Einstein time for stationary objects.\par

The distance between the mirrors is $M.$ The distance traveled by a very short pulse is very nearly $2M.$ Within our observable world, the $M$ cannot be any positive small number. Is there an ``assumed'' physical process that uses a ``very small length $L,$'' as compared to a standard meter, that closely approximates this process?\par

Consider a theoretical reason why an electron and a proton are kept within a close range in the hydrogen atom. Photons are emitted by the proton, absorbed by the electron, and emitted by the electron and absorbed by the proton. One could assign a general approximate $L$ to the ``distance'' between an electron and the proton and, hence, the distance each interacting photon covers. Due to the linear requirements for the infinitesimal-world, a viewed from this world often requires idealization. (This does not mean that the NSPPM view is not what might actually be happening.) In this case, a linear photon path is used. Consider that just at the moment a photon is absorbed by the proton that the proton emits a photon and, when this is absorbed by the electron, the electron immediately emits another photon, etc. Although there is no counting-device, for this back-and-forth process one surmises that over a standard period of time, an ``extremely large number'' of interactions take place. \par

Since the $L$ involved is extremely small as compared to a standard meter, this ``smallness'' allows one to ``model'' an infinitesimal light-clock (inf-light-clock). The counting numbers that correspond to nonzero measures are all members of $\nat_\infty.$ The number $L$ uses the same unit as used to measure $c$ locally, say meters. Further, the number $L$ needs to be taken from an infinite set of special infinitesimals. Why? By Theorem 11.1.1 (11, p. 108) given such an $L$ and any nonzero real number $r$, then there is an $A\in \nat_\infty$ such that $2AL \approx r.$ When {\tt st}, is applied to A(2L) the result is r. Thus A(2L) gives a measure of how far the entire collection of interchanged photons has ``traveled'' linearly in the inf-light-clock. \par

Consider measuring a local distance between two fixed points $F(1)$ and $F(2)$. Since the velocity of light as measured locally will be $c$, such an inf-light-clock can calculate a measure for the ``light'' distance between $F(1)$ and $F(2)$. Let $A$ be the count at the moment a photon leaves point $F(1)$ and $B$ the count for the same or a synchronized inf-light-clock when the photon registers its presence at $F(2).$ The (light) distance from $F(1)$ to $F(2)$ is $2L(B - A).$ Applying {\tt st} yields $r(2) - r(1),$ which is a very accurate distance measurement.\par

The same inf-light-clock can be used to measure ``the (light) time'' between two local events using the time unit $u.$ If $L$ is in meters, then ``seconds'' can be used for $u.$ The time between two successive events $E(1)$ and $E(2)$ occurring at the same point, where inf-light-clock counts $B$ for event $E(2)$ and $A$ for event $E(1)$ is $(2L/c)(B - A).$ (As with ``ticks'' there can be portions of a counting number.) Applying {\tt st}, this yields the standard time measurement. Letting a photon have the ballistics property within an infinitesimal neighborhood, the basic derivation yields that, when it ``moves'' from one neighborhood to another in our physical universe, it acquires the wave property that the standard locally measured velocity c is not altered by the velocity of the source. Atomic clocks also function using photon properties. \par

The calculus is the most successful mathematical theory ever devised. But, for the question of whether something actually exists in some sort of reality that is akin to these infinitesimal entities and we use such analogue models because we can neither describe nor comprehend the infinitesimal-world in any other way, please consider the following as written by Robinson.\par

``For phenomena on a different scale, such as are considered in Modern Physics, the dimensions of a particular body or process may not be observable directly. Accordingly, the question whether or not the scale of non-standard analysis is appropriate to the physical world really amounts to asking whether or not such a system provides a better explanation of certain observable phenomena than the 
standard system. . . . The possibility that this is the case should be borne 
in mind.''\par

\noindent Fine Hall,\par
\noindent Princeton University. \pars

One of these better explanations might be a NSPPM process that gives photons particle properties and one wave property, even if the frequency  property is only a probabilistic statement.\parm

\noindent{\bf 4. Infinitesimal Light-clocks and Gravitational Fields.}\parm
In what follows, the $\Pi$ objects are members of $\nat_\infty.$ Let $F_s$ be a position where the gravitational potential is not zero and $\Pi_s$ an inf-light-clock count at $F_s$. Considered an identical inf-light-clock located at position $F_m$ with inf-light-clock count $\Pi_m$, where there is no gravitational potential. In the usual manner when compared, usually, $\Pi_s < \Pi_m.$ How is this possible? \par

The basic assumption used here and within modern physical science is that ``length'' $L$ is not altered. The infinitesimal $u= 1/(c10^\omega),\ \omega \in \nat_\infty,$ used at $F_s$ and at $F_m$ does not vary when the velocity of light is ``measured'' at $F_s,$ it will measure to be $c$. The reason for this is that due to the change in energy at $F_s$ (compared to $F_m$) produced by a gravitational potential any form of ``timing'' device used to measure the velocity of light has also been affected by the gravitational potential $F_s.$
To make such a comparisons physically, consider information a as propagated by ``light'' from $F_s$ to $F_m.$ If during propagation the slowing of light by gravitational potentials is assumed, than as the light propagates through the gravitational field and arrives at $F_m$ it would regain the original velocity $c$. Under this assumption, comparatively, $c_s$ used at 
$F_s$, when ``viewed'' from $F_m,$ is less than the $c$ without a gravitational field. Does this comparative ``slowing'' of the velocity of light follow from the theory of inf-light-clocks? \par

Consider $\sqrt {g_1}\,u(\Pi_m - \Pi_m')$ and nonzero $\st {u\Pi_m} = r,\ \st {u\Pi_m'} = r'.$ Then, by Theorem 11.1.1 (11, p. 108), there exist $\Gamma_m,\ \Gamma'_m \in \nat_\infty $ such that $\sqrt{g_1}\,u\Pi_m \approx \sqrt {g_1}\, r\approx u\Gamma_m,\ \sqrt{g_1}\,u\Pi_m' \approx \sqrt {g_1}\, r\approx u\Gamma_m'.$ Hence, 
$$\sqrt{g_1}\, u(\Pi_m - \Pi_m')=u(\Gamma_m - \Gamma_m') + \eps,\eqno (B13)$$
where $\eps \in \mu(0).$ Technically, the $\eps$ cannot be removed from (B13). But, equation (B1), if written in inf-light-clock form, appears to lead to a contradiction for the expression $\sqrt {g_1}u\Delta t^m = u\Delta t^s$ when $u$ is divided and the result is viewed from the infinitesimal world. This occurs since $\sqrt {g_1}$ need not be a member of the nonstandard rational numbers $\Hyper Q$, while $\Delta t^s/\Delta t^m \in \Hyper Q.$ This difficulty does not occur if the clock being used is assumed to vary over an interval $[a,b]$ or the entire nonnegative real numbers.
\par
This ``contradiction'' is eliminated, when $\Delta t^m,\ \Delta t^s$ are translated into inf-light-clock notation, by including the count units in the notation. Recall that $L = 1/10^\omega$ is in meters and $u$ is in seconds. Let $T_m$ be interpreted as count ``ticks'' at $F_m$ and $T_s$ the count ``ticks'' at $F_s.$ The translations are 
$$\sqrt{g_1}\, u(\Pi_m - \Pi_m')({\rm m-sec.}) =  u(\Pi_s - \Pi_s')({\rm s-sec.}),$$ 
$$\sqrt{g_1}\, u(\Pi_m - \Pi_m')T_m=u(\Pi_s - \Pi_s')T_s.\eqno (B14)$$
As an example, let $\sqrt {g_1}= 1/\pi.$ Then diving by non-zero $u$ yields

$$(\Pi_m - \Pi_m')({\rm in\ m-ticks})= \pi\, (\Pi_s - \Pi_s')({\rm in
 \ s-ticks})\eqno (B15),$$
where $(\Pi_m - \Pi_m'),\ (\Pi_s - \Pi_s') \in \hypernat.$ \par

Since the $L$ is invariant and identical inf-light-clocks are used, there is only one way that equations (B14), (B15) can be interpreted to avoid a contradiction. The irrational $\pi$ implies that there can be ``partial'' ticks as well as ``partial'' seconds. For identical non-infinitesimal light-clocks, a partial tick comes about when the photons that produce the tick have not traversed the entire $2M$ distance. A basic reason why this can occur is that, in comparison and for a gravitational field as a propagation medium, the velocity of light has been altered. Assuming that the mathematical model faithfully represents such physical aspects and that the process of ``counting'' is a universal process in that the human ``concept'' of counting is not somehow of other altered by the field, then making an informal *-transfer of this yields   
$$\sqrt{g_1}\, c_m = \sqrt{g_1}c = c_s. \eqno (B16)$$
Equation (B16) also explains (B3) and the atomic-clock measures of time in terms of portions of a sec.\par
 
\noindent{\bf 5. Inf-light-clocks and continuity.}\parm

Suppose that the variations in a static field potential are continuous in terms of the distance $r$ from a center of mass. Then, for the field being considered, $\sqrt {g_1}$ is a continuous function in $r$. From (B16), this implies that $c_s$ is a continuous function in $r$. Let $r \in [a,b],\ b \not=a\geq 0.$ If an ideal physical light-clock is employed, then,  
for each member of $[a,b],$ there would exist distinct counting numbers registered by the light-clock. Hence, this gives a one-to-one function  $f\colon [a,b] \to \nat.$ Assuming $r$ is modeled in this ``trivial'' fashion, then this is impossible. But, this ``impossibility'' is removed when inf-light-clocks are utilized. \par

The basic Theorem 11.1.1 [12, p. 108] shows that for every $0 \not= r \in [a,b]$ there exists an $\Gamma_r\in \nat_\infty$ such that $\st {u\Gamma_r}= r.$ There are, however, infinitely many $x\in \nat_\infty$ that have this same property since, at the least, for each $a \in \nat,\ \st {u(\Gamma_r + a)}= r.$ For $q \in \mu (p),$ let $m(q,p) = \{q +x \mid x \in \mu(0)\}.$ If $t\in \mu(q,p),$ then $t \approx q \approx p$ yields $t \approx p$ and $t \in \mu (p).$  In like manner, if $t \in \mu(p),$ then $t\in \mu(q,p)$ implies that $\mu(q,p) = \mu(p).$ Let $M(\Gamma_r)=\{\Gamma_r/(c10^\omega) + x\mid x \in \mu (0)\}.$ Then $M(\Gamma_r) = \mu(r).$ If $p \in [a,b],\ r \not= p,$ then $\mu(r)\cap \mu(p) = M(\Gamma_r) \cap M(\Gamma_p) =\emptyset.$ Thus, by the Axiom of Choice, there is a one-to-one map $g\colon [a,b] \to \nat_\infty.$ Significantly, there is a many-to-one surjection $k\colon {\cal B}=\bigcup \{M(\Gamma_x)\cap \Hyper [a,b] \mid x \in [a,b]\} = \bigcup \{\mu(x)\cap \Hyper [a,b]\mid x \in [a,b]\}\to  [a,b],$ where $k[M(\Gamma_x)\cap \Hyper [a,b]] = \{x\}.$ The function $k$ is actually the restriction of the standard part operator, {\tt st}, that is defined on the set $G(0)\subset\hyperreal$, where $G(0)= \bigcup\{M(\Gamma_x)\mid x\in \real\}.$ In this new notation, if an inf-light-clock count is $\Pi_r$ and $\st {u\Pi_r} =r,$ then $M(\Pi_r) = \mu (r).$\par

 The set of all ``open'' sets $\tau ,$ where each is contained in $\real$, is called a ``topology'' for (on) $\real$. The set $\Hyper \tau$ does not, in general, form a topology for $\hyperreal.$ There is a topology $\cal T$ for $\hyperreal$ called the Q-topology, where $\Hyper \tau \subset {\cal T}$ and the set $\Hyper \tau$ is used as a base for $\cal T.$ A member of $\cal T$ is called a Q-open set. For the Q-topology, the {\tt st} operator is a Q-continuous mapping on $G(0)$ onto $\real$ and $\mu(r)$ and $G(0)$ are Q-open sets [13]. In general, $\mu(r) \notin \Hyper \tau$. Any function $g$ defined on $\Hyper [a,b]$ is ``microcontinuous'' if and only if for each $p \in [a,b]$ and $q\in \Hyper [a,b],$ where $q \approx p,$ $g(q) \approx g(p)$ [11]. The operator {\tt st} is also ``microcontinuous'' on $\Hyper [a,b]$ for if $p\in [a,b],$ $q \in \Hyper [a,b]$ and $q\approx p$, then $q \in G(0)$ and $\st {q} = \st {p} = p.$ Further, as expected, the restriction of {\tt st} to $\Hyper [a,b],$ in the induced Q-topology, is a Q-continuous operator.\par

 For comparison using the ``$M$'' notation, a function $f\colon [a,b] \to \real,$ is (standard real number) continuous on $[a,b]$ if and only if for each $r\in [a,b],\ \hyper f[M(\Gamma_r)\cap \Hyper [a,b]]\subset M(\Gamma_{f(r)}).$ The function $\hyper f$ is also microcontinuous and Q-continuous on $\Hyper [a,b].$ Further, for compact $[a,b]$, it follows that $\Hyper f[M(\Gamma_r)\cap \Hyper [a,b]]= M(\Gamma_{f(r)})\cap [c,d],$ where, since the image is compact and connected, $f[[a,b]] = [c,d].$ But, $\st {M(\Gamma_r)\cap \Hyper [a,b]]} = \{r\}.$ Thus, based upon the compact and connected properties of $[a,b]$ and image $f[[a,b]]=[c,d]$, in general, the Q-continuity of {\tt st} is a very specific and a ``stronger'' type of continuity than standard continuity. This follows since $\Hyper f[M(\Gamma_r)\cap \Hyper [a,b]]= M(\Gamma_{f(r)})\cap [c,d]$, while the Q-continuity of {\tt st} requires that {\tt st}$[M(\Gamma_r)\cap \Hyper [a,b]] = \st {\{r\}}=\{r\}.$ The Q-continuity of {\tt st} on $G(0)$ is not related to a standard field continuity that might be a property of $\sqrt {g_1}.$ Q-continuity on $G(0)$ applies to any form of alteration in $\sqrt {g_1},$ even an abrupt quantum physical alteration. This follows since, obviously, for any function $f\colon [a,b] \to \real,$ $\st {\hyper f(p)}= f(p).$ Hence, the operator {\tt st} has no affect upon the values of $f(p)$ for any $p\in [a,b]$ in the sense that {\tt st} merely mimics many of the standard $f$-characteristics.   \parm

\centerline{\bf References}\smallskip
\noindent [1] Bergmann, G., Introduction to the Theory of Relativity, Dover, New York, 1976.\par
\noindent [2] Craig, H. V., Vector and Tensor Analysis, McGraw-Hill, New York, 1943.\par
\noindent [3] Cranshaw, T. E., J. P. Schiffer and P. A. Egelstaff, Measurement of the red shift using the M\"{o}ssbauer effect in $\rm Fe^{57},$ Phys. Rev. Letters 4(4)(1960):163-164,\par
\noindent [4] Fokker, Albert Einstein, inventor of chronogeometry, Synth\'{e} se 9:442-444.\par
\noindent [5] Herrmann, R. A., Nonstandard Analysis Applied to Special and General Relativity - The Theory of Infinitesimal Light-Clocks, 1992,93,94,95.\hfil\break http://arxiv/abs/math/0312189 \par
\noindent [6] Lawden, D. F., An Introduction to Tensor Calculus, Relativity and Cosmology, John Wiley \& Sons, New York, 1982.\par
\noindent [7] McConnell, A. J., Applications of Tensor Analysis, Dover, New York, 1957.\par
\noindent [8] Ohanian, H. and R. Ruffini, Gravitation and Spacetime, W. W. Nortin Co., New York, 1994.\par
\noindent [9] Poncelet, V., Trait\'{e} des propri\'{e}t\'{e}s projective des figures, 1822.\par
\noindent [10] Rindler, W., Essential Relativity, Springer-Verlag, New York, 1977. \parm

\centerline{\bf Additional Additional References}\parm
\noindent [11] Davis, M., Applied Nonstandard Analysis (Wiley-Interscience), New York, 1977. \par 
\noindent [12] Herrmann, R. A., The Theory of Ultralogics, 1993, \hfil\break http://www.arxiv.org/abs/math.GM/9903081\hfil\break 
http://www.arxiv.org/abs/math.GM/9903082.\par
\noindent [13] Herrmann, R. A. 1976. The Q-topology, Whyburn type filters and the cluster set map, Proc. Amer. Math. Soc., 59(1976):161--166.\eject
\centerline{NOTES}\vfil\eject
     %This also has the same size font as megstep 1
 \catcode`@=11 % from plain.tex
\newdimen\pagewidth \newdimen\pageheight \newdimen\ruleht
% These values were modified:
   \hsize=6.5in  \vsize=54pc  \maxdepth=2.2pt  \parindent=2pc
   \hoffset=.35in
\pagewidth=\hsize  \pageheight=\vsize  \ruleht=.5pt

%This routine is used by \output; this is different from 
%  the one found in App. E since some are not needed here.
\def\onepageout#1{\shipout\vbox{\offinterlineskip
  \vbox to \pageheight {\makeheadline
          #1 % the content of page
          \makefootline \boxmaxdepth=\maxdepth}}
  \advancepageno}

\output{\onepageout{\unvbox255}}
\newbox\partialpage
\def\begind{\begingroup
  \output={\global\setbox\partialpage=\vbox{\unvbox255\bigskip}}\eject
  \output={\doublecolumnout} \hsize=2.85in \vsize=109pc} %2.85
  % Again the sizes have been changed  !!

\def\doublecolumnout{\splittopskip=\topskip \splitmaxdepth=\maxdepth
  \dimen@=54pc \advance\dimen@ by-\ht\partialpage
  % Need to change the value of \dim@ also...
  \setbox0=\vsplit255 to\dimen@ \setbox2=\vsplit255 to\dimen@
  \onepageout\pagesofar
  \unvbox255 \penalty\outputpenalty}
\def\pagesofar{\unvbox\partialpage
  \wd0=\hsize \wd2=\hsize \hbox to\pagewidth{\box0\hfil\hfil\box2}}
\def\balancecolumns{\setbox0=\vbox{\unvbox255} \dimen@=\ht0
  \advance\dimen@ by\topskip \advance\dimen@ by-\baselineskip
  \divide\dimen@ by2 \splittopskip=\topskip
  {\vbadness=10000 \loop \global\setbox3=\copy0
    \global\setbox1=\vsplit3 to\dimen@
    \ifdim\ht3>\dimen@ \global\advance\dimen@ by1pt \repeat}
  \setbox0=\vbox to\dimen@{\unvbox1}
  \setbox2=\vbox to\dimen@{\dimen2=\dp3 \unvbox3\kern-\dimen2 \vfil}
  \pagesofar}

\begind
                       %but it is slight more spread out.
\def\sub#1{{\leftskip=0.2in \noindent #1 \par}\par}
\noindent Due to revisions, page\par 
\noindent locations are only\par 
\noindent approximations. One or two\par 
\noindent may need to be added to the\par
\noindent page number.\par
\noindent \par
 
\indent\indent{\bf A}\par
\vskip 12pt 
\noindent Abel and Cauchy 10.\par
\noindent absolute\par
\sub {length, no meaning 12.}\par
\sub {realism 8.}\par 
\sub {time  10, 19, 20}\par
\sub {time, no meaning 10.}\par
\sub {only known one 26.}\par
\noindent alterations derived\par
\sub {decay rates 63.}\par
\sub {energy shifts 62.}\par
\sub {gravitational 65.}\par
\sub {mass effects 63.}\par
\sub {transverse Doppler 62.}\par
\noindent altered by $P$-process 41.\par
\noindent analogue model 39.\par
\sub {light-clock counts 42, 43.}\par
\sub {Riemannian geometry 16, 76.}\par
\noindent approximate, continuum 19.\par
\noindent atomic clocks 63.\par
\noindent atomic, electromagnetic\par
 radiation 62.\par
\noindent {\quad}\par
\noindent {\quad} 
\indent\indent{\bf B}\par 
\noindent {\quad}\par 
\noindent Barnes 61.\par
\noindent black hole\par
\sub {diverse scenarios 75.}\par
\sub {formation 74.}\par
\sub {halo effect possible 75.}\par
\sub {leads to a quasi-white hole 75.}\par
\sub {possible transitional zone 75.}\par
\sub {spherical shell effect possible 75.}\par
\noindent bookkeeping technique 43.\par
\noindent bounded, cosmos 76.\par
\noindent bounded, finite hyperreal\par
 numbers 20.\par
\noindent Breitner 68.\par
\noindent Builder 61.\par
\noindent {\quad}\par 
\noindent {\quad}\par
\noindent {\quad}\par
\noindent {\quad}\par
\indent\indent{\bf C}\par
\noindent {\quad}\par
\noindent c, velocity of light \par
\sub {possible not fixed for NSPPM time 40.}\par
\noindent Cartesian coordinate system,\par
 inertial 18.\par
\noindent catalyst, time 24.\par 
\noindent Cauchy, error  9.\par
\noindent chronogeometry 16.\par
\noindent {\quad}\par 
\noindent {\quad}\par 
\noindent chronotopic interval 17, 56.\par
\noindent clock, its many definitions 10.\par
\noindent close to, infinitesimally 20.\par
\noindent collapse\par
\sub {optical appearance 55.}\par
\sub {restricted 69.}\par
\noindent comparisons to standard only\par
has human meaning 73.\par
\noindent conceptual observer 10.\par
\sub {cannot reject 10.}\par
\sub {Einstein 10.}\par
\noindent constancy, velocity of light 20.\par
\noindent content, descriptive 19.\par
\noindent continuity, S 28.\par
\noindent continuum\par
\sub {approximating 20.}\par
\sub {time 10.}\par
\noindent contraction, length 25.\par
\noindent contradiction\par
\sub {Einstein's postulates and the derivation 24.}\par
\noindent coordinate\par
\sub {change acceptable 16.}\par
\sub {gravitation field, alteration in 16.}\par
\sub {systems and Riemannian geometry 12.}\par
\sub {transformations, differentiable with nonvanishing Jocabians 68.}\par
\noindent Copernican principle = \par
Cosmological principle 68.\par
\noindent cosmological expansion line-\par
 element derived 68.\par
\noindent Cosmological principle =\par
 Copernican principle 68.\par
\noindent cosmological redshift 69.\par
\noindent cosmos, bounded 76.\par
\noindent count, infinite 20.\par
\noindent counting mechanism 16.\par
\noindent ``creation'' (formation) \par
\sub {white holes at 74.}\par
\sub {explosive effect 74.}\par  
\sub {pseudo-white holes 74.}\par
\noindent criticism, Fock 14.\par 
\noindent {\quad}\par
\noindent {\quad}\par
\indent\indent{\bf D}\par
\noindent {\quad}\par
\noindent de Sitter line-element 66.\par
\noindent decay rates\par
\sub {alterations derived 63.}\par
\sub {gravitational derived 66.}\par
\noindent deceleration parameter 62.\par
\noindent derivations from fundamental\par
 properties, 9.\par
\noindent descriptive content 18.\par
\noindent dilation, time 25.\par 
\noindent Dingle\par
\sub {no absolute motion 20.}\par
\sub {only known absolute 26.}\par
\noindent directed numbers, not modeled\par
 by 61.\par
\noindent distance function  27.\par
\noindent {\quad}\par
\noindent {\quad}\par
\indent\indent{\bf E}\par
\noindent {\quad}\par
\noindent Eddington-Finkelstein\par
 transformation 70.\par
\noindent Einstein 10.\par
\sub {logical errors and Fock 14.}\par
\sub {measures 20, 25.}\par
\sub {original paper 10.}\par
\sub {hypotheses 11.}\par
\noindent Einstein--Rosen bridge\par
 (wormholes), none yet 69.\par
\noindent electromagnetic propagation 8.\par
\sub {Galilean, infinitesimal 20.}\par
\sub {Euclidean neighborhood 25.}\par
\noindent electromagnetic radiation,\par
 atomic 62.\par
\noindent {\quad}\par
\noindent {\quad}\par
\noindent {\quad}\par

\noindent emis, effects 17, 61.\par
\sub {defined 40.}\par
\noindent empty, space-time 24.\par
\noindent energy shift, Schr\"odinger equation\par
 approach 62.\par
\noindent equilinear 40.\par
\noindent Equivalence Principle\par
\sub {does not generally hold 15.}\par
\sub {effects infinitesimal and local 15.}\par
\noindent error of generalization 11.\par
\noindent errors, logical  11.\par
\noindent \ae{ther} = medium 8.\par
\noindent \ae{ther}\par
\sub {calculations 8.}\par
\sub {removed by postulating 12.}\par 
\noindent Euclidean neighborhood, \par
electromagnetic propagation\par
25.\par 
\noindent evidence, indirect 10.\par
\noindent expansion of universe\par
\sub {and Special Theory 25.}\par
\sub {rate, extreme at formation 76.}\par
\sub {NSPPM velocity effect 68.}\par
\noindent explosive effects, at formation 76.\par
\noindent {\quad}\par 
\noindent {\quad}\par 
\indent\indent{\bf F}\par
\noindent {\quad}\par
\noindent finite = bounded = limited 20.\par
\noindent first approximation, Newtonian 67.\par
\noindent Fock 14, 69.\par
\sub {comparison with human intuition only has meaning 73.}\par
\sub {equivalence principle 15.}\par
\sub {harmonic coordinates 15.}\par
\noindent Fokker, chronogeometry 16.\par
\noindent force-like, interaction 18.\par
\noindent fractal curve  12.\par
\noindent Friedmann\par
\sub {closed universe differential equation 69.}\par
\sub {model, positive curvature 69.}\par
\sub {open universe model 69.}\par
\noindent function, universal 62.\par
\noindent fundamental properties,\par
 derivations from 9.\par
\noindent {\quad}\par
\noindent {\quad}\par
\indent\indent{\bf G}\par
\noindent {\quad}\par
\noindent Galilean\par
\sub {electromagnetic propagation theory 8, 20.}\par
\sub {theory of average velocities (velocities) 27-28.}\par
\noindent General Theory,  logical errors 11.\par
\noindent generalization, error of 11.\par
\noindent geometry, human construct 17.\par
\noindent Gerber  58.\par
\noindent gravitational\par
\sub {alterations in the radioactive decay 66.}\par
\sub {field, space-time geometry not a physical cause 57-58.}\par
\sub {redshift derived  66.}\par
\noindent {\quad}\par
\noindent {\quad}\par
\indent\indent{\bf H}\par
\noindent {\quad}\par
\noindent halo effect with some black holes\par
 75.\par
\noindent harmonic coordinates and Fock 15.\par
\noindent Heaviside 68.\par
\noindent Hubble Law 69.\par
\noindent human comprehension\par
\sub {and geometry 8.}\par
\sub {Planck statement on 8.}\par
\noindent human intuition, only meaning is\par 
by comparison with 73.\par 
\noindent human mind and imaginary\par
 entities  9.\par
\noindent Humphreys 75.\par
\noindent Huygens, medium 8.\par
\noindent hypotheses, Einstein 11.\par
\noindent {\quad}\par

\indent\indent{\bf I}\par
\noindent {\quad}\par
\noindent imaginary entities 9.\par
\noindent indirect evidence  10.\par
\noindent indistinguishable \par
\sub {effects, 57.}\par
\sub {for dt 27.}\par
\sub {first level 27.}\par
\noindent {\quad}\par
\noindent {\quad}\par
\noindent {\quad}\par
\noindent inertial 19.\par
\noindent infinite, Robinson numbers 20.\par
\noindent infinitesimal\par
\sub {light-clock 20, 58.}\par
\noindent infinitesimalizing 25.\par
\sub {and the calculus rules 26.}\par
\sub {simple behavior 17.}\par
\noindent infinitesimally,  close to 20.\par
\noindent infinitesimals, \par
\sub {Cauchy 9.}\par
\sub {Einstein error 11.}\par
\sub {Newton 10.}\par
\sub {Robinson 10.}\par
\noindent instantaneous velocity = ultimate\par
 velocity  10.\par
\noindent instantaneous, snapshot effects 61.\par
\noindent interaction, force-like only 18.\par
\noindent invariant\par
\sub {forms 58.}\par
\sub {solution methods 58.}\par
\sub {statement 16.}\par
\noindent Ives 56.\par
\noindent Ives-Stillwell 61.\par
\noindent {\quad}\par
\noindent {\quad}\par
\indent\indent{\bf J}\par
\noindent {\quad}\par
\noindent Jefferson 68.\par
\noindent {\quad}\par
\noindent {\quad}\par
\indent\indent{\bf K}\par
\noindent {\quad}\par
\noindent Kennedy-Thorndike 41.\par
\noindent Kerr transformation 73.\par
\sub {as science-fiction 73.}\par
\noindent Kruskal space, insurmountable\par
 difficulties with 73.\par
\noindent {\quad}\par
\noindent {\quad}\par
\indent\indent{\bf L}\par
\noindent {\quad}\par
\noindent language, convenient,\par
 Riemannian geometry 55.\par
\noindent language, corresponds to math.\par structure 56.\par
\noindent Laplacian 67.\par
\noindent length contraction\par
\sub {(emis) effect 44.}\par
\sub {modern approach 44.}\par 
\sub {not absolute effect 56.}\par
\noindent length, no alteration in 56.\par
\noindent light-clock, infinitesimalized 56.\par
\noindent light propagation 10.\par
\sub {Milne 26.}\par
\sub {only known absolute 26.}\par
\sub {principles 25.}\par
\noindent light velocity measurement, how\par
 made 30.\par
\noindent light-clock\par
\sub {count change 42.}\par
\sub {counting mechanism 16.}\par
\sub {counts as an analogue model 42.}\par
\sub {diagram 35.}\par
\sub {infinitesimal 19.}\par
\sub {ticks 19.}\par
\sub {timing, orientation 58.}\par
\noindent limited =  finite 20, 28.\par
\noindent line-element\par
\sub {cosmological expansion derived 68.}\par
\sub {de Sitter 67.}\par
\sub {Galilean  15.}\par
\sub {linear effect 62.}\par
\sub {Minkowski 12, 15, 55.}\par
\sub {modified Schwarzschild 67.}\par
\sub {partial 64.}\par
\sub {proper time-like 10.}\par
\sub {quasi-Schwarzschild 66-67.}\par 
\sub {quasi-time-like 57.}\par
\sub {Robertson-Walker, derived 68.}\par
\sub {Schwarzschild 66.}\par
\noindent linear effect line-element 60.\par
\noindent linear light propagation,\par
 to-and-fro 41.\par
\noindent local measure of the velocity 55.\par
\noindent location, fixed in a NSPPM  16.\par
\noindent logical error\par
\sub {in Special Theory arguments 42.}\par
\sub {General Theory 11.}\par
\sub {modeling error 11.}\par
\sub {predicate errors 10.}\par
\noindent {\quad}\par
\noindent {\quad}\par
\noindent {\quad}\par
\sub {time 16.}\par
\noindent Lorentz 8.\par
\sub {transformation  11.}\par
\sub {altering realism  9.}\par
\noindent luminiferous \ae{ther} 8.\par
\noindent {\quad}\par
\noindent {\quad}\par
\indent\indent{\bf M}\par
\noindent {\quad}\par
\noindent m superscript, relative motion with\par
 respect to
 stationary,\par standard, or
 observer, altered\par 41.\par
\noindent M-M = Michelson-Morley 41.\par
\noindent MA-model scenario 74.\par
\noindent mass alterations derived 64.\par 
\noindent material particle 74.\par
\noindent math. models and Newton 10.\par           
\noindent math. structure, corresponding to\par
 physical language 44.\par
\noindent Maxwell \ae{ther}, 8. \par
\sub {removed by postulating 12.}\par
\noindent measure\par
\sub {Einstein 11, 19.}\par
\sub {velocity of light 19.}\par
\noindent mechanical behavior, the\par
 calculus 26.\par
\noindent medium\par
\sub {ether, Maxwell 8, 12.}\par
\sub {Huygens, 8.}\par
\sub {Newton 8.}\par
\sub {NSPPM 10, 19.}\par
\sub {the NSPPM 19.}\par
\sub {Thomson 8.}\par
\noindent metamorphic, (i.e. sudden)\par 
structured  change 68.\par
\noindent Milne, light propagation theory 25.\par
\noindent Minkowski-type interval =\par
 chronotopic 15.\par
\noindent Minkowski-type line-element 15,\par
 55.\par
\noindent missing physical quantities 43.\par
\noindent model theoretic error of\par
 generalization  11.\par
\noindent model, nature required to follow\par
 13.\par
\noindent modeling, mathematical, schism 10.\par
\noindent modified Schwarzschild\par
 line-element 67.\par
\noindent monad 27, 58.\par
\noindent monadic cluster 58.\par
\noindent monadic neighborhood 58.\par
\noindent motion, relative, superscript m, with\par
 respect to stationary, standard, or\par
observer, 41.\par
\noindent {\quad}\par 
\noindent {\quad}\par 
\indent\indent{\bf N}\par
\noindent {\quad}\par
\noindent N-world relative velocity, N-world, \par
nonderived 32.\par
\noindent nature, required to behave as model\par
 dictates  14.\par
\noindent near to, infinitesimally 20.\par
\noindent Newton\par
\sub {calculus and mechanical behavior 26.}\par 
\sub {ether, medium 8.}\par
\sub {infinitesimals 9.}\par
\sub {natural world implies math. 9.}\par
\sub {ultimate velocity 9.}\par
\sub {velocities 26.}\par
\noindent Newtonian\par
\sub {first approximation 67.}\par
\sub {gravitational potential 65.}\par
\sub {time 10, 19.}\par 
\noindent nonstandard electromagnetic field\par
-- NSPPM 10, 19.\par
\noindent nonstandard physical world model\par
= NSP-world 9.\par 
\noindent NSPPM = nonstandard\par
 photon-particle medium 10, 19.\par
\noindent NSP-world = nonstandard physical\par
 world  10.\par
\sub {linear ruler 27.}\par
\sub {time and a stationary subparticle field 27.}\par
\noindent numbers\par
\sub {directed, not modeled by 61.}\par
\sub {real 11.}\par
\noindent {\quad}\par
\noindent {\quad}\par
\indent\indent{\bf O}\par
\noindent {\quad}\par
\noindent observer\par
\sub {conceptual, cannot reject 13.}\par
\sub {privileged, rejection 13, 16.}\par
\sub {standard, stationary, s, 41.}\par
\noindent Ohanian 68.\par
\noindent operational definition 10.\par
\noindent operator, separating 62.\par
\noindent {\quad}\par
\noindent {\quad}\par
\indent\indent{\bf P}\par
\noindent {\quad}\par
\noindent partial differential\par
 calculus and Einstein 25.\par
\noindent partial line-elements 74.\par
\noindent particle, material 74.\par
\noindent Patton and Wheeler 18.\par
\noindent Phillps 58.\par
\noindent philosophy\par
\sub {realism 8.}\par
\sub {scientism 8.}\par
\noindent photon, language 27.\par
\noindent physical\par
\sub {meaning to contraction, none 36-37.}\par
\sub {quantities, missing 44.}\par
\sub {language and math. structure 43.}\par
\noindent Planck, human comprehension 8.\par
\noindent positive curvature Friedmann\par
 model 69.\par
\noindent postulate away the existence of real\par
 Maxwellian substratum 12.\par
\noindent potential velocity 19, 55-56.\par
\noindent predicate, logical  10.\par
\noindent pregeometry, Patton and\par
 Wheeler 18.\par
\noindent privileged observer\par
\sub {fixed in NSPPM 15.}\par
\sub {inertial Cartesian coordinate system 19.}\par
\sub {rejection of 3.}\par
\noindent processes, certain ones\par
\sub {slowing down 75.}\par
\noindent Prokhovnik 20.\par
\noindent proper-time-like, line-element 56.\par
\noindent pseudo-white hole  76.\par
\noindent {\quad}\par
\noindent {\quad}\par
\indent\indent{\bf Q}\par
\noindent {\quad}\par
\noindent quasi-Schwarzschild line-element\par
 66-67.\par
\noindent quasi-time-like, line-element 57.\par
\noindent {\quad}\par
\noindent {\quad}\par
\indent\indent{\bf R}\par
\noindent {\quad}\par
\noindent radar or reflected light pulse \par
method 73.\par
\noindent radiation, atomic,\par
 electromagnetic 62.\par
\noindent radioactive decay rates\par
\sub {Special theory alterations derived 44, 64.}\par
\sub {gravitational alterations derived 66.}\par
\noindent radius, Schwarzschild 69, 73.\par
\noindent real numbers  11.\par
\noindent realism, absolute 8.\par
\noindent redshift\par
\sub {cosmological 55, 68.}\par 
\sub {gravitational derived 66, 69.}\par
\sub {transverse (Doppler) derived 63.}\par
\noindent reflected light pulse method 73.\par 
\noindent relative motion, m superscript with\par
 respect to 
stationary,\par
 standard, or
observer, altered 41.\par
\noindent relative velocity\par
\sub {measured of 36.}\par 
\sub {nonderived N-world 32.}\par
\noindent relativistic redshift, transverse 63.\par
\noindent Riemannian geometry 14.\par
\noindent Riemannian geometry, 
\sub {an analogue model 16.}\par
\sub {convenient language 55.}\par
\sub {coordinate systems 14.}\par 
\noindent Robertson-Walker line-element\par
 derived 68.\par
\noindent Robinson, infinitesimals 10.\par
\noindent {\quad}\par
\noindent {\quad}\par
\indent\indent{\bf S}\par
\noindent {\quad}\par
\noindent s superscript for stationary,\par
 standard or, sometimes,
observer\par 
41.\par
\noindent S-continuity 28.\par
\noindent Sagnac type of experiment 42.\par
\noindent schism\par
\sub {in mathematical modeling 10.}\par
\sub {and Newton  9.}\par
\noindent Schr\"odinger equation,\par
 time-dependent 63.\par
\sub {applied to any energy shift 63.}\par
\noindent Schwarzschild\par
\sub {line-element, derived 66.}\par
\sub {radius, reduced 69, 73.}\par
\sub {surface, collapse through 69, 72.}\par
\noindent scientism 8.\par
\noindent separating operator 62.\par
\noindent simple behavior, in the small 17.\par
\noindent simplicity, rule 25.\par
\noindent slow down, certain processes 76.\par
\noindent small, in the 17.\par
\noindent snapshot effect 61.\par
\noindent space-time\par
\sub {geometry, not a physical cause 55.}\par
\sub {empty 24.}\par
\noindent Special theory, logical error in\par
arguments 42.\par
\noindent spherical shell effect, black holes\par
75.\par
\noindent spherical wavefront (light)\par
 concepts 11.\par
\noindent standard part operator 20.\par
\noindent standard, stationary, observer, s\par
superscript 41.\par
\noindent statement, invariant 16.\par
\noindent stationary, standard, observer, s\par
superscript 41.\par
\noindent subparticle NSPPM, removing \par
energy 75.\par
\noindent subparticles, ultimate 58.\par
\noindent NSPPM\par
\sub {expansion, contraction 68.}\par
\sub {fixed observer in 16.}\par 
\sub {Galilean rules for velocity composition 24, 27-28.}\par
\sub {light propagation principles 27-28.}\par
\sub {medium 27}\par
\noindent Surdin 58.\par
\noindent surface, Schwarzschild 65, 70, 74.\par
\noindent {\quad}\par
\noindent {\quad}\par
\indent\indent{\bf T}\par
\noindent {\quad}\par
\noindent textural expansion not relative to\par
the Special Theory 25.\par
\noindent Thomson, medium 8.\par
\noindent ticks, light-clock 20.\par
\noindent time continuum 10.\par
\noindent time dilation, no such\par
\sub {effect 61.}\par
\sub {modern approach 25.}\par
\noindent time\par
\sub {absolute 10.}\par
\sub {and logical error 16.}\par
\sub {as a catalyst 25.}\par
\sub {continuum, returning to 20.}\par
\sub {light-clock ticks 20.}\par
\sub {NSP-world 20.}\par
\sub {universal 10.}\par
\noindent time-dependent indistinguishable \par
effects 57.\par
\noindent time-dependent Schr\"odinger\par
 equation 63.\par
\noindent timing infinitesimal light-clock,\par
 orientations 58.\par
\noindent to-and-fro\par
\sub {linear light propagation 38.}\par
\sub {natural world measurement, light 19.}\par
\noindent transitional phase, the apparent\par
 turbulent physical behavior\par
 75.\par
\noindent transitional zone 73.\par
\noindent transverse\par
\sub {Doppler, redshift derived 63.}\par
\noindent turbulent physical behavior in\par
 transitional zone 75.\par
{\quad}\par
{\quad}\par

\indent\indent{\bf U}\par
\noindent {\quad}\par
\noindent ultimate\par
\sub {subparticles  58.}\par
\sub {ultrawords 68.}\par
\sub {velocity, Newton  9.}\par
\noindent ultralogic 68.\par
\noindent ultrawords, ultimate 68.\par
\noindent universal\par
\sub {function 62.}\par
\sub {time 10.}\par 
\noindent Upham 58. \par
\noindent {\quad}\par
\noindent {\quad}\par
\indent\indent{\bf V}\par
\noindent {\quad}\par
\noindent velocity\par
\sub {different composition rules 25.}\par
\sub {Galilean 27-28.}\par
\sub {light, constancy 18.}\par
\sub {local measure 55-56.}\par
\sub {Newton's concepts based upon 26.}\par
\sub {potential 18, 55-56.} \par
\sub {relative, measure of 36.}\par
\sub {NSPPM, expansion, contraction 68.}\par
\noindent {\quad}\par
\noindent {\quad}\par
\indent\indent{\bf W}\par
\noindent {\quad}\par
\noindent wormholes, none.\par 69.\par

\end